\documentclass[12pt, a4paper]{scrartcl}  
\usepackage[T1]{fontenc}
\usepackage[utf8]{inputenc}
\usepackage[english]{babel}
\usepackage{amsmath, amsfonts}
\usepackage{amssymb}
\usepackage{amsthm}
\usepackage[all]{xy}
\usepackage{cancel}
\usepackage{mathtools}
\usepackage[inline]{enumitem}
\usepackage{float}
\usepackage{tikz}
\usetikzlibrary{arrows}
\usetikzlibrary{arrows.meta, calc}
\usepackage{hyperref}
\hypersetup{
    colorlinks=true,
    linkcolor=blue,
    citecolor=blue,
    urlcolor=blue,
    pdfpagemode=UseOutlines
}
\usepackage{aliascnt}
\usepackage{cleveref}

\theoremstyle{plain}

\newaliascnt{Thm}{sa}\newtheorem{Thm}[Thm]{Theorem}\aliascntresetthe{Thm}
\newaliascnt{Lem}{sa}\newtheorem{Lem}[Lem]{Lemma}\aliascntresetthe{Lem}
\newaliascnt{Prp}{sa}\newtheorem{Prp}[Prp]{Proposition}\aliascntresetthe{Prp}
\newaliascnt{Cor}{sa}\newtheorem{Cor}[Cor]{Corollary}\aliascntresetthe{Cor}
\newaliascnt{Con}{sa}\aliascntresetthe{Con}
\newaliascnt{Fct}{sa}\aliascntresetthe{Fct}
\newaliascnt{Def}{sa}\newtheorem{Def}[Def]{Definition}\aliascntresetthe{Def}
\newaliascnt{DefRem}{sa}\aliascntresetthe{DefRem}
\newaliascnt{DefThm}{sa}\aliascntresetthe{DefThm}
\newaliascnt{DefLem}{sa}\aliascntresetthe{DefLem}
\newaliascnt{Axm}{sa}\aliascntresetthe{Axm}
\newaliascnt{Not}{sa}\newtheorem{Not}[Not]{Notation}\aliascntresetthe{Not}
\newaliascnt{Rem}{sa}\newtheorem{Rem}[Rem]{Remark}\aliascntresetthe{Rem}
\newaliascnt{Cau}{sa}\aliascntresetthe{Cau}
\newaliascnt{Eg}{sa}\newtheorem{Eg}[Eg]{Example}\aliascntresetthe{Eg}
\newaliascnt{Tho}{sa}\aliascntresetthe{Tho}

\crefname{sa}{Theorem}{Theorems}                 \Crefname{sa}{Theorem}{Theorems}
\crefname{Thm}{Theorem}{Theorems}                \Crefname{Thm}{Theorem}{Theorems}
\crefname{Lem}{Lemma}{Lemmas}                    \Crefname{Lem}{Lemma}{Lemmas}
\crefname{Prp}{Proposition}{Propositions}        \Crefname{Prp}{Proposition}{Propositions}
\crefname{Cor}{Corollary}{Corollaries}           \Crefname{Cor}{Corollary}{Corollaries}
\crefname{Con}{Conjecture}{Conjectures}          \Crefname{Con}{Conjecture}{Conjectures}
\crefname{Fct}{Facts}{Facts}                     \Crefname{Fct}{Facts}{Facts}
\crefname{Def}{Definition}{Definitions}          \Crefname{Def}{Definition}{Definitions}
\crefname{DefRem}{Definition}{Definitions}       \Crefname{DefRem}{Definition}{Definitions}
\crefname{DefThm}{Definition}{Definitions}       \Crefname{DefThm}{Definition}{Definitions}
\crefname{DefLem}{Definition}{Definitions}       \Crefname{DefLem}{Definition}{Definitions}
\crefname{Axm}{Axiom}{Axioms}                    \Crefname{Axm}{Axiom}{Axioms}
\crefname{Not}{Notation}{Notations}              \Crefname{Not}{Notation}{Notations}
\crefname{Rem}{Remark}{Remarks}                  \Crefname{Rem}{Remark}{Remarks}
\crefname{Cau}{Caution}{Cautions}                \Crefname{Cau}{Caution}{Cautions}
\crefname{Eg}{Example}{Examples}                 \Crefname{Eg}{Example}{Examples}
\crefname{Tho}{Thoughts}{Thoughts}               \Crefname{Tho}{Thoughts}{Thoughts}

\newcommand{\N}{\mathbb{N}}
\newcommand{\Z}{\mathbb{Z}}
\newcommand{\R}{\mathbb{R}}
\newcommand{\Q}{\mathbb{Q}}
\newcommand{\id}{\mathrm{id}}
\newcommand{\Acal}{\mathcal{A}}
\newcommand{\Bcal}{\mathcal{B}}
\newcommand{\Ccal}{\mathcal{C}}
\newcommand{\Dcal}{\mathcal{D}}
\newcommand{\Ecal}{\mathcal{E}}
\newcommand{\Fcal}{\mathcal{F}}
\newcommand{\Hcal}{\mathcal{H}}
\newcommand{\Ical}{\mathcal{I}}
\newcommand{\Mcal}{\mathcal{M}}
\newcommand{\Ncal}{\mathcal{N}}
\newcommand{\Pcal}{\mathcal{P}}
\newcommand{\Qcal}{\mathcal{Q}}
\newcommand{\Rcal}{\mathcal{R}}
\newcommand{\Scal}{\mathcal{S}}
\newcommand{\Tcal}{\mathcal{T}}
\newcommand{\Ucal}{\mathcal{U}}
\newcommand{\Wcal}{\mathcal{W}}
\newcommand{\Xcal}{\mathcal{X}}
\newcommand{\Ycal}{\mathcal{Y}}
\newcommand{\Zcal}{\mathcal{Z}}
\newcommand{\kerstyle}[1]{\mathbf{#1}}
\newcommand{\Ek}{\kerstyle{E}}
\newcommand{\Gk}{\kerstyle{G}}
\newcommand{\Kk}{\kerstyle{K}}
\newcommand{\Mk}{\kerstyle{M}}
\newcommand{\Pk}{\kerstyle{P}}
\newcommand{\Qk}{\kerstyle{Q}}
\newcommand{\Nk}{\kerstyle{N}}
\newcommand{\Rk}{\kerstyle{R}}
\newcommand{\Tk}{\kerstyle{T}}
\newcommand{\Uk}{\kerstyle{U}}
\newcommand{\Xk}{\kerstyle{X}}
\newcommand{\Yk}{\kerstyle{Y}}
\newcommand{\Zk}{\kerstyle{Z}}
\newcommand{\deltabf}{{\boldsymbol \delta}}
\newcommand{\mubf}{{\boldsymbol \mu}}
\newcommand{\nubf}{{\boldsymbol \nu}}
\newcommand{\lambdabf}{{\boldsymbol \lambda}}
\newcommand{\rhobf}{{\boldsymbol \rho}}
\newcommand{\ismapof}{\precsim}
\newcommand{\sm}{\setminus}						
\newcommand{\ins}{\subseteq} 					
\newcommand{\sni}{\supseteq} 					
\newcommand{\cmpl}{\mathsf{c}}
\newcommand{\pr}{\mathrm{pr}} 	
\newcommand{\ev}{\mathrm{ev}} 
\newcommand{\inj}{\hookrightarrow}
\newcommand{\dshto}{\dashrightarrow}
\newcommand{\E}{\mathbb{E}}
\newcommand{\I}{\mathbf{1}}
\newcommand{\Pa}{\mathrm{Pa}} 		
\newcommand{\Ch}{\mathrm{Ch}} 	
\newcommand{\Anc}{\mathrm{Anc}} 		
\newcommand{\Desc}{\mathrm{Desc}} 	
\newcommand{\Pred}{\mathrm{Pred}} 
\DeclareMathOperator*{\Indep}{\perp\!\!\!\perp}
\DeclareMathOperator*{\Perp}{\perp}
\DeclareMathOperator*{\nPerp}{\cancel\perp}
\DeclareMathOperator*{\given}{|}
\DeclareMathOperator{\doit}{do}
\newcommand{\dcup}{\,\dot{\cup}\,}
\newcommand{\lp}{\left ( }
\newcommand{\rp}{\right ) }
\newcommand{\lB}{\left [ }
\newcommand{\rB}{\right ] }
\newcommand{\lC}{\left \{ }
\newcommand{\rC}{\right \} }
\newcommand{\arrhead}{{latex}}
\newcommand{\arrtail}{{stealth reversed}}
\newcommand{\arrstar}{{Rays[n=6]}}
\newcommand*{\hut}[1][]{\mathrel{\tikz [baseline=-0.25ex,\arrhead-\arrtail, #1] \draw [#1] (0pt,0.5ex) -- (1.3em,0.5ex);}}
\newcommand*{\tuh}[1][]{\mathrel{\tikz [baseline=-0.25ex,\arrtail-\arrhead, #1] \draw [#1] (0pt,0.5ex) -- (1.3em,0.5ex);}}
\newcommand*{\sus}[1][]{\mathrel{\tikz [baseline=-0.25ex,\arrstar-\arrstar, #1] \draw [#1] (0pt,0.5ex) -- (1.3em,0.5ex);}}
\let\oldFootnote\footnote
\newcommand\nextToken\relax
\renewcommand\footnote[1]{%
    \oldFootnote{#1}\futurelet\nextToken\isFootnote}
\newcommand\isFootnote{%
    \ifx\footnote\nextToken\textsuperscript{,}\fi}
\newcommand{\Asterisk}{\mathord{\ast}}
\newcommand{\Rimp}{\nobreak\hspace{4pt}\penalty0\hspace{-4pt}\ensuremath{\implies}\nobreak\hspace{4pt}}
\newcommand{\Riff}{\nobreak\hspace{4pt}\penalty0\hspace{-4pt}\ensuremath{\iff}\nobreak\hspace{4pt}}
\newcommand{\Thetacal}{\mathit{\Theta}}
\newcommand{\Gamcal}{\mathit{\Gamma}}
\newcommand{\Pws}{\mathbf{2}}
\newcommand{\st}{\,\middle|\,}
\newcommand{\syd}{\triangle}
\providecommand{\notni}{\not\owns}
\begin{document}

\begin{titlepage}
\title{Transitional Conditional Independence}
\author{Patrick Forr\'e\\[+10pt]
    \small{AI4Science Lab}\\[-5pt]
 \small{Korteweg-de Vries Institute for Mathematics}\\[-5pt]
    \small{University of Amsterdam}\\[-5pt]
    \small{\texttt{p.d.forre@uva.nl}}}%
\date{}%
\maketitle
\begin{abstract}
Statistical models contain variables that are not random: parameters, treatments, environments, design points.
Ordinary conditional independence cannot express relations involving such variables. To apply it one must first
put a distribution on them, and that changes the meaning of the statement.

This paper introduces \emph{transitional conditional independence}. It relates three variables on a Markov kernel
$\Kk(W|T)$ with non-stochastic input $T$, and is defined by a single factorization:
\[ \Xk \Indep_{\Kk(W|T)} \Yk \given \Zk \quad :\iff \quad
   \exists\, \Qk(X|Z):\; \Kk(X,Y,Z|T) = \Qk(X|Z) \otimes \Kk(Y,Z|T). \]
The relation asserts a Markov kernel $\Qk(X|Z)$ that is the same for every input $t$. It therefore yields a
factorization rather than an almost-sure identity between conditional expectations, and it needs no distribution
on the input space.

The relation is asymmetric. We show that the asymmetry is essential: symmetrizing it destroys the statements it
was built to make. We prove left and right versions of all separoid rules except Symmetry. Ten of them hold on
arbitrary measurable spaces, the remaining ones under one condition on the spaces involved, and we give criteria
for when Symmetry itself holds. We axiomatize the resulting structure and show that it arises from any symmetric
separoid by a shift.

We give several applications. Ancillarity, sufficiency and adequacy become factorizations that hold pointwise in the
parameter, without a prior and without null sets; the theorems of Fisher--Neyman and of Basu take this form. The
invariance hypothesis of invariant prediction, ``$Y \Indep E \given X_S$'', receives its intended meaning: one
kernel predicts $Y$ from $X_S$ in every environment $E$. And Bayesian networks with non-stochastic input nodes satisfy
a directed global Markov property whose graphical criterion returns a kernel and a factorization, on arbitrary
input spaces.%
\end{abstract}

{
\small
\tableofcontents
}

\thispagestyle{empty}
\end{titlepage}%

\pagestyle{headings}
\pagenumbering{arabic}
\pagestyle{headings}
\setcounter{page}{5}
\setcounter{section}{0}

\section{Introduction}
\label{sec:introduction}

Conditional independence is a basic tool of statistics, probability theory and machine learning, see
\cite{Bis06, Mur12}, and in particular of probabilistic graphical models and causality, see
\cite{Daw93, Lau96, SGS2000, Pearl09, KF09, PJS17, Daw02, Ric02, Ric03, Ali09, Col12, Eva14, Eva16, Eva18, MMC20, Bon21, Gei90}
and many more.
Already in its invention paper, see \cite{Daw79}, a strong motivation for (further) developing conditional independence was the ability to express statistical concepts like sufficiency, adequacy and ancillarity, etc., in terms of conditional independence. 
For example, an ancillary statistic $S(X)$ w.r.t.\ model $\Pk(X|\Theta)$, see \cite{Fis25, Bas64},
is a function $S$ of the data $X$ that has the same probability distribution $\Pk(S(X))$ under any chosen model parameters $\Theta=\theta$. 
The goal is then to formalize this equivalently as a (conditional) independence relation: $S(X) \Indep_{\Pk(X|\Theta)} \Theta$.
This comes with two challenges.
First, in the non-Bayesian setting, the parameters $\Theta$ of the model $\Pk(X|\Theta)$ are not considered random variables, and thus the usual
stochastic (conditional) independence cannot express such concepts in its vanilla form. 
Second, the dependence relation between $X,S(X)$ and $\Theta$ then becomes asymmetric, with deterministic input variables $\Theta$ and stochastic output variables $X$ and $S(X)$; \Cref{fig:cadmg-1} shows the same asymmetry in a graphical model, where the square nodes are the non-stochastic inputs.
These two points similarly hold true for the notions of sufficiency, adequacy, etc.
So any extension of conditional independence that aims at capturing such concepts equivalently 
needs to embrace and incorporate the discussed asymmetry.

A second example, and one where the difficulty is not one of foundations but of daily practice, is
\emph{invariant prediction}, see \cite{PBM16, PBP19}. One observes a response $Y$ and covariates
$X=(X_1,\dots,X_p)$ in several \emph{environments} (different laboratories, different time periods, different
experimental conditions) and looks for subsets $S$ of the covariates that are \emph{invariant}, i.e.\ for which
the conditional distribution of $Y$ given $X_S$ is the same in every environment. The statement one wants to write
down is ``$Y \Indep E \given X_S$'' with the environment $E$ in the middle slot. Ordinary conditional independence
can only read this by treating $E$ as a random variable, which requires a distribution over the environments; and
that changes the meaning, turning a statement about \emph{every} environment into a statement about a mixture of
them. If the environment index is time, or a continuum of interventions, or simply a label chosen by the
experimenter, no such distribution is available at all.
Transitional conditional independence gives the intended statement its exact meaning without inventing one: with
the environment as the non-stochastic input, ``$Y \Indep E \given X_S$'' says precisely that there \emph{exists}
one Markov kernel $\Qk(Y|X_S)$, the same for all environments, that reproduces the conditional distribution of $Y$
given $X_S$ in each of them. We make this precise in \Cref{sec:invariant-prediction}, where we also point out which
steps of the invariant-prediction methodology do and do not follow from the separoid calculus.

\begin{figure}[ht]
    \centering
    \begin{tikzpicture}[scale=1, transform shape]
      \tikzstyle{every node} = [draw,shape=rectangle,color=teal, minimum size = 0.75cm]
       \node (v1) at (3,0) {$v_1$};
       \node (v2) at (8,0) {$v_2$};
        \tikzstyle{every node} = [draw,shape=circle,color=blue]
        \node (v3) at (11,0) {$v_3$};
        \node (v4) at (4.75,-1.5) {$v_4$};
        \node (v5) at (3,-3) {$v_5$};
        \node (v6) at (6.5,-3) {$v_6$};
        \node (v7) at (9.5,-3) {$v_7$};
        \node (v8) at (11.5,-3) {$v_8$};
   \foreach \from/\to in {v1/v4, v2/v6, v2/v7}
    \draw[-{Latex[length=3mm,width=2mm]}, color=teal] (\from) -- (\to); 
   \foreach \from/\to in {v4/v5, v4/v6, v3/v8, v3/v7, v6/v5}
    \draw[-{Latex[length=3mm,width=2mm]}, color=blue] (\from) -- (\to); 
     \end{tikzpicture}
     \caption{Bayesian network with input nodes $v_1,v_2$ and output nodes $v_3,\dots,v_8$.
         The graph allows us to read off the transitional conditional independencies:
    $X_7 \Indep X_1 \given X_2$ and $X_7 \Indep X_1,X_5 \given X_2,X_4,X_6$, etc.
    }
    \label{fig:cadmg-1}
\end{figure}

Over time several extensions of conditional independence have been proposed and studied, see
\cite{Daw79, Daw79b, Daw80, Daw98, Daw01, Gil01, CD17, RERS17, Cho17, FM20, Fri20}, each coming with a different
focus and motivation. The two examples above, together with the graphical models taken up in
\Cref{sec:applications-causal_models}, force three requirements on any such extension, and no existing notion
meets all three at once. \emph{It must be asymmetric}, since ancillarity, sufficiency, adequacy and the invariance
statement all distinguish the two sides of the bar, whereas the notions of \cite{RERS17, Cho17, Fri20} are
symmetric and can express them only after committing to one particular spelling, see the explicit model in
\Cref{eg:symmetrization-loses}. \emph{It must yield a Markov kernel}, not merely an almost-sure identity between
conditional expectations: the kernel is what lets one state ancillarity pointwise in the parameter, assert that
one predictor works in every environment, and read a factorization off a graph, and the $\Qcal$-extended notion of
\cite{FM20} asserts none. And \emph{it must satisfy enough relevance rules} (as many of the separoid rules of
\cite{Daw01} as possible, in left and right versions) for an induction over a graph to go through; for the
extended conditional independence of \cite{CD17} the one rule that such an induction uses at every step, the
mirror weak union, is left open in general. \Cref{sec:main:comparison} opens with a longer list of properties and
\Cref{tab:comparison} records where each notion stands on them. The categorical notions of \cite{Cho17, Fri20} are
best seen as complementary rather than competing, since the disintegration that we have to construct is there
either assumed or bypassed by an abstract factorization, see \Cref{sec:main:comparison:cat}.

It is instructive to see why the two most obvious repairs do not work.
Turning the non-stochastic parameter $\Theta$ into a random variable requires a prior $\Pk(\Theta)$, which the
non-Bayesian setting does not provide; and without it there is no conditional distribution $\Pk(\Theta|S)$ on which
one could impose a condition in the first place.
Alternatively, one may keep the family $\lp\Pk_t(W)\rp_{t \in \Tcal}$ and simply require ordinary conditional
independence for every $t \in \Tcal$ separately. This, however, silently conditions on \emph{all} of $\Tcal$, so that
ancillarity, a statement that compares different parameter values, can no longer be expressed at all.
Our proposed notion of \emph{transitional conditional independence} avoids both traps and meets the three
requirements above in arguably most possible generality. One qualification belongs here rather than in a footnote:
like ordinary conditional independence for general random variables, it satisfies neither Composition nor
Intersection, and the absence of Intersection is precisely why the classical route from a pairwise to a global
Markov property is not available, see \Cref{rem:pairwise-markov}.

\paragraph*{Contributions}
The notion we propose is a single factorization. Fix a Markov kernel $\Pk(W|T)$ from a space $\Tcal$ of
non-stochastic inputs to a space $\Wcal$ of outcomes, and let $X$, $Y$, $Z$ be three variables on it. We define
that $X$ is \emph{transitionally conditionally independent} from $Y$ given $Z$ with respect to $\Pk(W|T)$,
in symbols $X \Indep_{\Pk(W|T)} Y \given Z$, if
\[ \exists\, \Qk(X|Z):\qquad \Pk(X,Y,Z|T) = \Qk(X|Z) \otimes \Pk(Y,Z|T), \]
where $\Pk(X,Y,Z|T)$ is the joint push-forward Markov kernel of $\Pk(W|T)$, where $\Pk(Y,Z|T)$ is its marginal,
and where $\otimes$ denotes the product of Markov kernels. Every object in this display is a Markov kernel. No
joint distribution over $\Tcal$ occurs, and none is invented.

The definition is short, it comes with a factorization by construction, and it can be written down over arbitrary
measurable spaces. Its entire content, and the entire difficulty of this paper, sits in the existential
quantifier. One has to produce a kernel $\Qk(X|Z)$ that is a probability measure in $X$ for every value of $Z$,
that is at the same time measurable in $Z$ simultaneously for all events, and that is one and the same for every
input $t \in \Tcal$. The last of these three demands is the source of the asymmetry: for $\Tcal = \Asterisk$ the
relation is symmetric on standard spaces, see \Cref{cor:sep-ci:symmetry} n). The first two are what make three of the separoid rules below depend on the underlying measurable spaces. The measure theoretic part and the
separoid theoretic part of this paper are therefore not independent of each other.

The asymmetry is the feature a reader is most likely to query, so two things about it are worth knowing early.
First, the relation always ties its second argument to the input: $\Xk \Indep_\Kk \Yk \given \Zk$ is equivalent to
$\Xk \Indep_\Kk \Tk \otimes \Yk \given \Zk$, so it reads as ``$\Xk$ is produced by $\Zk$ alone, while $\Yk$ and
whatever is left of the input are free''. Symmetry can therefore be expected only once the conditioning variable
has exhausted the input. \Cref{rem:reading-the-asymmetry} turns this into a reading rule for the notation, and
\Cref{rem:Q-no-T} adds its companion: moving $\Tk$ from the second into the third argument recovers the weaker,
``for every $t$ separately'' notion, so both readings live inside the same calculus and are told apart by which
slot $\Tk$ occupies. Second, symmetrizing is not a harmless simplification:
\Cref{eg:symmetrization-loses} exhibits a three-coin model in which the disjunctive symmetrization $\Indep^\lor$
holds while the statistical property it is meant to express fails, and in which it assigns different truth values
to two statements that transitional conditional independence itself proves equivalent. That the same asymmetry
reappears in the categorical treatment of \cite{FK23}, introduced there as the categorical counterpart of the
notion defined here, is evidence that it is forced by the problem rather than chosen.

The display above is written in a language that has to be built first. Rather than giving an ad hoc definition of
extended conditional independence, we go back to the roots of measure theoretic probability and develop
``conditional'' versions of its basic objects in \Cref{sec:transitional_probability_theory}, before defining the
relation itself in \Cref{sec:transitional_conditional_independence}: \emph{transition probability spaces}
$\Wcal \times \Tcal$, the codomain and domain of a fixed Markov kernel $\Pk(W|T)$, and \emph{transitional random
variables} $\Xk:\, \Wcal \times \Tcal \dshto \Xcal$ on them, which generalize random variables and already unify
random and non-stochastic variables in a single object.

Furthermore, transitional conditional independence satisfies left and right versions of the \emph{separoid rules}, see \cite{Daw01}, except Symmetry.
Exactly three of them ($\Tk$-Restricted Right Redundancy, Left Weak Union and $\Tk$-Restricted Symmetry, see
b), f) and m) of \Cref{thm:separoid_axioms-tci} and \Cref{cor:sep-ci:symmetry}) require that the codomains
involved form a \emph{disintegration triple}, i.e.\ that a certain conditional Markov kernel exists; that is for
instance the case when the first space is standard and the second one countably generated, and all the other rules
hold on arbitrary measurable spaces. We prove the required disintegration theorem in
\Cref{sec:conditional_markov_kernels}, extending the classical statements for probability measures, see
\cite{Kle20, Rao05}, and we give criteria that dispense with the assumption on the spaces altogether, see
\Cref{lem:trivial-right-argument} and \Cref{sep:tci:l-uni}, as well as criteria for when Symmetry holds.

Transitional conditional independence implies the usual notion of conditional independence for random variables
(the corner case $\Tcal=\Asterisk$) and the two notions of extended conditional independence whose definitions we
reproduce, see \cite{CD17} and \cite{FM20}; for the symmetric proposals of \cite{RERS17, Cho17, Fri20} we prove no
implication and compare formally, through the symmetrization $\Indep^\lor$, see \Cref{sec:main:comparison:sym}. It
also encodes deterministic functional relations exactly, with no null sets involved, see \Cref{thm:variation-ci},
and although it captures only corner cases of variation conditional independence, a formal analogy between the two
lets one combine them by a logical ``and'' without losing any separoid rule, see \Cref{sec:main:comparison:vci}.

Transitional conditional independence can express the statistical concepts of ancillary, sufficient and adequate statistic,
see \cite{Fis22, Fis25}, $S(X)$ for statistical model $\Pk(W|\Theta)$ and transitional random variables $X$ (and $Y$) via:
\begin{enumerate}
    \item Ancillarity: $\displaystyle S(X) \Indep_{\Pk(W|\Theta)} \Theta.$
    \item Sufficiency: $\displaystyle X \Indep_{\Pk(W|\Theta)} \Theta \given S(X).$
    \item Adequacy: $\displaystyle X\Indep_{\Pk(W|\Theta)} \Theta,Y \given S(X).$
\end{enumerate}

As an application of \emph{transitional conditional independence} and its separoid rules we show that 
``conditional''/transitional graphical models like Bayesian networks with (non-stochastic) input variables will follow a 
\emph{directed global Markov property}, i.e.\ they entail transitional
conditional independence relations that are graphically encoded by \emph{id-separation}, the $J$-shift of ordinary
d-separation, see \cite{Pearl09, Lau90, Gei90, Ver93} for the latter. The proof relies on the fact that id-separation
and transitional conditional independence follow the very same asymmetric separoid rules. The same argument applies
to graphs with cycles and latent confounders: $\sigma$-separation satisfies the same asymmetric separoid rules,
obtained from its symmetric ones by the same shift, and the chaining is then verbatim the one given here; what
changes is the model class, from Bayesian networks to (uniquely solvable) structural causal models, together with
the factorization the induction starts from, see \cite{Ric03, FM17, FM18, FM20}. We restrict ourselves to the acyclic
case here for clarity of exposition. Since the criterion hands one factorizations and Markov kernels directly, it
is also the natural starting point for the identification of causal effects, e.g.\ for the rules of $\doit$-calculus,
see \cite{Pearl09, FM20}, a direction we do not pursue here.

The mechanism behind this is worth isolating, and we do so in \Cref{sec:sym-sep-asym-sep}: we introduce
\emph{$\tau$-$\kappa$-separoids}, the asymmetric counterpart of \cite{Daw01}'s separoids, and show that the
$\tau$-shift of any symmetric separoid is one. Id-separation is the $J$-shift of d-separation, so its fourteen
asymmetric rules are a formal consequence of the five classical symmetric ones, and five further rules of the
same kind follow once Composition and Intersection for d-separation are added, with not a single walk to be
inspected again. We also study when the usual notion of a symmetric separoid is recovered.

\paragraph*{Overview}
\Cref{sec:transitional_probability_theory} develops the transitional probability theory that the definition needs,
up to the existence of conditional Markov kernels; \Cref{sec:transitional_conditional_independence} defines the
relation and establishes its separoid rules; \Cref{sec:applications-statistics} carries the statistical
applications; \Cref{sec:applications-causal_models} proves the global Markov property, where what is added over
existing formulations is not the graphical criterion but the \emph{conclusion}, an actual Markov kernel and a
factorization rather than a family of almost-sure identities, and with no assumption whatsoever on the input
spaces, see \Cref{rem:gmp-what-is-gained}; and \Cref{sec:main:comparison} compares the relation to the other
notions in the literature, with \Cref{eg:symmetrization-loses} making precise what a symmetric notion cannot do
and \Cref{tab:comparison} summarizing the outcome. \Cref{sec:discussion} discusses the findings and gives an
outlook. All proofs are in the corresponding appendices, which also contain the detailed statements behind
\Cref{sec:main:comparison}, see \Cref{sec:supp:comparison}.

\section{Transitional Probability Theory}
\label{sec:transitional_probability_theory}

\paragraph*{Notations}
We will use curly letters like $\Wcal$, $\Tcal$, $\Xcal$, $\Ycal$, $\Zcal$ to indicate \emph{measurable spaces}. 
We implicitly assume that they are endowed with a fixed \emph{$\sigma$-algebra}, 
which we will denote by $\Bcal_\Wcal$, $\Bcal_\Tcal$, etc.,
if needed. If we say that $D \ins \Wcal$ is a \emph{measurable subset} we will mean $D \in \Bcal_\Wcal$.
We will, unless stated otherwise, always assume that \emph{topological spaces} like $\R^D$, $[0,1]$, $\bar\R := [-\infty,+\infty]$, etc., are endowed with their 
\emph{Borel $\sigma$-algebra}. 
Similarly, we will assume that \emph{product spaces} like $\Wcal \times \Tcal$ carry the \emph{product $\sigma$-algebra}. 
For the \emph{space of probability measures} $\Pcal(\Wcal)$ on $\Wcal$ we will use the smallest $\sigma$-algebra 
such that all evaluations maps
$j_D:\, \Pcal(\Wcal) \to [0,1]$ given by $j_D(\Pk):=\Pk(D)$ for $D \in \Bcal_\Wcal$ are measurable.
Maps will usually be denoted by capital letters $X$, $Y$, $Z$ in correspondence to their codomains $\Xcal$, $\Ycal$, $\Zcal$, resp. We use bold font letters $\Kk$, $\Pk$, $\Xk$, etc., to indicate probability distributions or Markov kernels. Later we will use $\Gk$ to denote graphs.
If we say that $X:\, \Wcal \to \Xcal$ is a \emph{measurable map} we implicitly assume that $\Wcal$ and $\Xcal$ are 
measurable spaces and that 
$\Bcal_\Wcal \sni X^*\Bcal_\Xcal :=\lC X^{-1}(A)\,|\,A \in \Bcal_\Xcal \rC$.
We call $\Xcal$ (or $\Bcal_\Xcal$) \emph{countably generated} if $\Bcal_\Xcal=\sigma(\Ecal)$ for a countable subset 
$\Ecal \ins \Bcal_\Xcal$, and \emph{countably separated} if there is a countable $\Ecal \ins \Bcal_\Xcal$ that separates the 
points of $\Xcal$.
A measurable space $\Xcal$ is called \emph{standard} (or \emph{standard Borel}) if it is measurably isomorphic to a Borel 
subset of $[0,1]$; equivalently, if $\Bcal_\Xcal$ is the Borel $\sigma$-algebra of a Polish topology on $\Xcal$, e.g.\ 
$\Xcal = \R^D$, $\Z$, $[0,1]^\N$ or any Borel subset thereof.
Every standard measurable space is countably generated and countably separated, and countable products of standard 
measurable spaces are again standard; see \cite{Kec95} 12.B, \cite{Bog07} 6.5.5--6.5.8 and \cite{Fremlin} 424B.
Throughout, $\Asterisk$ denotes the \emph{one-point space} $\lC \ast \rC$, which is standard and discrete and
carries only one $\sigma$-algebra; we write the same symbol $\ast$ for its unique element and identify the two.
So $\Tcal = \Asterisk$ says that there is exactly one input value, and the probability measures on $\Wcal$ are
exactly the Markov kernels $\Asterisk \dshto \Wcal$.

\subsection{Transition Probabilities/Markov Kernels}

Here we will review the notion of transition probabilities, also known as Markov kernels. We mainly introduce our
suggestive notations, which make the later theory more intuitive. In more abstract terms, we give here an explicit description of many constructions that also appear in the Kleisli category of the Giry monad,
see \cite{Law62, Gir82, Kle65}.

\begin{Def}[Markov kernel]
    \label{def-markov-kernel}
    Let $\Tcal$, $\Wcal$ be measurable spaces.
    A \emph{Markov kernel}  or \emph{transition probability} from $\Tcal$ to $\Wcal$
    is - per definition - a map:
    \[   \Kk:\, \Bcal_\Wcal \times \Tcal \to [0,1],\quad (D,t) \mapsto \Kk(D|t),\]
    such that:
    \begin{enumerate}
        \item For each $t \in \Tcal$ the mapping:
            \[ \Bcal_\Wcal \to [0,1],\quad D \mapsto \Kk(D|t)\]
            is a probability measure (i.e.\ normalized and countably additive).
        \item For each $D \in \Bcal_\Wcal$ the mapping:
            \[ \Tcal \to [0,1],\quad t \mapsto \Kk(D|t) \]
            is measurable.
    \end{enumerate}
\end{Def}

\begin{Not}[Markov kernel]
    We will most of the time use the dashed arrow $\dshto$ to $\Wcal$ instead of a usual arrow $\to$ on other spaces 
    to indicate Markov kernels:
    \[\Kk:\, \Tcal \dshto \Wcal,\quad (D,t) \mapsto \Kk(D|t).  \]
    Furthermore, we will often use suggestive notations as follows:
    \[\Kk(W|T):\, \Tcal \dshto \Wcal,\quad (D,t) \mapsto \Kk(W \in D|T=t):=\Kk(D|t).  \]
    We also use the following notations. For fixed $D \in \Bcal_\Wcal$ the map:
    \[ \Kk(W \in D|T):\, \Tcal \to [0,1],\quad t \mapsto \Kk(W \in D|T=t),  \]
    and for fixed $t \in \Tcal$ the map:
    \[ \Kk(W|T=t):\, \Bcal_\Wcal \to [0,1],\quad D \mapsto \Kk(W \in D|T=t).  \]
    We might also use the same notation as above to represent the Markov kernel as a measurable probabilistic map:
    \[ \Kk(W|T):\, \Tcal \to \Pcal(\Wcal),\quad t \mapsto \Kk(W|T=t).\]
    Here $W$ and $T$ are considered suggestive symbols only, but one could give $W$ the meaning 
    of the (identity or) projection map 
    $\pr_\Wcal$ onto $\Wcal$. 
    From the moment a map $T$ into $\Tcal$ is also present, the notation becomes ambiguous: $\Kk(W|T)$ could also mean $\Kk(W|T)$ where we plugged in $T$ for $t$ in ``$T=t$'',
     similar to conditional expectations $\E[W|T]$, but the meaning should become clear from the context.
\end{Not}

\begin{Rem}[Markov kernels generalize probability distributions]\hfill
    \begin{enumerate}
        \item Every probability distribution $\Pk(W) \in \Pcal(\Wcal)$ can be considered as a \emph{constant Markov kernel} from $\Tcal$ to $\Wcal$ via:
            \[  \Kk(W|T):\, \Tcal \dshto \Wcal,\quad (D,t) \mapsto \Kk(W \in D|T=t):=\Pk(W \in D).   \]
        \item Every Markov kernel from the one-point space $\Tcal=\Asterisk$ to $\Wcal$:
            \[ \Kk(W|T):\, \Asterisk \dshto \Wcal,\quad (D,\ast) \mapsto \Kk(W \in D|T=\ast),\]
            defines a unique probability distribution $\Pk(W) \in \Pcal(\Wcal)$ given via:
            \[ \Pk(W \in D):=\Kk(W \in D|T=\ast).\]
            So we can identify probability distributions on $\Wcal$ with Markov kernels $\Asterisk \dshto \Wcal$.
    \end{enumerate}
\end{Rem}

\begin{Rem}[Markov kernels generalize deterministic maps]
    Consider a measurable map $\tilde X: \Tcal \to \Xcal$. Then we can turn $\tilde X$ into a Markov kernel 
    $\deltabf_{\tilde X}(X|T)$ via:
    \[  \deltabf_{\tilde X}(X|T):\, \Tcal \dshto \Xcal, \quad (A,t) \mapsto \deltabf_{\tilde X}(X \in A|T=t):= \I_A(\tilde X(t)),\]
    which puts $100\%$ of the probability mass onto the point $\tilde X(t)$.
    We will often also use the notation without the dummy variable $X$: $\deltabf(\tilde X|T):=\deltabf_{\tilde X}(X|T)$.
\end{Rem}

\subsection{Constructing Transition Probabilities from Others}
\label{sec:other-markov-kernels}

In the following we will demonstrate how one can construct new Markov kernels from others. 
The constructions include marginalization, product, composition, push-forward.
Later an own subsection is dedicated to conditioning.
Note that the measurability of those constructions is either clear or can be proven using Dynkin's $\pi$-$\lambda$ theorem, 
see \cite{Kle20} Thm.\ 1.19, also see \cite{Bog07} Thm.\ 1.9.3.

\begin{Def}[Marginalizing Markov kernels]
    \label{def-marginal-markov-kernels}
    Let 
    \[\Kk(X,Y|T):\,\Tcal \dshto  \Xcal \times \Ycal \]
    be a Markov kernel in two variables.
    We can then define the \emph{marginal Markov kernels} as follows:
    \[ \Kk(X|T):\, \Tcal \dshto \Xcal,\quad (A,t) \mapsto \Kk(X \in A, Y \in \Ycal|T=t),  \]
    and:
    \[ \Kk(Y|T):\, \Tcal \dshto \Ycal,\quad (B,t) \mapsto \Kk(X \in \Xcal, Y \in B|T=t).  \]
\end{Def}

\begin{Def}[Product of Markov kernels, general form]
    \label{def-product-markov-kernels-general}
    Let $I$ be a finite set of \emph{variables}, each variable $v \in I$ coming with a measurable space $\Xcal_v$,
    and write $\Xcal_S := \prod_{v \in S} \Xcal_v$ and $X_S := (X_v)_{v \in S}$ for $S \ins I$.
    Consider two Markov kernels:
    \[ \Qk(X_{O_1}|X_{I_1}):\, \Xcal_{I_1} \dshto \Xcal_{O_1}, \qquad\qquad
       \Kk(X_{O_2}|X_{I_2}):\, \Xcal_{I_2} \dshto \Xcal_{O_2}, \]
    with disjoint output variables, i.e.\ $O_1 \cap O_2 = \emptyset$, and such that no variable is an input and an
    output of the same or of the other kernel, i.e.\ $I_1 \cap O_1 = I_2 \cap O_2 = I_2 \cap O_1 = \emptyset$. Then their \emph{product} is the Markov kernel:
    \[ \Qk(X_{O_1}|X_{I_1}) \otimes \Kk(X_{O_2}|X_{I_2}):\,
       \Xcal_{(I_1 \cup I_2) \sm O_2} \dshto \Xcal_{O_1 \cup O_2}, \]
    given for measurable $E \ins \Xcal_{O_1 \cup O_2}$ and $x \in \Xcal_{(I_1 \cup I_2)\sm O_2}$ by:
    \[ \int \!\! \int \I_E\lp x_{O_1},x_{O_2} \rp \,
        \Qk\lp X_{O_1} \in dx_{O_1} \,\big|\, X_{I_1} = x_{I_1} \rp \,
        \Kk\lp X_{O_2} \in dx_{O_2} \,\big|\, X_{I_2} = x_{I_2} \rp, \]
    where the outer integration is over $x_{O_2} \in \Xcal_{O_2}$ and the inner one over
    $x_{O_1} \in \Xcal_{O_1}$, and where the arguments $x_{I_1 \cap O_2}$ of $\Qk$ are the ones supplied by the
    outer integration.
    In words: one integrates the left kernel against the right one over exactly those arguments of the left kernel
    that are \emph{outputs} of the right one; all remaining arguments of both kernels become inputs of the product
    and all outputs of both become outputs.
\end{Def}

\begin{Rem}[Special cases of the product]
    \label{rem:product-wirings}
    All the wirings that occur in this paper are instances of \Cref{def-product-markov-kernels-general}:
    \begin{enumerate}
        \item \emph{Shared inputs}, e.g.\ $\Xk(X|W,T) \otimes \Yk(Y|W,T)$ with
            $O_1=\{X\}$, $O_2=\{Y\}$, $I_1=I_2=\{W,T\}$. Here $I_1 \cap O_2 = \emptyset$, so nothing is
            integrated out and the product is $(w,t) \mapsto \Xk(X|W=w,T=t) \otimes \Yk(Y|W=w,T=t)$, the
            \emph{product measure}. In other words, in such an expression $X$ and $Y$ are conditionally independent
            given $(w,t)$; see also \Cref{rem:independent-copies}.
        \item \emph{Left input $=$ right output}, e.g.\ $\Qk(X|Z) \otimes \Kk(Y,Z|T)$ with $O_1=\{X\}$,
            $I_1=\{Z\}$, $O_2=\{Y,Z\}$, $I_2=\{T\}$, which is a Markov kernel $\Tcal \dshto \Xcal \times \Ycal
            \times \Zcal$. This is the wiring of the definition of transitional conditional independence.
        \item \emph{No shared variables}, e.g.\ $\Kk(U) \otimes \Kk(W|Z)$, a Markov kernel
            $\Zcal \dshto \Ucal \times \Wcal$.
        \item \emph{Adjoining an input to the output}: $\deltabf(T|T) \otimes \Kk(W|T)$ denotes the Markov kernel
            $\Tcal \dshto \Tcal \times \Wcal$, $t \mapsto \deltabf_t \otimes \Kk(W|T=t)$, obtained from the
            Dirac kernel of the identity of $\Tcal$. (Strictly speaking this leaves the setting above, where a
            variable cannot be an input and an output at the same time; we allow it as a convenient abbreviation.)
    \end{enumerate}
    Associativity and the commutation rule of \Cref{rem:markov-kernels-products} hold in this generality, with
    the same proofs.
\end{Rem}

\begin{Def}[Product of Markov kernels, the basic wiring]
    \label{def-product-markov-kernels}
    Consider two Markov kernels:
    \[ \Qk(Z|Y,W,T):\, \Ycal \times \Wcal\times \Tcal \dshto \Zcal,\qquad  
    \Kk(W,U|T,X):\, \Tcal \times \Xcal \dshto \Wcal \times \Ucal. \]
    Then we define the \emph{product Markov kernel}:
    \[ \Qk(Z|Y,W,T) \otimes \Kk(W,U|T,X) :\, \Ycal \times  \Tcal \times \Xcal \dshto \Zcal \times \Wcal \times \Ucal, \]
    using  measurable sets $E \ins \Zcal \times \Wcal \times \Ucal$ via: $(E,(y,t,x)) \mapsto $
    \[ \int\int \I_E(z,w,u)\, \Qk(Z \in dz|Y=y,W=w,T=t) \, \Kk((W,U) \in d(w,u)|T=t,X=x),\]
where the inner integration is over $z \in \Zcal$ and the outer integration over $(w,u) \in \Wcal\times \Ucal$.
\end{Def}

\begin{Def}[Composition of Markov kernels]
        \label{def-composition-markov-kernels}
    Consider two Markov kernels:
    \[ \Qk(Z|Y,W,T):\, \Ycal \times \Wcal\times \Tcal \dshto \Zcal,\qquad  \Kk(W,U|T,X):\, 
    \Tcal \times \Xcal \dshto \Wcal\times\Ucal. \]
    Then we define their \emph{composition}:
    \[ \Qk(Z|Y,W,T) \circ \Kk(W,U|T,X) :\, \Ycal \times  \Tcal \times \Xcal \dshto \Zcal, \]
    using  measurable sets $C \ins \Zcal$ via: $(C,(y,t,x))    \mapsto  $
    \[
     \int \Qk(Z \in C|Y=y,W=w,T=t) \, \Kk(W \in dw|T=t,X=x).
     \]
    Note that we implicitly marginalized $U$ out, i.e.\ in the composition we integrate over all variables (here: $W$ and $U$) 
    from the right hand Markov kernel.
    This is on purpose: $\circ$ is the composition of the Kleisli category of the Giry monad, the arrows of which
    are Markov kernels, and it is the operation for which $\deltabf$ of a measurable map is functorial, i.e.\
    $\deltabf_\varphi \circ \deltabf_\psi = \deltabf_{\varphi \circ \psi}$. It is therefore \emph{not} the
    composition of multi-output maps: a right hand kernel with two outputs $W$ and $U$ loses $U$, exactly as
    $\Kk(Z|W) \circ \Kk(W,U|X)$ has no place to put a $U$. Whenever the further outputs are to be kept one uses the
    product $\otimes$ instead, which retains all outputs; by \Cref{rem:markov-kernels-products} 1.\ the
    composition is precisely the corresponding marginal of the product, so nothing is lost by having both
    operations.
\end{Def}

\begin{Rem}
    \label{rem:markov-kernels-products}
    \begin{enumerate}
        \item
It is clear from the \Cref{def-composition-markov-kernels}, \ref{def-product-markov-kernels} and \ref{def-marginal-markov-kernels} that the composition:
    \[\Qk(Z|Y,W,T) \circ \Kk(W,U|T,X)\]
    is the $Z$-marginal of the product:
    \[ \Qk(Z|Y,W,T) \otimes \Kk(W,U|T,X).\]
    \item Both, products and compositions, are each associative, but clearly not commutative in general.
    \item If the left Markov kernel $\Qk(Z|Y,T)$ has no dependence in the second arguments w.r.t.\ to a first argument of the right Markov kernel $\Kk(W,U|T,X)$, i.e.\ no $W$ in the above terms, then they commute by Fubini's theorem:
        \[ \Qk(Z|Y,T) \otimes \Kk(W,U|T,X)=  \Kk(W,U|T,X) \otimes \Qk(Z|Y,T). \]
\end{enumerate}
\end{Rem}

\begin{Rem}[Composition of deterministic Markov kernels]
    Consider measurable maps:
    \[ \tilde X:\, \Tcal \to \Xcal,\qquad \tilde Z:\, \Xcal \to \Zcal, \]
    and their composition $\tilde Z \circ \tilde X$. 
    Then the composition of the corresponding Markov kernels satisfies:
    \[  \deltabf_{\tilde Z\circ \tilde X}(Z|T)  = \deltabf_{\tilde Z}(Z|X) \circ \deltabf_{\tilde X}(X|T), \]
    where $\deltabf_{\tilde Z}(Z \in C|X=x) := \I_C(\tilde Z(x))$ and $\deltabf_{\tilde X}(X \in A|T=t) := \I_A(\tilde X(t))$.\\
    So the composition of Markov kernels extends the composition of maps.
\end{Rem}

\begin{Def}[Push-forward Markov kernel w.r.t.\ measurable maps]
    Consider a Markov kernel $\Kk(W|T):\, \Tcal \dshto \Wcal$ and a measurable map:
    $ X:\, \Wcal \times \Tcal \to \Xcal $.
    Then we define the \emph{push-forward Markov kernel}: 
    \[ X_*\Kk(W|T) := \Kk(X(W,T)|T):= \Kk(X|T):\, \Tcal \dshto \Xcal,\] 
    of $\Kk(W|T)$ w.r.t.\ $X$ 
    via:  $(A,t) \mapsto$
    \[  \Kk(X \in A|T=t) := \Kk(W \in X^{-1}_t(A)|T=t), \]
    where:
    \[ X^{-1}_t(A) = X^{-1}(A)_t:= \{ w\in \Wcal\,|\,X(w,t) \in A \}.\]
\end{Def}

\begin{Rem}
    We can also write push-forwards as compositions:
    \[\Kk(X|T) = \deltabf(X|W,T) \circ \Kk(W|T),\]
    where:
    $\deltabf(X \in A|W=w,T=t) := \I_A(X(w,t)) = \I_{X^{-1}(A)}(w,t).$
    In this sense compositions of Markov kernels generalize push-forward Markov kernels.
\end{Rem}

\begin{Def}[Push-forward Markov kernel w.r.t.\ another Markov kernel]
    Consider  Markov kernels $\Kk(W|T):\, \Tcal \dshto \Wcal$ and 
    $ \Xk(X|W,T):\, \Wcal \times \Tcal \dshto \Xcal $.
    Then we define the \emph{push-forward Markov kernel} as the composition: 
    \[ \Kk(X|T):= \Xk(X|W,T) \circ \Kk(W|T):  \, \Tcal \dshto \Xcal.\] 
\end{Def}

\begin{Rem}\label{rem:extending_markov_kernel}
    Any Markov kernel
    $ \Kk(W|T):\, \Tcal \dshto \Wcal$
    can always be extended to include the canonical projection map $T=\pr_\Tcal:\, \Wcal \times \Tcal \to \Tcal$ via: 
    \[\Kk(W,T|T):\, \Tcal \dshto \Wcal \times \Tcal,\quad (E,t) \mapsto \]
    \[ \Kk( (W,T)  \in E| T=t) = \Kk( W \in E_t|T=t), \]
    where $E_t = \{ w \in \Wcal\,|\,(w,t) \in E\}$.
    Using Definition~\ref{def-product-markov-kernels}, we can also write this as:
    \[ \Kk(W,T|T) = \Kk(W|T) \otimes \deltabf(T|T) = \deltabf(T|T) \otimes \Kk(W|T),\]
    where $\deltabf(T \in D|T=t) := \I_D(t)$ for measurable $D \ins \Tcal$ and $t \in \Tcal$.
\end{Rem}

\subsection{Null Sets w.r.t.\ Transition Probabilities}

\begin{Def}[Null sets w.r.t.\ transition probabilities]
Let $\Kk(W|T) :\, \Tcal \dshto \Wcal$  be a transition probability.
A subset $M \ins \Wcal \times \Tcal$ will be called a \emph{$\Kk(W|T)$-null set} 
if every section/fibre $M_t:=\{ w \in \Wcal\,|\, (w,t) \in M\}$ is a $\Kk(W|T=t)$-null set, i.e.\ 
there exist measurable $N_t \ins \Wcal$ with $M_t \ins N_t$ and $\Kk(W \in N_t|T=t)=0$, for every $t \in \Tcal$.\\
\end{Def}

We are usually interested in measurable null sets. 
The notion of null sets w.r.t.\ transition probabilities generalizes the notion of null sets in probability spaces, which can be recovered by taking $\Tcal=\Asterisk$.

\subsection{Transition Probability Spaces}

We will now give the definition of a transition probability space, which will generalize the notion of probability spaces.

\begin{Def}[Transition probability space]
    \label{def:transition_probability_space-main}
    Consider measurable spaces $\Tcal$ and $\Wcal$ and a Markov kernel/transition probability:
    \[\Kk(W|T):\, \Tcal \dshto \Wcal,\quad (D,t) \mapsto \Kk(W \in D|T=t).  \]
    We then call the tuple: 
    $(\Wcal \times \Tcal, \Kk(W|T))$ a \emph{transition probability space}.
    It naturally comes with the canonical projection map: 
    \[T:\,\Wcal \times \Tcal \to \Tcal,\quad T(w,t):=t,\]
    and the Markov kernel: $\Kk(W,T|T) := \Kk(W|T) \otimes \deltabf(T|T)$, which then satisfies $\Kk(T|T) = \deltabf(T|T)$.
    Throughout we assume $\Tcal \neq \emptyset$, so that at least one transition probability is actually given;
    the case $\Tcal = \Asterisk$ recovers ordinary probability spaces.
\end{Def}

As for null sets above, the notion of transition probability space generalizes the notion of probability spaces, recovered by $\Tcal=\Asterisk$.

\subsection{Transitional Random Variables}

In this subsection we will introduce the notion of transitional random variables, which will generalize the usual notion
of random variables, formalizes what one could call ``conditional'' random variables. Furthermore, we start from a bit more
general point of view as we not only allow for (deterministic) measurable maps, but also for stochastic maps, which again will be formalized as Markov kernels.

\begin{Rem}[On the name]
    We call these objects \emph{transitional} rather than \emph{conditional} random variables. The reason is that the
    dependence on $T$ is not a conditioning on an event or a $\sigma$-algebra ($T$ carries no distribution at all)
but an \emph{input} dependence, exactly as for a transition probability. The name also matches the ambient
    \emph{transition probability space} and keeps ``conditional'' free for the conditional Markov kernels of
    \Cref{sec:conditional_markov_kernels}, where genuine conditioning does take place.
\end{Rem}

\begin{Def}[Transitional random variables]
    \label{def:transitional_random_variables}
    If  $(\Wcal\times \Tcal, \Kk(W|T))$ is a transition probability space then a \emph{transitional random variable}
    is a Markov kernel:
    \[ \Xk =\Xk(X|W,T):\, \Wcal \times \Tcal \dshto \Xcal  \]
    to any other measurable space $\Xcal$. It will come with its push-forward Markov kernel:
    \[ \Kk(X|T) := \Xk(X|W,T) \circ \Kk(W|T). \]
\end{Def}

\begin{Rem}[How to picture them, and one warning]
    \label{rem:trv-picture}
    Transitional random variables can be thought of as ``conditional'' random variables, or as stochastic processes
    or random fields $(X_t)_{t \in \Tcal}$ of which we care not about a joint law $\Pk((X_t)_{t \in \Tcal})$ but
    about how the transition probability $\Pk(X_t|T=t)$ depends on the ``parameter'' $t$. They already unify random
    variables and deterministic, non-stochastic variables in one object.
    The extra generality of \emph{probabilistic} maps is not free, and we flag its price once: two occurrences of
    the same genuinely stochastic $\Xk$ in one expression are \emph{independent copies} given $(w,t)$, so that
    $\Xk$ need not even be a map of itself, see \Cref{rem:ismapof-properties} item \ref{rem:no-reflexivity}. All
    separoid rules hold regardless; only the order-theoretic statements are restricted to the deterministic
    transitional random variables.
\end{Rem}

\begin{Rem}
\begin{enumerate}
\item 
    If  $(\Wcal\times \Tcal, \Kk(W|T))$ is a transition probability space then any measurable map
    $X:\, \Wcal \times \Tcal \to \Xcal$ induces a transitional random variable $\deltabf(X|W,T)$ given by:
    \[ \deltabf(X \in A|W=w,T=t) := \I_A(X(w,t)).\]
    By slight abuse of notation we will call $X$ itself a transitional random variable as well, by actually referring to 
    $\deltabf(X|W,T)$. Transitional random variables of this form will be of the main focus in the following.
\item
 A transitional random variable $X:\, \Wcal \times \Tcal \to \Xcal$ can be considered as a family of random variables 
 measurably parameterized by $t \in \Tcal$.
    For $t \in \Tcal$ we have the measurable maps:
    \[ X_t:\, \Wcal \to \Xcal,\quad w \mapsto X_t(w):= X(w,t),\]
    each of which can be considered a random variable on the probability space $(\Wcal,\Kk(W|T=t))$.
    Note that in this setting we are not modelling the joint distribution of $(X_t)_{t \in \Tcal}$, but rather 
    how the individual distribution of $X_t$ depends on and varies with $t\in \Tcal$.
\item Note that by going from transition probability space $\lp \Wcal \times \Tcal, \Kk(W|T) \rp$ to the one 
     $\lp \Xcal \times \Tcal, \Kk(X|T) \rp$ for a transitional random variable $\Xk(X|W,T)$ the projection map:
    \[ X: \, \Xcal \times \Tcal \to \Xcal,\quad (x,t) \mapsto x, \]
    can be considered a transitional random variable of the form $\deltabf(X|X,T)$. 
    So with only slight loss of generality one can replace a general transitional random variable $\Xk$ 
    by one of the form $\deltabf(X|W,T)$. More will be said in \Cref{thm:reparamererizing-trv}.
\item 
The notion of transitional random variables generalizes the notion of random variables and formalizes 
what one could call a \emph{(probabilistic) ``conditional'' random variable}. 
Note that the Markov kernel can be given without any conditioning operation. 
\item
Transitional random variables can model \emph{probabilistic programs}. For each user chosen input $T=t$ 
a random input $w \sim \Kk(W|T=t)$ 
is drawn. Then the input $(w,t)$ is presented to the probabilistic program $\Xk$ and an output is sampled $x_1 \sim \Xk(X|W=w,T=t)$.
Using Markov kernels to represent transitional random variables allows for random noise inside the program that generates the output 
$x_1$. So even when presented with the same input $(w,t)$ again another output $x_2 \neq x_1$ might be drawn. So $\Xk(X|W=w,T=t)$ models
the output distribution for fixed input $(w,t)$. Certainly, if one has no insight into the input sampling procedure $\Kk(W|T)$ one might only be interested in the push-forward: $\Kk(X|T)$, which directly describes the output distribution for each user chosen input $T=t$.
\item
If we want to model a deterministic variable with no stochasticity we could consider transitional random variables of the form 
$\deltabf_\varphi(X|T)$ that do not depend on the $W$-argument.
\end{enumerate}
\end{Rem}

\begin{Eg}[Special transitional random variables of importance]
        Let $(\Wcal\times \Tcal, \Kk(W|T))$ be a transitional probability space.
        Then we denote by:
        \begin{enumerate}
            \item $T$ the \emph{canonical projection} onto $\Tcal$:
                \[T:= \pr_\Tcal:\,\Wcal \times \Tcal \to \Tcal,\quad (w,t) \mapsto T(w,t):= t.\]
                We also put:
                \[ \Tk = \Tk(T|W,T)=\deltabf(T|W,T) : \, \Wcal \times \Tcal \dshto \Tcal, \quad (D,(w,t)) \mapsto \I_D(t).  \]
            \item $\ast$ the \emph{constant transitional random variable}:
                \[ \ast :\,\Wcal \times \Tcal \to \Asterisk,\quad (w,t) \mapsto \ast,\]
                We also use the same symbol $\ast$ to denote the Markov kernel:
                \[ \deltabf_\ast:\, \Wcal \times \Tcal \dshto \Asterisk,\quad (D,\ast) \mapsto \I_D(\ast).\]
        \end{enumerate}
\end{Eg}

\subsection{Ordering the Class of Transitional Random Variables}

We now introduce several comparison relations between transitional random variables. 
All them model to some degree that one variable $\Xk$ is a (deterministic or stochastic) 
measurable function of another one $\Yk$
(up to some form of null set).

\begin{Not} 
    \label{not:ismapof}
    Let $(\Wcal\times \Tcal, \Kk(W|T))$ be a transitional probability space and $\Xk$, $\Yk$, $\Zk$
    be transitional random variables with joint Markov kernel:
    \[ \Kk(X,Y,Z|T) := \lp  \Xk(X|W,T) \otimes \Yk(Y|W,T)  \otimes \Zk(Z|W,T)\rp \circ \Kk(W|T).  \]
    We put:
    \begin{enumerate}
        \item $X \ismapof Y$, for transitional random variables of the form $\Xk=\deltabf(X|W,T)$, $\Yk=\deltabf(Y|W,T)$,
            if there exists a measurable map $\varphi:\, \Ycal \to \Xcal$ such that $X = \varphi \circ Y$.
        \item $\Xk \ismapof_\Kk \Yk$ if there exists a measurable map $\varphi:\, \Ycal \to \Xcal$ such that: 
            \[ \Kk(X,Y|T) = \deltabf_\varphi(X|Y) \otimes \Kk(Y|T),  \]
            where $\Kk(Y|T)$ is the marginal of $\Kk(X,Y|T)$.
        \item $\Xk \ismapof_\Kk^\ast \Yk$ if there exists a Markov kernel  $\Qk(X|Y):\, \Ycal \dshto \Xcal$ such that: 
            \[ \Kk(X,Y|T) = \Qk(X|Y) \otimes \Kk(Y|T).  \]
    \end{enumerate}
\end{Not}

\begin{Rem}
    \label{rem:ismapof-properties}
    \begin{enumerate}
        \item 
        We have the implications:
    \[ X \ismapof Y \implies \Xk \ismapof_\Kk \Yk \implies \Xk \ismapof_\Kk^\ast \Yk.\]
        \item The relation $\ismapof_\Kk$ will be the most crucial one in the following.
        \item \label{rem:no-reflexivity} Note that for general $\Xk$ we do not even have reflexivity:
         $\Xk \ismapof_\Kk \Xk$.
         Indeed, take $\Tcal = \Asterisk$ and $\Xk(X|W) := \Ncal(0,1)$ for every $w \in \Wcal$. Writing
         $\Xk \otimes \Xk$ for the joint kernel of two occurrences of $\Xk$ we get
         $\Kk(X_1,X_2) = \Ncal(0,1) \otimes \Ncal(0,1)$, which is \emph{not} of the form
         $\deltabf_\varphi(X_1|X_2) \otimes \Kk(X_2)$ for any measurable $\varphi$, because the latter is carried by
         the graph of $\varphi$, a Lebesgue null subset of $\R^2$.
         The general phenomenon behind this is worth stating: whenever a transitional random variable occurs twice in
         an expression, the two occurrences are \emph{conditionally independent copies} given $(w,t)$, by the very
         definition of the joint Markov kernel above. So $\Xk \Indep_\Kk \Xk \given \deltabf_\ast$ does
         \emph{not} say what a reader trained on random variables will expect. On \emph{deterministic}
         transitional random variables, i.e.\ those of the form $\deltabf(X|W,T)$, reflexivity does hold, see
         \Cref{lem:restricted-reflexivity}.
     \item We also do not have anti-symmetry, i.e.\ we can not conclude from: $\Xk \ismapof_\Kk \Yk \ismapof_\Kk \Xk$ 
         that then $\Xk = \Yk$ holds, since such variables might differ on some null-set.
     \item We can fix the anti-symmetry by going over to almost-sure anti-symmetry, i.e.\ by defining: 
         \[ \Xk \approx_\Kk \Yk :\iff \Xk \ismapof_\Kk \Yk \ismapof_\Kk \Xk.\]
    \end{enumerate}
\end{Rem}

The relation $\ismapof_\Kk$ satisfies the following rules, of which ``product extension'' is the most important
one: it is exactly the rule that $\ismapof_\Kk^\ast$ lacks, which is why we work with $\ismapof_\Kk$
throughout and introduced $\ismapof_\Kk^\ast$ only for the comparison in \Cref{sec:main:comparison}.

\begin{Thm}
    Let $(\Wcal\times \Tcal, \Kk(W|T))$ be a transitional probability space and  $\Xk$, $\Yk$, $\Zk$, $\Uk$
    be transitional random variables.
    The relation $\ismapof_\Kk$ satisfies the following rules:
    \begin{enumerate}
        \item Almost-sure anti-symmetry: $\Xk \ismapof_\Kk \Yk \ismapof_\Kk \Xk \implies: \Xk \approx_\Kk \Yk $.
        \item Transitivity: $\Xk \ismapof_\Kk \Yk \ismapof_\Kk \Zk \implies \Xk \ismapof_\Kk \Zk.$
        \item Bottom element: $\deltabf_\ast \ismapof_\Kk \Xk$.
        \item Product stays bounded: $(\Xk \ismapof_\Kk \Zk)  \land  (\Yk \ismapof_\Kk \Zk)  \implies \Xk \otimes \Yk \ismapof_\Kk \Zk.$
        \item Product extension: $\Xk \ismapof_\Kk \Yk \implies \Xk \ismapof_\Kk \Yk \otimes \Zk.$
        \item Product compatibility: $ (\Xk \ismapof_\Kk \Zk) \land (\Yk \ismapof_\Kk \Uk) \implies  \Xk \otimes \Yk \ismapof_\Kk \Zk \otimes \Uk. $
    \end{enumerate}
    Furthermore, the relation  $\ismapof_\Kk$ turns the sub-class of all transitional random variables on
    $(\Wcal\times \Tcal, \Kk(W|T))$ of the form
    $\Xk=\deltabf(X|W,T)$, where $X:\, \Wcal \times \Tcal \to \Xcal$ is a measurable map, and where $\Xcal$ may vary,          
           into a \emph{bounded join-semi-lattice} up to almost-sure anti-symmetry $\approx_\Kk$ with join $\otimes$ and bottom element $\deltabf_\ast$. 
\end{Thm}

Note that the mentioned sub-class might not be a set, but all the properties of a bounded join-semi-lattice can be proven.
The proofs will be given in \Cref{sec:join-semi-lattice}: rule 1.\ holds by the definition of $\approx_\Kk$,
rule 2.\ is \Cref{lem:transitivity}, rule 3.\ is \Cref{lem:bottom-element}, rule 4.\ is \Cref{lem:prod-bounded},
rule 5.\ is \Cref{lem:join-upper-bound-2} and rule 6.\ is \Cref{lem:products-bound}; all six hold for arbitrary
transitional random variables. \Cref{thm:join-semi-lattice} collects them for \emph{deterministic} transitional
random variables, where in addition reflexivity, see \Cref{lem:restricted-reflexivity}, and idempotency, see
\Cref{lem:idempotency}, hold; the join-semi-lattice statement is \Cref{thm:join-semi-lattice} together with
\Cref{cor:join-semi-lattice}.

\subsection{Disintegration of Transition Probabilities}
\label{sec:conditional_markov_kernels}

In this subsection we present results about the existence of conditional Markov kernels. Since we will factorize a joint Markov kernel 
into a marginal part and a conditional part such procedures are also called disintegration.
First, we will talk about the essential uniqueness of such factorizations and then existence. 
For proofs see \Cref{sec:supp-disintegration}.
For probability measures the disintegration theorem over a standard first factor is classical, see \cite{Kal17} Thm.\ 1.25, \cite{Kle20} Ch.\ 8.3 and \cite{Rao05}; a version for analytic measurable spaces is given in \cite{Bog20}. What is proven here is the corresponding statement for \emph{transition} probabilities, i.e.\ with an additional parameter, together with the joint measurability in the conditioning and the parameter variable; we give an elementary and self-contained proof in \Cref{sec:supp-disintegration}.

\begin{Def}[Conditional Markov kernels]
     Consider a Markov kernel
    \[ \Kk(X,Y|Z):\, \Zcal \dshto \Xcal \times \Ycal\]
    and its marginal $\Kk(Y|Z)$. A \emph{conditional Markov kernel of $\Kk(X,Y|Z)$ conditioned on $Y$ given $Z$}
    is a Markov kernel:
    \[ \Kk(X|Y,Z):\, \Ycal \times \Zcal \dshto \Xcal\]
    such that:
    \[   \Kk(X,Y|Z) = \Kk(X|Y,Z) \otimes \Kk(Y|Z).\]
\end{Def}

\begin{Def}[Disintegration triple]
    \label{def:disintegration-triple}
    A triple $(\Xcal,\Ycal,\Zcal)$ of measurable spaces is called a \emph{disintegration triple} if
    \emph{every} Markov kernel:
    \[ \Kk(X,Y|Z):\, \Zcal \dshto \Xcal \times \Ycal\]
    admits a conditional Markov kernel:
    \[ \Kk(X|Y,Z):\, \Ycal \times \Zcal \dshto \Xcal,\]
    i.e.\ such that:
    \[   \Kk(X,Y|Z) = \Kk(X|Y,Z) \otimes \Kk(Y|Z).\]
\end{Def}

Being a disintegration triple is exactly the property that is needed to invoke the disintegration
theorem, and it is the only thing that the separoid rules of \Cref{sec:separoid-rules-tci} will
require of the underlying measurable spaces. It will therefore be convenient to phrase all such
assumptions in these terms rather than in terms of concrete conditions like ``standard'' or
``countably generated''.

\begin{Lem}[Essential uniqueness of conditional Markov kernels]
    \label{lem:ess-unique}
    If we have Markov kernels:    
    \[\Pk(X|Y,Z),\;\Qk(X|Y,Z):\, \Ycal\times \Zcal \dshto \Xcal,\qquad \text{ and }\qquad \Kk(Y|Z):\,\Zcal \dshto \Ycal,\]
    between any measurable spaces $\Xcal$, $\Ycal$, $\Zcal$ such that:
    \[ \Pk(X|Y,Z) \otimes \Kk(Y|Z) = \Qk(X|Y,Z) \otimes \Kk(Y|Z), \]
    then for every $A \in \Bcal_\Xcal$ the set:
    \[N_A := \lC (y,z) \in \Ycal \times \Zcal \,|\, \Pk(X \in A|Y=y,Z=z) \neq \Qk(X \in A|Y=y,Z=z)\rC\]
     is a measurable $\Kk(Y|Z)$-null set. \\
     If, furthermore, $\Bcal_\Xcal$ is countably generated
     then also $N := \bigcup_{A \in \Bcal_\Xcal} N_A$ is a measurable $\Kk(Y|Z)$-null set.
 \end{Lem}

\begin{Thm}[Existence of conditional Markov kernels]
    \label{thm-regular-conditional-Markov-kernel}
    Let
    $ \Kk(X,Y|Z):\, \Zcal \dshto \Xcal \times \Ycal$
    be a Markov kernel. 
    Furthermore, assume one of the following:
    \begin{enumerate}
        \item $\Xcal$ standard, $\Ycal$ countably generated, $\Zcal$ arbitrary;
        \item $\Xcal$ standard, $\Ycal$ arbitrary, $\Zcal$ discrete\footnote{\label{fn:discrete}A measurable
            space $\Ucal$ is called \emph{discrete} if it is countable with $\Bcal_\Ucal = \Pws^\Ucal$.};
        \item $\Xcal$ arbitrary, $\Ycal$ discrete, $\Zcal$ arbitrary;
        \item $\Kk(X,Y|Z=z) \ll \mubf \otimes \nubf$ for all $z \in \Zcal$, for $\sigma$-finite measures
            $\mubf$ on $\Xcal$ and $\nubf$ on $\Ycal$, with a density $k(x,y|z) \in [0,\infty)$ that is jointly
            measurable in $(x,y,z)$;\footnote{\label{fn:density-weaker}The $\sigma$-finiteness of $\nubf$ is not
            needed, see \Cref{rem:density-hypothesis} in \Cref{sec:supp-disintegration}.}
        \item $X \ismapof_\Kk (Y,Z)$;\footnote{\label{fn:thm-ismapof}In points 5.\ and 6.\ the relation
            $\ismapof_\Kk$ is read as in \Cref{not:ismapof} for the transition probability space
            $\lp \lp \Xcal \times \Ycal \rp \times \Zcal, \Kk(X,Y|Z) \rp$, i.e.\ with $\Zcal$ in the role of the
            input space; see \Cref{sec:supp-disintegration} for the two conditions written out.}
        \item $Y \ismapof_\Kk Z$.
    \end{enumerate}
    Then a conditional Markov kernel $\Kk(X|Y,Z)$ of $\Kk(X,Y|Z)$ exists.
    In case 4.\ it is given by the familiar quotient of densities:
    \[ k(x|y,z) = \frac{k(x,y|z)}{k(y|z)}, \qquad k(y|z) := \int_\Xcal k(x,y|z)\,\mubf(dx),\]
    wherever $0 < k(y|z) < \infty$, and by an arbitrary fixed probability measure on $\Xcal$ elsewhere.
    The first three points make no reference to the kernel and thus establish disintegration triples
    $(\Xcal,\Ycal,\Zcal)$, see \Cref{def:disintegration-triple}; the last three are conditions on the given
    $\Kk(X,Y|Z)$.
\end{Thm}

The following corollary is a well known result for probability measures (see \cite{Kle20, Rao05}):

\begin{Cor}[Conditional probability distributions]
    Let $X$ and $Y$ be random variables on probability space $(\Wcal,\Pk(W))$ 
    with standard measurable spaces $\Xcal$, $\Ycal$, resp., as codomains.
    Then there always exist \emph{regular\footnote{The word ``regular'' refers to the fact that the conditional probabilities are 
    Markov kernels as defined above.}
     conditional probability distributions} $\Pk(X|Y)$ and $\Pk(Y|X)$ satisfying:
    \[ \Pk(X,Y) = \Pk(X|Y) \otimes \Pk(Y),\qquad \Pk(X,Y) = \Pk(Y|X) \otimes \Pk(X).  \]
    Furthermore, these conditional probability distributions are essentially unique in the strong sense of the
    second statement of \Cref{lem:ess-unique}, which applies since $\Xcal$ and $\Ycal$ are standard and hence
    countably generated.
\end{Cor}

\section{Transitional Conditional Independence}
\label{sec:transitional_conditional_independence}

\subsection{Definition of Transitional Conditional Independence}

In this section we will introduce the notion of transitional conditional independence for transitional random variables. 
It strengthens the extended conditional independence of \cite{CD17} and the $\Qcal$-extended conditional
independence of \cite{FM20}, see \Cref{lem:tci-implies-eci} and \Cref{thm:supp:ci-equivalences}, and it specializes
to ordinary stochastic conditional independence in the corner case $\Tcal = \Asterisk$; for the other proposals in
the literature, see \cite{Daw79, Daw80, Daw01, RERS17, Cho17, Fri20, FK23}, we record how they relate without
proving implications.
A comparison with other notions of (extended) conditional independence from the literature will be done in 
\Cref{sec:supp:comparison}. 

\begin{Rem}[Repeated occurrences are independent copies]
    \label{rem:independent-copies}
    The joint Markov kernel above is a \emph{product} of the three kernels, see
    \Cref{rem:product-wirings} 1. This encodes a modelling decision that is worth making explicit:
    given $(w,t)$, the three variables $X$, $Y$, $Z$ are drawn \emph{independently} from
    $\Xk(X|W=w,T=t)$, $\Yk(Y|W=w,T=t)$ and $\Zk(Z|W=w,T=t)$. For \emph{deterministic} transitional random
    variables, i.e.\ for $\Xk = \deltabf(X|W,T)$ etc., this is no restriction at all: the product is then the Dirac
    kernel of the joint map $(X,Y,Z)$, see \Cref{thm:join-semi-lattice} point 4. For genuinely stochastic ones,
    however, two occurrences of the \emph{same} $\Xk$ in one expression are two independent copies given $(w,t)$;
    this is exactly why $\ismapof_\Kk$ fails to be reflexive, see \Cref{rem:ismapof-properties} item \ref{rem:no-reflexivity}.
\end{Rem}

\begin{Def}[Transitional conditional independence]
\label{def:transitional_conditional_independence}
Let $(\Wcal \times \Tcal, \Kk(W|T))$ be a transition probability space with Markov kernel:
\[\Kk(W|T): \, \Tcal \dashrightarrow  \Wcal.\]
Consider transitional random variables:
$\Xk:\, \Wcal \times \Tcal \dshto \Xcal$ and 
    $\Yk:\, \Wcal \times \Tcal \dshto \Ycal$ and $\Zk:\, \Wcal \times \Tcal \dshto \Zcal$.
The joint push-forward Markov kernel is then given by:
\[ \Kk(X,Y,Z|T) := \lp\Xk(X|W,T) \otimes \Yk(Y|W,T) \otimes \Zk(Z|W,T) \rp \circ \Kk(W|T). \]
We say that \emph{$\Xk$ is (transitionally) independent of $\Yk$ conditioned on $\Zk$ w.r.t.\ $\Kk=\Kk(W|T)$}, in symbols:
\[\Xk \Indep_{\Kk(W|T)} \Yk \given \Zk, \]
if there \underline{\emph{exists}} a Markov kernel: 
\[ \Qk(X|Z):\; \Zcal \dashrightarrow \Xcal,\]
such that:
\begin{align} 
    \Kk(X,Y,Z|T) = \Qk(X|Z) \otimes \Kk(Y,Z|T), \label{eq:ci}
\end{align}
where $\Kk(Y,Z|T)$ is the marginal of $\Kk(X,Y,Z|T)$.\\
We use the following notations for the following special case:
\[\Xk \Indep_{\Kk(W|T)} \Yk \qquad :\iff \qquad \Xk \Indep_{\Kk(W|T)} \Yk \given \deltabf_\ast. \]
For transitional random variables of the forms $\Xk=\deltabf(X|W,T)$, $\Yk=\deltabf(Y|W,T)$, $\Zk=\deltabf(Z|W,T)$,  where 
$X:\,\Wcal\times \Tcal \to \Xcal$, $Y:\,\Wcal\times \Tcal \to \Ycal$, $Z:\,\Wcal \times \Tcal \to \Zcal$, etc.,
are measurable maps, we might also just write $X$, $Y$, $Z$, instead of $\Xk$, $\Yk$, $\Zk$, in those relations $\Indep$. 
E.g. we would write: 
\[ X\Indep_{\Kk(W|T)} Y \given Z \qquad  :\iff \qquad \deltabf(X|W,T) \Indep_{\Kk(W|T)} \deltabf(Y|W,T) \given \deltabf(Z|W,T). \]
\end{Def}

Several remarks are in order.

\begin{Rem}
    If a candidate $\Qk(X|Z)$ is found then
    for the Equation \ref{eq:ci} to hold it is sufficient to check that for all $t \in \Tcal$ and all measurable $A \ins \Xcal$, $B \ins \Ycal$, $C \ins \Zcal$ one has that:
\begin{align*}
    \Kk(X \in A, Y \in B, Z \in C|T=t) & = \int_C \int_B  \Qk(X \in A|Z=z)\,\Kk(Y \in dy,Z \in dz|T=t).
\end{align*}
\end{Rem}

\begin{Rem}[Essential uniqueness]
    \label{ess-unique-ci}
    The Markov kernel $\Qk(X|Z)$ appearing in the conditional independence $\Xk \Indep_{\Kk(W|T)} \Yk \given \Zk$ in \Cref{def:transitional_conditional_independence} is then a version of a conditional Markov kernel $\Kk(X|Y,Z,T)$ and is thus essentially unique (up to $\Kk(Z|T)$-null set)
    in the sense of \Cref{lem:ess-unique}.
    To be precise: two such kernels $\Qk(X|Z)$ and $\Qk'(X|Z)$ agree for $\Kk(Z|T=t)$-almost all $z$, \emph{for
    every} $t \in \Tcal$; there is in general no single $\Zcal$-null set outside of which they agree. Note that this
    concerns only the \emph{null set}: the kernels themselves do not depend on $t$, see \Cref{rem:Q-no-T}. 
\end{Rem}

\begin{Not}
    The Markov kernel $\Qk(X|Z)$ appearing in the conditional independence $\Xk\Indep_{\Kk(W|T)}\Yk\given\Zk$ is essentially unique
    as remarked in \ref{ess-unique-ci} and we can suggestively write it as:
    \[\Kk(X|\cancel{Y},Z,\cancel{T}):= \Kk(X|\cancel{T,Y},Z) := \Qk(X|Z),\]
    or similarly with crossed variables in different order.
    So we have in case of $\Xk\Indep_{\Kk(W|T)}\Yk\given\Zk$:
    \[ \Kk(X,Y,Z|T) = \Kk(X|\cancel{Y},Z,\cancel{T})\otimes \Kk(Y,Z|T).\]
This notation indicates that $\Kk(X|\cancel{Y},Z,\cancel{T})$ is a version of the conditional Markov kernel $\Kk(X|Y,Z,T)$, but does not (directly)
depend on the arguments of $Y$ and $T$.
\end{Not}

\begin{Rem}[Conditional independence includes conditional independence from $\Tk$]
    \label{add-T}
 By $\Tk$-Inverted Right Decomposition \ref{sep:tci:inv-r-dec} and Right Decomposition \ref{sep:tci:r-dec} we have the equivalence:
 \[\Xk\Indep_{\Kk(W|T)}\Yk\given\Zk \quad\iff \quad \Xk\Indep_{\Kk(W|T)}\Tk\otimes\Yk\given\Zk.\]
 So independence from $\Yk$ automatically comes with independence from the input $\Tk$.
 The reason we can afford this is that we do \emph{not} require the conditioning variable $\Zk$ to be
 ``orthogonal'' to, or ``functionally independent'' of, $\Tk$, as other notions of extended conditional
 independence do. On the contrary, $Z$ may be a direct function of $T$.
 Dependence on $T$ is thus not forbidden; it is re-introduced through the third argument.
 It is this interplay between the second and the third argument of $\Indep_\Kk$ that makes the asymmetric notion
 flexible, and \Cref{rem:Q-no-T} below explains what the two arguments do differently.
\end{Rem}

\begin{Rem}[Why $\Qk$ does not depend on $T$]
    \label{rem:Q-no-T}
    The kernel $\Qk(X|Z)$ in \Cref{def:transitional_conditional_independence} is a Markov kernel on $\Zcal$ alone:
    it is \emph{not} allowed to depend on $t \in \Tcal$. This is deliberate, and it is where the whole content of
    the notion sits. To see the difference, consider a Markov kernel $\Kk(X,Y,Z|T)$ that factorizes as
    \[ \Kk(X,Y,Z|T) = \Kk(X|Z,T) \otimes \Kk(Y|Z,T) \otimes \Kk(Z|T), \]
    i.e.\ $X$ and $Y$ are conditionally independent given $Z$ \emph{and} $T$, in the ordinary sense, for every
    fixed $t$ separately. This is in general \emph{weaker} than $\Xk \Indep_\Kk \Yk \given \Zk$, which demands
    one single kernel $\Qk(X|Z)$ serving all $t$ at once.
    Both statements are expressible in our language, and it is the third argument that distinguishes them:
    \[ \Xk \Indep_\Kk \Yk \given \Zk \otimes \Tk \]
    is the ``for every $t$ separately'' version (there $\Tk$ sits in the conditioning slot, so the kernel
    $\Qk(X|Z,T)$ may depend on $t$) while
    \[ \Xk \Indep_\Kk \Yk \given \Zk \]
    is the uniform one. Under mild hypotheses the weaker version is exactly the family of ordinary conditional
    independences: if $\Xcal$ is standard and $\Zcal$ countably generated then
    $\Xk \Indep_\Kk \Yk \given \Zk \otimes \Tk$ holds if and only if $X \Indep^\omega Y \given Z$ holds under
    $\Kk(X,Y,Z|T=t)$ for every $t \in \Tcal$, see \Cref{thm:supp:ci-equivalences} and the corollary following it.
    The weaker version is the one that the $\Qcal$-extended conditional independence of \cite{FM20} formalizes
    (for $\Qcal \sni \lC \deltabf_t \,|\, t \in \Tcal\rC$), see \Cref{sec:main:comparison}; the uniform one is
    what expresses invariance
    (\Cref{prp:invariance-tci}), ancillarity (\Cref{sec:applications-statistics}) and the global Markov property
    (\Cref{thm-gmp-mI-CBN}); and that is why $\Qk$ carries no $T$.
\end{Rem}

\begin{Rem}[How to read the asymmetry]
    \label{rem:reading-the-asymmetry}
    \Cref{add-T} also says where the asymmetry of $\Indep_\Kk$ comes from. Unfolding both sides through it, the
    two statements
    \[ \Xk \Indep_\Kk \Yk \given \Zk \qquad \text{and} \qquad \Yk \Indep_\Kk \Xk \given \Zk \]
    read
    \[ \Xk \Indep_\Kk \Tk\otimes\Yk \given \Zk \qquad \text{and} \qquad
       \Yk \Indep_\Kk \Tk\otimes\Xk \given \Zk, \]
    and there is no reason for these to be equivalent. Each of them asks its \emph{left} variable to be produced
    from $\Zk$ alone, by one kernel serving all inputs at once, while leaving the right variable free to use the
    input. The asymmetry is therefore not an artefact of the definition; it is the asymmetry between ``is produced
    by $\Zk$'' and ``may depend on $T$''.
    A useful way to read the notation is thus
    \[ \Xk \Indep_\Kk \Yk \given \Zk \qquad \text{as} \qquad \Xk \Indep_\Kk \Rk \otimes \Yk \given \Zk, \]
    where $\Rk$ stands for whatever is left of the input once $\Yk$ and $\Zk$ have been accounted for, informally
    ``$\Rk= \Tk \setminus (\Yk,\Zk)$'': the whole input always sits on the right of the bar, and $\Zk$ is the only thing
    $\Xk$ is allowed to use. Symmetry can then be expected only once the input has been exhausted: either
    because the conditioning variable already contains it, which is $\Tk$-Restricted Symmetry
    \ref{sep:tci:t-res-sym}, or because there is no input at all, $\Tcal = \Asterisk$, which is Symmetry
    \ref{sep:tci:sym}. Both still ask for a disintegration triple; exhausting the input is necessary for symmetry,
    not sufficient.
\end{Rem}

\begin{Rem}[Existence of conditional Markov kernels expressed as conditional independence]
    \label{conditional-markov-kernel-as-ci}
    Let $\Xk$, $\Yk$ be transitional random variables on transition probability space $(\Wcal \times \Tcal, \Kk(W|T))$.
    Then we can express the existence of a conditional Markov kernel $\Kk(X|Y,T)$ of $\Kk(X,Y|T)$ equivalently in
    one of the following equivalent statements: 
    \[ \Xk \ismapof_\Kk^\ast \Yk \otimes \Tk \qquad \iff \qquad \Xk\Indep_{\Kk(W|T)}\deltabf_\ast\given\Yk\otimes\Tk  
     \qquad \iff \qquad \Xk\Indep_{\Kk(W|T)}\Tk\given\Yk\otimes\Tk.\]
    Note that for standard measurable space $\Xcal$ and countably generated $\Ycal$ the above statements always hold by 
    \Cref{thm-regular-conditional-Markov-kernel}.
\end{Rem}

\subsection{Transitional Conditional Independence for Random Variables}

\begin{Rem}[Transitional conditional independence for random variables]
    \label{ci-random-variables}
    If we translate \emph{transitional conditional independence} to random variables $X,Y,Z$
on a probability space $(\Wcal, \Pk(W))$, i.e.\ taking $\Tcal=\Asterisk$, then we arrive at:
\begin{align} 
    X \Indep_{\Pk(W)} Y \given Z \qquad \iff \qquad \exists \Qk(X|Z):\, \Pk(X,Y,Z) = \Qk(X|Z) \otimes \Pk(Y,Z).
\end{align}
Such a $\Qk(X|Z)$ would then clearly be a regular version of both, $\Pk(X|Y,Z)$ and $\Pk(X|Z)$. 
It thus directly implies what we will call \emph{weak conditional independence}:
\begin{align} 
    X \Indep_{\Pk(W)}^\omega Y \given Z  &\quad:\iff \quad 
    \forall A \in \Bcal_\Xcal:\, \E[\I_A(X)|Y,Z] = \E[\I_A(X)|Z] \quad \Pk(W)\text{-a.s.},
\end{align}
which makes use of the conditional expectations for each $A$, which exist for all measurable spaces, 
in contrast to regular conditional probability distributions. Both notions of conditional independence 
can be defined for all measurable spaces. Transitional conditional independence incorporates 
the existence of such a $\Pk(X|Z)$ and a factorization of the joint directly into its definition. 
Certainly, if a regular version of $\Pk(X|Z)$ does not even exist the variables 
are declared (transitionally) conditionally dependent. 
But in case $\Pk(X|Z)$ exists, e.g.\ for standard measurable $\Xcal$ and $\Zcal$ 
by \Cref{thm-regular-conditional-Markov-kernel}, 
then both notions of conditional independence are equivalent. So the choice of which notion to pick
depends on how much meaning one finds in the existence of such a regular version $\Pk(X|Z)$ and 
a factorization. In the applications to causal graphical models, where one wants to connect and work with many different
subsystems, the existence of such conditional Markov kernels is crucial, because otherwise those subsystems might 
not even be well-defined.

Even though, one might argue that asking to check for the existence of a regular versions of $\Pk(X|Z)$ seems like
an unnecessary burden, from the point on we prove how the existence of such Markov kernels 
can be inherited through the (asymmetric) \emph{separoid rules}, see \Cref{thm:separoid_axioms-tci},
or can be guaranteed just through graphical criteria, see the \emph{global Markov property} in \Cref{thm-gmp-mI-CBN},
it will turn out to be very useful to get such conditional Markov kernels (almost) for free.

\end{Rem}

\subsection{Transitional Conditional Independence for Deterministic Variables}

We now demonstrate how transitional conditional independence behaves when applied to the other corner case of deterministic functions
that contain no stochasticity.

\begin{Thm}[Transitional conditional independence for deterministic variables]
    \label{thm:variation-ci}
    Let $F:\, \Tcal \to \Fcal$ and $H:\,\Tcal \to \Hcal$ be measurable maps and $\Fcal$ a standard measurable space.
    We now consider them as (deterministic) transitional random variables on the transition probability space $(\Wcal \times \Tcal,\Kk(W|T))$.
    Let $\Yk:\, \Wcal \times \Tcal \dshto \Ycal$ be another transitional random variable.\\
    Then the following statements are equivalent:
    \begin{enumerate}
        \item $\displaystyle F \Indep_{\Kk(W|T)} \Yk \given H$.
        \item There exists a measurable function $\varphi:\, \Hcal \to \Fcal$ such that $F=\varphi \circ H$.
    \end{enumerate}
\end{Thm}

\begin{Rem} 
    \begin{enumerate}
    \item    Note that the second statement is  independent of $Y$.
    \item It is worth stressing that the identity $F = \varphi \circ H$ holds \emph{pointwise}, for every
        $t \in \Tcal$, and not merely almost surely: the proof evaluates Markov kernels at each $t$ separately, so
        no null set appears anywhere. The same is true for the propensity score in \Cref{propensity}. This is a
        genuine gain of the transitional set-up over the classical one, where the corresponding statements hold only
        up to null sets.
    \item \label{rem:variation-ci-hypotheses}
        In the direction 1.\ $\implies$ 2.\ one may replace the standardness of $\Fcal$ by the weaker
        assumption that $\Bcal_\Fcal$ separates the points of $\Fcal$, at the price of obtaining $\varphi$ only on
        $H(\Tcal)$, i.e.\ a measurable map $\varphi:\, H(\Tcal) \to \Fcal$ with $F=\varphi \circ H$.
        Standardness of $\Fcal$ is used twice and for two different purposes: for the point separation just
        mentioned, and for Kuratowski's extension theorem, which extends $\varphi$ from $H(\Tcal)$ to all of
        $\Hcal$. 
    \item  \Cref{thm:variation-ci} shows  how transitional conditional independence can express certain 
        functional conditional (in)dependences. It also shows its (restricted) relation to 
        variation conditional independence, see \cite{CD17, CD17-supp}.            
    \item  The full equivalence in \Cref{thm:variation-ci} for standard $\Fcal$ 
        needs Kuratowski's extension theorem for standard measurable spaces (see \cite{Kec95} 12.2). 
    See \Cref{sec:kuratowski-extension}.
    The proof of \Cref{thm:variation-ci} can be found in \Cref{thm:supp:variation-ci} in \Cref{sec:tci-det}.
    \end{enumerate}
\end{Rem}

\begin{Eg}
    If, for example, $\Tcal = \Tcal_1 \times \Tcal_2$ and $T_i:\,\Wcal \times \Tcal_1 \times \Tcal_2 \to \Tcal_i$ the canonical projection onto $\Tcal_i$,
    then $F$ is a function in two variables $(t_1,t_2)$.
    We then  have:
    \[ F \Indep_{\Kk(W|T)} T_1 \given T_2,\]
    if and only if $F$ - as a function - is only dependent on the argument $t_2$ (and not on $t_1$).
\end{Eg}

\subsection{Criteria for Checking Transitional Conditional Independence}
\label{sec:checking-tci}

The whole content of \Cref{def:transitional_conditional_independence} sits in the existential quantifier, so it is
worth recording how one produces the kernel $\Qk(X|Z)$, or rules out that it exists. The first observation says
what a witness must look like.

\begin{Rem}[Narrowing the search for $\Qk(X|Z)$]
    \label{rem:narrowing-Q}
    By \Cref{def:transitional_conditional_independence} the relation $\Xk \Indep_\Kk \Yk \given \Zk$ says that:
    \[ \Kk(X,Y,Z|T) = \Qk(X|Z) \otimes \Kk(Y,Z|T) \]
    for some Markov kernel $\Qk(X|Z)$. Multiplying this identity by any probability measure $\Pk(T)$ on $\Tcal$
    from the right and marginalizing $\Tcal$ and $\Ycal$ out gives:
    \[ \Pk(X,Z) = \Qk(X|Z) \otimes \Pk(Z), \]
    where $\Pk(X,Z)$ and $\Pk(Z)$ are the marginals of $\Pk(X,Y,Z) := \Kk(X,Y,Z|T) \circ \Pk(T)$; for
    $\Pk(T) := \deltabf_{t_0}$ they are $\Kk(X,Z|T=t_0)$ and $\Kk(Z|T=t_0)$. So a witness $\Qk(X|Z)$ is itself a
    conditional Markov kernel of $\Pk(X,Z)$, for \emph{every} prior $\Pk(T)$ at once, and any version $\Pk(X|Z)$
    of that conditional therefore satisfies $\Qk(X|Z) = \Pk(X|Z)$ $\Pk(Z)$-almost surely, for each event
    separately, see \Cref{lem:ess-unique}. The converse fails: for a given $\Pk(T)$ one may \emph{not} assume
    that $\Pk(X|Z)$ is a witness for $\Xk \Indep_\Kk \Yk \given \Zk$, since it is left unconstrained wherever
    $\Pk(Z)$ puts no mass. Requiring in addition that $\Kk(Z|T=t) \ll \Pk(Z)$ for every $t \in \Tcal$ repairs
    this, as the next theorem shows.
\end{Rem}

No distribution on the input is being asserted here: $\Pk(T)$ is a device for producing one candidate for the
witness $\Qk(X|Z)$, and \Cref{rem:Q-no-T} applies to it verbatim. The second hypothesis below costs nothing
beyond a disintegration triple, and holds already when $\Xcal$ is standard, with $\Zcal$ arbitrary, by point 2.\
of \Cref{thm-regular-conditional-Markov-kernel} with $\Asterisk$ in its third slot.

\begin{Thm}[A witness from a prior]
    \label{thm:tci-prior}
    Let $\Pk(T)$ be a probability measure on $\Tcal$ and put:
    \[ \Pk(X,Y,Z) := \Kk(X,Y,Z|T) \circ \Pk(T). \]
    Assume that:
    \begin{enumerate}
        \item $\Kk(Z|T=t) \ll \Pk(Z)$ for every $t \in \Tcal$;
        \item $\Pk(X,Z)$ admits a conditional Markov kernel $\Pk(X|Z):\, \Zcal \dshto \Xcal$, i.e.\
            $\Pk(X,Z) = \Pk(X|Z) \otimes \Pk(Z)$.
    \end{enumerate}
    Then, for every version of $\Pk(X|Z)$:
    \[ \Xk \Indep_\Kk \Yk \given \Zk \qquad \iff \qquad
       \Kk(X,Y,Z|T) = \Pk(X|Z) \otimes \Kk(Y,Z|T). \]
    The proof is given in \Cref{thm:supp:tci-prior} in \Cref{sec:supp:checking}.
\end{Thm}

\begin{Rem}[Hypothesis 1.\ cannot be dropped]
    \label{rem:prior-hypotheses}
    Without 1.\ the kernel $\Pk(X|Z)$ is left unconstrained wherever $\Pk(Z)$ puts no mass, and a version chosen
    badly there breaks the criterion. Take $\Tcal = \lC 1,2 \rC$, $\Xcal = \Zcal = \lC 0,1 \rC$,
    $\Ycal = \Asterisk$ and:
    \[ \Kk(Z|T=1) := \deltabf_0, \qquad \Kk(Z|T=2) := \deltabf_1, \qquad
       \Qk(X|Z=0) := \Qk(X|Z=1) := \deltabf_0, \]
    with $\Kk(X,Y,Z|T) := \Qk(X|Z) \otimes \Kk(Y,Z|T)$, so that $\Xk \Indep_\Kk \Yk \given \Zk$ holds by
    construction. For the prior $\Pk(T) := \deltabf_1$ we get $\Pk(Z) = \deltabf_0$, so hypothesis 2.\ holds for
    every version of $\Pk(X|Z)$ with $\Pk(X|Z=0) = \deltabf_0$, while $\Pk(X|Z=1)$ is left arbitrary; choosing
    $\Pk(X|Z=1) := \deltabf_1$ makes the right hand side of the criterion fail at $t=2$. Hypothesis 1.\ fails
    here, since $\Kk(Z|T=2) = \deltabf_1 \not\ll \deltabf_0 = \Pk(Z)$, and it is precisely what rules this out:
    for each fixed event, two versions of $\Pk(X|Z)$ agree at $\Pk(Z)$-almost all $z$ by \Cref{lem:ess-unique},
    and 1.\ makes that difference invisible to every $\Kk(Y,Z|T=t)$.
\end{Rem}

Once the kernel is dominated by a product of $\sigma$-finite measures, the conditional $\Pk(X|Z)$ of
\Cref{thm:tci-prior} 2.\ can be written down from densities and the criterion becomes an identity between
densities. One half of this needs no prior at all: any measurable $q:\, \Xcal \times \Zcal \to [0,\infty]$ with
$\int_\Xcal q(x|z)\,\mubf_\Xcal(dx) = 1$ for every $z \in \Zcal$ and
$k(x,y,z|t) = q(x|z) \cdot k(y,z|t)$ $\mubf$-almost everywhere, for every $t \in \Tcal$, defines a witness
$\Qk(X \in A|Z=z) := \int_A q(x|z)\, \mubf_\Xcal(dx)$, by the same computation that identifies the density of
$\Pk(X|Z) \otimes \Kk(Y,Z|T=t)$ in the proof below. The content is the converse.

\begin{Cor}[A witness from a prior in terms of densities]
    \label{cor:tci-density}
    Let $\mubf_\Xcal$, $\mubf_\Ycal$, $\mubf_\Zcal$ be $\sigma$-finite measures on $\Xcal$, $\Ycal$, $\Zcal$, put
    $\mubf := \mubf_\Xcal \otimes \mubf_\Ycal \otimes \mubf_\Zcal$ and assume that:
    \[ \Kk(X,Y,Z|T=t) \ll \mubf \qquad \text{ for every } t \in \Tcal. \]
    For every $t \in \Tcal$ fix a density $k(x,y,z|t)$ of $\Kk(X,Y,Z|T=t)$ w.r.t.\ $\mubf$, not necessarily
    jointly measurable in $t$, and write $k(y,z|t)$ for the corresponding marginal density of
    $\Kk(Y,Z|T=t)$. Let $\Pk(T)$ be a probability measure on $\Tcal$ such that:
    \[ \Kk(Y,Z|T=t) \ll \Pk(Y,Z) \qquad \text{ for every } t \in \Tcal, \]
    where we have put $\Pk(X,Y,Z) := \Kk(X,Y,Z|T) \circ \Pk(T)$. Then $\Pk(X,Y,Z) \ll \mubf$; let $p(x,y,z)$ be a
    density of it w.r.t.\ $\mubf$, with marginal densities $p(x,z)$ and $p(z)$.\footnote{\label{fn:hs-version}No
    statement below depends on these choices: two densities of the same measure w.r.t.\ $\mubf$ agree
    $\mubf$-almost everywhere, and altering $k(x,y,z|t)$ or $p(x,y,z)$ on $\mubf$-null sets alters $p(x|z)$ and
    $\Pk(X|Z)$ only on null sets of the corresponding reference measures.} Finally fix any probability density $w$
    on $\lp \Xcal, \mubf_\Xcal \rp$ and put:
    \[ p(x|z) := \begin{cases} p(x,z) \big/ p(z), & \text{ if } 0 < p(z) < \infty,\\
                               w(x), & \text{ otherwise,}\end{cases} \]
    \[ \Pk(X \in A|Z=z) := \int_A p(x|z)\, \mubf_\Xcal(dx). \]
    Then $\Pk(X|Z):\, \Zcal \dshto \Xcal$ is a conditional Markov kernel of $\Pk(X,Z)$ as in
    \Cref{thm:tci-prior} 2., so that \Cref{thm:tci-prior} applies, and the following are equivalent:
    \begin{enumerate}
        \item $\Xk \Indep_\Kk \Yk \given \Zk$;
        \item $\Kk(X,Y,Z|T) = \Pk(X|Z) \otimes \Kk(Y,Z|T)$;
        \item for every $t \in \Tcal$: $k(x,y,z|t) \cdot p(z) = p(x,z) \cdot k(y,z|t)$ for $\mubf$-almost all
            $(x,y,z)$, where the exceptional set may depend on $t$.
    \end{enumerate}
    The proof is given in \Cref{cor:supp:tci-density} in \Cref{sec:supp:checking}.
\end{Cor}

\begin{Rem}[The lemma of Halmos and Savage]
    \label{rem:halmos-savage}
    A prior $\Pk(T)$ as in \Cref{cor:tci-density} always exists under the single assumption:
    \[ \Kk(X,Y,Z|T=t) \ll \mubf \qquad \text{ for every } t \in \Tcal, \]
    for a $\sigma$-finite measure $\mubf$ on $\Xcal \times \Ycal \times \Zcal$, which for this purpose need not
    be a product. Indeed, by a classical lemma of \cite{HS49}, see \Cref{lem:supp:halmos-savage}, a family of
    measures that is dominated by a $\sigma$-finite measure is already dominated by one of its own countable
    mixtures. Applied to $\lp \Kk(X,Y,Z|T=t) \rp_{t \in \Tcal}$ and $\mubf$, it provides a countable subset
    $\Ccal \ins \Tcal$ and weights $p(t)>0$ with $\sum_{t \in \Ccal} p(t) = 1$ such that the \emph{discrete}
    prior $\Pk(T) := \sum_{t \in \Ccal} \deltabf_t \cdot p(t)$, a probability measure on $\Tcal$ whatever its
    $\sigma$-algebra, satisfies:
    \[ \Kk(X,Y,Z|T=t) \ll \Pk(X,Y,Z) = \sum_{t \in \Ccal} \Kk(X,Y,Z|T=t) \cdot p(t)
       \qquad \text{ for every } t \in \Tcal, \]
    and hence, by marginalization, also $\Kk(Y,Z|T=t) \ll \Pk(Y,Z)$, which is the hypothesis of
    \Cref{cor:tci-density}. For a discrete prior the density there is explicit,
    \[ p(x,y,z) = \sum_{t \in \Ccal} k(x,y,z|t) \cdot p(t), \]
    a countable sum of the given densities; for a general prior the Radon-Nikodym theorem supplies $p(x,y,z)$ only
    abstractly, since writing it as $\int_\Tcal k(x,y,z|t)\, \Pk(T \in dt)$ would need the joint measurability in
    $t$ that \Cref{cor:tci-density} avoids.
\end{Rem}

\begin{Rem}[How to get a prior]
    \label{rem:prior-cases}
    Three cases give a prior $\Pk(T)$ as in \Cref{cor:tci-density} directly, without
    \Cref{rem:halmos-savage}, and they cover most of what one meets in practice.
    \begin{enumerate}
        \item If $\Tcal$ is countable, then \emph{any} strictly positive probability weight function $p$ on
            $\Tcal$ gives one, namely $\Pk(T) := \sum_{t \in \Tcal} \deltabf_t \cdot p(t)$: since all weights are
            strictly positive, $\Pk(E) = \sum_{t \in \Tcal} \Kk(E|T=t) \cdot p(t)$ vanishes if and only if
            $\Kk(E|T=t)$ vanishes for every $t$, so $\Kk(X,Y,Z|T=t) \ll \Pk(X,Y,Z)$ for every $t$.
        \item If some $t_0 \in \Tcal$ satisfies $\Kk(X,Y,Z|T=t) \ll \Kk(X,Y,Z|T=t_0)$ for every $t \in \Tcal$
            (for instance if $\Tcal = \Asterisk$), then the prior $\Pk(T) := \deltabf_{t_0}$ does,
            $p(x,y,z) = k(x,y,z|t_0)$, and the kernel of \Cref{cor:tci-density} is the ordinary conditional at
            that one input value, $p(x|z) = k(x,z|t_0)/k(z|t_0)$ wherever $0 < k(z|t_0) < \infty$, as in
            \Cref{rem:narrowing-Q}.
        \item If $\Xcal$, $\Ycal$, $\Zcal$ and $\Tcal$ are countable and discrete, take counting measures for
            $\mubf_\Xcal$, $\mubf_\Ycal$, $\mubf_\Zcal$ and the prior $\Pk(T)$ of 1.; then
            $k(x,y,z|t)$ is the mass function of $\Kk(X,Y,Z|T)$, there are no non-empty null sets, and
            \Cref{cor:tci-density} 3.\ reads:
            \[ \Xk \Indep_\Kk \Yk \given \Zk \quad \iff \quad
               \forall x,y,z,t:\; k(x,y,z|t) \cdot p(z) = p(x,z) \cdot k(y,z|t), \]
            a finite system of polynomial identities in the entries of $k$ when all four spaces are finite.
    \end{enumerate}
\end{Rem}

\begin{Rem}[What \Cref{cor:tci-density} tests]
    \label{rem:hs-reading}
    In \Cref{cor:tci-density} the existential quantifier of
    \Cref{def:transitional_conditional_independence} has disappeared: its 1.\ asserts that \emph{some} witness
    $\Qk(X|Z)$ exists, its 2.\ that the kernel $\Pk(X|Z)$ built there is one, and there is nothing left to search
    for. What its condition 3.\ tests is best seen by splitting it into two classical statements, which
    \Cref{cor:supp:tci-density} shows to be jointly equivalent to it:
    \begin{enumerate}[label=(\alph*)]
        \item $X \Indep^\omega_{\Pk(X,Y,Z)} Y \given Z$, ordinary conditional independence under the single
            probability measure $\Pk(X,Y,Z)$; and
        \item $\Xk \Indep_\Kk \Tk \given \Yk \otimes \Zk$, i.e.\ $(Y,Z)$ is \emph{sufficient} for the input in
            the sense of \Cref{sec:ancillarity-sufficiency}: a conditional Markov kernel $\Kk(X|Y,Z)$ exists and
            does not depend on $t$. In the density form of \cite{HS49} this says that for every $t \in \Tcal$ the
            Radon-Nikodym derivative:
            \[ \frac{\Kk(X \in dx, Y \in dy, Z \in dz|T=t)}{\Pk(X \in dx, Y \in dy, Z \in dz)} \]
            has a version that does not vary with $x$, i.e.\ is a function of $(y,z)$ alone. Equivalently, the
            $\sigma$-algebra generated by $(Y,Z)$ is a sufficient $\sigma$-algebra for the family
            $\lp \Kk(X,Y,Z|T=t) \rp_{t \in \Tcal}$.
    \end{enumerate}
    This is the informal reading ``the conditional of $X$ given $(Y,Z,T)$ depends neither on $Y$ nor on $T$'',
    with each ``neither'' isolated into one classical condition. Two caveats. Half (b) is a condition on all of
    $\Tcal$ and cannot be read off the prior alone: the prior pins down the null sets of the family, not its
    densities, and for a discrete prior only half (a) is a statement about countably many input values. And
    $\Pk(X|Z)$ itself depends on the choice of $\Pk(T)$ and $w$, while \Cref{cor:tci-density} does not, since it
    holds for every admissible choice; where $\Xk \Indep_\Kk \Yk \given \Zk$ does hold, any two choices agree
    $\Kk(Z|T=t)$-almost surely for every $t$, by \Cref{lem:ess-unique}.
\end{Rem}

For linear Gaussian Markov kernels the criterion can be solved explicitly, and the answer separates cleanly into a
classical part and a transitional one.

\begin{Thm}[Linear Gaussian Markov kernels]
    \label{thm:tci-gaussian}
    Let $\Tcal = \R^{d_T}$, $\Xcal = \R^{d_X}$, $\Ycal = \R^{d_Y}$, $\Zcal = \R^{d_Z}$ and let $\Pk(X,Y,Z|T)$ be
    the linear Gaussian Markov kernel:
    \[ \Pk(X,Y,Z|T=t) = \Ncal \lp
        \begin{bmatrix} \Gamma_X \\ \Gamma_Y \\ \Gamma_Z \end{bmatrix} t +
        \begin{bmatrix} \gamma_X \\ \gamma_Y \\ \gamma_Z \end{bmatrix} ,\;
        \begin{bmatrix} \Sigma_{XX} & \Sigma_{XY} & \Sigma_{XZ} \\
                        \Sigma_{YX} & \Sigma_{YY} & \Sigma_{YZ} \\
                        \Sigma_{ZX} & \Sigma_{ZY} & \Sigma_{ZZ} \end{bmatrix} \rp,  \]
    with a positive definite covariance matrix. Then:
    \[ X \Indep_{\Pk(X,Y,Z|T)} Y \given Z \qquad \iff \qquad
       \Sigma_{XY} = \Sigma_{XZ}\Sigma_{ZZ}^{-1}\Sigma_{ZY} \quad \land \quad
       \Gamma_X = \Sigma_{XZ}\Sigma_{ZZ}^{-1}\Gamma_Z. \]
    Here the prior $\Pk(T) := \deltabf_0$ is as in \Cref{rem:prior-cases} 2., since any two Gaussian distributions
    with non-degenerate covariance matrices are equivalent, both being equivalent to the Lebesgue measure, so the
    kernel of \Cref{cor:tci-density} is the
    ordinary conditional at $t=0$; in the affirmative case it is a witness, namely the linear Gaussian Markov
    kernel:
    \begin{align*}
        \Qk(X|Z=z) &= \Pk(X|Z=z,T=0) \\
        &= \Ncal \lp \gamma_X + \Sigma_{XZ}\Sigma_{ZZ}^{-1}(z-\gamma_Z),\;
          \Sigma_{XX} - \Sigma_{XZ}\Sigma_{ZZ}^{-1}\Sigma_{ZX} \rp,
    \end{align*}
    so that $\Pk(X,Y,Z|T) = \Pk(X|Z,T=0) \otimes \Pk(Y,Z|T)$.
    The proof is given in \Cref{thm:supp:tci-gaussian} in \Cref{sec:supp:checking}.
\end{Thm}

\begin{Rem}[What the second condition says]
    \label{rem:gaussian-reading}
    The first condition is the classical criterion for conditional independence of jointly Gaussian variables and
    constrains the covariance alone. The second constrains only the coefficients through which the input enters,
    and it is what transitional conditional independence adds: the input may act on $X$, but only through $Z$ and
    in exactly the proportion in which $Z$ already carries it. Two special cases make this concrete. If
    $\Gamma_Z = 0$, i.e.\ the conditioning variable does not react to the input at all, the condition forces
    $\Gamma_X = 0$: nothing that reacts to the input may stand on the left of the bar. And for $d_T = 0$, i.e.\
    $\Tcal = \Asterisk$, the second condition is empty and one recovers the classical Gaussian criterion, in
    accordance with Symmetry \ref{sep:tci:sym}. Both conditions are closed algebraic conditions on the parameters,
    so in this model class the relation is decided by linear algebra.
\end{Rem}

\subsection{Separoid Rules for Transitional Conditional Independence}
\label{sec:separoid-rules-tci}

In the following we will list all the left and right versions of the separoid rules (see \cite{Daw01} for the symmetric versions or
\Cref{sec:sym-sep-asym-sep}) that hold for transitional conditional independence.
Note that almost all of these work for all measurable spaces. Some of the rules, especially Left Weak Union, require
the existence of a conditional Markov kernel, i.e.\ that the codomains of the transitional random variables
involved form a \emph{disintegration triple} in the sense of \Cref{def:disintegration-triple}.
By \Cref{thm-regular-conditional-Markov-kernel} this is in particular the case when the first space is standard
and the second one is countably generated, when the first space is standard and the parameter space is discrete,
or when the second space is discrete.

Formally we will show that the class of transitional random variables whose codomains form disintegration triples
together with transitional conditional independence will form what we will call a \emph{$T$-$\ast$-separoid},
i.e.\ an asymmetric analogue of the separoid rules of \cite{Daw01}, see \Cref{def:t-k-separoid} below;
\Cref{sec:sym-sep-asym-sep} shows how such asymmetric rules arise from symmetric ones.
Note that these rules have not been proven in this amplitude for other versions of extended conditional independence, see \cite{CD17}.

The proofs for these separoid rules for transitional random variables will be given in
\Cref{sec:proofs-separoid-rules-for-transitional-conditional-independence}.

\begin{Thm}[Separoid rules for transitional conditional independence]
    \label{thm:separoid_axioms-tci}
    Consider a transition probability space $\lp \Wcal \times \Tcal, \Kk(W|T) \rp$ and transitional random variables
    $\Xk:\, \Wcal \times \Tcal \dshto \Xcal$ and
    $\Yk:\, \Wcal \times \Tcal \dshto \Ycal$ and $\Zk:\, \Wcal \times \Tcal \dshto \Zcal$ and
    $\Uk:\, \Wcal \times \Tcal \dshto \Ucal$.
   Then the ternary relation $\Indep = \Indep_{\Kk(W|T)}$ satisfies the following rules:
{\setlength{\itemsep}{1pt}\setlength{\parskip}{0pt}
\begin{enumerate}[label=\alph*)]
    \item Extended Left Redundancy \ref{sep:tci:ext-l-red}:
    \item[] \mbox{$\Xk \ismapof_\Kk \Zk$}\Rimp\mbox{$\Xk \Indep  \Yk \given \Zk$}.
    \item $\Tk$-Restricted Right Redundancy \ref{sep:tci:r-red} (for $(\Xcal,\Zcal,\Tcal)$ a disintegration triple)\footnote{\label{fn:sep-ax-std}The disintegration triple is assumed only to guarantee the existence of a conditional Markov kernel, see \Cref{def:disintegration-triple} and \Cref{thm-regular-conditional-Markov-kernel}. It can be traded for criteria checked on the kernels at hand, several of which restrict the spaces not at all, see \Cref{sep:tci:r-red} for b), m) and n), and \Cref{lem:l-uni-criteria} for f).}:
    \item[] \mbox{$\Xk\Indep \deltabf_\ast \given \Zk\otimes\Tk$} always holds.
    \item Left Decomposition \ref{sep:tci:l-dec}:
	\item[]\mbox{$\Xk\otimes\Uk \Indep  \Yk \given \Zk$}\Rimp\mbox{$\Uk \Indep  \Yk \given \Zk$}.
    \item Right Decomposition \ref{sep:tci:r-dec}:
    \item[] \mbox{$\Xk \Indep  \Yk\otimes\Uk \given \Zk$}\Rimp\mbox{$\Xk \Indep  \Uk \given \Zk$}.
    \item $\Tk$-Inverted Right Decomposition \ref{sep:tci:inv-r-dec}:
    \item[] \mbox{$\Xk \Indep  \Yk \given \Zk$}\Rimp\mbox{$\Xk \Indep  \Tk\otimes\Yk \given \Zk$}.
    \item Left Weak Union \ref{sep:tci:l-uni} (for $(\Xcal,\Ucal,\Zcal)$ a disintegration triple)\footref{fn:sep-ax-std}:
    \item[] \mbox{$\Xk\otimes\Uk \Indep  \Yk  \given \Zk$}\Rimp\mbox{$\Xk \Indep  \Yk \given \Uk\otimes\Zk$}.
    \item Right Weak Union \ref{sep:tci:r-uni}:
  \item[] \mbox{$\Xk \Indep  \Yk\otimes\Uk \given \Zk$}\Rimp\mbox{$\Xk \Indep  \Yk \given \Uk\otimes\Zk$}.
  \item Left Contraction \ref{sep:tci:l-con}:\item[]
 \mbox{$(\Xk \Indep  \Yk \given \Uk\otimes\Zk) \land (\Uk \Indep  \Yk \given \Zk)$}\Rimp\mbox{$\Xk\otimes\Uk \Indep  \Yk \given \Zk$}.
\item Right Contraction \ref{sep:tci:r-con}:
\item[] \mbox{$(\Xk \Indep  \Yk \given \Uk\otimes\Zk) \land (\Xk \Indep\Uk \given \Zk)$}\Rimp\mbox{$\Xk \Indep \Yk\otimes\Uk \given \Zk$}.
\item Right Cross Contraction \ref{sep:tci:rc-con}:
\item[] \mbox{$(\Xk \Indep  \Yk \given \Uk\otimes\Zk) \land (\Uk \Indep  \Xk \given \Zk)$}\Rimp\mbox{$\Xk  \Indep  \Yk\otimes\Uk \given \Zk$}.
\item Flipped Left Cross Contraction \ref{sep:tci:flc-con}:
\item[] \mbox{$(\Xk \Indep  \Yk \given \Uk\otimes\Zk) \land (\Yk \Indep  \Uk \given \Zk)$}\Rimp\mbox{$\Yk \Indep  \Xk\otimes\Uk \given \Zk$}.
\end{enumerate}}
\end{Thm}%

One may wonder why \emph{four} contraction rules h)--k) are listed. The answer is that they are exactly the four
cases that arise when Contraction is unfolded for the symmetrized relation $\Indep^\lor$, see
\Cref{thm:sym-t-k-separoid}; dropping any one of them would break that argument.

\begin{Rem}
In particular, we have the equivalence:
\[ (\Xk \Indep  \Yk\otimes\Uk \given \Zk) \quad\iff\quad (\Xk \Indep  \Yk \given \Uk\otimes\Zk) \quad\land\quad (\Xk \Indep\Uk \given \Zk)  ,\]
 and, if $(\Xcal,\Ucal,\Zcal)$ is a disintegration triple:
\[(\Xk\otimes\Uk \Indep  \Yk \given \Zk) \quad\iff\quad (\Xk \Indep  \Yk \given \Uk\otimes\Zk) \quad\land\quad (\Uk \Indep  \Yk \given \Zk).\]
\end{Rem}

\Cref{tab:rules-hypotheses} at the end of this section records, for each rule, whether it costs a hypothesis and
where it is proven, and pairs it with its graphical counterpart from \Cref{d-sep-separoid-axioms}; the inductive
proof of the global Markov property in \Cref{sec:main:global_markov} does nothing but move along the rows of that
table.

\begin{Cor}[Symmetry]
    \label{cor:sep-ci:symmetry}
    Let the setting be like in \Cref{thm:separoid_axioms-tci}. We then have:
{\setlength{\itemsep}{1pt}\setlength{\parskip}{0pt}
\begin{enumerate}[resume,label=\alph*)]
    \setcounter{enumi}{11}
\item Restricted Symmetry \ref{sep:tci:res-sym}:
    \item[] \mbox{$(\Xk \Indep  \Yk \given \Zk) \land (\Yk \Indep \deltabf_\ast \given \Zk)$}\Rimp\mbox{$\Yk \Indep \Xk \given \Zk$}.
    \item $\Tk$-Restricted Symmetry \ref{sep:tci:t-res-sym} (for $(\Ycal,\Zcal,\Tcal)$ a disintegration triple)\footref{fn:sep-ax-std}:
    \item[] \mbox{$\Xk \Indep  \Yk \given \Zk\otimes\Tk$}\Rimp\mbox{$\Yk \Indep  \Xk \given \Zk\otimes\Tk$}.
    \item Symmetry \ref{sep:tci:sym} (for $\Tcal=\Asterisk$ and $(\Ycal,\Zcal,\Asterisk)$ a disintegration triple)\footref{fn:sep-ax-std}:
    \item[] \mbox{$\Xk \Indep  \Yk \given \Zk$}\Rimp\mbox{$\Yk \Indep  \Xk \given \Zk$}.
\end{enumerate}}
\end{Cor}

\begin{Cor}
    Let the setting be like in \Cref{thm:separoid_axioms-tci}. We then have:
{\setlength{\itemsep}{1pt}\setlength{\parskip}{0pt}
\begin{enumerate}[resume,label=\alph*)]
    \setcounter{enumi}{14}
    \item Inverted Left Decomposition \ref{sep:tci:inv-l-dec}:
    \item[] \mbox{$(\Xk \Indep  \Yk \given \Zk) \land \lp \Uk \ismapof_\Kk \Xk\otimes\Zk\rp$}\Rimp\mbox{$\Xk\otimes\Uk \Indep \Yk \given \Zk$}.
    \item $\Tk$-Extended Inverted Right Decomposition \ref{sep:tci:ext-inv-r-dec}:
    \item[] \mbox{$(\Xk \Indep  \Yk \given \Zk) \land  \lp \Uk \ismapof_\Kk \Tk \otimes \Yk\otimes\Zk \rp$}\Rimp\mbox{$\Xk \Indep \Tk\otimes\Yk\otimes\Uk \given \Zk$}.
    \item Equivalent Exchange \ref{sep:tci:eq-ex}:
    \item[] \mbox{$\lp \Xk \Indep_\Kk  \Yk \given \Zk \rp \land \lp \Zk \approx_\Kk \Zk'
        \rp$}\Rimp\mbox{$\Xk \Indep_\Kk\Yk \given \Zk'.$}
\item Full Equivalent Exchange \ref{sep:tci:full-eq-ex}: If
    $(\Xk' \approx_\Kk \Xk)\,\land\, (\Yk' \approx_\Kk \Yk)\,\land\,(\Zk' \approx_\Kk \Zk)$ then:
\item[] \mbox{$\Xk \Indep_\Kk  \Yk \given \Zk \qquad$}\Riff\mbox{$\qquad
        \Xk' \Indep_\Kk \Yk' \given \Zk'.$}
\end{enumerate}}

\end{Cor}
The rules a)--k) above are the ones we will need throughout, and they are worth a name of their own; the three
symmetry rules l)--n) need not be postulated, as they already follow from them, see \Cref{rem:symmetry-derived}.

Only b), f), m) and n) cost a hypothesis, and in each case what is needed is the existence of one conditional
Markov kernel, which a disintegration triple supplies uniformly but which may also be established for the kernels
at hand. For b), m) and n) the required statement is of the form $\Xk \Indep_\Kk \deltabf_\ast \given \Zk \otimes \Tk$,
which \Cref{sep:tci:r-red} shows to be \emph{equivalent} to the existence of that kernel and for which it lists
four criteria. For f), the rule on which the induction behind the
global Markov property \Cref{thm-gmp-mI-CBN} rests, the full list is the following, and its last case shows that
the list leaves nothing out.

\begin{Lem}[Left Weak Union, criteria for the hypothesis]
    \label{lem:l-uni-criteria}
    Assume $\Xk \otimes \Uk \Indep_\Kk \Yk \given \Zk$. Then:
    \[ \Xk \Indep_\Kk \Yk \given \Uk \otimes \Zk \]
    holds in each of the following cases:
{\setlength{\itemsep}{1pt}\setlength{\parskip}{0pt}
    \begin{enumerate}
        \item $(\Xcal,\Ucal,\Zcal)$ is a disintegration triple, see \Cref{def:disintegration-triple};
        \item $\Xk \ismapof_\Kk \Uk \otimes \Zk$;
        \item $\Uk \ismapof_\Kk \Zk$;
        \item there is a probability measure $\Pk(T)$ on $\Tcal$ such that the mixture
            $\Pk(X,U,Z) := \Kk(X,U,Z|T) \circ \Pk(T)$ admits a conditional Markov kernel $\Pk(X|U,Z)$ and has a
            $\Zcal$-marginal $\Pk(Z)$ with $\Kk(Z|T=t) \ll \Pk(Z)$ for all $t \in \Tcal$, which is automatic for
            $\Tcal = \Asterisk$;
        \item $\Kk(X,U,Z|T)$ admits a conditional Markov kernel $\Kk(X|U,Z,T)$ and
            $\Tk \ismapof_\Kk \Uk \otimes \Zk$;
        \item $\Xk \Indep_\Kk \deltabf_\ast \given \Uk \otimes \Zk$.
    \end{enumerate}}
    Only case 1.\ is a condition on the underlying measurable spaces; cases 2.--6.\ hold on arbitrary ones. Case
    6.\ is moreover \emph{necessary}, so that, given the premise, it is equivalent to the conclusion and
    1.--5.\ are criteria for it. The proof is given in \Cref{sep:tci:l-uni} in
    \Cref{sec:proofs-separoid-rules-for-transitional-conditional-independence}.
\end{Lem}

The following definition is stated so that it applies verbatim to the graphical relation of
\Cref{sec:graph-theory} as well; the general mechanism producing such rules is explained in
\Cref{sec:sym-sep-asym-sep}.

\begin{Def}[$\tau$-$\kappa$-separoid]
    \label{def:t-k-separoid}
    Let $\Omega$ be a class equipped with an associative and commutative operation $\lor$ (up to a fixed notion of
    isomorphism $\cong$), a neutral element $\varnothing$ and a transitive relation $\ll$ that is compatible with
    $\cong$ and satisfies \emph{product extension}, i.e.\ $\alpha \ll \beta \implies \alpha \ll \beta \lor \gamma$,
    and write $\alpha \approx \beta$ for
    $\lp \alpha \ll \beta \rp \land \lp \beta \ll \alpha \rp$.
    Let $\tau, \kappa \in \Omega$ be two
    distinguished elements such that $\ll$ is reflexive at $\tau$ and $\tau \lor \tau \approx \tau$.
    An (in general asymmetric) ternary relation $\Indep$ on $\Omega$ is called a
    \emph{$\tau$-$\kappa$-separoid} if it is invariant under $\cong$ and under $\approx$ in all three arguments
    and satisfies the eleven rules a)--k) of \Cref{thm:separoid_axioms-tci}, where throughout $\Tk$ is replaced by $\tau$,
    $\deltabf_\ast$ by $\kappa$, $\otimes$ by $\lor$ and $\ismapof_\Kk$ by $\ll$.
    We write \emph{$T$-$\ast$-separoid} for a $\Tk$-$\deltabf_\ast$-separoid; and for a CDAG $\Gk=(J,V,E)$ we call
    the corresponding structure on the subsets of $J \dcup V$ (with $\lor = \cup$,
    $\varnothing = \kappa = \emptyset$, $\ll \,=\, \ins$ and $\tau = J$) a \emph{$J$-$\emptyset$-separoid}.
    For transitional conditional independence the data are $\lor = \otimes$, $\cong$ the isomorphism of codomains,
    $\varnothing = \kappa = \deltabf_\ast$, $\ll \,=\, \ismapof_\Kk$ and $\tau = \Tk$.
    The disintegration triple hypotheses attached to b) and f) in \Cref{thm:separoid_axioms-tci} are hypotheses on
    \emph{measurable spaces} and are simply void in the abstract setting.
\end{Def}

\begin{Rem}[The symmetry rules are derived, not imposed]
    \label{rem:symmetry-derived}
    The three symmetry rules l)--n) of \Cref{cor:sep-ci:symmetry} are deliberately \emph{not} part of
    \Cref{def:t-k-separoid}: they are consequences of a)--k).
    Indeed, the premise $\beta \Indep \kappa \given \gamma$ of l) becomes $\beta \Indep \varnothing \given \gamma$
    by d) Right Decomposition, since $\kappa \cong \kappa \lor \varnothing$, and k) Flipped Left Cross
    Contraction with $\varnothing$ in its fourth slot then yields l) Restricted Symmetry, using
    $\varnothing \lor \gamma \cong \gamma$ and $\alpha \lor \varnothing \cong \alpha$;
    m) $\tau$-Restricted Symmetry is l) applied with $\gamma \lor \tau$ in place of $\gamma$, its
    second premise being b) $\tau$-Restricted Right Redundancy; and n) Symmetry is the case
    $\tau \cong \varnothing$ of m), where $\gamma \lor \tau \cong \gamma$ by neutrality; for transitional
    conditional independence this is $\Tcal = \Asterisk$, for id-separation it is $J = \emptyset$.
    Beyond the rules the derivations use only the associativity, commutativity and neutrality of $\lor$ up to
    $\cong$ and the invariance of $\Indep$; they are written out abstractly in \Cref{thm:plus-t-k-separoid} and
    are exactly the ones used in \Cref{cor:sep-ci:symmetry} for transitional conditional independence itself.
    They also account for the hypotheses in \Cref{cor:sep-ci:symmetry}: l) is unconditional because k) and d) are,
    while m) inherits the disintegration triple hypothesis of b), and n) inherits it from m).
\end{Rem}

\begin{Rem}
    \label{rem:t-k-separoid-instances}
    Two comments on the hypotheses of \Cref{def:t-k-separoid}.
    \begin{enumerate}
        \item The two conditions on the distinguished element $\tau$, reflexivity of $\ll$ at $\tau$ and
            $\tau \lor \tau \approx \tau$, are automatically satisfied in both instances and impose \emph{no}
            restriction on the class $\Omega$. For transitional conditional independence
            $\tau = \Tk = \deltabf(T|W,T)$ is by construction the Dirac kernel of the canonical projection and
            hence \emph{deterministic}, so $\Tk \ismapof_\Kk \Tk$ holds by \Cref{lem:restricted-reflexivity} and
            $\Tk \otimes \Tk \approx_\Kk \Tk$ by \Cref{lem:idempotency}, on every transition probability space
            and whatever the ambient class; and for id-separation $\tau = J$ with $J \cup J = J$.
            Note that $\Tk \otimes \Tk \cong \Tk$ would be \emph{false}: the codomains
            $\Tcal \times \Tcal$ and $\Tcal$ need not be measurably isomorphic. So the coarser equivalence
            $\approx$ is the right one here, and it is for this reason that $\approx$-invariance is part of the
            definition; for transitional conditional independence it is Full Equivalent Exchange
            \ref{sep:tci:full-eq-ex}, which again holds for arbitrary transitional random variables.
        \item Accordingly, transitional conditional independence forms a $T$-$\ast$-separoid on the class of
            \emph{all} transitional random variables whose codomains form disintegration triples, e.g.\ all of
            those with standard codomains; see \Cref{cor:t-star-separoid}.
            Determinism enters only in the finer statement that this class is a bounded join-semi-lattice: the
            relation $\ismapof_\Kk$ is \emph{not} reflexive on genuinely stochastic transitional random
            variables, see \Cref{rem:ismapof-properties} item \ref{rem:no-reflexivity}, so $\approx_\Kk$ is an equivalence relation only on the
            deterministic ones, see \Cref{lem:restricted-reflexivity}.
            The global Markov property of \Cref{sec:applications-causal_models} uses only the rules a)--k), whose
            individual hypotheses are recorded in \Cref{tab:rules-hypotheses}, and therefore does not need the full
            separoid structure; in particular the input spaces $\Xcal_j$, $j \in J$, may be arbitrary there.
    \end{enumerate}
\end{Rem}

\begin{table}[ht]
    \centering
    \footnotesize
    \setlength{\tabcolsep}{5pt}
    \renewcommand{\arraystretch}{1.15}
    \begin{tabular}{@{}lllll@{}}
        \hline
        & rule & transitional & hypothesis & graphical \\
        \hline
        a) & Extended Left Redundancy       & \ref{sep:tci:ext-l-red} & none                                & \ref{sep:sig:l-red} \\
        b) & $\Tk$-Restr.\ Right Redundancy & \ref{sep:tci:r-red}     & $(\Xcal,\Zcal,\Tcal)$ a d.t.        & \ref{sep:sig:r-red} \\
        c) & Left Decomposition             & \ref{sep:tci:l-dec}     & none                                & \ref{sep:sig:l-dec} \\
        d) & Right Decomposition            & \ref{sep:tci:r-dec}     & none                                & \ref{sep:sig:r-dec} \\
        e) & $\Tk$-Inverted Right Decomp.   & \ref{sep:tci:inv-r-dec} & none                                & \ref{sep:sig:inv-r-dec} \\
        f) & Left Weak Union                & \ref{sep:tci:l-uni}     & $(\Xcal,\Ucal,\Zcal)$ a d.t.        & \ref{sep:sig:l-uni} \\
        g) & Right Weak Union               & \ref{sep:tci:r-uni}     & none                                & \ref{sep:sig:r-uni} \\
        h) & Left Contraction               & \ref{sep:tci:l-con}     & none                                & \ref{sep:sig:l-con} \\
        i) & Right Contraction              & \ref{sep:tci:r-con}     & none                                & \ref{sep:sig:r-con} \\
        j) & Right Cross Contraction        & \ref{sep:tci:rc-con}    & none                                & \ref{sep:sig:rc-con} \\
        k) & Flipped Left Cross Contraction & \ref{sep:tci:flc-con}   & none                                & \ref{sep:sig:flc-con} \\
        \hline
        l) & Restricted Symmetry            & \ref{sep:tci:res-sym}   & none                                & \ref{sep:sig:res-sym} \\
        m) & $\Tk$-Restricted Symmetry      & \ref{sep:tci:t-res-sym} & $(\Ycal,\Zcal,\Tcal)$ a d.t.        & \ref{sep:sig:j-res-sym} \\
        n) & Symmetry                       & \ref{sep:tci:sym}       & $\Tcal=\Asterisk$, $(\Ycal,\Zcal,\Asterisk)$ a d.t. & \ref{sep:sig:sym} \\
        \hline
    \end{tabular}
    \caption{The fourteen rules, their cost and their graphical counterparts. Each row is one rule, in its
    transitional form for $\Indep_\Kk$, see \Cref{thm:separoid_axioms-tci} and \Cref{cor:sep-ci:symmetry}, and in
    its graphical form for id-separation $\Perp^{\id}_\Gk$, see \Cref{d-sep-separoid-axioms} and
    \Cref{rem:d-sep-symmetry}, with $\Tk$ read as $J$, $\otimes$ as $\cup$ and $\ismapof_\Kk$ as $\ins$; the
    numbers point at the proofs in the appendices. ``d.t.'' abbreviates \emph{disintegration triple}, see
    \Cref{def:disintegration-triple}; by \Cref{thm-regular-conditional-Markov-kernel} it suffices, for instance, that the
    first space is standard and the second one countably generated. Ten of the fourteen rules hold on
    \emph{arbitrary} measurable spaces; the exceptions are b), f), m) and the special case n) of m).
    Both columns of rules are instances of one abstract structure, the $\tau$-$\kappa$-separoid of
    \Cref{def:t-k-separoid}; for the graphical column this is moreover a formal consequence of
    \Cref{thm:plus-t-k-separoid}, whereas the transitional column is proved directly in
    \Cref{sec:proofs-separoid-rules-for-transitional-conditional-independence}.}
    \label{tab:rules-hypotheses}
\end{table}

\section{Applications to Statistical Theory}
\label{sec:applications-statistics}

In the following we will collect some illustrative applications of transitional conditional independence.

\subsection{Ancillarity, Sufficiency, Adequacy}
\label{sec:ancillarity-sufficiency}

In this subsection we want to relate the concepts of ancillarity, sufficiency and adequacy, 
see \cite{Fis22, Fis25, Bas59, Bas64, Daw75, Gho10}, to transitional conditional independence.

\begin{Eg}[Certain statistics expressed as conditional independence] Let $\Pk(W|\Theta)$ be a statistical model, considered as a Markov kernel $\Thetacal \dshto \Wcal$. Let $X$ and $Y$ be two transitional random variables w.r.t.\ $\Pk(W|\Theta)$.
    A \emph{statistic} of $X$ is a measurable map $S:\,\Xcal \to \Scal$, which we consider as the transitional random variable $S\ismapof X$ given via:
    \[ S:\,\Wcal \times \Thetacal \to \Scal,\quad (w,\theta) \mapsto S(X(w,\theta)). \]
\begin{enumerate}
\item \emph{Ancillarity}. $S$ is an \emph{ancillary statistic} of $X$ w.r.t.\ $\Theta$ if and only if:
\[ S \Indep_{\Pk(W|\Theta)} \Theta.\]
This means that every parameter $\Theta=\theta$ induces the same distribution for $S$: 
\[\Pk(S|\Theta=\theta) = \Pk(S|\cancel{\Theta}).\]
\item \emph{Sufficiency}. $S$ is a \emph{sufficient statistic} of $X$ w.r.t.\ $\Theta$ if and only if:
\[ X \Indep_{\Pk(W|\Theta)} \Theta \given S.\]
This means that there is a Markov kernel $\Pk(X|S,\cancel{\Theta})$, not dependent on $\Theta$, such that:
\[ \Pk(X,S|\Theta) = \Pk(X|S,\cancel{\Theta}) \otimes \Pk(S|\Theta).\]
So $X$ only ``interacts'' with the parameters $\Theta$ through $S$.
\item \emph{Adequacy}. $S$ is an \emph{adequate statistic} of $X$ for $Y$ w.r.t.\ $\Theta$ if and only if:
\[ X \Indep_{\Pk(W|\Theta)} \Theta,Y \given S.\]
This means we have a factorization:
\[ \Pk(X,Y,S|\Theta) = \Pk(X|\cancel{Y},S,\cancel{\Theta}) \otimes \Pk(Y,S|\Theta),\]
for some Markov kernel $\Pk(X|\cancel{Y},S,\cancel{\Theta})$, only dependent on $S$.
This means that all information of $X$ about the (parameters and/or) labels $Y$ are fully captured already by $S$.
\end{enumerate}
Two words on the status of these equivalences. Ancillarity is an equivalence on \emph{arbitrary} measurable
spaces: the witnessing object is a Markov kernel out of $\Asterisk$, i.e.\ a probability measure, so no
measurability question can arise. Sufficiency and adequacy are equivalences with the classical notions read in the
\emph{kernel} sense, i.e.\ as the existence of the Markov kernel displayed above. On a standard $\Xcal$ this
agrees with the classical, per-event formulation, which for sufficiency asks only that
$\E_\theta\lB \I_A(X) \given S \rB$ admit a $\theta$-free version for each $A$ separately, see \cite{HS49}, and
for adequacy that $\E_\theta\lB \I_A(X) \given Y,S \rB$ admit a version that is both $\theta$-free and a function
of $S$ alone. On a general $\Xcal$ the kernel version is strictly stronger, since it must produce one countably
additive kernel rather than one version per event, see the corresponding remark after \Cref{cor:basu-classical}.
\end{Eg}

Now we want to show that the classical Fisher-Neyman factorization criterion for  sufficiency (see \cite{Fis22, Ney35, HS49, Bur61})
is in line with the reformulation of sufficiency as a transitional conditional independence. The proof is given in \Cref{thm:supp:fisher-neyman} in \Cref{sec:supp:appl-statistics}. 

\begin{Thm}[Fisher-Neyman]
    \label{thm:fisher-neyman-main}
    Let $\Xcal$, $\Scal$, $\Thetacal$ be measurable spaces with $\Xcal$ standard. 
    Let $\mubf$ be a $\sigma$-finite measure on $\Xcal$ and $S:\, \Xcal \to \Scal$ a measurable map. 
    Let $\Pk(X|\Theta):\, \Thetacal \dshto \Xcal$ be a statistical model that is absolutely continuous w.r.t.\ $\mubf$: $\Pk(X|\Theta) \ll \mubf$. Then the following two statements are equivalent:
    \begin{enumerate}
        \item $\Pk(X|\Theta)$ has a Radon-Nikodym derivative\footnote{\label{fn:fisher-neyman-joint}It is not necessary to assume joint measurability for the equivalence to hold; the proof \emph{produces} versions of $p_\theta$, $g_\theta$ and $f$ for which the maps $(x,\theta) \mapsto p_\theta(x)$, $(s,\theta) \mapsto g_\theta(s)$ and $x \mapsto f(x)$ are jointly measurable, so that one may always assume this w.l.o.g.; see the last paragraph of the proof of \Cref{thm:supp:fisher-neyman}. This joint measurability is what the likelihood principle, \Cref{thm:proper-likelihood-principle}, uses.} $p_\theta$ w.r.t.\ $\mubf$ of the form:
            \begin{align}  
             p_\theta(x) = h(x) \cdot g_\theta(S(x)), \label{eq:fisher-neyman}
            \end{align} 
            with measurable maps $h:\, \Xcal \to \R_{\ge 0}$ and $g_\theta:\, \Scal \to \R_{\ge0}$ for $\theta \in \Thetacal$.
        \item $S$ is a sufficient statistic for $\Pk(X|\Theta)$, i.e.\ we have the transitional conditional independence: \begin{align*}  X \Indep_{\Pk(X|\Theta)} \Theta \given S. \end{align*} 
    \end{enumerate}
\end{Thm}

Note that \Cref{thm:fisher-neyman-main} requires the existence of a Radon-Nikodym derivative w.r.t.\ a reference measure.
Our definition of conditional independence generalizes the factorization theorem to Markov kernels (per definition) without the necessity of densities and/or reference measures.

\subsection{Basu's Theorem as a Rule for Transitional Conditional Independence}
\label{sec:basu}

Once ancillarity and sufficiency are conditional independence statements, Basu's theorem, see \cite{Bas55, Bas58},
becomes a \emph{rule}: it takes two transitional conditional independences as input and returns a third. It is not
a rule of the separoid calculus, since it needs one genuinely statistical hypothesis in addition, namely bounded
completeness. The theorem is \cite{Bas55}, stated there for a \emph{complete} sufficient statistic; bounded
completeness is the standard weakening of that hypothesis under which the conclusion still holds, and \cite{Bas58}
is the sequel on partial converses; see also \cite{KT75, Leh81}.

\begin{Def}[Boundedly complete given $\Zk$ w.r.t.\ $\Kk$]
    \label{def:boundedly-complete}
    Let $\Uk$ and $\Zk$ be transitional random variables on a transition probability space
    $\lp \Wcal \times \Tcal, \Kk(W|T) \rp$. We call $\Uk$ \emph{boundedly complete given $\Zk$ w.r.t.\
    $\Kk = \Kk(W|T)$} if for every bounded measurable map $g:\, \Ucal \times \Zcal \to \R$ we have the
    implication:
    \begin{align*}
        &\lp \forall t \in \Tcal\; \forall C \in \Bcal_\Zcal:\;
            \int_{\Ucal \times C} g(u,z)\, \Kk(du,dz|T=t) = 0 \rp\\
        \implies \quad &\lp \forall t \in \Tcal:\; g = 0 \;\; \Kk(U,Z|T=t)\text{-a.s.} \rp.
    \end{align*}
\end{Def}

\begin{Rem}[Reading the definition]
    \label{rem:boundedly-complete-reading}
    For $\Zk = \deltabf_\ast$ the sets $C$ are $\emptyset$ and $\lC \ast \rC$ and this is the classical notion:
    $\Uk$ is \emph{boundedly complete w.r.t.\ $\Kk$} if the only bounded $g$ with
    $\E\lB g(U) \given T=t \rB = 0$ for all $t$ is the one that vanishes almost surely. Whenever a conditional
    Markov kernel $\Kk(U|Z,T)$ exists the condition takes the familiar form: if
    $\E\lB g(U,Z) \given Z=z, T=t \rB = 0$ for every $t$ and $\Kk(Z|T=t)$-almost every $z$, then $g=0$ holds
    $\Kk(U,Z|T=t)$-almost surely for every $t$. The two ``almost every'' clauses cannot be dropped: a pointwise
    reading in $z$ would be a strictly weaker hypothesis, under which \Cref{thm:basu} fails.

    Apart from the domination assumption of \Cref{thm:fisher-neyman-main}, this is the only hypothesis in this
    paper that constrains the transition probabilities themselves rather than the measurable spaces involved. It
    is not a separoid condition, and it is not preserved by $\otimes$: if $\Uk' \ismapof_\Kk \Uk \otimes \Zk$ and
    $\Uk$ is boundedly complete given $\Zk$ w.r.t.\ $\Kk$ then so is $\Uk'$, but refinements need not be:
$X_1+X_2$ is boundedly complete for $\Ncal(\theta,1)^{\otimes 2}$ and $(X_1,X_2)$ is not. Note also that, like
    everything else here, the condition is stated for \emph{every} $t$ separately and needs no distribution on
    $\Tcal$; and that it refers only to the joint $\Kk(U,Z|T)$, so no disintegration is presupposed anywhere.
\end{Rem}

\begin{Thm}[Conditional Basu, as a rule for transitional conditional independence]
    \label{thm:basu}
    Let $\Xk$, $\Yk$, $\Uk$, $\Zk$ be transitional random variables on a transition probability space
    $\lp \Wcal \times \Tcal, \Kk(W|T) \rp$ such that:
{\setlength{\itemsep}{1pt}\setlength{\parskip}{0pt}
    \begin{enumerate}
        \item $\displaystyle \Xk \Indep_\Kk \Yk \given \Uk \otimes \Zk$;
        \item $\displaystyle \Xk \Indep_\Kk \Yk \given \Zk$;
        \item $\Uk$ is \emph{boundedly complete given $\Zk$ w.r.t.\ $\Kk$}, see \Cref{def:boundedly-complete}.
    \end{enumerate}}
    Then we have:
    \[ \Xk \Indep_\Kk \Yk \otimes \Uk \given \Zk. \]
    No assumption on the measurable spaces is needed; in particular no disintegration triple is required, since
    the only conditional Markov kernels in the argument are the ones asserted by 1.\ and 2.
    The proof is given in \Cref{thm:supp:basu} in \Cref{sec:supp:appl-statistics}.
\end{Thm}

\begin{Rem}[What bounded completeness buys]
    \label{rem:basu-buys}
    Both hypotheses 1.\ and 2.\ are \emph{necessary}: hypothesis 1.\ follows from the conclusion by Right Weak
    Union \ref{sep:tci:r-uni}, and hypothesis 2.\ by Right Decomposition \ref{sep:tci:r-dec} after commuting the
    two right factors, see \Cref{lem:separoid-compatibility-1} points 2.\ and 6. What is missing without 3.\ is exactly one
    step. Given 1., Right Contraction \ref{sep:tci:r-con} would already give the conclusion from
    $\Xk \Indep_\Kk \Uk \given \Zk$; so the content of the theorem is that bounded completeness lets one replace
    $\Yk$ by $\Uk$ in the second slot of 2. That is precisely the Basu move, from independence of the
    \emph{parameter} to independence of the \emph{statistic}, and it is a substitute for a Composition-type
    step, which the calculus does not have, see \Cref{thm:separoid_axioms-tci}, bought with a statistical
    hypothesis rather than a separoid one.
    The proof uses hypothesis 2.\ only through the weaker statement $\Xk \Indep_\Kk \deltabf_\ast \given \Zk$, so
    the theorem holds in that sharper form as well; we state it with 2.\ because that is the shape in which it is
    applied and because it makes the necessity above visible.
\end{Rem}

\begin{Cor}[Basu]
    \label{cor:basu-classical}
    Let $\Xcal$ and $\Thetacal$ be measurable spaces, let $\Pk(X|\Theta):\, \Thetacal \dshto \Xcal$ be a
    statistical model and let $R$ and $S$ be statistics on $\Xcal$, with arbitrary measurable codomains, such
    that:
{\setlength{\itemsep}{1pt}\setlength{\parskip}{0pt}
    \begin{enumerate}
        \item $S$ is \emph{sufficient}: $\displaystyle X \Indep_{\Pk(X|\Theta)} \Theta \given S$;
        \item $R$ is \emph{ancillary}: $\displaystyle R \Indep_{\Pk(X|\Theta)} \Theta$;
        \item $S$ is \emph{boundedly complete} w.r.t.\ $\Pk(X|\Theta)$, i.e.\ boundedly complete given
            $\deltabf_\ast$ in the sense of \Cref{def:boundedly-complete}.
    \end{enumerate}}
    Then we have:
    \[ R \Indep_{\Pk(X|\Theta)} \Theta, S. \]
\begin{proof}
    Apply \Cref{thm:basu} with $R$, $\Theta$, $S$ and $\deltabf_\ast$ in the roles of $\Xk$, $\Yk$, $\Uk$ and
    $\Zk$; the three hypotheses then correspond in the order listed. Hypothesis 1.\ there is
    $R \Indep_{\Pk} \Theta \given S$, which follows from sufficiency: $R$ is a map of $X$, hence
    $R \ismapof_\Pk X \otimes S$, so Inverted Left Decomposition \ref{sep:tci:inv-l-dec} turns
    $X \Indep_\Pk \Theta \given S$ into $X \otimes R \Indep_\Pk \Theta \given S$ and Left Decomposition
    \ref{sep:tci:l-dec} drops $X$. Hypothesis 2.\ there is ancillarity, since
    $\Zk = \deltabf_\ast$, and 3.\ is bounded completeness.
\end{proof}
\end{Cor}

\noindent
As in \Cref{thm:basu}, no assumption on the measurable spaces is needed, though \Cref{cor:basu-classical} is not
a generalization of the classical statement: hypothesis 2.\ already \emph{asserts} the conditional Markov kernel
$\Qk(X|S)$, which on a non-standard $\Xcal$ is more than classical sufficiency of $\sigma(S)$ demands.

\begin{Rem}[What the conclusion says]
    \label{rem:basu-reading}
    Unfolding \Cref{def:transitional_conditional_independence}, the conclusion of \Cref{thm:basu} asserts one
    single Markov kernel $\Qk(X|Z)$ with
    \[ \Kk(X,Y,U,Z|T) = \Qk(X|Z) \otimes \Kk(Y,U,Z|T), \]
    i.e.\ the same kernel for every value of the input. In the situation of \Cref{cor:basu-classical}, where
    $\Zk = \deltabf_\ast$, it is a single probability measure $\Qk(R)$ with
    \[ \Pk(R,S|\Theta=\theta) = \Qk(R) \otimes \Pk(S|\Theta=\theta) \qquad \text{ for every } \theta \in \Thetacal, \]
    i.e.\ $R$ and $S$ are independent under every $\Pk(R,S|\Theta=\theta)$ \emph{and} the law of $R$ is the same
    one for all of them. The classical conclusion of Basu's theorem and the ancillarity of $R$ are thus packed into
    a single relation. Two further readings are available. By \Cref{add-T} the conclusion may equally be spelled
    \[ R \Indep_{\Pk(X|\Theta)} S, \]
    the input variable being carried along by the relation whether it is written or not; and in the terminology
    of the previous subsection it says that the trivial statistic is \emph{adequate} for $R$ with respect to $S$.
    It is worth noting that this is exactly the shape of statement on which a disjunctive symmetrization goes
    wrong: by \Cref{eg:symmetrization-loses} 1.\ the relation $R \Indep^\lor_\Pk S$ can hold while
    $R \Indep_\Pk S$ fails, because a symmetric relation cannot record \emph{which} of $R$ and $S$ is the
    ancillary one, and that is the whole content of Basu's theorem.
\end{Rem}

\begin{Eg}[Sample mean and sample variance]
    \label{eg:basu-gauss}
    Let $\Thetacal = \R \times \R_{>0}$ with coordinates $\mu$ and $\sigma^2$, and let
    $\Pk(X|\Theta) = \Ncal(\mu,\sigma^2)^{\otimes n}$ be the model of $n \geq 2$ independent observations
    $X = (X_1,\dots,X_n)$ on $\Xcal = \R^n$. Put:
    \[ \bar X := \frac1n \sum_{i=1}^n X_i, \qquad\qquad
       V := \frac{1}{n-1}\sum_{i=1}^n \lp X_i - \bar X \rp^2, \]
    and apply \Cref{thm:basu} with $V$, $\Theta$, $\bar X$ and $\sigma^2$ in the roles of $\Xk$, $\Yk$, $\Uk$ and
    $\Zk$, where $\mu$ and $\sigma^2$ also denote the two coordinate projections of $\Theta$, deterministic
    transitional random variables. The three hypotheses hold.
    \emph{1.\ Sufficiency.} $\bar X$ is sufficient for $\mu$ once $\sigma^2$ is conditioned on, i.e.\
    $X \Indep_\Pk \Theta \given \bar X, \sigma^2$: the conditional law of $X$ given $\bar X = m$ is that of
    $m \mathbf{1} + \varepsilon$ with $\varepsilon \sim \Ncal\lp 0, \sigma^2(I - \tfrac1n \mathbf{1}\mathbf{1}^\top)\rp$,
    which depends only on $(m,\sigma^2)$ and measurably so. Since $V$ is a map of $X$, the argument of
    \Cref{cor:basu-classical} gives $V \Indep_\Pk \Theta \given \bar X, \sigma^2$.
    \emph{2.\ Conditional ancillarity.} $\tfrac{n-1}{\sigma^2}V \sim \chi^2_{n-1}$, so the law of $V$ depends on
    $\Theta$ only through $\sigma^2$, and $\sigma^2 \mapsto \tfrac{\sigma^2}{n-1}\chi^2_{n-1}$ is measurable; so
    it really is a Markov kernel $\Qk(V|\sigma^2)$ and $V \Indep_\Pk \Theta \given \sigma^2$.
    \emph{3.\ Conditional bounded completeness.} Because $\sigma^2$ is a deterministic function of $\Theta$, we
    have $\Pk(\bar X, \sigma^2|\Theta=(\mu,v)) = \Ncal(\mu, v/n) \otimes \deltabf_v$, so choosing $C \ni v$ in
    \Cref{def:boundedly-complete} collapses the condition to $\int g(m,v)\, \Ncal(dm;\mu,v/n) = 0$ for all $\mu$,
    one $v$ at a time; and $\lp \Ncal(\mu,v/n) \rp_{\mu \in \R}$ is a complete exponential family for each $v$.
    (In general, completeness for each $z$ separately does \emph{not} give \Cref{def:boundedly-complete}; here it
    does, through this collapse.)

    \Cref{thm:basu} therefore gives:
    \[ V \;\Indep_{\Pk(X|\Theta)}\; \mu, \bar X \given \sigma^2, \]
    where we have written $\mu$ rather than $\Theta$ on the right: the two are interchangeable there, since
    $\Theta \approx_\Pk \mu \otimes \sigma^2$ and $\sigma^2$ already occupies the conditioning slot, see Full
    Equivalent Exchange \ref{sep:tci:full-eq-ex}, $\Tk$-Inverted Right Decomposition \ref{sep:tci:inv-r-dec} and
    Right Decomposition \ref{sep:tci:r-dec}. Spelled out, the relation says that there is \emph{one} Markov kernel
    $\Qk(V|\sigma^2)$, namely the law of $\tfrac{\sigma^2}{n-1}\chi^2_{n-1}$, with
    \[ \Pk\lp V, \bar X \given \Theta = (\mu,\sigma^2) \rp
       = \Qk(V|\sigma^2) \otimes \Pk\lp \bar X \given \Theta=(\mu,\sigma^2)\rp
       \qquad \text{ for every } (\mu,\sigma^2). \]
    So the sample mean and the sample variance are independent, for every value of both parameters, with the law
    of the sample variance depending on $\sigma^2$ alone: the classical statement, obtained without ever
    putting a distribution on $\Thetacal$ and without declaring $\sigma^2$ known. Had we fixed $\sigma^2$ instead,
    i.e.\ taken $\Thetacal = \R$ and $\Zk = \deltabf_\ast$, \Cref{cor:basu-classical} would give the same
    conclusion for that one sub-model; the point of the conditioning slot is that the single kernel
    $\Qk(V|\sigma^2)$ serves all of them at once. Which of the two statements one gets is decided by the slot
    $\sigma^2$ occupies, in the spirit of \Cref{rem:Q-no-T}.
\end{Eg}

\subsection{Comparison of Experiments}
\label{sec:blackwell}

Sufficiency compares a statistic with the data it is computed from. Blackwell's \emph{comparison of experiments},
see \cite{Bla51, Bla53, LeC64, Tor91}, compares two \emph{different} experiments on the same parameter, and its
definition has exactly the shape transitional conditional independence was built for: it asks for one Markov kernel
that serves every parameter value at once.

\begin{Def}[Experiments and garblings]
    \label{def:experiment}
    An \emph{experiment} with parameter space $\Thetacal$ is a Markov kernel
    $\Pk(X|\Theta):\, \Thetacal \dshto \Xcal$. Given two of them, $\Ek_1 = \Pk_1(X_1|\Theta)$ and
    $\Ek_2 = \Pk_2(X_2|\Theta)$, we say that $\Ek_1$ is \emph{at least as informative} as $\Ek_2$, in symbols
    $\Ek_1 \succeq \Ek_2$, if there is a Markov kernel (a \emph{garbling})
$\Qk(X_2|X_1):\, \Xcal_1 \dshto \Xcal_2$, the \emph{same} for all parameters, with:
    \[ \Pk_2(X_2|\Theta=\theta) = \Qk(X_2|X_1) \circ \Pk_1(X_1|\Theta=\theta)
       \qquad \text{ for every } \theta \in \Thetacal. \]
\end{Def}

\begin{Thm}[Blackwell's order is a transitional conditional independence]
    \label{thm:blackwell}
    Let $\Ek_1$ and $\Ek_2$ be experiments with parameter space $\Thetacal$. Then the following are equivalent:
    \begin{enumerate}
        \item $\Ek_1 \succeq \Ek_2$;
        \item there is a Markov kernel $\Kk(X_1,X_2|\Theta):\, \Thetacal \dshto \Xcal_1 \times \Xcal_2$ with
            marginals $\Pk_1(X_1|\Theta)$ and $\Pk_2(X_2|\Theta)$ such that, on the transition probability space
            $\lp \lp \Xcal_1 \times \Xcal_2 \rp \times \Thetacal, \Kk(X_1,X_2|\Theta) \rp$:
            \[ X_2 \Indep_{\Kk(X_1,X_2|\Theta)} \Theta \given X_1. \]
    \end{enumerate}
    No assumption on the measurable spaces $\Xcal_1$, $\Xcal_2$, $\Thetacal$ is needed.
    The proof is given in \Cref{thm:supp:blackwell} in \Cref{sec:supp:appl-statistics}.
\end{Thm}

\begin{Rem}
    \label{rem:blackwell}
    \begin{enumerate}
        \item The quantifier is the whole point. \Cref{def:experiment} asks for one garbling valid at every
            $\theta$, and \Cref{def:transitional_conditional_independence} asks for one kernel $\Qk(X_2|X_1)$ valid
            at every input; the two existential quantifiers are the same one. Classically the same statement is
            reached only indirectly, by quantifying over \emph{all} priors and all bounded loss functions and
            comparing Bayes risks, an equivalence that itself needs hypotheses, a finite parameter space in
            \cite{Bla53} and regularity conditions in general, see \cite{Tor91}. That detour exists precisely
            because no single prior on $\Thetacal$ is available, which is the first of the two failure modes of
            \Cref{sec:main:comparison:classical}.
        \item The setting of sufficiency is the special case in which one experiment arises from the other by a
            \emph{deterministic} garbling. Let $\Ek_1 = \Pk(X|\Theta)$ and $\Ek_2 = \Pk(S|\Theta)$ for a statistic
            $S:\, \Xcal \to \Scal$. Then $\Ek_1 \succeq \Ek_2$ always holds, witnessed by the joint law of $(X,S)$
            together with $S \Indep_\Pk \Theta \given X$, which is Extended Left Redundancy
            \ref{sep:tci:ext-l-red} applied to $S \ismapof X$. In the other direction, sufficiency of $S$ (i.e.\
            $X \Indep_\Pk \Theta \given S$ in the sense of \Cref{sec:ancillarity-sufficiency}) always implies
            $\Ek_2 \succeq \Ek_1$: marginalizing $S$ out of the sufficiency factorization leaves a garbling
            $\Qk(X|S)$.
            The converse of that last implication is \emph{not} automatic, and the reason is instructive.
            $\Ek_2 \succeq \Ek_1$ only asks for \emph{some} garbling, whereas sufficiency asks the \emph{canonical}
            coupling, the joint law of $(X,S)$, to factorize. If, say, $\Pk(X|\Theta)$ does not depend on $\Theta$
            at all, then $\Ek_2 \succeq \Ek_1$ holds with a constant garbling for every statistic $S$ whatsoever,
            while $X \Indep_\Pk \Theta \given S$ still demands a conditional Markov kernel $\Qk(X|S)$, which on a
            non-standard $\Xcal$ need not exist. This is the same gap between the kernel-sense and the per-event
            reading of sufficiency that \Cref{sec:ancillarity-sufficiency} and \Cref{cor:basu-classical} record; on standard
            $\Xcal$ the two directions do match up.
        \item Statement 2.\ is existentially quantified over the coupling as well, and it has to be: the two
            experiments are given by their marginals only, and Blackwell's condition is exactly the assertion that
            some joint realizing them exhibits the required factorization. The proof shows that the garbled coupling
            $\Kk(X_1,X_2|\Theta) := \Qk(X_2|X_1) \otimes \Pk_1(X_1|\Theta)$ always works, so nothing is lost by
            restricting attention to it.
        \item Blackwell's order is the \emph{exact} case of a theory whose working content is quantitative: Le
            Cam's deficiency, see \cite{LeC64, Tor91}, measures by how much $\Ek_1 \succeq \Ek_2$ fails. An
            approximate version of transitional conditional independence, in which the separoid rules degrade with
            an additive budget, would be its natural companion; see the outlook in \Cref{sec:discussion}.
    \end{enumerate}
\end{Rem}

\subsection{Invariant Prediction and Environments}
\label{sec:invariant-prediction}

Our next application is not a reformulation of a classical concept but a problem that, in its usual formulation,
has no formal meaning at all in terms of ordinary conditional independence.
In \emph{invariant prediction}, see \cite{PBM16, PBP19}, one observes a response $Y$ and covariates
$X = (X_1,\dots,X_p)$ in several \emph{environments} (different laboratories, different time
periods, different experimental conditions) and searches for subsets $S \ins \{1,\dots,p\}$ that are
\emph{invariant}: the conditional distribution of $Y$ given $X_S$ should be one and the same in every environment.
The environment is an index chosen by the experimenter, not a random draw; if it is time, or a continuum of
interventions, there is no distribution on the set of environments at hand and none should be invented.

We therefore take the environment as the non-stochastic input: let $\lp \Wcal \times \Tcal, \Kk(W|T)\rp$ be a
transition probability space whose input space $\Tcal$ is the set of environments, and let $Y$ and
$X_1,\dots,X_p$ be transitional random variables on it, so that $\Kk(Y,X|T=t)$ is the joint distribution of
$(Y,X)$ in the environment $t \in \Tcal$. To stress the reading we write $E := T$ for the deterministic
transitional random variable given by the input, i.e.\ ``the environment''.

\begin{Def}[Invariant set]
    \label{def:invariant-set}
    A subset $S \ins \{1,\dots,p\}$ is called \emph{invariant} (for $Y$) if there exists a Markov kernel
    \[ \Qk(Y|X_S):\, \Xcal_S \dshto \Ycal \]
(the \emph{same} for all environments) such that:
    \[ \Kk(Y,X_S|T) = \Qk(Y|X_S) \otimes \Kk(X_S|T). \]
    The kernel $\Qk(Y|X_S)$ of course depends on $S$: it is a kernel out of $\Xcal_S$, and a different subset
    calls for a different kernel. What ``the same'' quantifies over is the \emph{environment}, and only that. For
    each candidate $S$ the definition asks whether such a kernel exists at all, and invariant prediction is the
    search over $S$ for a subset where it does.
\end{Def}

\begin{Prp}[Invariance is transitional conditional independence]
    \label{prp:invariance-tci}
    For every $S \ins \{1,\dots,p\}$ the following statements are equivalent:
    \begin{enumerate}
        \item $S$ is invariant in the sense of \Cref{def:invariant-set};
        \item $Y \ismapof^\ast_\Kk X_S$;
        \item $Y \Indep_{\Kk(W|T)} \deltabf_\ast \given X_S$;
        \item $Y \Indep_{\Kk(W|T)} E \given X_S$.
    \end{enumerate}
\begin{proof}
    1.\ and 2.\ are the same statement, by \Cref{not:ismapof} 3.
    2.\ $\iff$ 3.: unfolding \Cref{def:transitional_conditional_independence} with $\deltabf_\ast$ in the second
    slot gives $\Kk(Y,\ast,X_S|T) = \Qk(Y|X_S)\otimes\Kk(\ast,X_S|T)$, and the one-point factor
    can be dropped on both sides.
    3.\ $\iff$ 4.\ is \Cref{add-T} together with $\deltabf_\ast \otimes \Tk \cong \Tk$.
\end{proof}
\end{Prp}

\begin{Rem}
    \label{rem:invariance-discussion}
    Three comments.
    \begin{enumerate}
        \item Statement 4.\ is the way invariance is usually written, namely as ``$Y \Indep E \given X_S$''. Under
            ordinary conditional independence that expression is only meaningful once one puts a distribution on
            $\Tcal$, and it then changes its meaning: it becomes a statement about a \emph{mixture} of environments
            rather than about every environment. \Cref{prp:invariance-tci} shows that transitional conditional
            independence gives the expression its intended meaning verbatim, with no distribution on $\Tcal$ and no
            null sets in $\Tcal$.
        \item Two of the separoid rules have a direct reading here. Left Decomposition \ref{sep:tci:l-dec} says
            that if the pair $(Y,U)$ is invariant given $X_S$ then so are $Y$ and $U$: the rule drops the first
            factor of a product, and the second one is dropped after reordering, see
            \Cref{lem:separoid-compatibility-1} 2.\ and 6.; and Left Contraction \ref{sep:tci:l-con}
            says that if $U$ is invariant given $X_S$ and $Y$ is invariant given $(X_S,U)$, then $(Y,U)$ is
            invariant given $X_S$, so invariance can be built up in stages.
        \item It is worth being explicit about what does \emph{not} follow. Invariance of $S_1$ and of $S_2$ does
            \emph{not} imply invariance of $S_1 \cap S_2$: none of the rules a)--k) of
            \Cref{thm:separoid_axioms-tci} concludes anything about the \emph{intersection} of two conditioning
            variables, and the statement already fails for the graphical
            relation of \Cref{sec:graph-theory}. For the CDAG with $J=\{e\}$, $V=\{x_1,x_2,y\}$ and edges
            $e \tuh x_2 \tuh x_1 \tuh y$ we have $y \Perp^{\id}_\Gk \emptyset \given \{x_1\}$ and
            $y \Perp^{\id}_\Gk \emptyset \given \{x_2\}$, since every walk from $y$ to $e$ passes through
            $x_1$ or $x_2$ as a non-collider, but $y \nPerp^{\id}_\Gk \emptyset \given \emptyset$, because with
            an empty conditioning set the directed walk $e \tuh x_2 \tuh x_1 \tuh y$ is open. Accordingly the identification results of \cite{PBM16}, which read off a causal
            set from the \emph{intersection of all} invariant sets, do not rest on such a closure property but on
            the assumption that the causal set is itself invariant; that is a modelling assumption and not a
            consequence of the calculus.
    \end{enumerate}
    Deciding statement 3.\ from data is exactly the testing problem that the procedures of \cite{PBM16, PBP19}
    address; see the outlook in \Cref{sec:discussion}.
\end{Rem}

\begin{Rem}[Invariant representations]
    \label{rem:invariant-representations}
    In machine learning one usually does not restrict the search to subsets $S$ of the given features. Instead one
    \emph{learns} a representation $Z = g(X)$, with $g$ ranging over a parameterized class of maps, and asks for
    \[ Y \Indep_{\Kk(W|T)} E \given Z \]
    in place of $Y \Indep_{\Kk(W|T)} E \given X_S$; see
    \cite{BCV13, CKNH20, LBL19, ABGL19, RSTP18, MOC18, MMC20, FDF20, FTF21, SLB21}.
    Nothing in \Cref{def:invariant-set} or \Cref{prp:invariance-tci} uses that $X_S$ is a coordinate projection:
    both hold verbatim for an arbitrary transitional random variable $Z$ in the conditioning slot, so the statement
    ``there is one Markov kernel $\Qk(Y|Z)$ valid in every environment'' is available for learned representations
    exactly as it is for subsets. The subset case is the special case $g = \pr_S$.
\end{Rem}

\subsection{Invariant Reductions}

\begin{Eg}[Invariant reduction]
Let $\Pk(X|\Theta)$ be a statistical Model, given as a Markov kernel.
Assume that we are only interested in a certain quantity of the parameters $\Gamma=\Gamma(\Theta)$,
considered as a measurable function $\Gamma:\, \Thetacal \to \Gamcal$ into a measurable space $\Gamcal$.
For the estimation of $\Gamma$ we then might only need parts of the information encoded in the data $X$.
An \emph{invariant reduction of $\Pk(X|\Theta)$ w.r.t.\ $\Gamma$}, see \cite{Hal65}, is then a measurable
function $U:\, \Xcal \to \Ucal$ such that $\Pk(U(X)|\Theta)$ depends on $\Theta$ only through $\Gamma$, in the
sense that there is a Markov kernel $\Qk(U|\Gamma):\, \Gamcal \dshto \Ucal$ with
$\Pk(U|\Theta=\theta) = \Qk(U|\Gamma=\Gamma(\theta))$ for every $\theta \in \Thetacal$. This is precisely the
transitional conditional independence:
\[  U \Indep_{\Pk(X|\Theta)} \Theta \given \Gamma,   \]
and the occurring Markov kernel $\Pk(U|\Gamma,\cancel{\Theta})$ then gives the correct model to further work with.
\end{Eg}

\subsection{Reparameterizing Transitional Random Variables}

We want to generalize two somewhat related well-known results, see \cite{Cen82}, from random variables to transitional random variables:

1.) Since the paper \cite{Dar53} it was developed that for a real-valued random variable $X$ that has a continuous cumulative distribution function $F$ and quantile function $R=F^{-1}$ that $E:=F(X)$ is uniformly distributed on $[0,1]$ and $R(E) = X$ a.s. 

2.) It is known that for random variables $X$ and $Z$ with a well-behaved joint distribution $\Pk(X,Z)$ there exists a random variable $E$ that is independent of $Z$ and a measurable map $g$ such that $X=g(E,Z)$ a.s. 

To establish such results for transitional random variables we will use the following constructions.

\begin{Def}
    \label{def:itcdf-tqf}
    Let $X$ be a transitional random variable with values in $\Xcal = \bar \R=[-\infty,+\infty]$ on a transition probability space 
    $\lp \Wcal \times \Zcal, \Kk(W|Z)\rp$. We then define the \emph{interpolated transitional cumulative distribution function (itcdf)} of $X$ as:
    \begin{align*}
        F(x;u|z) &:= \Kk(X < x|Z=z) + u \cdot \Kk(X=x|Z=z),
    \end{align*}
    with $u \in [0,1]$,
    and the (``transitional'') \emph{quantile function (tqf)} of $X$ as:
    \begin{align*}
        R(e|z) &:=  \inf\lC \tilde{x} \in \bar \R\,|\, F(\tilde{x};1|z) \ge e  \rC,
    \end{align*}
    for $e \in [0,1]$.
\end{Def}

\begin{Thm}
    \label{thm:reparamererizing-trv}
    Let $\lp \Wcal \times \Zcal, \Kk(W|Z)\rp$ be any transition probability space and $X$ be a transitional random variable with values in
    a standard measurable space $\Xcal$ and $\iota:\, \Xcal \inj \bar \R$ a fixed embedding onto a Borel subset of $\bar \R$ 
    (i.e.\ w.l.o.g.\ $\Xcal=\bar\R$). 
    Let $\Kk(U)$ be the uniform distribution on $\Ucal:=[0,1]$. We put $\bar\Wcal = \Ucal \times \Wcal$, $\bar W = (U,W)$ and $\Kk(\bar W|Z) = \Kk(U) \otimes \Kk(W|Z)$. We then consider the transitional random variables $X$, $U$, $Z$, $E$ on the transition probability space 
    $\lp \bar\Wcal \times \Zcal, \Kk(\bar W |Z) \rp$ where:
    \begin{align*}
        E &:=F(X;U|Z):\, \bar \Wcal \times \Zcal \to \Ecal:=[0,1],
    \end{align*}
    and $F$ is the itcdf of $X$ from \Cref{def:itcdf-tqf}.
    Then we have the transitional independence:
    \[ E \Indep_{\Kk(\bar W|Z)}Z, \qquad \text{ with } \qquad \Kk(E|\cancel{Z}) \text{ the uniform distribution on } \Ecal=[0,1],   \]
    and:
    \[ X = R(E|Z) \quad \Kk(\bar W|Z)\text{-a.s.},\]
    where $R$ is the tqf of $X$ from \Cref{def:itcdf-tqf}.
\end{Thm}

The proof of this theorem can be found in \Cref{ccdf2} in \Cref{sec:proofs-reparameterizing}.

\subsection{Propensity Score}
\label{sec:propensity}

For a random variable $X$ and binary random variable $Y \in \{0,1\}$ the \emph{propensity score} is $e(x):=\Pk(Y=1|X=x)$. It is the ``smallest'' statistic of $X$ such that $Y \Indep X \given e(X)$ (see \cite{RR83}). This is one of the core concepts of causal inference
using the potential outcome formulation.
We now claim that the above can be generalized to arbitrary (non-binary) $Y$ with Markov kernel $\Pk(Y|X)$, 
even when no distribution for $X$ is specified.

\begin{Thm}[Propensity score]
\label{propensity}
Let $\Pk(Y|X)$ be a Markov kernel.
We define the \emph{propensity} of $x \in \Xcal$ w.r.t.\ $\Pk(Y|X)$ as:
\[E(x):= \Pk(Y|X=x) \; \in \; \Pcal(\Ycal).\]
Note that the map $E: \,\Xcal \to \Ecal:=\Pcal(\Ycal)$ is measurable and $E \ismapof X$.
Now let $S: \, \Xcal \to \Scal$ be another measurable map ($S \ismapof X$).
Then we have the equivalence:
\[ Y \Indep_{\Pk(Y|X) } X \given S \qquad \iff \qquad E \ismapof S.\]
In particular, as $E \ismapof E$, we have:
\[ Y \Indep_{\Pk(Y|X) } X \given E.\]
So $E$ is in this sense the smallest statistic of $X$ such that the above conditional independence holds.
\end{Thm}

The proof is given in \Cref{thm:supp:propensity} in \Cref{sec:supp:appl-statistics}.

\subsection{Likelihood Principles}
\label{sec:likelihood}

In this section we prove the \emph{likelihood principle}, which states that the likelihood function, considered as a random variable, is a sufficient statistic of the data.
Throughout this subsection we fix, for every $\theta \in \Thetacal$, a version $p_\theta$ of the density of
$\Pk(X|\Theta=\theta)$ w.r.t.\ the reference measure $\mubf$, jointly measurable in $(x,\theta)$, which is
possible by \Cref{thm:fisher-neyman-main} and its footnote; every statement about the
likelihood function below is relative to this choice, and likewise for the measure $\nubf$ and the function
$L_\nubf$ appearing in the quasi-minimality clause.
For the history of the likelihood principle and discussions see \cite{SBC62, Fis22, Hac65, Edw74, Edw92, Roy97, Bir62, Jay03, May14, EvaM13, Gan15}. The proof can be found in \Cref{thm:supp:proper-likelihood-principle} in \Cref{sec:supp:appl-statistics}.

\begin{Thm}[The likelihood principle]
    \label{thm:proper-likelihood-principle}
    Let $\Xcal$, $\Thetacal$ be measurable spaces with $\Xcal$ standard and let $\mubf$ be a $\sigma$-finite measure on $\Xcal$.
    Consider a statistical model, written as the Markov kernel: $\Pk(X|\Theta): \, \Thetacal \dshto \Xcal$.
    For each $\theta \in \Thetacal$ assume that the Radon-Nikodym derivative $p_\theta$ exists:
    \[ p_\theta(x) := \frac{\Pk(X \in dx|\Theta=\theta)}{\mubf(dx)}(x).\]
    Then consider the \emph{likelihood function}\footnote{We endow $\Rcal$ with the smallest $\sigma$-algebra $\Bcal_\Rcal$ such that the evaluation map $\ev_\theta:\, \Rcal \to \R_{\ge0}$, $r \mapsto r(\theta)$, is measurable for every $\theta \in \Thetacal$.}:
    \[ L_\mubf:\, \Xcal \to \Rcal:=\R_{\ge 0}^\Thetacal, \qquad x \mapsto \lp \theta \mapsto p_\theta(x) \rp. \]
    Then $L_\mubf$ is measurable, the Fisher-Neyman \Cref{thm:fisher-neyman-main} applies with $p_\theta(x)=\ev_\theta(L_\mubf(x))$, and we get the transitional conditional independence:
     \[ X \Indep_{\Pk(X|\Theta) } \Theta \given L_\mubf. \]
   Furthermore, let $S$ be any other measurable map of $X$, i.e.\ $S \ismapof X$. Then we have:
   \begin{enumerate}
       \item \emph{Sufficiency}: If $L_\mubf \ismapof S$ then also: $\displaystyle X \Indep_{\Pk(X|\Theta) } \Theta \given S$.
       \item \emph{Quasi-minimality}: If $S$ satisfies: $\displaystyle X \Indep_{\Pk(X|\Theta) } \Theta \given S$, then there exists a measure $\nubf$ such that $\nubf$ has a density w.r.t.\ $\mubf$, $\Pk(X|\Theta)$ has a density w.r.t.\ $\nubf$ and the corresponding likelihood function $L_\nubf$ satisfies: $L_\nubf \ismapof S$.
   \end{enumerate}
\end{Thm}

In this sense the likelihood captures all information of the parameters $\Theta$ about the data $X$ 
and it does so most efficiently (modulo the multiplicative factor $h(x)=\frac{d\nubf}{d\mubf}(x)$)\footnote{
To remove that factor $h(x)$ to arrive at proper minimality conditions one could consider a \emph{likelihood ratio principle}, e.g.\ by replacing $p_\theta(x)$ by: $\frac{p_\theta(x)}{p_{\tilde \theta}(x)}$, $\frac{p_\theta(x)}{\sup_{\tilde \theta \sim \pi}p_{\tilde \theta}(x)}$, $\frac{p_\theta(x)}{\E_{\tilde\theta \sim \pi}[p_{\tilde \theta}(x)]}$, etc., also see \Cref{cor:bay-stat-ci}.}.

A \emph{dual version of the likelihood principle}, where the roles of data point $x$ and parameter $\theta$ are swapped in a certain sense, can be derived from \Cref{propensity} applied to a statistical model $\Pk(X|\Theta)$ as follows.

\begin{Thm}[A dual likelihood principle]
    \label{thm:likelihood-principle}
    Let $\Pk(X|\Theta)$ be a statistical model. Define the \emph{dual likelihood function} as: 
    $R(\theta):=\Pk(X|\Theta=\theta) \in \Pcal(\Xcal)$.
    We then have the transitional conditional independence:
 \[ X \Indep_{\Pk(X|\Theta) } \Theta \given R. \]
 Furthermore, for any other measurable map $S$ of the parameters $\Theta$, i.e.\ $S \ismapof \Theta$, we 
 have the equivalence:
 \[ X \Indep_{\Pk(X|\Theta) } \Theta \given S \qquad \iff \qquad R \ismapof S. \]
\end{Thm}

In this sense the dual likelihood function captures all information of the parameters $\Theta$ about the data $X$  and it does so most efficiently.

\subsection{Bayesian Statistics}

    Let $\Pk(X|\Theta)$ be a statistical model between standard measurable spaces, $\Xcal$ and $\Thetacal$, 
    and $\Pk(\Theta|\Pi)$ be a prior with hyperparameters $\Pi=\pi$.
Then by the standard Bayesian setting we have a joint (transition) probability distribution:
\[\Pk(X,\Theta|\Pi) := \Pk(X|\Theta) \otimes \Pk(\Theta|\Pi). \]
A conditional Markov kernel gives us the posterior (transitional) probability distributions:
\[\Pk(\Theta|X,\Pi),\]
which is unique up to $\Pk(X|\Pi)$-null set. 
We now define the transitional random variable: 
\[ Z(x,\pi):= \Pk(\Theta|X=x,\Pi=\pi), \] 
which gives us a joint (transition) probability distribution: $\Pk(X,\Theta,Z|\Pi)$.

The following result then formalizes the basic idea that the posterior (given via $Z$) most efficiently incorporates all 
information from the data ($X$) about the state of the
parameters ($\Theta$) as soon as a prior ($\Pi$) is specified.

\begin{Thm}[Bayesian statistics]
    \label{thm:bay-stat-ci}
With the above notations we have the conditional independence:
\[ \Theta \Indep_{\Pk(X,\Theta|\Pi) } X \given Z.\]
Now let $S$ be another deterministic measurable function in $(X,\Pi)$, i.e.\ $S \ismapof (X,\Pi)$. 
Then we have the equivalence:
\[ \Theta \Indep_{\Pk(X,\Theta|\Pi) } X \given S \qquad \iff \qquad Z \ismapof_{\Pk(X,\Theta|\Pi)} S  \]
\end{Thm}

The proof of this theorem is given in \Cref{thm:supp:bayesian-statistics} in \Cref{sec:supp:appl-statistics}.

\Cref{thm:bay-stat-ci} states both assertions with $\Theta$ on the \emph{left}. Turning them around gives the
\emph{likelihood principle for Bayesian statistics}: the posterior $Z$ is a minimal sufficient statistic for
$\Pk(X|\Theta)$. The passage uses Symmetry and therefore deserves two comments.
First, $\Tk$-Restricted Symmetry \ref{sep:tci:t-res-sym} needs its premise in the form
$\cdot \given \Zk \otimes \Tk$, here with $\Tk = \Pi$: from $\Theta \Indep X \given Z$ we obtain
$\Theta \Indep \Pi \otimes X \given Z$ by $\Tk$-Inverted Right Decomposition \ref{sep:tci:inv-r-dec} and then
$\Theta \Indep X \given Z,\Pi$ by Right Weak Union \ref{sep:tci:r-uni}.
Second, the disintegration triples have to be checked. For the first assertion the relevant triple is
$(\Xcal,\Pcal(\Thetacal),\cdot)$ with an arbitrary third component, see
\Cref{thm-regular-conditional-Markov-kernel} 1., and it is one: $\Xcal$ is standard, and $\Pcal(\Thetacal)$, the
codomain of the posterior $Z$, is a standard measurable space, hence countably generated, because
$\Thetacal$ is standard, see \cite{Kec95} 17.24 and \cite{Sch74} Appendix \S5 Thm.\ 7+8.
For the second assertion the relevant triples are $(\Xcal,\Scal,\cdot)$ and $(\Thetacal,\Scal,\cdot)$ (the two
directions of the equivalence need Symmetry in opposite directions), so there we have to require in addition that
the codomain $\Scal$ of the statistic $S$ is countably generated. This is no restriction in practice, but it cannot
be dropped from the argument.
Note that in contrast to \Cref{thm:proper-likelihood-principle}, where we used the likelihood function $L_\mubf$ w.r.t.\ some reference measure $\mubf$, the posterior $Z$ does not require a density w.r.t.\ a reference measure and also provides proper minimality conditions:

\begin{Cor}[The likelihood principle for Bayesian statistics]
    \label{cor:bay-stat-ci}
We have the transitional conditional independence:
\[ X \Indep_{\Pk(X,\Theta|\Pi) } \Theta \given Z,\Pi.\]
Furthermore, if $S$ is another measurable map in $(X,\Pi)$, i.e.\ $S \ismapof (X,\Pi)$, whose codomain
$\Scal$ is countably generated (e.g.\ standard), then for every fixed prior $\Pi=\pi$ we have the equivalence:
\[ X \Indep_{\Pk(X,\Theta|\Pi=\pi) } \Theta \given S \qquad \iff \qquad Z \ismapof_{\Pk(X,\Theta|\Pi=\pi)} S  \]
 \end{Cor}

The \emph{dual likelihood principle}, \Cref{thm:likelihood-principle}, also holds analogously:

\begin{Thm}[A dual likelihood principle for Bayesian statistics]
    We also have the transitional conditional independence with $R(\theta):=\Pk(X|\Theta=\theta)$:
    \[ X \Indep_{\Pk(X,\Theta|\Pi) } \Theta,\Pi \given R. \]
    If $\Xcal$ is countably generated and $S$ is another measurable map in $\Theta$, $S \ismapof \Theta$, 
    then we get the equivalence:
    \[X \Indep_{\Pk(X,\Theta|\Pi) } \Theta,\Pi \given S \qquad \iff \qquad R \ismapof_{\Pk(X,\Theta|\Pi)} S\] 
\end{Thm}

The proof of this theorem is given in \Cref{thm:supp:bayesian-likelihood} in \Cref{sec:supp:appl-statistics}.
\section{Applications to Graphical Models}
\label{sec:applications-causal_models}

In this section we relate transitional conditional independence to d-separation in graphs; the goal is the global
Markov property for Bayesian networks with input nodes, \Cref{thm-gmp-mI-CBN}.
We first introduce the few graph theoretic notions that this needs. Of importance, also on its own, is the notion of
\emph{d-separation}, see \cite{Pearl09, Lau90, Gei90, Ver93, Lau96, KF09, SGS2000}, which we will use in its
ordinary, \emph{symmetric} form, ignoring the input nodes entirely.
On top of it we introduce a variant that is adapted to graphs with input nodes, \emph{input-d-separation}
(\emph{id-separation} for short): it is d-separation from the set $B \cup J$ instead of from $B$, and it is
\emph{asymmetric}.
We set up id-separation in such a way that it forms a \emph{$J$-$\emptyset$-separoid},
see \Cref{def:t-k-separoid}.
The reason is that we want to match those separoid rules to the ones for transitional conditional independence
in \Cref{thm:separoid_axioms-tci}.
The gain of this two-step approach is that \emph{all} separoid rules for id-separation, including the more exotic
looking ones like Flipped Left Cross Contraction, will be \emph{derived} from the classical symmetric rules for
d-separation by purely formal manipulations, without inspecting a single walk.
The separoid calculus of this section generalizes to graphs with cycles and bi-directed edges once d-separation
is replaced by $\sigma$-separation: its symmetric rules are available in \cite{Ric03, FM17, FM18, FM20}, and the
asymmetric ones then follow by the same shift, see \Cref{thm:plus-t-k-separoid}. For clarity of exposition we
restrict ourselves to the acyclic case.

\subsection{Conditional Directed Acyclic Graphs (CDAGs)}
\label{sec:graph-theory}

\begin{figure}[ht]
    \centering
    \begin{tikzpicture}[scale=0.85, transform shape]
        \tikzstyle{every node} = [draw,shape=circle,color=blue]
        \node (v1) at (5,5) {$v_4$};
        \node (v2) at (3,3) {$v_5$};
        \node (v3) at (7,3) {$v_6$};
        \node (v7) at (11,3) {$v_7$};
        \node (v5) at (5,1) {$v_8$};
        \node (v6) at (7,-1) {$v_9$};
        \node (v8) at (9,5) {$v_3$};
        \tikzstyle{every node} = [draw,shape=rectangle,color=teal, minimum size = 0.75cm]
        \node (v4) at (1,5) {$v_1$};
        \node (v11) at (13,5) {$v_2$};
        \draw[-{Latex[length=3mm, width=2mm]}, color=teal] (v4) to (v2);
        \draw[-{Latex[length=3mm, width=2mm]}, color=teal] (v11) to (v7);
        \draw[-{Latex[length=3mm, width=2mm]}, color=blue] (v1) to (v2);
        \draw[-{Latex[length=3mm, width=2mm]}, color=blue] (v1) to (v3);
        \draw[-{Latex[length=3mm, width=2mm]}, color=blue] (v2) to (v5);
        \draw[-{Latex[length=3mm, width=2mm]}, color=blue] (v3) to (v5);
        \draw[-{Latex[length=3mm, width=2mm]}, color=blue] (v5) to (v6);
        \draw[-{Latex[length=3mm, width=2mm]}, color=blue] (v7) to (v6);
        \draw[-{Latex[length=3mm, width=2mm]}, color=blue] (v8) to (v1);
        \draw[-{Latex[length=3mm, width=2mm]}, color=blue] (v8) to (v7);
    \end{tikzpicture}
    \caption{Conditional Directed Acyclic Graph (CDAG). Input nodes are drawn as squares, output nodes as circles.}
    \label{fig:cdag}
\end{figure}

\begin{Def}[Conditional directed acyclic graphs (CDAGs)]
    \label{def-cdag}
    A \emph{conditional directed acyclic graph (CDAG)} $\Gk=(J,V,E)$ consists of two disjoint, finite sets of
    vertices/nodes: the set of \emph{input nodes} $J$, the set of \emph{output nodes} $V$;
    and a set of \emph{directed edges} $E \ins \lC w \tuh v\,|\, w \in J \cup V,\, v \in V \rC$,
    such that the directed graph $(J \dcup V, E)$ is \emph{acyclic}.
    So - per definition - there won't be any arrow heads pointing to input nodes $j \in J$.
    We drop ``conditional'' from the name if $J=\emptyset$ (DAG).
\end{Def}

An example of a CDAG is given in \Cref{fig:cdag}.

\begin{Rem}
    \label{rem:cdag-not-new}
    A CDAG is not a new class of graphs: it is a finite directed acyclic graph on $J \dcup V$ together with the
    requirement that no node of the designated subset $J$ has a parent. What is new is not the graph but its
    \emph{semantics}: the nodes of $J$ will carry no distribution, they index a family of Markov kernels. In
    particular all graph theoretic notions below are the classical ones, and the only genuinely new object is the
    asymmetric relation $\Perp^{\id}_\Gk$ of \Cref{def:d-separation}, which is d-separation composed with the
    fixed set $J$.
    For asymmetric independence models attached to graphs in a different way we refer to the local independence
    graphs of \cite{Did08} and \cite{MH20} with their $\delta$- resp.\ $\mu$-separation, discussed in
    \Cref{sec:main:comparison-local}.
\end{Rem}

\begin{Not}
    \label{not-cdag}
    Let $\Gk=(J,V,E)$ be a CDAG. We use the following shorthands:
    \begin{enumerate*}[label=(\roman*)]
        \item $v \in \Gk$ means $v \in J \dcup V$;
        \item $v \tuh w \in \Gk$ means $v \tuh w \in E$;
        \item $w \hut v \in \Gk$ means $v \tuh w \in E$;
        \item $v \sus w \in \Gk$ means $v \tuh w \in \Gk \,\lor\, v \hut w \in \Gk$.
    \end{enumerate*}
    A \emph{walk} from $v$ to $w$ in $\Gk$ is a finite sequence of nodes and edges
    $\pi = \lp v = v_0 \sus v_1 \sus \cdots \sus v_{n-1} \sus v_n = w \rp$ in $\Gk$ for some $n \ge 0$, i.e.\ such
    that $v_{k-1} \sus v_k \in \Gk$ for every $k=1,\dots,n$; the repeated appearance of the same nodes and edges
    is allowed, and so is the \emph{trivial walk} $\pi=\lp v_0 \rp$ consisting of a single node (in case $v=w$).
    It is a \emph{directed walk} if all arrow heads point in the direction of $w$ and none point back, i.e.\ if it
    is of the form $v=v_0 \tuh \cdots \tuh v_n=w$. We further put:
    \begin{enumerate*}[label=(\roman*),resume]
        \item $\Pa^\Gk(v) := \lC w \in \Gk \,|\, w \tuh v \in \Gk \rC$, the \emph{parents} of $v$;
        \item $\Ch^\Gk(v) := \lC w \in \Gk \,|\, v \tuh w \in \Gk \rC$, the \emph{children} of $v$;
        \item $\Anc^\Gk(v)$, the \emph{ancestors} of $v$, i.e.\ all $w \in \Gk$ admitting a directed walk
            $w \tuh \cdots \tuh v$ in $\Gk$ (in particular $v \in \Anc^\Gk(v)$, via the trivial walk);
        \item $\Desc^\Gk(v)$, the \emph{descendants} of $v$, i.e.\ all $w \in \Gk$ admitting a directed walk
            $v \tuh \cdots \tuh w$ in $\Gk$ (in particular $v \in \Desc^\Gk(v)$);
    \end{enumerate*}
    and we extend these notions to sets $A \ins J \dcup V$ by taking unions, e.g.\ $\Anc^\Gk(A) = \bigcup_{v \in A}\Anc^\Gk(v)$.
    Finally, a \emph{topological order} of $\Gk$ is a total order $<$ of $J \dcup V$ such that
    $v \in \Pa^\Gk(w) \implies v<w$ for all $v,w \in \Gk$; we write
    $\Pred^\Gk_<(v):= \lC w \in J \dcup V \,|\, w<v \rC$ for the set of \emph{predecessors} of $v$ w.r.t.\ $<$.
    Note that $\Gk$, being acyclic, always has a topological order, and that, conversely, a finite directed graph
    with a topological order is acyclic. Note also that
    $\Pa^\Gk(j)=\emptyset$ for every input node $j \in J$ and that $\Pa^\Gk(v) \notni v$ for every $v \in \Gk$.
\end{Not}

\subsection{d-Separation and id-Separation in CDAGs}

We first recall the classical notion of d-separation, see \cite{Pearl09, Lau90, Gei90, Ver93, Lau96, KF09, SGS2000},
which does not refer to the input nodes at all, and then adapt it to graphs that also have input nodes.

\begin{Def}[Colliders and d-blocked walks]
    \label{def:d-blocked}
    Let $\Gk=(J,V,E)$ be a CDAG, $C \ins J \dcup V$ a subset of nodes and $\pi$ a walk in $\Gk$:
    $\pi = \lp v_0 \sus \cdots \sus v_n \rp$, $n \ge 0$.
    We call the node $v_k$ of $\pi$ a \emph{collider of $\pi$} if two arrow heads of $\pi$ point at it, i.e.\ if
    $0<k<n$ and $v_{k-1} \tuh v_k \hut v_{k+1}$, and a \emph{non-collider of $\pi$} if at most one arrow head of $\pi$
    points at it. Since every edge of $\Gk$ carries exactly one arrow head, the non-colliders of $\pi$ are precisely
    the two \emph{end nodes} $v_0$, $v_n$ and those \emph{inner nodes} $v_k$, $0<k<n$, that form a chain
    ($v_{k-1} \tuh v_k \tuh v_{k+1}$ or $v_{k-1} \hut v_k \hut v_{k+1}$) or a fork ($v_{k-1} \hut v_k \tuh v_{k+1}$).
    We then say that the walk $\pi$ is \emph{$C$-d-blocked}, or \emph{d-blocked by $C$}, if
    $\pi$ has a non-collider in $C$ or a collider outside of $\Anc^\Gk(C)$;
    and that $\pi$ is \emph{$C$-d-open} otherwise, i.e.\ if every non-collider of $\pi$ lies outside of $C$ and every
    collider of $\pi$ lies in $\Anc^\Gk(C)$.
    Note the asymmetry between the two clauses: a non-collider blocks if it lies in $C$ itself, whereas a collider
    blocks only if neither it nor any of its descendants lies in $C$.
    Note that being a collider is a property of a \emph{position} in $\pi$: the same node of $\Gk$ may occur at
    several positions of $\pi$ and be a collider at some of them and a non-collider at others.
\end{Def}

\begin{Def}[d-separation and id-separation]
    \label{def:d-separation}
    Let $\Gk=(J,V,E)$ be a CDAG and $A,B,C \ins J \dcup V$ (not necessarily disjoint) subsets of nodes.
    \begin{enumerate}
        \item[1.] We say that \emph{$A$ is d-separated from $B$ given $C$ in $\Gk$}, in symbols:
            \[ A \Perp^d_\Gk B \given C, \]
            if every walk in $\Gk$ from a node in $A$ to a node in $B$ is d-blocked by $C$.
        \item[2.] We say that \emph{$A$ is input-d-separated (id-separated) from $B$ given $C$ in $\Gk$}, in symbols:
            \[ A \Perp^{\id}_\Gk B \given C \qquad :\iff \qquad A \Perp^{d}_\Gk B \cup J \given C, \]
            i.e.\ if every walk in $\Gk$ from a node in $A$ to a node in $J \cup B$ is d-blocked by $C$.
    \end{enumerate}
    If the corresponding condition fails we write: $A \nPerp^d_\Gk B \given C$, $A \nPerp^{\id}_\Gk B \given C$, resp.
    As special cases we abbreviate:
    $\displaystyle \quad A \Perp^d_\Gk B \; :\iff \; A \Perp^d_\Gk B \given \emptyset$, and
    $\displaystyle \quad A \Perp^{\id}_\Gk B \; :\iff \; A \Perp^{\id}_\Gk B \given \emptyset$.
\end{Def}

\begin{Rem}
    \label{rem:d-sep-conventions}
    \begin{enumerate}
        \item[a)] If $A,B,C$ are pairwise disjoint then $\Perp^d_\Gk$ is exactly the classical notion of d-separation
            in DAGs, see \cite{Pearl09, Lau90, Gei90, Ver93, Lau96, KF09, SGS2000}: the collider clause
            $v_k \notin \Anc^\Gk(C)$ is the usual requirement that neither $v_k$ nor any of its descendants lies in
            $C$. The only difference is that we quantify over all \emph{walks} and not only over all paths, which
            defines the same relation, since a $C$-open walk between two nodes exists if and only if a $C$-open path
            between them does. Walks are what makes the surgeries in the proofs elementary, and they also allow one
            to replace $\Anc^\Gk(C)$ by $C$ itself in the collider clause (a walk may travel from a collider down
            to a node of $C$ and back again), which is the form in which the collider condition is used throughout
            \Cref{sec:supp:separoid-rules-sigma} and \Cref{sec:supp:markov_property}, see \Cref{lem:strictly-open}.
        \item[b)] Our convention that the two end nodes of a walk count as non-colliders means that a walk with an end
            node in $C$ is always d-blocked. This only matters for non-disjoint $A,B,C$, where it provides
            Redundancy \ref{sep:d:red}, i.e.\ $A \ins C \implies A \Perp^d_\Gk B \given C$, and with it Extended Left
            Redundancy \ref{sep:sig:l-red}.
        \item[c)] The relation $\Perp^d_\Gk$ is symmetric, whereas $\Perp^{\id}_\Gk$ is not: id-separation depends on
            its right argument $B$ only through $B \cup J$, i.e.\ the input nodes always need to be separated from $A$
            as well. It is precisely this asymmetry that lets the separoid rules of id-separation match the asymmetric
            separoid rules of transitional conditional independence from \Cref{thm:separoid_axioms-tci}, see
            \Cref{d-sep-separoid-axioms} below, which in turn is what makes the global Markov property in the next
            section work in the presence of non-stochastic input variables.
        \item[d)] For $J = \emptyset$ the two relations coincide: $\Perp^{\id}_\Gk \,=\, \Perp^{d}_\Gk$.
    \end{enumerate}
\end{Rem}

\subsection{Separoid Rules for d-Separation and id-Separation}

Here we collect the formal rules that d-separation and id-separation satisfy. Short, self-contained walk-based
proofs of all of them are given in \Cref{sec:supp:separoid-rules-sigma}.
We start with the classical, symmetric rules for d-separation. The five rules of \Cref{thm:d-sep-sym-rules} are the
(symmetric) separoid, i.e.\ semi-graphoid, axioms, see \cite{PP85, Spo94, Daw01, Lau96, Pearl09}; they are the only
ones that are needed to derive the asymmetric rules a)--n) for id-separation below.

\begin{Thm}[Symmetric separoid rules for d-separation]
    \label{thm:d-sep-sym-rules}
    Let $\Gk=(J,V,E)$ be a CDAG and $A,B,C,D \ins J \dcup V$ (not necessarily disjoint) subsets of nodes.
    Then the ternary relation $\Perp \,=\, \Perp^d_\Gk$ satisfies the following rules:
{\setlength{\itemsep}{1pt}\setlength{\parskip}{0pt}
\begin{enumerate}[label=\arabic*)]
    \item Symmetry \ref{sep:d:sym}:
    \item[] \mbox{$A \Perp B \given C$}\Rimp\mbox{$B \Perp A \given C$}.
    \item Redundancy \ref{sep:d:red}:
    \item[] \mbox{$A \ins C$}\Rimp\mbox{$A \Perp B \given C$}.
    \item Decomposition \ref{sep:d:dec}:
    \item[] \mbox{$A \Perp B \cup D \given C$}\Rimp\mbox{$A \Perp B \given C$}.
    \item Weak Union \ref{sep:d:uni}:
    \item[] \mbox{$A \Perp B \cup D \given C$}\Rimp\mbox{$A \Perp B \given D \cup C$}.
    \item Contraction \ref{sep:d:con}:
    \item[] \mbox{$(A \Perp B \given D \cup C) \land (A \Perp D \given C)$}\Rimp\mbox{$A \Perp B \cup D \given C$}.
\end{enumerate}}
\end{Thm}

For d-separation, but not for general (stochastic) conditional independence, two further rules hold, which turn
$\Perp^d_\Gk$ into a \emph{compositional graphoid}, see
\cite{Gei90, Ver93, Lau90, Lau96, SGS2000, Pearl09}.

\begin{Lem}[Additional rules for d-separation]
    \label{thm:d-sep-extra-rules}
    In the situation of \Cref{thm:d-sep-sym-rules} the relation $\Perp \,=\, \Perp^d_\Gk$ satisfies in addition:
{\setlength{\itemsep}{1pt}\setlength{\parskip}{0pt}
\begin{enumerate}[label=\arabic*),start=6]
    \item Composition \ref{sep:d:com}:
    \item[] \mbox{$(A \Perp B \given C) \land (A \Perp D \given C)$}\Rimp\mbox{$A \Perp B \cup D \given C$}.
    \item Intersection \ref{sep:d:int}: If $B \cap D = \emptyset$ then:
    \item[] \mbox{$(A \Perp B \given D \cup C) \land (A \Perp D \given B \cup C)$}\Rimp\mbox{$A \Perp B \cup D \given C$}.
\end{enumerate}}
\end{Lem}

\begin{Rem}
    By Symmetry \ref{sep:d:sym} every one of the rules 2)-7) also holds in its ``left'' version, e.g.\
    $A \cup D \Perp B \given C \implies A \Perp B \given D \cup C$ (left weak union) or
    $(A \Perp B \given C) \land (D \Perp B \given C) \implies A \cup D \Perp B \given C$ (left composition).
    A convenient consequence of Redundancy \ref{sep:d:red}, Symmetry \ref{sep:d:sym}, Decomposition \ref{sep:d:dec}
    and Composition \ref{sep:d:com} is that neither argument sees its intersection with the conditioning set, see
    More Redundancies \ref{sep:d:m-red}: for all $E_1,E_2 \ins C$ we have:
    \[ A \Perp^d_\Gk B \given C \quad\iff\quad (A \cup E_1) \Perp^d_\Gk (B \cup E_2) \given C \quad\iff\quad
       (A\sm C) \Perp^d_\Gk (B \sm C) \given C. \]
\end{Rem}

We now turn to id-separation. Note that these rules match the rules of transitional conditional independence in
\Cref{thm:separoid_axioms-tci}. We formally show that the subsets of $J \dcup V$ of a CDAG $\Gk$ together with the
relations $=$, $\ins$, $\Perp^{\id}$, operation $\cup$ and element $\emptyset$ form a $J$-$\emptyset$-separoid,
see \Cref{def:t-k-separoid}.

\begin{Thm}[Separoid rules for id-separation]
    \label{d-sep-separoid-axioms}
    Let $\Gk=(J,V,E)$ be a CDAG and $A,B,C,D \ins J \dcup V$ subsets of nodes.
    Then the ternary relation $\Perp \,=\, \Perp^{\id}_\Gk$ satisfies the following rules:
{\setlength{\itemsep}{1pt}\setlength{\parskip}{0pt}
\begin{enumerate}[label=\alph*)]
    \item Extended Left Redundancy \ref{sep:sig:l-red}:
    \item[] \mbox{$A \ins C$}\Rimp\mbox{$A \Perp  B \given C$}.
     \item $J$-Restricted Right Redundancy \ref{sep:sig:r-red}:
     \item[] \mbox{$A \Perp  \emptyset \given C \cup J$} always holds.
	\item Left Decomposition \ref{sep:sig:l-dec}:
	\item[]\mbox{$A \cup D \Perp  B \given C$}\Rimp\mbox{$D \Perp  B \given C$}.
  \item Right Decomposition \ref{sep:sig:r-dec}:
    \item[] \mbox{$A \Perp  B \cup D \given C$}\Rimp\mbox{$A \Perp  D \given C$}.
    \item $J$-Inverted Right Decomposition \ref{sep:sig:inv-r-dec}:
    \item[] \mbox{$A \Perp  B \given C$}\Rimp\mbox{$A \Perp  J \cup B \given C$}.
    \item Left Weak Union \ref{sep:sig:l-uni}:
    \item[] \mbox{$A \cup D \Perp  B  \given C$}\Rimp\mbox{$A \Perp  B \given D \cup C$}.
  \item Right Weak Union \ref{sep:sig:r-uni}:\item[] \mbox{$A \Perp  B \cup D \given C$}\Rimp\mbox{$A \Perp  B \given D \cup C$}.
\item Left Contraction \ref{sep:sig:l-con}:\item[]
 \mbox{$(A \Perp  B \given D \cup C) \land (D \Perp  B \given C)$}\Rimp\mbox{$A \cup D \Perp  B \given C$}.
\item Right Contraction \ref{sep:sig:r-con}:
\item[] \mbox{$(A \Perp  B \given D \cup C) \land (A \Perp  D \given C)$}\Rimp\mbox{$A \Perp  B \cup D \given C$}.
\item Right Cross Contraction \ref{sep:sig:rc-con}:
\item[] \mbox{$(A \Perp  B \given D \cup C) \land (D \Perp  A \given C)$}\Rimp\mbox{$A  \Perp  B \cup D \given C$}.
\item Flipped Left Cross Contraction \ref{sep:sig:flc-con}:
\item[] \mbox{$(A \Perp  B \given D \cup C) \land (B \Perp  D \given C)$}\Rimp\mbox{$B \Perp  A \cup D \given C$}.
\end{enumerate}}
\end{Thm}

\begin{Rem}
In particular, we have the equivalences:
$$ (A \Perp  B \cup D \given C) \quad\iff\quad (A \Perp  B \given D \cup C) \quad\land\quad (A \Perp  D \given C)  ,$$
$$ (A \cup D \Perp  B \given C) \quad\iff\quad (A \Perp  B \given D \cup C) \quad\land\quad (D \Perp  B \given C).$$
\end{Rem}

\begin{Rem}[Symmetry]
    \label{rem:d-sep-symmetry}
Let the assumptions be like in \Cref{d-sep-separoid-axioms}.
    We also have the following rules:
    {\setlength{\itemsep}{1pt}\setlength{\parskip}{0pt}
\begin{enumerate}[label=\alph*)]
        \setcounter{enumi}{11}
    \item Restricted Symmetry \ref{sep:sig:res-sym}:
    \item[] \mbox{$(A \Perp  B \given C) \land (B \Perp  \emptyset \given C)$}\Rimp\mbox{$B \Perp  A \given C$}.
    \item $J$-Restricted Symmetry \ref{sep:sig:j-res-sym}:
    \item[] \mbox{$A \Perp  B \given C \cup J$}\Rimp\mbox{$B \Perp  A \given C \cup J$}.
    \item Symmetry \ref{sep:sig:sym}: If $J=\emptyset$ then:
    \item[] \mbox{$A \Perp  B \given C$}\Rimp\mbox{$B \Perp  A \given C$}.
    \end{enumerate}}
\end{Rem}

\begin{Lem}[More separoid like rules]
Let the assumptions be like in \Cref{d-sep-separoid-axioms}.
    We also have the following rules:
    {\setlength{\itemsep}{1pt}\setlength{\parskip}{0pt}
\begin{enumerate}[label=\alph*)]
        \setcounter{enumi}{14}
        \item Left Composition \ref{sep:sig:l-com}:
        \item[] \mbox{$(A \Perp  B \given  C) \land (D \Perp  B \given C)$}\Rimp\mbox{$A \cup D \Perp  B \given C.$}
        \item Right Composition \ref{sep:sig:r-com}:
        \item[] \mbox{$(A \Perp  B \given  C) \land (A \Perp  D \given C)$}\Rimp\mbox{$A \Perp  B \cup  D \given C.$}
         \item Left Intersection \ref{sep:sig:l-int}: If $A \cap D = \emptyset$ then:
        \item[] \mbox{$(A \Perp  B \given D \cup  C) \land (D \Perp  B \given A \cup C)$}\Rimp\mbox{$A \cup D \Perp  B \given C.$}
        \item Right Intersection \ref{sep:sig:r-int}: If $B \cap D = \emptyset$ then:
        \item[] \mbox{$(A \Perp  B \given D \cup  C) \land (A \Perp  D \given B \cup C)$}\Rimp\mbox{$A \Perp  B \cup D \given C.$}
        \item More Redundancies \ref{sep:sig:m-red}:
        \item[] \mbox{$A \Perp B \given C \iff (A\sm C) \Perp (B \sm C) \given C$}\Riff\mbox{$A \cup C \Perp J \cup B \cup C \given C$}.
    \end{enumerate}}
\end{Lem}

\begin{Rem}[How the rules for id-separation are obtained]
    \label{rem:d-sep-further-rules-sketch}
    All nineteen rules above are derived from the symmetric rules of \Cref{thm:d-sep-sym-rules} and
    \Cref{thm:d-sep-extra-rules} by purely formal
    arguments; not a single walk has to be inspected again. The recipe is always the same: unfold the definition
    $A \Perp^{\id}_\Gk B \given C \iff A \Perp^d_\Gk B \cup J \given C$ on both sides, and then apply the symmetric
    rules to the enlarged sets. For instance, Right Weak Union \ref{sep:sig:r-uni} reads
    $A \Perp^d (B \cup J) \cup D \given C \implies A \Perp^d B \cup J \given D \cup C$ and is thus an instance of
    Weak Union \ref{sep:d:uni}, and $J$-Inverted Right Decomposition \ref{sep:sig:inv-r-dec} becomes a tautology,
    because $(J \cup B)\cup J = B \cup J$.
    It is worth emphasising that the eleven rules a)--k), and with them the three symmetry rules l)--n), only use the
    five core rules of \Cref{thm:d-sep-sym-rules}; the two extra rules of \Cref{thm:d-sep-extra-rules} enter only in
    o)--s). The one derivation that looks as if it needed Composition is \emph{Flipped Left Cross Contraction}
    \ref{sep:sig:flc-con}, but it does not: the trick is to move $J$, rather than $B$, into the conditioning set,
    which turns the first assumption into $B \Perp^d A \given (D \cup J) \cup C$, so that Contraction
    \ref{sep:d:con} with the \emph{single} set $D \cup J$ from the second assumption already gives the conclusion.
    \emph{Right Intersection} \ref{sep:sig:r-int} cannot use Intersection \ref{sep:d:int} for the pair
    $(B \cup J, D)$, since these two sets need not be disjoint ($D$ may contain input nodes); one first has to
    pass to the set $(B \cup J)\sm D$, which is legitimate since the removed part $J \cap D$ lies in the conditioning
    set $D \cup C$.
    Finally, More Redundancies \ref{sep:sig:m-red} shows that one may always assume $A$ and $B$ to be disjoint from
    the conditioning set $C$, or alternatively that $C \ins A$ and $J \cup C \ins B$.
\end{Rem}
\subsection{Bayesian Networks with Input Nodes}
\label{sec:causal_bayesian_networks}

We now introduce a definition of Bayesian networks that allows for non-stochastic input variables.

\begin{Def}[Bayesian network with input nodes]
    \label{def:cbn}
    A \emph{Bayesian network with input nodes} $\Mk$ consists of:
    \begin{enumerate}
        \item a (finite) conditional directed acyclic graph (CDAG) $\Gk=(J,V,E)$, see \Cref{def-cdag},
        \item input variables $X_j$, $j \in J$, and
            (stochastic) output variables $X_v$, $v \in V$,
        \item a measurable space $\Xcal_v$ for every $v \in J \dcup V$,
            where $\Xcal_v$ is standard\footnote{\label{fn:bn-std}What the proof actually uses is only that
            $(\Xcal_A,\Xcal_v,\Xcal_C)$ is a disintegration triple, see \Cref{def:disintegration-triple}, for all
            $A \ins V$, $v \in V$ and $C \ins J \dcup V$; by \Cref{thm-regular-conditional-Markov-kernel} this holds whenever
            every $\Xcal_v$, $v \in V$, is standard. No assumption on the input spaces $\Xcal_j$, $j \in J$, is
            needed.}
            if $v \in V$,
        \item a Markov kernel, suggestively written as: $\displaystyle  \Pk_v\lp X_v\,\Vert\, X_{\Pa^{\Gk}(v)}\rp$:
        \[\begin{array}{ccrcl}
           & & \Xcal_{\Pa^{\Gk}(v)} & \dshto & \Xcal_v, \\
           && \qquad (A,x_{\Pa^{\Gk}(v)}) & \mapsto &  \Pk_v\lp X_v \in A\,\Vert\, X_{\Pa^{\Gk}(v)} = x_{\Pa^{\Gk}(v)} \rp,
        \end{array}\]
        for every $v \in V$, where we write for $D \ins J \dcup V$:
               \begin{align*}
                   \Xcal_D &:= \prod_{v \in D} \Xcal_v,&
                   \Xcal_\emptyset &:= \Asterisk,\\
                   X_D &:= (X_v)_{v \in D},&
                   X_\emptyset &:= \ast,\\
                   x_D &:= (x_v)_{v \in D},&
                   x_\emptyset &:= \ast.
              \end{align*}
    \end{enumerate}
    By abuse of notation, we denote the Bayesian network as:
    \[  \Mk= \lp \Gk, \lp  \Pk_v\lp X_v\,\Vert\, X_{\Pa^{\Gk}(v)}\rp  \rp_{v \in V}  \rp. \]
    We drop ``with input nodes'' if $J=\emptyset$.
\end{Def}

\begin{Def}
    \label{def:joint-markov-kernel-bn}
    Any Bayesian network $\Mk$ with input nodes comes with its \emph{joint Markov kernel}:
            \[\Pk(X_V\Vert X_J) :\, \Xcal_J \dshto \Xcal_{V},   \]
            given by:
       \[ \Pk(X_V\Vert X_J) :=  \bigotimes_{v \in V}^> \Pk_v\lp X_v\,\Vert\, X_{\Pa^{\Gk}(v)}\rp,\]
       where the product $\otimes^>$ is taken in reverse order of a fixed topological order $<$, i.e.\
       children appear only on the left of their parents in the product.
     Note that by \Cref{rem:markov-kernels-products} about associativity and (restricted) commutativity of the product the
     joint Markov kernel does actually not depend on the topological order.
\end{Def}

\begin{Rem}
    \label{rem:g-formula}
    The joint Markov kernel of \Cref{def:joint-markov-kernel-bn} is the kernel version of what is known in the
    literature as the \emph{truncated factorization} or \emph{$g$-formula}, see \cite{Rob86, Pearl09, FM20}: the
    factors belonging to the input nodes are absent, and $\Xcal_J$ is a parameter space rather than a sample space.
    Every statement of this section is a statement about that kernel and involves no further causal assumptions.
\end{Rem}

\subsection{Global Markov Property for Bayesian Networks}
\label{sec:main:global_markov}

We now turn to probably the most striking application of transitional conditional independence:
the \emph{global Markov property} for Bayesian networks that allow for (non-stochastic) input variables.
The global Markov property relates the graphical structure $\Gk$ to transitional
conditional independence relations between the corresponding transitional random variables $X_A$.
So checking the graph for id-separation relations, see \Cref{def:d-separation}, will then automatically imply
the existence of a conditional Markov kernel that does not depend on the specified variables,
which can be stochastic or not.

This will be the first time the global Markov property will be proven in this generality of measure theoretic probability, in the presence of input variables and with such a strong notion of conditional independence.
Graphical models that allow for latent confounders, cycles or selection bias can be treated by the same
strategy, see \cite{FM17, FM18, FM20, Ric03, Eva16, Eva18, RERS17}: one replaces d-separation by
$\sigma$-separation, which satisfies the same asymmetric separoid rules by the same shift, and the chaining
argument below is then reused verbatim; what has to be supplied is the factorization the induction starts from,
i.e.\ the analogue of \Cref{cor:local-markov} for the model class in question. For clarity of exposition we
restrict ourselves here to the acyclic case without latent confounders.

The proof of the global Markov property follows similar arguments as used in \cite{Lau90, Ver93, Ric03, FM17, FM18, RERS17},
namely chaining the separoid rules for transitional conditional independence (see \Cref{thm:separoid_axioms-tci}) and the ones for
id-separation (see \Cref{d-sep-separoid-axioms}) together in an inductive way.
The main difference here is that we never rely on the Symmetry property but instead use
the left and right versions of the separoid rules separately. Note, again, that the validity of those separoid rules
in this vast generality is a non-trivial result,
see \Cref{thm:separoid_axioms-tci},
and was only known for corner cases for other notions of extended conditional independence, rendering those less useful
for the applications to graphical models.
The full proof of the global Markov property is given in \Cref{sec:supp:markov_property}.

\begin{Thm}[Global Markov property for Bayesian networks with input nodes]
\label{thm-gmp-mI-CBN}
Consider a Bayesian network $\Mk$ with input nodes, with CDAG $\Gk=(J,V,E)$
and joint Markov kernel $\Pk(X_V\Vert X_J)$.
Then for all $A, B, C \ins J \dcup V$ (not-necessarily disjoint) we have the implication:
\[ A \Perp^{\id}_{\Gk} B \given C \qquad \implies \qquad X_A \Indep_{\Pk(X_V\Vert X_J)} X_B \given X_C.   \]
Recall that we have - per definition - an implicit dependence on $J$, $X_J$, resp., in the second argument on each
side; this is why the graphical premise is id-separation and not plain d-separation, see \Cref{def:d-separation}.
\end{Thm}

\begin{Rem}[The global Markov property is a statement about independence models]
    \label{rem:gmp-independence-models}
    Structurally, \Cref{thm-gmp-mI-CBN} is not a statement about Markov kernels at all but about
    \emph{independence models}. The CDAG $\Gk$ induces the asymmetric independence model
    \[ \Ical(\Gk) := \lC (A,B,C) \,\big|\, A \Perp^{\id}_\Gk B \given C \rC \]
    on the subsets of $J \dcup V$, and the Bayesian network $\Mk$ induces the asymmetric independence model
    \[ \Ical(\Mk) := \lC (A,B,C) \,\big|\, X_A \Indep_{\Pk(X_V\Vert X_J)} X_B \given X_C \rC; \]
    the global Markov property is exactly the inclusion $\Ical(\Gk) \ins \Ical(\Mk)$.
Both models satisfy the same rules: $\Ical(\Gk)$ is a $J$-$\emptyset$-separoid, see
    \Cref{def:t-k-separoid}, and $\Ical(\Mk)$ satisfies the corresponding rules a)--k) of
    \Cref{thm:separoid_axioms-tci}, with the hypotheses recorded in \Cref{tab:rules-hypotheses}. Of the
    conditional rules the proof invokes only f) Left Weak Union, and only with output-node spaces in its first two
    slots, which is why the input spaces $\Xcal_j$ may be arbitrary here; if in addition all $\Xcal_v$ and $\Xcal_J$ are standard then
    $\Ical(\Mk)$ is a full $T$-$\ast$-separoid, see \Cref{cor:t-star-separoid}.
    The proof consists of nothing but transporting the rules of the first model into the second along an induction
    over $\#V$.
    Neither the graph nor its separation criterion refers to probability or measure theory; this is why we keep the
    notation for $\Gk$ purely graph theoretic and reserve the double bar $\Vert$ for the kernels.
\end{Rem}

\begin{Rem}[What is gained over earlier formulations]
    \label{rem:gmp-what-is-gained}
    Global Markov properties for graphs with input or intervention nodes have been proven before. What
    \Cref{thm-gmp-mI-CBN} adds is not the graphical criterion but the \emph{conclusion}.
    For the $\Qcal$-extended conditional independence of \cite{FM20} the conclusion reads, in the notation of
    \Cref{sec:main:comparison:qeci},
    \[ X_A \Indep^\omega_{\Pk(X_V\Vert X_J) \otimes \Qcal} X_J, X_B \given X_C, \]
    i.e.\ a family of almost-sure conditional-expectation identities, one for every distribution on the input
    nodes; it does not assert the existence of any Markov kernel, so the object $\Qk(X_A|X_C)$ has to be
    constructed by hand afterwards, which is where \cite{FM20} run into measurability questions and settle them
    only under additional hypotheses, see their Appendix C. For the extended conditional independence of
    \cite{CD17} the obstruction is a different one: the mirror weak union that an induction over a graph uses at
    every step is exactly the axiom they leave open in general, see \Cref{sec:main:comparison:eci}.
    \Cref{thm-gmp-mI-CBN} hands one the kernel and the factorization directly, and this is the practical
    difference; see \Cref{tab:comparison}.
\end{Rem}

\begin{Cor}[Ordered local Markov property]
    \label{cor:local-markov}
    Let $\Mk$ be a Bayesian network with input nodes and CDAG $\Gk=(J,V,E)$, and let $<$ be a topological order of
    $\Gk$. Then for every $v \in V$ we have:
    \[ X_v \Indep_{\Pk(X_V\Vert X_J)} X_{\Pred^\Gk_<(v)} \given X_{\Pa^\Gk(v)}. \]
\begin{proof}
    Order the product of \Cref{def:joint-markov-kernel-bn} so that the factors of the nodes $w > v$ stand to the
    left, which is possible by \Cref{rem:markov-kernels-products}. Since every $\Pk_w$ is a probability (and not merely a sub-probability) kernel, marginalizing out these
    leftmost factors one after the other gives:
    \[ \Pk\lp X_v, X_{\Pred^\Gk_<(v)} \,\Vert\, X_J \rp
       = \Pk_v\lp X_v \,\Vert\, X_{\Pa^\Gk(v)} \rp \otimes \Pk\lp X_{\Pred^\Gk_<(v)} \,\Vert\, X_J \rp, \]
    where we used $\Pa^\Gk(v) \ins \Pred^\Gk_<(v)$. The same inclusion makes $X_{\Pa^\Gk(v)}$ a coordinate
    projection of $X_{\Pred^\Gk_<(v)}$, so the displayed identity is the defining factorization of
    \Cref{def:transitional_conditional_independence} with $\Qk\lp X_v|X_{\Pa^\Gk(v)}\rp := \Pk_v$.
\end{proof}
\end{Cor}

\begin{Rem}[On a pairwise Markov property]
    \label{rem:pairwise-markov}
    \Cref{cor:local-markov} is the strongest ``local'' statement that comes for free from the factorization, and
    it is the form in which the proof of \Cref{thm-gmp-mI-CBN} uses the model, see the display $(\dagger)$ there.
    One may ask for a \emph{pairwise} Markov property instead, i.e.\ for a set of statements about single
    non-adjacent nodes that already implies the global one. For ordinary conditional independence such implications
    rest on the Intersection property, which is available as soon as the joint distribution has a strictly
    positive density w.r.t.\ a product of $\sigma$-finite measures on the single-variable spaces, see
    \cite{Lau96, Pearl09}.
    Transitional conditional independence does not satisfy Intersection in general\footnote{Already in the corner
    case $\Tcal=\Asterisk$ of ordinary random variables: let $V$ be a random variable whose distribution is not a
    Dirac measure and let $X$, $Y$, $Z$ be three variables that are all almost surely equal to $V$; then
    $X \Indep Y \given Z$ and $X \Indep Z \given Y$ but not $X \Indep Y \otimes Z$.} (in contrast to
    id-separation, see \Cref{thm:d-sep-extra-rules} and Right Intersection \ref{sep:sig:r-int}), so the usual
    route from pairwise to global is not available here. Under which additional hypotheses on the Markov kernels an
    Intersection rule, and with it a pairwise Markov property, can be recovered is an interesting open question.
\end{Rem}

\begin{Rem}
    The global Markov property says that already checking the graphical criterion $\displaystyle A \Perp^{\id}_{\Gk} B \given C$ is
    enough to get
    the existence of a Markov kernel, suggestively written as:
    $\displaystyle \Pk\lp X_A\,|\,\cancel{X_B},X_{C \cap V}\,\Vert\, X_{C \cap J},\cancel{X_J} \rp,$
    such that:
    \[ \Pk(X_A,X_B,X_C\Vert X_J) = \Pk\lp X_A\,|\,\cancel{X_B},X_{C \cap V}\,\Vert\, X_{C \cap J},\cancel{X_J} \rp
    \otimes \Pk(X_B,X_C\Vert X_J). \]
    Note that the Markov kernel on the right hand side does not depend on $X_B$ and depends on the input variables
    only through $X_{C \cap J}$. It is exactly this extra conclusion, and not merely a numerical independence
    statement, that makes transitional conditional independence useful in this context.
\end{Rem}

\begin{Eg}
    Let $A \ins J \dcup V$ be an ancestral subset of $\Gk$, i.e.\ $A = \bigcup_{v \in A} \Anc^\Gk(v)$.
    Then we have:
    $\displaystyle (V \cap A) \Perp^{\id}_{\Gk} (J \sm A) \given J \cap A$,
    and thus, by the global Markov property, a Markov kernel $\Pk(X_{V \cap A}\Vert X_{J \cap A})$ such that:
    \[ \Pk(X_A\Vert X_J) = \Pk(X_{V \cap A}\Vert X_{J \cap A}) \otimes \Pk(X_{J \cap A}\Vert X_J). \]
    So for ancestral subsets we can only work with input variables from $J \cap A$ and ignore the ones from $J \sm A$,
    which is in correspondence with our expectations about ancestral relations.
\end{Eg}
\section{Comparison to Other Notions of Conditional Independence}
\label{sec:main:comparison}

Now that the theory is in place we can say precisely how transitional conditional independence relates to the other
notions of (extended) conditional independence in the literature, and in which respect it improves on them.
Throughout this section we fix a transition probability space $\lp \Wcal \times \Tcal, \Kk(W|T)\rp$ and transitional
random variables $\Xk$, $\Yk$, $\Zk$ with joint Markov kernel $\Kk(X,Y,Z|T)$, and we recall the definition:
\[ \Xk \Indep_{\Kk(W|T)} \Yk \given \Zk \qquad :\iff \qquad
   \exists \Qk(X|Z):\; \Kk(X,Y,Z|T) = \Qk(X|Z) \otimes \Kk(Y,Z|T). \]
It is useful to fix the axes along which these notions can be compared. Arguably, an extended notion of
conditional independence should:
{\setlength{\itemsep}{1pt}\setlength{\parskip}{0pt}
\begin{enumerate}[label=(\roman*)]
    \item \emph{Expressiveness.} Be able to express classical statistical concepts like sufficiency, adequacy and
        ancillarity.
    \item \emph{Asymmetry.} Embrace and anticipate the inherent asymmetry in the dependency between the stochastic
        (output) and the non-stochastic (input) parts of the variables.
    \item \emph{Generality.} Work for large classes of measurable spaces, of random and non-stochastic variables
        and of (transition) probability distributions, e.g.\ for variables that do not even have densities.
    \item \emph{Rules.} Satisfy reasonable relevance relations, e.g.\ as many of the separoid rules of
        \cite{Daw01} as possible.
    \item \emph{Factorization.} Give rise to meaningful factorizations of the distributions involved.
    \item \emph{Balance.} Be \emph{strong} enough that establishing it yields something (a factorization, a
        Markov kernel, an identification result), and at the same time \emph{weak} enough that there are usable
        criteria for establishing it, such as a directed global Markov property on the theoretical side or a
        statistical test on the empirical side.
\end{enumerate}}
The first five property columns of \Cref{tab:comparison} at the end of this section record (ii), (iv), (v), (iii)
and (i), in that order, and the last two record the other half of (vi): whether the relation can actually be obtained:
empirically, by a conditional independence test, and theoretically, by a global Markov property.
On (i), (ii), (iv) and (v) the notions genuinely differ; (iii) is listed because it is where the difficulty
of the present paper lies, even though most notions considered score well on it. Points (ii), (iv) and (v) are
jointly what a global Markov property for conditional probabilistic (causal) graphical models needs, while (iii)
is what makes that Markov property available on arbitrary input spaces.
Point (vi) is the trade-off the whole paper negotiates, and it is the reason a table cannot settle the comparison
on its own: \Cref{eg:symmetrization-loses} shows what is lost when a relation is made too weak, while
\Cref{thm-gmp-mI-CBN} is the criterion that keeps the present one usable. Conditional independence \emph{testing},
the other such criterion, is discussed in \Cref{sec:discussion}.
How the notions relate to each other formally is the subject of the present section; simplicity, finally, is a
matter of taste, on which we let the reader compare \Cref{def:transitional_conditional_independence} with
\Cref{def:ext-ci}.
The detailed statements and all proofs of this section are collected in \Cref{sec:supp:comparison}.

\subsection{Conditional Independence of Random Variables and the Two Failure Modes}
\label{sec:main:comparison:classical}

For ordinary random variables $X,Y,Z$ on a probability space $(\Wcal,\Pk(W))$, i.e.\ in the corner case
$\Tcal=\Asterisk$, there are two classical ways of writing conditional independence. The \emph{factorization form}
\begin{align}
    \Pk(X,Y|Z) &= \Pk(X|Z) \otimes \Pk(Y|Z) & \Pk(Z)\text{-a.s.} \label{eq:main:cond-ind-fact}
\end{align}
presupposes that the conditional distributions exist, which on general measurable spaces they need not; and where they
exist they are only unique up to null sets. The \emph{weak form}
\begin{align}
    \forall A \in \Bcal_\Xcal:\quad \E\lB \I_A(X)\,\big|\,Y,Z\rB &= \E\lB \I_A(X)\,\big|\,Z\rB & \Pk(W)\text{-a.s.},
    \label{eq:main:cond-ind-weak}
\end{align}
written $X \Indep^\omega_{\Pk(W)} Y \given Z$, avoids all existence questions, but pays for it: the conditional
expectations are defined separately for each event $A$ and need not be countably additive in $A$, so
\eqref{eq:main:cond-ind-weak} hands one no object to compute with. Transitional conditional independence keeps the
factorization of \eqref{eq:main:cond-ind-fact} and turns the existence of $\Pk(X|Z)$ from a hypothesis into part of
the assertion. The first two formulations agree as soon as $\Xcal$ is standard and $\Zcal$ countably generated,
with $\Ycal$ arbitrary; the factorization form needs in addition that $\Pk(X,Y|Z)$ and $\Pk(Y|Z)$ exist, so all
three agree once $\Xcal$ and $\Ycal$ are standard and $\Zcal$ is countably generated, see
\Cref{thm:supp:ci-equivalences}.

We can now make the two failure modes announced in \Cref{sec:introduction} precise. They appear as soon as one
leaves this corner case and asks, as \cite{Daw79} did, that conditional independence express statistical concepts
such as ancillarity, sufficiency and adequacy for a model
$\Pk(X|\Theta)$. First, the theory of \eqref{eq:main:cond-ind-fact}--\eqref{eq:main:cond-ind-weak} is of a purely
probabilistic nature, so the non-random parameter $\Theta$ has to be turned into a random variable, which requires a
prior $\Pk(\Theta)$ that the non-Bayesian setting does not provide, and with it disappear the conditional
distributions $\Pk(\Theta|S)$ on which one would want to impose conditions.
Second, one might try to repair this by declaring
$X \Indep Y \given Z$ to mean $X \Indep^\omega_{\Kk(W|T=t)} Y \given Z$ for every $t \in \Tcal$ separately. This
\emph{naive extension} silently conditions on \emph{all} of $\Tcal$, and it does not merely fail to express
ancillarity, it trivializes it: under a single $\Kk(W|T=t)$ the parameter is a constant, and a constant is
independent of everything, so ``$S$ is ancillary'' would come out true for every statistic $S$ whatsoever. The
statement ``$S(X)$ has the same distribution under every $\theta$'' compares different values of the parameter and
is not a statement about any single $\Kk(W|T=t)$.
Both observations point in the same direction: an extension of conditional independence that is to capture these
concepts must be \emph{asymmetric}, and it must treat the input variable $\Tk$ as a variable that may appear on the
right of the bar without being conditioned upon. This is exactly what the defining factorization above does, and it
is why \Cref{sec:applications-statistics} can characterize ancillarity, sufficiency and adequacy as
\emph{equivalences} rather than as implications.

\subsection{Variation Conditional Independence}
\label{sec:main:comparison:vci}

Variation conditional independence, see \cite{Daw01b, CD17}, is a non-probabilistic, set-theoretic relation: for maps
$X,Y,Z$ on a set $\Wcal$ one puts $X \Indep_v Y \given Z$ if $\Rcal(X|Y=y,Z=z)=\Rcal(X|Z=z)$ for all attainable
$(y,z)$, where $\Rcal(\,\cdot\,|\,\cdot\,)$ denotes the attainable range. It is in fact a \emph{strong} separoid, see
\cite{CD17} Thm.\ 2.7, a property that stochastic conditional independence does not have.
Its relation to transitional conditional independence is a \emph{formal} one: replacing the space of probability
measures $\Pcal(\Xcal)$ by the power set $2^\Xcal$ and measurable maps by arbitrary maps turns Markov kernels
$\Zcal \dshto \Xcal$ into maps $\Zcal \to 2^\Xcal$, and the defining factorization of transitional conditional
independence into the defining factorization of variation conditional independence, see \Cref{rem:supp:var-tci} and
\Cref{eg:supp:var-tci}. Beyond this analogy the two notions only meet in corner cases: transitional conditional
independence captures exactly the \emph{deterministic} ones, where it is equivalent to a functional dependence:
for measurable maps $F:\, \Tcal \to \Fcal$ and $H:\, \Tcal \to \Hcal$ with $\Fcal$ standard we have
\[ F \Indep_{\Kk(W|T)} \Yk \given H \qquad \iff \qquad \exists \text{ measurable } \varphi:\; F = \varphi \circ H, \]
see \Cref{thm:variation-ci}. Since the $T$-shift of a symmetric separoid is again an asymmetric separoid,
\Cref{thm:plus-t-k-separoid}, the two relations can be combined with a logical ``and'' without losing any of the
separoid rules, which is the practically useful way of using them together.

\subsection{Extended Conditional Independence}
\label{sec:main:comparison:eci}

\emph{Extended conditional independence} was introduced in \cite{CD17} for a family $\Ecal = (\Pk_t(W))_{t \in \Tcal}$
of probability measures and variables of a restricted shape: $X,Y,Z$ live on $\Wcal$, $\Phi,\Theta$ live on $\Tcal$,
the joint map $(\Phi,\Theta)$ is required to be injective, and one asks that
$\E_t\lB h(X)|Y,Z\rB = g_{h,\phi}(Z)$ hold $\Pk_t(W)$-almost surely for all $t \in \Phi^{-1}(\phi)$, see
\Cref{def:ext-ci}. Reading $(Y,\Theta)$ and $(Z,\Phi)$ as transitional random variables, transitional conditional
independence is the stronger notion:
\[ X \Indep_{\Pk(W|T)} (Y,\Theta) \given (Z,\Phi) \qquad \implies \qquad X \Indep_\Ecal (Y,\Theta) \given (Z,\Phi), \]
see \Cref{lem:tci-implies-eci}.

The instructive difference is less in what the two relations assert than in how the variables are organized.
\cite{CD17} keep the stochastic and the non-stochastic ones in separate compartments: $X,Y,Z$ are maps on the
sample space, $\Phi,\Theta$ are \emph{decision variables} on the regime space, the first argument of the relation
has to be fully stochastic, and the second right-hand argument $\Theta$ enters the definition only through the
requirement that $(\Phi,\Theta)$ be injective; the condition itself never mentions it. Every statement thus
carries a bookkeeping of which variable is of which kind, and the definition is technically involved besides,
quantifying over all bounded measurable $h$ and over all fibres $\Phi^{-1}(\phi)$. Here there is a single class of
objects, the transitional random variables, i.e.\ Markov kernels $\Wcal \times \Tcal \dshto \Xcal$: the
deterministic ones, the genuinely stochastic ones, the input $\Tk$ and the trivial $\deltabf_\ast$ are all members
of it, no argument of $\Indep_\Kk$ is restricted in shape, and the rules quantify over the whole class at once.

The same difference shows up in the rules and in what can be built on them. \cite{CD17} establish their axioms
P2$'$--P5$'$ together with the mirror axioms P3$''$ and P5$''$ in full generality; symmetry is available to them
only in a restricted form, and one axiom, the mirror weak union P4$''$, is left open in general and recovered only
under additional hypotheses: a discrete regime space, discrete stochastic variables, or a dominating regime,
see \cite{CD17} \S4 and \S5. Their \S6 takes the other way out and offers a \emph{pairwise} weakening of the
relation which satisfies the axioms unconditionally. Two of these observations deserve to be taken seriously
rather than scored off. First, the axiom that resists is exactly the one that resists here: P4$''$ is Left Weak
Union \ref{sep:tci:l-uni}, and their restriction of symmetry is our \Cref{sep:tci:sym}, the case
$\Tcal = \Asterisk$ of $\Tk$-Restricted Symmetry. Second, their enabling hypotheses reappear here, but as
\emph{criteria} rather than as standing assumptions: a discrete regime space is the discrete input space of point
2.\ of \Cref{thm-regular-conditional-Markov-kernel}, and a dominating regime is covered by case 4.\ of
\Cref{sep:tci:l-uni},
each one condition among several that one may check for the triple at hand.

The upshot is a difference in shape rather than in degree. There the theory advances behind a succession of
restrictions (on which variables may stand where, on which axiom is available, on which spaces or which regime
families make the missing axiom true), and no global Markov property is proven. Here the relation is defined once
for one class of objects, the rules are established for that class in one piece, and exactly three of them ask for
a disintegration triple; that uniformity is what an induction over a graph consumes, and it is what lets
\Cref{thm-gmp-mI-CBN} deliver a Markov kernel with no assumption whatsoever on the input spaces.

\subsection{Symmetric Extended Conditional Independence}
\label{sec:main:comparison:sym}

A symmetric notion of extended conditional independence was proposed in \cite{RERS17}, for \emph{kernels}
$q_V(x_V|x_W)$ in which the variables $x_W$ have been fixed and carry no distribution, the same situation that
$\Kk(W|T)$ describes here, with $x_W$ in the role of $T$. Their Definition 7 declares
$X_A \Indep X_B \given X_C \, [q_V]$ to hold if \emph{either} $A \cap W = \emptyset$ and
$q_V(x_A|x_B,x_C,x_{W \sm (B \cup C)})$ is a function of $(x_A,x_C)$ alone, \emph{or} the same statement with the
roles of $A$ and $B$ exchanged. Each of the two clauses is asymmetric, and asymmetric for exactly the reason
discussed in \Cref{rem:reading-the-asymmetry}: it requires the \emph{left} set to avoid the fixed variables. Their notion is
the \emph{disjunction} of the two clauses, and it is symmetric for that reason. It is therefore not merely an
example of a symmetric relation but is literally of the shape $\Indep^\lor$ below. They establish the
semi-graphoid axioms (their Proposition 9) and, on that basis, global Markov properties for conditional acyclic
directed mixed graphs and, recursively over the fixing operation, the nested Markov property. Their state spaces are finite discrete sets or
finite-dimensional real vector spaces rather than general measurable spaces, and their relation asserts a
functional dependence rather than producing a Markov kernel.
We do not reproduce it here and, in contrast to \Cref{sec:main:comparison:eci} and
\Cref{sec:main:comparison:qeci}, we prove no implication for it; the comparison below is the formal one, through
the symmetrization $\Indep^\lor$. Symmetry has a price, and \Cref{eg:symmetrization-loses} below makes it
explicit. The symmetrized
relation agrees with $\Indep_\Kk$ only as long as the input variable $\Tk$ is kept in the second argument; in the
equivalent spellings without it (and $\Indep_\Kk$ itself certifies them as equivalent, see \Cref{add-T}), it
is strictly weaker, and with $\deltabf_\ast$ in the second argument it is vacuous. A symmetric relation therefore
cannot carry the statements of \Cref{sec:applications-statistics} in a spelling-independent way.
Transitional conditional independence contains a symmetric notion as a special case rather than the other way round:
its symmetrization
\[ \Xk \Indep^\lor_{\Kk(W|T)} \Yk \given \Zk \qquad :\iff \qquad
   \lp \Xk \Indep_{\Kk(W|T)} \Yk \given \Zk \rp \;\lor\; \lp \Yk \Indep_{\Kk(W|T)} \Xk \given \Zk \rp \]
satisfies the symmetric separoid rules, see \Cref{thm:sym-t-k-separoid}, and is implied by, but does not imply,
transitional conditional independence. On the graph side the same phenomenon occurs: id-separation becomes symmetric
as soon as one conditions on all input nodes, $J$-Restricted Symmetry \ref{sep:sig:j-res-sym}, so the global Markov
property \Cref{thm-gmp-mI-CBN} immediately yields
\[ A \Perp^{\id}_\Gk B \given C \cup J \qquad \implies \qquad X_A \Indep^\lor_{\Pk(X_V\Vert X_J)} X_B \given X_C, X_J, \]
which recovers and strengthens the corresponding results of \cite{RERS17} and \cite{FM20}.

\begin{Eg}[What is lost by symmetrizing]
    \label{eg:symmetrization-loses}
    Let $\Thetacal = (0,1)$ and consider on $\Wcal = \lC 0,1 \rC^3$ the model
    \[ \Pk(X_1,X_2,Y|\Theta=\theta) = \mathbf{Ber}(X_1|\theta) \otimes \mathbf{Ber}(X_2|\theta) \otimes
       \mathbf{Ber}(Y|1/2), \]
    i.e.\ $X_1,X_2$ are two independent coin flips with unknown bias $\theta$ and $Y$ is an independent fair coin.
    \begin{enumerate}
        \item \emph{$\Indep^\lor$ cannot see which side carries the parameter.}
            We have $Y \Indep_\Pk X_1 \given \deltabf_\ast$, witnessed by $\Qk(Y) = \mathbf{Ber}(Y|1/2)$, and
            therefore $X_1 \Indep^\lor_\Pk Y \given \deltabf_\ast$. The unsymmetrized statement
            $X_1 \Indep_\Pk Y \given \deltabf_\ast$ is \emph{false}: it would require a distribution
            $\Qk(X_1)$ with $\Pk(X_1,Y|\Theta) = \Qk(X_1) \otimes \Pk(Y|\Theta)$ and hence
            $\Pk(X_1|\Theta=\theta) = \Qk(X_1)$ for every $\theta$.
            By \Cref{add-T} that false statement is $X_1 \Indep_\Pk \Theta,Y \given \deltabf_\ast$, i.e.\
            exactly the statement that the trivial statistic is adequate for $X_1$ with respect to $Y$;
equivalently, that $X_1$ is ancillary and independent of $Y$, see \Cref{sec:applications-statistics}.
            So the symmetrized relation holds while the statistical property it is meant to express fails, and it
            fails for the one reason $\Indep^\lor$ is blind to: it is $X_1$, not $Y$, that carries the parameter
            dependence. Note that this is adequacy in its $\Tk$-free spelling; written as
            $X_1 \Indep_\Pk \Theta,Y \given \deltabf_\ast$ the symmetrized statement is correctly false, since
            $\Theta,Y \Indep_\Pk X_1 \given \deltabf_\ast$ would make $\Theta$ constant. That the two spellings
            disagree is the content of point 2.
        \item \emph{$\Indep^\lor$ does not respect the equivalences of $\Indep_\Kk$.}
            Write $X := (X_1,X_2)$ and put $S := X_1$. By \Cref{add-T} the two statements
            \[ X \Indep_\Pk \Theta \given S \qquad \text{and} \qquad X \Indep_\Pk \deltabf_\ast \given S \]
            are \emph{equivalent}, and both say that $S$ is a sufficient statistic. Their symmetrizations are not
            equivalent. The second one is vacuous: $\deltabf_\ast \Indep_\Pk X \given S$ holds always, by Left
            Redundancy \ref{sep:tci:l-red}, so $X \Indep^\lor_\Pk \deltabf_\ast \given S$ holds for every model
            and every statistic whatsoever. The first one is false, as it should be: $S=X_1$ is not sufficient,
            since a kernel $\Qk(X|S)$ would have to reproduce $\Pk(X_2|\Theta=\theta) = \mathbf{Ber}(X_2|\theta)$
            without knowing $\theta$, and $\Theta \Indep_\Pk X \given S$ fails as well, since $S \in \lC 0,1 \rC$
            cannot determine $\theta \in (0,1)$.
    \end{enumerate}
    Points 1.\ and 2.\ are two instances of one phenomenon.
    If the second argument contains $\Tk$, the symmetrization costs nothing: the flipped statement
    $\Tk \otimes \Yk \Indep_\Kk \Xk \given \Zk$ forces its kernel to be of the form
    $\deltabf_t \otimes \Qk_0(Y|Z=z)$ for $\Kk(Z|T=t)$-almost every $z$ and every $t$, so
    $\Tk \ismapof_\Kk \Zk$ when $\Tcal$ is countably separated\footnote{\label{fn:dirac-selection}Both steps
    use countable separation and nothing else. Let $\Ecal = \lC A_n \rC_{n \in \N}$ be a countable separating
    family and put $D := \lC z \,|\, \forall n:\, \Qk(A_n|Z=z) \in \lC 0,1 \rC \rC$, which is measurable. For
    $z \in D$ the set $B_z := \bigcap_{n:\,\Qk(A_n|z)=1} A_n \cap \bigcap_{n:\,\Qk(A_n|z)=0} A_n^\cmpl$ is
    measurable, being a countable intersection, and has $\Qk(B_z|z)=1$; since $\Qk(\cdot|z)$ is a \emph{probability}
    measure this forces $B_z \neq \emptyset$, and since $\Ecal$ separates points $B_z$ is a singleton
    $\lC \varphi(z) \rC$. Hence $\Qk(A|Z=z) = \I_A(\varphi(z))$ for \emph{every} $A \in \Bcal_\Tcal$, and
    $\varphi$ is measurable, since $\varphi^{-1}(A) \cap D = \lC z \in D \,|\, \Qk(A|Z=z) = 1 \rC$. Extending
    $\varphi$ by a constant off $D$, whose complement $D^\cmpl$ is $\Kk(Z|T=t)$-null for every $t$, gives
    $\Tk \ismapof_\Kk \Zk$.
    Neither standardness of $\Tcal$ nor countable generation of $\Bcal_\Tcal$ is needed.},
    and $\Xk \Indep_\Kk \Tk \otimes \Yk \given \Zk$
    follows as soon as $(\Xcal,\Zcal,\Tcal)$ is a disintegration triple. That is the spelling in which
    \Cref{sec:applications-statistics} writes ancillarity, sufficiency and adequacy, and there the two relations
    agree.
    But $\Indep_\Kk$ identifies $\Xk \Indep_\Kk \Tk\otimes\Yk \given \Zk$ with
    $\Xk \Indep_\Kk \Yk \given \Zk$, see \Cref{add-T}, and $\Indep^\lor_\Kk$ does not: without $\Tk$ in the
    second argument it is strictly weaker, as point 1.\ shows, and for $\Yk = \deltabf_\ast$ it is vacuous, by
    Left Redundancy \ref{sep:tci:l-red}. That is how ancillarity and sufficiency, the existence of a conditional
    Markov kernel, \Cref{conditional-markov-kernel-as-ci}, and the invariance of a predictor across environments,
    \Cref{prp:invariance-tci} 3., all disappear.
    So the symmetrized relation assigns different truth values to statements that transitional conditional
    independence itself proves equivalent, and that is what makes it unusable as a carrier of the theory.
\end{Eg}

\subsection{Categorical Conditional Independence}
\label{sec:main:comparison:cat}

Conditional independence has also been formulated inside \emph{categorical probability}, where a Markov kernel is
an abstract morphism of a Markov category rather than a map into a space of measures. Two developments have to be
kept apart, because they differ on exactly the point that occupies \Cref{sec:conditional_markov_kernels}.

\cite{Cho17} define conditional independence \emph{through} disintegrations: for a fixed joint state,
$X \Indep Y \given Z$ holds if the disintegration on $X \otimes Y$ given $Z$ agrees, almost surely, with the
tensor of the two individual disintegrations (their Definition 6.6). The section carries the standing hypothesis
that the ambient category admits disintegration, so the conditionals are available by assumption; the relation is
symmetric, and they record this together with decomposition, weak union and contraction (their Proposition 6.10).

\cite{Fri20} drops that hypothesis. There, a morphism $f:\, A \to X \otimes Y$ carries the conditional
independence of $X$ and $Y$ given $A$ if it factors through the copy map, $f = (g \otimes h) \circ \mathrm{copy}_A$,
a condition that can be stated without conditionals existing at all, which is described there as the main
innovation over \cite{Cho17}. Several variants are given, according to whether one conditions on an input, on an
output, or on both (his Definitions 12.1, 12.12, 12.16 and 12.19), together with the corresponding semi-graphoid
properties (his Lemma 12.5 and its companions), and the framework then carries abstract versions of the
Fisher-Neyman factorization theorem, of Basu's theorem and of Bahadur's theorem (his Sections 14--16). Note that
his morphisms do have an input object, which may well be a parameter carrying no distribution; what makes the
relation symmetric is not the absence of an input but the shape of the condition, an existential factorization
that is invariant under exchanging the two outer legs. What \Cref{sec:main:comparison:sym} says about symmetric
relations (the symmetrization $\Indep^\lor$ and \Cref{eg:symmetrization-loses}) therefore applies to it. In the other direction its generality is considerable: it applies to any Markov category,
including $\mathbf{Stoch}$ on arbitrary measurable spaces, where conditionals need not exist.

The more recent \cite{FK23}, which proves a d-separation criterion in that setting, works instead with an
\emph{asymmetric} conditional independence for morphisms with inputs (their Definition 16, with the failure of
symmetry recorded in their Remark 17), reserving the symmetric notion of \cite{Cho17, Fri20} for states without
inputs. They do not arrive at it independently: they introduce it as the categorical rendering of the transitional
conditional independence of \Cref{def:transitional_conditional_independence} and say so. That the notion carries
over to a setting with no measure theory in it, and that it is needed precisely once input objects appear, is
evidence that the asymmetry is forced by the problem rather than by our construction.

The two treatments differ, however, on \emph{where} the ability to condition is assumed. \cite{Cho17} place their
whole section under the hypothesis that the ambient category admits disintegration, and the d-separation criterion
of \cite{FK23} assumes a Markov category \emph{with conditionals}; in both cases the assumption is global (every
object can be conditioned on), so a measure-theoretic instance forces \emph{every} space to be of a kind on
which conditional Markov kernels exist, standard Borel say, including the spaces that play the role of input
nodes. \cite{Fri20} avoids the assumption altogether, which is what buys him arbitrary Markov categories, but then
there is no d-separation criterion to be had. Here conditioning is neither assumed globally nor avoided: it is
requested for the individual triples that a given rule actually uses, see \Cref{def:disintegration-triple} and
\Cref{lem:trivial-right-argument}, and the global Markov property
\Cref{thm-gmp-mI-CBN} accordingly needs the output spaces to be standard and imposes \emph{no} condition on the
input spaces $\Xcal_j$, $j \in J$, at all.
The categorical and the measure-theoretic developments are complementary rather than competing. A Markov category
with conditionals \emph{assumes}
the disintegration that we have to construct, so the questions occupying
\Cref{sec:conditional_markov_kernels} and \Cref{sec:supp-disintegration} (on which spaces does a conditional
Markov kernel exist, and how badly is it non-unique) do not arise there; conversely, the categorical formulation reaches models, such as possibilistic
ones, that are not measure-theoretic at all.

\subsection{Extended Conditional Independence for Families of Distributions}
\label{sec:main:comparison:qeci}

Finally, fix a set $\Qcal \ins \Pcal(\Tcal)$ of distributions on the input space and define
\emph{$\Qcal$-extended conditional independence} by
\[ X \Indep^\omega_{\Kk(W|T)\otimes\Qcal} Y \given Z \qquad :\iff \qquad
   \forall\, \Qk(T) \in \Qcal:\quad X \Indep^\omega_{\Kk(W|T)\otimes\Qk(T)} Y \given Z, \]
a variant of the notion used in \cite{FM20}. It is remarkably simple, it inherits all separoid rules from the weak
conditional independence $\Indep^\omega$ on arbitrary measurable spaces, and by \Cref{thm:supp:ci-equivalences} we have
the implications
\[ X \Indep_{\Kk(W|T)} Y \given Z \quad\implies\quad X \Indep^\omega_{\Kk(W|T)\otimes\Qcal} T,Y \given Z
  \quad\implies\quad X \Indep^\omega_{\Kk(W|T)\otimes\Qcal} Y \given Z, \]
where the middle relation satisfies the \emph{asymmetric} rules of \Cref{thm:separoid_axioms-tci} without any
assumption on the spaces. Its one drawback is the decisive one: it asserts no factorization and hence produces no
Markov kernels. This is visible in the global Markov property, where $\Qcal$-extended conditional independence yields
only $X_A \Indep^\omega_{\Pk(X_V\Vert X_J)\otimes\Qcal} X_J,X_B \given X_C$, whereas \Cref{thm-gmp-mI-CBN} yields
$X_A \Indep_{\Pk(X_V\Vert X_J)} X_B \given X_C$ and thus hands one the kernel $\Qk(X_A|X_C)$ for free. Having to
construct those kernels by hand afterwards is where the measurability questions of \cite{FM20} Appendix C arise,
settled there under additional hypotheses and explicitly left for future work otherwise, and it was one of the
main motivations for developing transitional conditional independence.

\begin{Rem}[When the first implication is an equivalence]
    \label{rem:qeci-equivalence}
    The first of the two implications can be reversed in one important regime. Assume that $\Tcal$ is countable and
    discrete, that $\Xcal$ is standard, with $\Zcal$ arbitrary, and that some $\Qk(T) \in \Qcal$ satisfies
    $\Qk(T=t)>0$ for every $t \in \Tcal$. Then:
    \[ \Xk \Indep_{\Kk(W|T)} \Yk \given \Zk \qquad \iff \qquad
       \Xk \Indep^\omega_{\Kk(W|T)\otimes\Qk(T)} \Tk \otimes \Yk \given \Zk. \]
    Indeed, ``$\implies$'' is the implication above. For ``$\Longleftarrow$'' point 2.\ of
    \Cref{thm-regular-conditional-Markov-kernel}, applied with $\Asterisk$ in its third slot, provides a
    Markov kernel $\Rk(X|Z)$ which is a version of the conditional distribution of $X$ given $Z$ under
    $\Pk := \Kk(X,Y,Z|T) \otimes \Qk(T)$. Fix a measurable
    $A \ins \Xcal$. The weak conditional independence gives
    $\E_\Pk\lB \I_A(X) \,\middle|\, Y,Z,T \rB = \Rk(X \in A|Z)$ $\Pk$-almost surely, and since $\lC T = t \rC$ is
    an atom of $\sigma(\Tk)$ with $\Pk(T=t) = \Qk(T=t)>0$, restricting to it shows that $\Rk(X \in A|Z)$ is a
    version of $\E_{\Kk(X,Y,Z|T=t)}\lB \I_A(X) \,\middle|\, Y,Z \rB$ for every $t \in \Tcal$. Testing against
    rectangles $B \times C$ and applying Dynkin's lemma (with only the single event $A$ involved, as in
    \Cref{thm:supp:ci-equivalences}) yields $\Kk(X,Y,Z|T=t) = \Rk(X|Z) \otimes \Kk(Y,Z|T=t)$ for every
    $t \in \Tcal$, which is transitional conditional independence. So under these hypotheses nothing is lost by
    passing to the $\Qcal$-extended relation; the point is that outside them no such reduction is available:
neither on uncountable input spaces, where no single prior charges every point, nor on spaces on which the
    required conditional Markov kernel fails to exist, see \Cref{sec:main:comparison-summary}. The second
    implication, by contrast, is strict as soon as some $\Qk(T) \in \Qcal$ charges two points $t_1 \neq t_2$: for
    trivial $\Ycal$ and $\Zcal$ its right hand side is then the empty statement, while its left hand side forces
    $\Kk(X|T=t_1) = \Kk(X|T=t_2)$. For $\Qcal = \lC \deltabf_t \mid t \in \Tcal \rC$, on the other hand, the two
    sides coincide, since $\Tk$ is $\Kk(W|T) \otimes \deltabf_t$-almost surely constant.
\end{Rem}

\subsection{Comparison to Local Conditional Independence}
\label{sec:main:comparison-local}

A different asymmetric notion of irrelevance appears in the theory of stochastic processes: \emph{local
independence}, introduced by \cite{Sch70} and \cite{Aal87} and developed into a graphical theory by
\cite{Did07, Did08} and \cite{MH20}. There, for a multivariate counting or jump process $X=(X_v)_{v \in V}$ adapted
to a filtration, one says that $X_B$ is \emph{locally independent} of $X_A$ given $X_C$ if the compensator
(intensity) of $X_B$ w.r.t.\ the large filtration generated by $X_{A \cup B \cup C}$ is already measurable
w.r.t.\ the smaller filtration generated by $X_{B \cup C}$; that is, the instantaneous evolution of $X_B$ does not
depend on the past of $X_A$ once the past of $X_{B\cup C}$ is known.

Local independence and transitional conditional independence are both \emph{directed} relevance relations that
satisfy left and right versions of the separoid rules but not Symmetry, and both come with a graphical calculus
and a global Markov property, for local independence with respect to $\delta$-separation in \cite{Did08} and to
$\mu$-separation in \cite{MH20}, in directed (mixed) graphs that may contain cycles. The sources of the asymmetry
are, however, different: for local independence it is the direction of time and the filtration, for transitional
conditional independence it is the presence of non-stochastic input variables. Accordingly the two notions are of a
different nature (local independence constrains intensities at each time point, whereas transitional conditional
independence asserts the existence of a Markov kernel), and neither implies the other. Making the analogy
precise, e.g.\ by exhibiting local independence as a transitional conditional independence for a suitable
family of transition kernels along the filtration, is an interesting open direction that we do not pursue here.

\subsection{Summary}
\label{sec:main:comparison-summary}

\begin{table}[H]
    \centering
    \footnotesize
    \setlength{\tabcolsep}{4pt}
    \renewcommand{\arraystretch}{1.2}
    \begin{tabular}{lccccccc}
        \hline
        notion & asym- & separoid & yields & any & statis- & test- & global \\
               & metric & rules & kernels & space & tics & able & Markov \\
        \hline
        weak c.i.\ $\Indep^\omega$, see \cite{Daw79}              & --  & $+$   & --    & $+$   & --    & $+$ & $+$   \\
        variation c.i.\ $\Indep_v$, see \cite{Daw01b, CD17}       & --  & $+$   & --    & $+$   & --    & --  & --    \\
        extended c.i.\ $\Indep_\Ecal$, see \cite{CD17}            & $+$ & $(-)$ & --    & $+$   & $(+)$ & $+$ & $(-)$ \\
        kernel c.i., see \cite{RERS17}                            & --  & $+$   & --    & --    & --    & $+$ & $+$   \\
        categorical c.i., see \cite{Cho17}                        & --  & $+$   & $(-)$ & $(-)$ & --    & $+$ & $(+)$ \\
        categorical c.i., see \cite{Fri20}                        & --  & $+$   & --    & $+$   & $(+)$ & $+$ & $(+)$ \\
        $\Qcal$-extended c.i.\ ($\Tk$-shifted), see \cite{FM20}   & $+$ & $+$   & --    & $+$   & $(+)$ & $+$ & $+$   \\
        \textbf{transitional c.i.} $\Indep_\Kk$                   & $+$ & $+$   & $+$   & $+$   & $+$   & $+$ & $+$   \\
        \hline
    \end{tabular}
    \caption{Overview of the properties discussed in this section (``c.i.'' abbreviates ``conditional
    independence''). The first five property columns record what a relation \emph{gives}, the last two what it
    takes to \emph{get} it.
    \emph{Separoid rules}: the set of rules appropriate for the notion, i.e.\ the symmetric rules
    for a symmetric relation and the full asymmetric set of \Cref{thm:separoid_axioms-tci} for an asymmetric one,
for transitional conditional independence three of them under the disintegration triple hypothesis.
    \emph{Yields kernels}: the relation asserts the existence of a conditional Markov kernel that is not already
    available by assumption, together with a factorization.
    \emph{Any space}: no topological or countability assumption on the underlying measurable spaces is built into
    the definition. \emph{Statistics}: ancillarity, sufficiency and adequacy are \emph{equivalent} to, not merely
    implied by, the corresponding relation. \emph{Testable}: the relation is a hypothesis that can be decided from
    data; in the corner case $\Tcal=\Asterisk$ it is ordinary conditional independence between random
    variables, and for finite spaces it is a submodel of the saturated model, see \Cref{sec:discussion}.
    \emph{Global Markov}: a directed global Markov property reads the relation off a graph; for \cite{Cho17} and
    \cite{Fri20} this is available through the criterion of \cite{FK23}, which assumes a Markov category with
    conditionals. A ``$(\cdot)$'' indicates that the property holds only partially or only in restricted
    settings. For \cite{Cho17} the kernels are \emph{assumed} rather than constructed, which is also why the space
    entry is restricted; \cite{RERS17} works with finite discrete or finite-dimensional real state spaces rather
    than with general measurable spaces. The row for \cite{FM20} refers to
    the $\Tk$-shifted spelling $X \Indep^\omega_{\Kk(W|T)\otimes\Qcal} T,Y \given Z$; the unshifted relation is
    symmetric. We neither reproduce the definitions of \cite{RERS17, Cho17, Fri20} nor prove implications for them.
    Throughout, the entries for notions other than transitional conditional independence record what is available
    for them in the literature; we prove nothing about them here. See \Cref{sec:main:comparison:sym} and
    \Cref{sec:main:comparison:cat}.}
    \label{tab:comparison}
\end{table}

In summary, for the two \emph{extended} notions for which we have given the definition, namely the extended
conditional independence of \cite{CD17} and the $\Qcal$-extended conditional independence of \cite{FM20},
transitional conditional independence is stronger: it implies both, see \Cref{lem:tci-implies-eci} and
\Cref{thm:supp:ci-equivalences}. The implications are strict, and for a reason that is not an artefact: on a space
on which the required conditional Markov kernel does not exist, both weaker notions can hold while transitional
conditional independence fails for want of the factorizing kernel, see \Cref{conditional-markov-kernel-as-ci}. For the symmetric notions of \cite{RERS17, Cho17, Fri20} we compare only
formally, via the symmetrization $\Indep^\lor$ of \Cref{thm:sym-t-k-separoid}; the asymmetric categorical notion
of \cite{FK23} is a counterpart of transitional conditional independence in a Markov category rather than an
alternative to it, see \Cref{sec:main:comparison:cat}, and is therefore not listed separately. Among the notions considered
here, transitional conditional independence is the only one that combines the asymmetry needed for statistics with
the factorization needed for graphical models, on arbitrary measurable spaces. Two of the notions discussed above sit apart from this comparison. Variation conditional independence is a
statement about attainable ranges rather than about distributions and neither implies nor is implied by
transitional conditional independence, see \Cref{sec:main:comparison:vci}; it is listed in
\Cref{tab:comparison} only to record that it scores on none of the columns that matter here. Local independence,
which constrains intensities of a stochastic process, is not on this scale at all and does not appear in the
table, see \Cref{sec:main:comparison-local}. Both are best used side by side with transitional conditional
independence.
\section{Discussion}
\label{sec:discussion}

The theory developed here is a theory of the \emph{population} relation $\Xk \Indep_\Kk \Yk \given \Zk$.
For it to become a tool of statistical practice one needs a way of deciding, from data, whether such a relation holds.
In the corner case $\Tcal = \Asterisk$ this is the much-studied problem of conditional independence testing.
In the presence of a non-stochastic input $\Tk$ the null hypothesis changes shape: it asserts the existence of
\emph{one} Markov kernel $\Qk(X|Z)$ that works simultaneously for all values of the input, i.e.\ it is a statement
about a whole family of distributions and not about a single one, and it is asymmetric in $\Xk$ and $\Yk$.
For finite spaces this null hypothesis is a submodel of the saturated model and can be tested with a
likelihood-ratio statistic. A systematic development of \emph{transitional conditional independence testing},
together with its consistency properties and its use in constraint-based structure learning for Bayesian networks
with input nodes, is a natural next step that we leave to future work.

Related to this is the question of which \emph{further} rules the two independence models satisfy. Beyond the
separoid rules, d-separation satisfies Composition and Intersection, see \Cref{thm:d-sep-extra-rules}, while
transitional conditional independence in general satisfies neither; and for special families of distributions
considerably more is known, e.g.\ the gaussoid axioms in the Gaussian case, see \cite{BDKS19}, and the conditions
entering the characterizations of faithfulness in \cite{Sad17}. Which asymmetric analogues of these hold for
transitional conditional independence, and under which hypotheses on the Markov kernels, is open; an answer would
in particular settle the question of a pairwise Markov property raised in \Cref{rem:pairwise-markov}.

A second direction concerns the graphs. \Cref{sec:applications-causal_models} was carried out for acyclic graphs
without latent confounders. Since $\sigma$-separation satisfies the same asymmetric separoid rules, obtained by
the same shift, the chaining argument of the global Markov property is reused verbatim for graphs with cycles and
latent confounders once d-separation is replaced by it, the remaining ingredient being the factorization the
induction starts from, see \cite{Ric03, FM17, FM18, FM20}; and the resulting Markov kernels are then the natural starting point for the
identification of causal effects and the rules of $\doit$-calculus, see \cite{Pearl09, FM20}.

A third question is raised by \Cref{rem:Q-no-T} and runs through the whole paper without ever being settled. The
uniform relation $\Xk \Indep_\Kk \Yk \given \Zk$ is in general strictly stronger than its pointwise companion
$\Xk \Indep_\Kk \Yk \given \Zk \otimes \Tk$, and every application in \Cref{sec:applications-statistics} and
\Cref{sec:applications-causal_models} turns on that gap. We illustrate it repeatedly but never delimit it: for
which classes of transition probability spaces, or of transitional random variables, do the two coincide? A
characterization would say exactly when the extra strength of $\Indep_\Kk$ is free, and would be the natural
companion to a theory of testing it.

A fourth is quantitative. Transitional conditional independence is an exact statement: the kernel $\Qk(X|Z)$ either
exists or it does not. Several theories in which it appears as the exact case have a well-developed approximate
version: Le Cam's deficiency relaxes the comparison of experiments of \Cref{sec:blackwell}, see \cite{LeC64, Tor91}, and the privacy and information-flow literature relaxes exact non-interference to a budget. Is there an
approximate transitional conditional independence, measuring by how much the required factorization fails, whose
separoid rules degrade additively along a derivation? Such a calculus would make the rules of
\Cref{thm:separoid_axioms-tci} usable on estimates rather than on population statements.

We developed the theory of transition probability spaces, transitional random variables 
and transitional conditional independence. These concepts are most well behaved if the
underlying spaces have similar properties to standard measurable spaces. We isolated the exact requirement in the
notion of a \emph{disintegration triple}.
Furthermore, we proved the disintegration of transition probabilities, i.e.\ the existence
of conditional Markov kernels, on such spaces.\\
Transitional conditional independence was defined as an asymmetric notion of (ir)relevance relations.
We developed the theory of asymmetric separoids and showed that transitional conditional
independence and the graphical notion of id-separation, the $J$-shift of d-separation, satisfy all those
asymmetric separoid rules. We then showed how to relate those notions in graphical models and proved
a global Markov property for Bayesian networks with non-stochastic input variables in measure theoretic
generality.\\
We then compared transitional conditional independence to other notions of extended conditional
independence and showed that it is stronger than the two whose definition we reproduce, \cite{CD17} and
\cite{FM20}; for the symmetric notions we gave an explicit example of what a symmetrization loses, see
\Cref{eg:symmetrization-loses}.\\
We also showed that
transitional conditional independence can express classical statistical concepts like
ancillarity, sufficiency, adequacy and invariant reductions, etc.
We also demonstrated what it can say about Bayesian statistics, the likelihood principle,
propensity scores, etc.\\
Finally, we want to stress the simplicity of the definition of \emph{transitional conditional independence}:
\[ X \Indep_{\Kk} Y \given Z \qquad :\iff \qquad 
\exists \Qk(X|Z):\quad \Kk(X,Y,Z|T) = \Qk(X|Z) \otimes \Kk(Y,Z|T).  \]

\section*{Acknowledgments}
\addcontentsline{toc}{section}{Acknowledgments}
This work was partially supported by the European Research Council (ERC) under the European Union's Horizon 2020 research and innovation programme (grant agreement n$^{\mathrm{o}}$ 639466).
The author wants to express his gratitude towards Joris M. Mooij for many inspiring discussions and his constant support.

\addcontentsline{toc}{section}{References}

\newcommand{\etalchar}[1]{$^{#1}$}
\providecommand{\bysame}{\leavevmode\hbox to3em{\hrulefill}\thinspace}
\providecommand{\MR}{\relax\ifhmode\unskip\space\fi MR }
\providecommand{\MRhref}[2]{%
  \href{http://www.ams.org/mathscinet-getitem?mr=#1}{#2}
}
\providecommand{\href}[2]{#2}

\newpage
\appendix
\section*{Appendix}
\addcontentsline{toc}{section}{Appendix}

\section{The Measurable Extension Theorem}
\label{sec:measure-theory}
\label{sec:kuratowski-extension}

All measure theoretic notions that we use (countably generated and countably separated $\sigma$-algebras and
standard (Borel) measurable spaces) are standard in the field and were recalled in
\Cref{sec:transitional_probability_theory}; for the classical background see \cite{Kec95}, \cite{Bog07} 6.5,
\cite{Kle20} and \cite{Fremlin} 42. The only classical result that we use repeatedly and that is less widely
quoted is the following extension theorem of Kuratowski, so we state it here for reference.

\begin{Thm}[Kuratowski extension theorem for standard measurable spaces, see \cite{Kec95} 12.2]
    \label{thm-kuratowski-extension}
    Let $(\Xcal,\Bcal_\Xcal)$ be any measurable space,
    $\Wcal \ins \Xcal$ any subset endowed with the subspace $\sigma$-algebra $\Bcal_{\Xcal|\Wcal}$
    and $(\Ycal,\Bcal_\Ycal)$ be a standard measurable space.
    Let $f:\, (\Wcal,\Bcal_{\Xcal|\Wcal}) \to (\Ycal,\Bcal_\Ycal)$ be a measurable map.
    Then there exists a measurable map:
    \[ F:\,(\Xcal,\Bcal_\Xcal) \to (\Ycal,\Bcal_\Ycal)\]
    such that the restriction to $\Wcal$ equals $f$, i.e.\ $F_{|\Wcal} = f$.\\
    In short: There exists $F$ such that the following diagram commutes:
    \[\xymatrix{ (\Wcal,\Bcal_{\Xcal|\Wcal}) \ar@{^(->}[d]  \ar^-f[r]& (\Ycal,\Bcal_\Ycal) \\
            (\Xcal,\Bcal_\Xcal). \ar^{\exists F}[ru]
    }\]
\end{Thm}

\begin{Rem}
    \label{rem:no-measurable-selection}
    For a non-injective measurable map $f:\,\Xcal \to \Ycal$ between standard measurable spaces the
    image $f(\Xcal)$ is in general not measurable and no measurable right inverse of $f$ needs to
    exist; measurable selection theorems then require either additional structural assumptions on the
    fibres or a weakening of the measurability of the selector. See the discussions in
    \cite{Bog07, Fremlin, Bla75, Kec95}.
    For \emph{injective} $f$, in contrast, the Lusin--Souslin theorem, see \cite{Kec95} 15.1, gives
    $f(\Xcal) \in \Bcal_\Ycal$ together with a measurable inverse on the image.
\end{Rem}
\section{Proofs - Disintegration of Transition Probabilities}
\label{sec:supp-disintegration}

Here we will prove the existence and essential uniqueness of conditional Markov kernels for
standard measurable spaces.
For the classical disintegration theorem for probability distributions see \cite{Kal17} Thm.\ 1.25 and
\cite{Kle20} Ch.\ 8.3; also see \cite{Rao05, Fremlin, Bog07}.
The results below are the corresponding statements for Markov kernels, i.e.\ conditionally on a further
variable $Z$, and they provide the proofs of \Cref{lem:ess-unique} and
\Cref{thm-regular-conditional-Markov-kernel} of the main paper.

\subsection{Definition of Conditional Markov Kernels}

\begin{Def}[Conditional Markov kernel]
\label{def:reg-cond-trans-kernel-def}
Let $(\Xcal, \Bcal_\Xcal)$, $(\Ycal,\Bcal_\Ycal)$, $(\Zcal,\Bcal_\Zcal)$ be measurable spaces and:
\[\Kk(X,Y|Z): \, (\Zcal, \Bcal_\Zcal) \dshto  (\Xcal \times \Ycal,\Bcal_\Xcal \otimes \Bcal_\Ycal),\]
be a Markov kernel in two variables, and
\[\Kk(Y|Z): \, (\Zcal, \Bcal_\Zcal) \dshto  (\Ycal,\Bcal_\Ycal),\quad (B,z) \mapsto \Kk(X \in \Xcal, Y \in B|Z=z),\]
the marginal Markov kernel.
A \emph{conditional Markov kernel} of $\Kk(X,Y|Z)$ conditioned on $Y$ given $Z$ is a Markov kernel:
\[ \Kk(X|Y,Z) :\;  (\Ycal \times \Zcal, \Bcal_\Ycal\otimes \Bcal_\Zcal) \dshto  (\Xcal,\Bcal_\Xcal),\]
such that:
\[ \Kk(X,Y|Z) = \Kk(X|Y,Z) \otimes \Kk(Y|Z).\]
\end{Def}

\subsection{Essential Uniqueness of Conditional Markov Kernels}

\begin{Lem}[Essential uniqueness]
    \label{ess-unique}
    Consider Markov kernels:
    \[\Pk(X|Y,Z),\,\Qk(X|Y,Z):\, \Ycal\times \Zcal \dshto \Xcal,\]
    and
    \[\Kk(Y|Z):\,\Zcal \dshto \Ycal\]
    with any measurable spaces $\Xcal$, $\Ycal$, $\Zcal$ such that:
    \[ \Pk(X|Y,Z) \otimes \Kk(Y|Z) = \Qk(X|Y,Z) \otimes \Kk(Y|Z). \]
    We then have the following statements.
    \begin{enumerate}
     \item For every $A \in \Bcal_\Xcal$ the set:
    \[N_A := \lC (y,z) \in \Ycal \times \Zcal \,|\, \Pk(X \in A|Y=y,Z=z) \neq \Qk(X \in A|Y=y,Z=z)\rC\]
     is a  $\Kk(Y|Z)$-null set with $N_A \in \Bcal_\Ycal \otimes \Bcal_\Zcal$.
    \item If $\Bcal_\Xcal$ is countably generated then $N := \bigcup_{A \in \Bcal_\Xcal} N_A$ is a
        $\Kk(Y|Z)$-null set with $N \in \Bcal_\Ycal \otimes \Bcal_\Zcal$.
     \end{enumerate}
     \begin{proof}
         For $A \in \Bcal_\Xcal$
         we put $N_A^>:=h^{-1}(D^>)$, where $D^> := \{(x_1,x_2)\,|\,x_1,x_2 \in \R, x_1 > x_2\} \in \Bcal_{\R^2}$ and $h$ is given
         as the following composition of maps:
         \[ \xymatrix{h:\; \Ycal \times \Zcal \ar^-{\Pk\times\Qk}[r] & \Pcal(\Xcal) \times \Pcal(\Xcal)
                 \ar^-{j_A \times j_A}[r] & \R \times \R.
         }\]
         We define $N_A^<:=h^{-1}(D^<)$ similarly.
         Then $N_A = N_A^> \cup N_A^<$. \\
        1.) If $A \in \Bcal_\Xcal$ then the evaluation map $j_A$ is measurable by definition of the $\sigma$-algebra on $\Pcal(\Xcal)$.
         So $N_A^>, N_A^< \in \Bcal_\Ycal \otimes \Bcal_\Zcal$.
         By assumption we have:
           \[ \Pk(X|Y,Z) \otimes \Kk(Y|Z) = \Qk(X|Y,Z) \otimes \Kk(Y|Z). \]
         Evaluating both sides at $z$ on the measurable rectangle $A \times (N_A^>)_z$ gives thus the same value. So their difference equals $0$:
         \begin{align*}
             0& =  \int \I_{N_A^>}(y,z)\cdot \lp \Pk(X\in A|Y=y,Z=z) - \Qk(X \in A|Y=y,Z=z) \rp \, \Kk(Y \in dy|Z=z),
         \end{align*}
         where the integrand is:
         \[\I_{N_A^>}(y,z)\cdot \lp \Pk(X\in A|Y=y,Z=z) - \Qk(X \in A|Y=y,Z=z) \rp  \ge 0.\]
         This implies that: $N_A^>$ must be a $\Kk(Y|Z)$-null set in $\Bcal_\Ycal \otimes \Bcal_\Zcal$.
         By symmetry we get that also $N_A^<$ is a $\Kk(Y|Z)$-null set in $\Bcal_\Ycal \otimes \Bcal_\Zcal$
         and thus $N_A = N_A^< \cup N_A^>$ is
         a $\Kk(Y|Z)$-null set in $\Bcal_\Ycal \otimes \Bcal_\Zcal$ as well.\\
         2.) If now $\Bcal_\Xcal$ is countably generated then $\Bcal_\Xcal =\sigma\lp \Acal\rp$
         with a countable set $\Acal$ that is closed under finite intersections.
         One then puts $M:=\bigcup_{A \in \Acal} N_A$, which is, as countable union of measurable $\Kk(Y|Z)$-null sets,
         a measurable $\Kk(Y|Z)$-null set.
         Then one can define:
         \[ \Dcal := \{ A \in \Bcal_\Xcal \,|\, \forall (y,z) \in M^\cmpl:\, \Pk(X \in A|Y=y,Z=z)=\Qk(X \in A|Y=y,Z=z) \}.\]
         One easily sees that $\Dcal$ is closed under complements, countable disjoint unions and contains $\Xcal \in \Dcal$.
         This shows that $\Dcal$ is a Dynkin system (aka $\lambda$-system).
         Furthermore, we have: $\Acal \ins \Dcal$ and that $\Acal$ is closed under finite intersections.
         By Dynkin's lemma we get that:
         \[\Bcal_\Xcal = \sigma(\Acal) \ins \Dcal.\]
         If now $(y,z) \in N = \bigcup_{A \in \Bcal_\Xcal}N_A$ then there is an $A \in \Bcal_\Xcal$ such that:
         \[ \Pk(X \in A|Y=y,Z=z) \neq \Qk(X \in A|Y=y,Z=z).\]
         This implies $(y,z) \in M$ since $A \in \Dcal$ (otherwise we had equality above). Since this holds for all $(y,z)$ we get:
         \[ N = \bigcup_{A \in \Bcal_\Xcal}N_A \ins M =\bigcup_{A \in \Acal}N_A   \ins N, \]
         thus equality.  This shows that $N=M$ is a measurable $\Kk(Y|Z)$-null set.
     \end{proof}
 \end{Lem}

\subsection{Existence of Conditional Markov Kernels}
\label{sec:supp:existence-cond-markov-kernel}

\begin{Rem}[Existence of conditional Markov kernels]
    If $\Kk(X,Y|Z)$ is a Markov kernel then we want $\Kk(X|Y,Z)$ such that:
    \[\Kk(X,Y|Z) = \Kk(X|Y,Z) \otimes \Kk(Y|Z)\]
    holds. The heuristic here is to find something like a Radon-Nikodym derivative:
    \[ \Kk(X \in A|Y=y,Z=z) = \frac{\Kk(X \in A,Y \in dy |Z=z)}{\phantom{X\in A,\,}\Kk(Y \in dy|Z=z)}(y),  \]
    in a way that it is still a probability measure in $X$ and jointly measurable in $(y,z)$.\\
    To achieve measurability from the start we could restrict to $\Xcal=\Ycal=\R$ (or $\bar \R$ or $[0,1]$, etc.) and
    make use of Besicovitch derivation theorem, see \cite{Fremlin} 472D.
    For $\Kk(Y|Z)$-almost-all $(y,z) \in \Ycal \times \Zcal$ and all $x\in \Xcal$ we could construct:
    \[  \Kk(X \le x|Y=y,Z=z) = \inf_{m \in \N} \liminf_{n \in \N}
        \frac{\Kk(X \le\lceil x\rceil_m,
    Y \in B_{1/n}(y) |Z=z)}{\phantom{X \le\lceil x\rceil_m,\,}\Kk(Y \in B_{1/n}(y) |Z=z)}.  \]
    An alternative, which we will follow below, is to use Doob's derivation theorem for $\Xcal=\R$ (or $\bar \R$ or $[0,1]$, etc.)
    and countably generated $\Ycal$,
    see \cite{Del83} Thm.\ 58, \cite{Kle20} Example 11.17, or \Cref{thm:doob-derivation} below.
    This would yield that $\Kk(Y|Z)$-almost-all $(y,z) \in \Ycal \times \Zcal$ and all $x\in \Xcal$:
    \[ \Kk(X \le x| Y=y,Z=z) = \]
            \[\inf_{m \in \N}
        \liminf_{n \in \N}\lp \sum_{\substack{B \in \Ecal_n\\\Kk(Y \in B|Z=z)>0}} \frac{\Kk(X\le \lceil x\rceil_m, Y \in B|Z=z)}{\phantom{X\le \lceil x\rceil_m,\,}\Kk(Y \in B|Z=z)} \cdot \I_B(y) \rp,
    \]
    where $x < \lceil x\rceil_m:=\frac{\lfloor mx+1\rfloor}{m} \le x + \frac{1}{m}$ for $m \in \N$ and
     $\Bcal_\Ycal=\sigma( (B_n)_{n \in \N})$ and:
     \[ \Ecal_0:=\{\Ycal\}, \qquad\quad \Ecal_n:=\{ D \sm B_n, D \cap B_n\,|\, D \in \Ecal_{n-1} \}\sm\{\emptyset\}, \quad n \in \N. \]
    In both approaches,
     on the remaining points $(y,z)$, which lie inside the $\Kk(Y|Z)$-null set,
    we may choose arbitrarily, e.g.\ we can put:
    \[\Kk(X \le x|Y=y,Z=z) := \Kk(X\le x|Z=z).\]
\end{Rem}

\begin{Thm}[Doob's derivation theorem, see \cite{Del83} Thm.\ 58, \cite{Kle20} Example 11.17]
    \label{thm:doob-derivation}
    Let $\Ycal$ and $\Zcal$ be measurable spaces with $\Bcal_\Ycal$ countably generated\footnote{The assumption of being countably generated
    cannot be dropped in general; see \cite{Del83} Thm.\ 58, referring to \cite{Yor76} Thm.\ 3.}. Consider two non-negative finite transition measures:
    \[ \Qk(Y|Z), \, \Pk(Y|Z) :\, \Zcal \to \Mcal_+(\Ycal),\]
    where $\Mcal_+(\Ycal)$ is the set of non-negative finite measures on $\Ycal$.
    Furthermore, we assume for every $z \in \Zcal$ that $\Qk(Y|Z=z)$ is absolutely continuous w.r.t.\ $\Pk(Y|Z=z)$.
    Then there exists a jointly measurable map:
    \[ g:\, \Ycal \times \Zcal \to \R_{\ge 0},\]
    such that $g$ is the Radon-Nikodym derivative of $\Qk(Y|Z)$ w.r.t.\ $\Pk(Y|Z)$ (for all $Z=z$ simultaneously):
    \[ g(y,z) = \frac{\Qk(Y \in dy|Z=z)}{\Pk(Y \in dy|Z=z)}(y).  \]
    The latter means that for all $B \in \Bcal_\Ycal$ and all $z \in \Zcal$ we have:
    \[ \Qk(Y \in B|Z=z) = \int_B g(y,z)\, \Pk(Y \in dy|Z=z).\]
    \begin{proof}[Proof sketch]
    Since $\Ycal$ is countably generated:
     $\Bcal_\Ycal=\sigma( (B_n)_{n \in \N})$. We inductively then put:
     \[ \Ecal_0:=\{\Ycal\}, \qquad\quad \Ecal_n:=\{ D \sm B_n, D \cap B_n\,|\, D \in \Ecal_{n-1} \} \sm \{\emptyset\}, \quad n \in \N. \]
     Then the following sequence of (jointly) measurable functions $g_n:\, \Ycal \times \Zcal \to \R_{\ge 0}$ defined by:
    \[g_n(y,z) :=  \sum_{\substack{B \in \Ecal_n\\\Pk(Y \in B|Z=z)>0}} \frac{\Qk(Y \in B|Z=z)}{\Pk(Y \in B|Z=z)} \cdot \I_B(y),  \]
    is a uniformly integrable martingale w.r.t.\ $\lp\sigma(\Ecal_n)\rp_{n \in \N}$ for $\Pk(Y|Z=z)$ for every fixed $z \in \Zcal$.
    The proof requires the absolute continuity of $\Qk(Y|Z)$ w.r.t.\ $\Pk(Y|Z)$.
    So the limit:
    $g:=\lim_{n\in\N} g_n$
    exists in $\R_{\ge0}$  $\Pk(Y|Z)$-almost-surely and the convergence is also in $L^1$
    by the martingale convergence theorem, see \cite{Kle20} Thm.\ 11.7.
    As a countable limit of jointly measurable functions, $g$ is jointly measurable.
    The $L^1$-convergence and martingale property then implies that $g$ is the wanted Radon-Nikodym derivative.
    \end{proof}
\end{Thm}

\begin{Prp}[Existence of conditional Markov kernels for the unit interval]
\label{reg-cond-trans-kernel-existence-I}
\label{reg-cond-trans-kernel-existence-II}
Let $\Xcal, \Ycal, \Zcal$ be measurable spaces where $\Xcal=[0,1]$ and $\Ycal$ countably generated.
Let
\[\Kk(X,Y|Z): \, \Zcal \dshto  \Xcal \times \Ycal,\]
be a Markov kernel in two variables.
Then a conditional Markov kernel
conditioned on $Y$ given $Z$:
\[\Kk(X|Y,Z): \, \Ycal \times \Zcal \dshto  \Xcal,\]
exists.
\begin{proof}
    For fixed $A \in \Bcal_\Xcal$ and all $z \in Z$
    we have a finite measure $\Kk(X\in A, Y|Z=z)$ in $Y$, which is
    absolutely continuous w.r.t.\ the marginal $\Kk(Y|Z=z)$.
    Since also $\Bcal_\Ycal$ is countably generated, by Doob's derivation theorem, see \Cref{thm:doob-derivation},
    we get a (jointly) measurable map:
    \[g_A:\, \Ycal \times \Zcal \to \R_{\ge 0},\]
    such that for all $z \in \Zcal$ and $B \in \Bcal_\Ycal$:
    \[ \Kk(X \in A, Y \in B|Z=z) = \int_B g_A(y,z)\, \Kk(Y \in dy|Z=z).\]
    For $x \in \Xcal$ we will define:
    \[ G(x|y,z) := g_{[0,x]}(y,z).\]
As a next step we want to modify  $G(x|y,z)$ such that it becomes a cumulative distribution function in $x$, i.e.\ it corresponds to a probability distribution on $\Xcal$.
For this define $\Xcal_\Q:= \Xcal \cap \Q$, which is countable and dense in $\Xcal$.
First note that:
\[ N_1:=\lC (y,z) \in \Ycal \times \Zcal \,|\, G(1|y,z) \neq 1 \rC  \]
is a measurable $\Kk(Y|Z)$-null set.
Then, for every pair $x_1<x_2$ in $\Xcal_\Q$ consider:
\[ N_{(x_1,x_2)} := \{ (y,z) \,|\, G(x_1|y,z) > G(x_2|y,z) \} \in \Bcal_\Ycal \otimes \Bcal_\Zcal. \]
Since we have the equations:
\[ \begin{array}{rcl}
 &&\int \I_{N_{(x_1,x_2),z}}(y) \cdot G(x_1|y,z) \, \Kk(Y \in dy|Z=z) \\
&=& \Kk(X \le x_1, Y \in N_{(x_1,x_2),z} |Z=z) \\
&\stackrel{x_1 < x_2}{\le}& \Kk(X \le x_2, Y \in N_{(x_1,x_2),z} |Z=z)  \\
&=& \int \I_{N_{(x_1,x_2),z}}(y) \cdot G(x_2|y,z) \, \Kk(Y \in dy|Z=z) \\
&\stackrel{G(x_2|y,z) < G(x_1|y,z)}{\le} & \int \I_{N_{(x_1,x_2),z}}(y) \cdot G(x_1|y,z) \, \Kk(Y \in dy|Z=z)
\end{array}\]
we necessarily have $\Kk(Y \in N_{(x_1,x_2),z} |Z=z) = 0$ for every $z \in \Zcal$.\\
Then $N_{\mathrm{mon}}:= N_1 \cup  \bigcup_{x_1 < x_2 \in \Xcal_\Q} N_{(x_1,x_2)}$ is also a $\Kk(Y|Z)$-null set in $\Bcal_\Ycal \otimes \Bcal_\Zcal$.\\
It remains to secure right-continuity in $x$ along the rationals. For this we use the truncated approximation
from above that will also be used below: for $x \in \Xcal=[0,1]$ and $m \in \N$ put
$\lceil x\rceil_m:=\min(1,\lfloor mx+1\rfloor/m) \in \Xcal_\Q$, so that $\lceil x\rceil_m \downarrow x$ for
$x<1$ and $\lceil 1\rceil_m = 1$ for all $m$.
Now for $x \in \Xcal_\Q$ we define:
\[ D_x := \lC (y,z) \in N_{\mathrm{mon}}^\cmpl \,\Big|\, G(x|y,z) <  \inf_{m \in \N} G(\lceil x\rceil_m|y,z)  \rC
   \;\in\; \Bcal_\Ycal \otimes \Bcal_\Zcal.\]
Note that $\lceil x \rceil_m \in \Xcal_\Q$, so all the values $G(\lceil x\rceil_m|y,z)$ are defined, and that
$D_1=\emptyset$. We then get, writing $D_{x,z}$ for the $z$-section of $D_x$:
\[ \begin{array}{rcl}
&&  \int \I_{D_{x,z}}(y) \cdot G(x|y,z) \, \Kk(Y \in dy|Z=z) \\
&\stackrel{(1)}{\le}& \int \I_{D_{x,z}}(y) \cdot \inf_{m \in \N} G(\lceil x\rceil_m|y,z) \, \Kk(Y \in dy|Z=z) \\
&\stackrel{(2)}{\le}& \inf_{m \in \N} \int \I_{D_{x,z}}(y) \cdot G(\lceil x\rceil_m|y,z) \, \Kk(Y \in dy|Z=z) \\
&\stackrel{(3)}{=}& \inf_{m \in \N} \Kk(X \le \lceil x\rceil_m, Y \in D_{x,z} |Z=z) \\
&\stackrel{(4)}{=}&  \Kk(X \le x, Y \in D_{x,z} |Z=z) \\
&\stackrel{(3)}{=}& \int \I_{D_{x,z}}(y) \cdot G(x|y,z) \, \Kk(Y \in dy|Z=z),
\end{array}\]
where (1) is the definition of $D_x$, (2) is the trivial estimate
$\I_{D_z} \cdot \inf_m G(\lceil x\rceil_m) \le \I_{D_z} \cdot G(\lceil x\rceil_m)$, valid pointwise for every
single $m$ (so that no convergence theorem is needed here), (3) is the defining property of
$G(q|\cdot)=g_{[0,q]}$ for $q \in \Xcal_\Q$, and (4) holds because $\bigcap_{m \in \N} [0,\lceil x\rceil_m] = [0,x]$: the sets
$[0,\min_{m \le M} \lceil x\rceil_m]$, $M \in \N$, decrease to $[0,x]$, so continuity from above of the measure
$\Kk(X,Y|Z=z)$ applies. (The sequence $\lceil x\rceil_m$ itself need not be monotone in $m$.)
So the first and the last term agree, which forces equality in (1). Since the integrand in (1) is \emph{strictly}
smaller on $D_{x,z}$, this shows that $\Kk(Y \in D_{x,z}|Z=z)=0$ for all $z \in \Zcal$.
So $D:= N_{\mathrm{mon}} \cup \bigcup_{x \in \Xcal_\Q} D_x$ is again a $\Kk(Y|Z)$-null set in $\Bcal_\Ycal \otimes \Bcal_\Zcal$.\\
So far, we got that $G$, when restricted to $\Xcal_\Q \times D^\cmpl$,
is jointly measurable in $(y,z)$ for fixed $x$ and monotone non-decreasing and continuous from above in $x$ for fixed $(y,z)$ with $G(1|y,z)=1$.
We now aim to extend $G$ to $\Xcal \times \Ycal \times \Zcal$.\\
Recall the approximation $\lceil x\rceil_m = \min(1,\lfloor mx+1\rfloor/m) \in \Xcal_\Q$ from above:
the map $x \mapsto \lceil x\rceil_m$ is measurable, $\lceil 1\rceil_m=1$, and for $x \in [0,1)$ we have
$x < \lceil x\rceil_m \le x + \frac{1}{m}$, so that $\lceil x\rceil_m$ converges to $x$ from above for
$m \to \infty$.\\
We then define for all $(x,y,z) \in \Xcal \times \Ycal \times \Zcal$:
\[ F(x|y,z) := \inf_{m \in \N} \left \{ G(\lceil x\rceil_m|y,z)\right\} \cdot \I_{D^\cmpl}(y,z) + \Kk(X \le x|Z=z) \cdot \I_D(y,z). \]
It is clear that $F$ is again jointly measurable in $(y,z)$ for fixed $x$ and agrees with $G$  on $\Xcal_\Q \times D^\cmpl$ by construction.
As a monotone approximation from above it is clearly continuous from above, monotone non-decreasing and satisfies $F(1|y,z)=1$ for all $(y,z)$.
So for fixed $(y,z)$ now $F(\cdot|y,z)$ corresponds to a probability distribution $\Kk(X|Y=y,Z=z)$ on $\Bcal_\Xcal$, uniquely given by the defining relations on sets $[0,x]$:
\[ F(x|y,z) =: \Kk(X \le x|Y=y,Z=z),\]
for all $x \in \Xcal$. \\
Now define $\Dcal \ins \Bcal_\Xcal$ as the set of all $A \in \Bcal_\Xcal$ that satisfy:
\begin{enumerate}
    \item the map $(y,z) \mapsto \Kk(X \in A|Y=y,Z=z)$ is $(\Bcal_\Ycal\otimes\Bcal_\Zcal)$-$\Bcal_\R$-measurable, and:
    \item for all $z \in \Zcal$ and $B \in \Bcal_\Ycal$ the following equation holds:
        \[ \Kk(X \in A, Y \in B |Z=z) = \int \I_B(y) \cdot \Kk(X \in A | Y=y,Z=z) \, \Kk(Y \in dy | Z=z). \]
\end{enumerate}
Since $\Kk(X,Y \in B|Z=z)$ and $\Kk(X|Y=y,Z=z)$ are probability measures in $X$ the system $\Dcal$ is closed under countable disjoint unions and
complements and contains $\Xcal=[0,1]$. So $\Dcal$ is a Dynkin system.
We already know that for $x \in \Xcal_\Q$ the map $(y,z) \mapsto \Kk(X \le x|Y=y,Z=z)=F(x|y,z)$ is measurable.
Since for $x \in \Xcal_\Q$ and every $B \in \Bcal_\Ycal$, $z \in \Zcal$, we have:
\[ \I_B(y) \cdot \Kk(X \le x | Y=y,Z=z) = \I_B(y) \cdot G(x|y,z) \]
up to the $\Kk(Y|Z=z)$-null set $D_z$ we already get for those $x \in \Xcal_\Q$:
\[ \Kk(X \le x, Y \in B |Z=z) = \int \I_B(y) \cdot \Kk(X \le x | Y=y,Z=z) \, \Kk(Y \in dy | Z=z). \]
This shows that $\Ecal:=\{ [0,x] \,|\, x \in \Xcal_\Q\} \ins \Dcal$. Since $\Ecal$ is closed under finite intersections
Dynkin's lemma (see \cite{Kle20} Thm. 1.19) implies:
\[\Bcal_\Xcal=\sigma(\Ecal) \ins \Dcal.\]
This shows that the two conditions hold for all $A \in \Bcal_\Xcal$ and thus that $\Kk(X|Y,Z)$ is the desired conditional Markov kernel.
\end{proof}
\end{Prp}

\begin{Thm}[Existence of conditional Markov kernels]
\label{reg-cond-trans-kernel-existence-V}
Let $\Xcal$ be a standard measurable space,
$\Ycal$ be a countably generated measurable space and $\Zcal$ be any measurable space.
Let
\[\Kk(X,Y|Z): \, \Zcal \dshto  \Xcal \times \Ycal,\]
be a Markov kernel in two variables.
Then there exists a conditional Markov kernel $\Kk(X|Y,Z)$ conditioned on $Y$ given $Z$.
For the unconditional case, i.e.\ $\Zcal=\Asterisk$, also see \cite{Kal17} Thm.\ 1.25 and \cite{Kle20} Ch.\ 8.3.
Note that in that unconditional case $\Ycal$ does not need to be countably generated, see
\Cref{reg-cond-trans-kernel-existence-VI} below.
\begin{proof}
Since $\Xcal$ is standard we find, by the definition of a standard measurable space, an injective  measurable map:
\[\varphi:\, (\Xcal, \Bcal_\Xcal) \inj ([0,1],\Bcal_{[0,1]})=:(\Xcal', \Bcal_{\Xcal'}) \]
 that induces  a measurable isomorphism $(\Xcal, \Bcal_\Xcal) \cong (\varphi(\Xcal), \Bcal_{\Xcal'|\varphi(\Xcal)})$
 with $\varphi(\Xcal) \in \Bcal_{\Xcal'}$.
 So we can consider the push-forward Markov kernel $\Kk(X',Y|Z):= \Kk(\varphi(X),Y|Z)$:
\[\Kk(X',Y|Z): \, (\Zcal, \Bcal_\Zcal) \stackrel{\Kk(X,Y|Z)}{\dshto}  (\Xcal \times \Ycal,\Bcal_{\Xcal} \otimes \Bcal_{\Ycal})
\stackrel{\varphi \times \id_\Ycal }{\longrightarrow}  (\Xcal' \times \Ycal,\Bcal_{\Xcal'} \otimes \Bcal_{\Ycal}).\]
Since $\Xcal'=[0,1]$ and $\Bcal_\Ycal$ is countably generated we can apply \Cref{reg-cond-trans-kernel-existence-II} and we then get the conditional Markov kernel  $\Kk(X'|Y,Z)$:
$$\Kk(X'|Y,Z): \, (\Ycal \times \Zcal, \Bcal_{\Ycal}\otimes \Bcal_\Zcal) \dshto  (\Xcal',\Bcal_{\Xcal'}).$$
If we put $A':=\Xcal' \sm \varphi(\Xcal)$ we have $A' \in \Bcal_{\Xcal'}$.
So $\Kk(X' \in A' |Y=y,Z=z)$ is well-defined for every $(y,z) \in \Ycal \times \Zcal$.
Consider the set:
$$D:=\{(y,z) \in \Ycal \times \Zcal \,|\, \Kk(X' \in A'|Y=y,Z=z) >0 \}.$$
We first show that $D \in \Bcal_\Ycal \otimes \Bcal_\Zcal$.
We now consider the Markov kernel $\Kk(X'|Y,Z)$ as the measurable map:
$$ \Kk(X'|Y,Z):\, (\Ycal \times \Zcal, \Bcal_{\Ycal}\otimes \Bcal_\Zcal) \longrightarrow (\Pcal(\Xcal'), \Bcal_{\Pcal(\Xcal')}).$$
Now consider the map:
$$ j_{A'}:\, \Pcal(\Xcal') \to [0,1],\quad \mu \mapsto \mu(A').$$
Since $A' \in \Bcal_{\Xcal'}$, the map $j_{A'}$ is $\Bcal_{\Pcal(\Xcal')}$-$\Bcal_{[0,1]}$-measurable
by the very definition of the $\sigma$-algebra $\Bcal_{\Pcal(\Xcal')}$.
Then the composition:
$$H:\, (\Ycal \times \Zcal, \Bcal_{\Ycal}\otimes \Bcal_\Zcal) \stackrel{\Kk(X'|Y,Z)}{\longrightarrow}
(\Pcal(\Xcal'), \Bcal_{\Pcal(\Xcal')}) \stackrel{j_{A'}}{\longrightarrow} ([0,1],\Bcal_{[0,1]}) $$
is $\Bcal_{\Ycal}\otimes \Bcal_\Zcal$-measurable.
It follows that $D = H^{-1}((0,1]) \in \Bcal_{\Ycal}\otimes \Bcal_\Zcal$.
Since $y \mapsto (y,z)$ is $\Bcal_\Ycal$-$(\Bcal_\Ycal\otimes\Bcal_\Zcal)$-measurable it follows that the section
$D_z \in \Bcal_\Ycal$ for every $z \in \Zcal$.
So we can evaluate $\Kk(Y \in D_z|Z=z)=\Kk((Y,Z) \in D|Z=z)$ for every $z \in \Zcal$ and the map:
$$ \Zcal \stackrel{\Kk(Y,Z|Z)}{\longrightarrow} \Pcal(\Ycal\times\Zcal) \stackrel{j_D}{\longrightarrow} [0,1],\quad z \mapsto \Kk(Y \in D_z|Z=z),$$
is $\Bcal_\Zcal$-measurable, again by the definition of the $\sigma$-algebra on $\Pcal(\Ycal\times\Zcal)$
and since $D \in \Bcal_\Ycal \otimes \Bcal_\Zcal$.
So we can integrate:
$$\begin{array}{rcl}
0 &=& \Kk(\varphi(X) \in A', Y \in D_z |Z=z) \\
 &=& \int \I_{D_z}(y) \cdot \Kk(X' \in A' | Y=y,Z=z) \, \Kk(dy | Z=z),
\end{array}$$
where the integrand is strictly positive on $D_z$.
It follows that for all $z \in \Zcal$ we have: $\Kk(Y \in D_z|Z=z)=0$. \\
For $A \in \Bcal_\Xcal$ let $\tilde A \in \Bcal_{\Xcal'}$ such that $A =\varphi^{-1}(\tilde A)$.
Since $\varphi$ is injective with $\varphi(\Xcal) \in \Bcal_{\Xcal'}$ we have $\varphi(A)=\tilde A \cap \varphi(\Xcal) \in \Bcal_{\Xcal'}$,
independently of the choice of $\tilde A$.
So we can define:
$$\begin{array}{rcl}
&& \Kk(X \in A|Y=y,Z=z) \\
&:=& \Kk(X' \in \varphi(A) |Y=y,Z=z) \cdot \I_{D^\cmpl}(y,z) + \Kk_0(X \in A) \cdot \I_{D}(y,z),
\end{array}$$
with any probability distribution $\Kk_0$ on $\Xcal$.
For fixed $(y,z)$ this is a probability measure on $\Bcal_\Xcal$: for $(y,z) \in D$ this is clear, and for
$(y,z) \in D^\cmpl$ we have $\Kk(X' \in \varphi(\Xcal)|Y=y,Z=z) = 1 - \Kk(X'\in A'|Y=y,Z=z) = 1$, while
$\sigma$-additivity follows from the injectivity of $\varphi$.
For fixed $A \in \Bcal_\Xcal$ the map $(y,z) \mapsto \Kk(X \in A|Y=y,Z=z)$ is
$\Bcal_{\Ycal}\otimes \Bcal_\Zcal$-measurable, since $\varphi(A) \in \Bcal_{\Xcal'}$, so that $j_{\varphi(A)}$
is measurable by definition of $\Bcal_{\Pcal(\Xcal')}$, and since $D \in \Bcal_{\Ycal}\otimes \Bcal_\Zcal$.
So we get the Markov kernel:
$$\Kk(X|Y,Z): \, (\Ycal \times \Zcal, \Bcal_{\Ycal}\otimes \Bcal_\Zcal) \dshto  (\Xcal,\Bcal_{\Xcal}).$$
Furthermore, since $D$ is a $\Kk(Y|Z)$-null set, we have for all $A \in \Bcal_\Xcal$, $B \in \Bcal_\Ycal$ and $z \in \Zcal$:
$$ \Kk(X \in A, Y \in B|Z=z)
 = \int \I_{B}(y) \cdot \Kk(X \in A | Y=y,Z=z) \, \Kk(Y \in dy | Z=z).$$
This shows the claim.
\end{proof}
\end{Thm}

\begin{Lem}
    \label{lem:conditional-markov-kernel-deterministic}
    Let $\Xcal$, $\Ycal$, $\Zcal$ be measurable spaces and:
    \[\Kk(X,Y|Z): \, \Zcal \dshto  \Xcal \times \Ycal,\]
     a Markov kernel in two variables.
     Assume one of the following:
     \begin{enumerate}
        \item $X \ismapof_\Kk (Y,Z)$.
        \item $Y \ismapof_\Kk Z$.
     \end{enumerate}
    then there exists a conditional Markov kernel conditioned on $Y$ given $Z$.\\
    \begin{proof}
        1.) If $X \ismapof_\Kk (Y,Z)$ then there exists a measurable map $\varphi:\, \Ycal \times \Zcal \to \Xcal$ such that:
        \[  \Kk(X,Y|Z) = \deltabf_\varphi(X|Y,Z) \otimes \Kk(Y|Z). \]
        So $\Kk(X|Y,Z) := \deltabf_\varphi(X|Y,Z)$ is a conditional Markov kernel.\\
        2.) If $Y \ismapof_\Kk Z$ then there exists a measurable map $\varphi:\, \Zcal \to \Ycal$ with
        $\Kk(Y|Z=z) = \deltabf_{\varphi(z)}$ for every $z \in \Zcal$. This already fixes the \emph{joint} kernel:
        for $A \in \Bcal_\Xcal$ and $B \in \Bcal_\Ycal$ we have
        \[ \Kk(X \in A, Y \in B|Z=z) = \Kk(X \in A|Z=z) \cdot \I_B(\varphi(z)), \]
        since for $\varphi(z) \notin B$ both sides vanish, $\lC Y \in B \rC$ being $\Kk(X,Y|Z=z)$-null, and for
        $\varphi(z) \in B$ the complement $\lC Y \in B^\cmpl \rC$ is $\Kk(X,Y|Z=z)$-null. So:
        \begin{align*}
            \Kk(X,Y|Z) &= \deltabf_\varphi(Y|Z) \otimes \Kk(X|Z) \\
                       &=  \Kk(X|Z) \otimes \deltabf_\varphi(Y|Z)\\
                       &=  \Kk(X|Z) \otimes \Kk(Y|Z).
        \end{align*}
        So $\Kk(X|Y,Z) := \Kk(X|Z)$ is a conditional Markov kernel.
    \end{proof}
\end{Lem}

\begin{Rem}
    \Cref{reg-cond-trans-kernel-existence-V} and \Cref{lem:conditional-markov-kernel-deterministic} above,
    together with \Cref{reg-cond-trans-kernel-existence-VI}, \Cref{reg-cond-trans-kernel-existence-VII} and
    \Cref{reg-cond-trans-kernel-existence-VIII} below, prove \Cref{thm-regular-conditional-Markov-kernel} of the
    main paper: points 1., 5.\ and 6.\ by the first two, and points 2., 3.\ and 4.\ by the latter three, in that
    order.
    Unfolding \Cref{not:ismapof} there, point 5.\ says that
    $\Kk(X,Y|Z=z) = \deltabf_{\varphi(\cdot\,,z)}(X|Y) \otimes \Kk(Y|Z=z)$ for a measurable
    $\varphi:\, \Ycal \times \Zcal \to \Xcal$ (the map may depend on $z$), and point 6.\ says that
    $\Kk(Y|Z=z) = \deltabf_{\varphi(z)}$ for a measurable $\varphi:\, \Zcal \to \Ycal$.
    Note that a standard measurable space is countably generated, see \cite{Bog07} 6.5.8, so point 1.\ in
    particular applies whenever both $\Xcal$ and $\Ycal$ are standard.
\end{Rem}
\begin{Cor}[Existence for discrete parameter spaces]
    \label{reg-cond-trans-kernel-existence-VI}
    Let $\Xcal$ be a standard measurable space, $\Ycal$ an \emph{arbitrary} measurable space and $\Zcal$ a
    \emph{discrete} measurable space, i.e.\ $\Zcal$ countable with $\Bcal_\Zcal = \Pws^\Zcal$. Let
    \[\Kk(X,Y|Z): \, \Zcal \dshto  \Xcal \times \Ycal\]
    be a Markov kernel in two variables. Then there exists a conditional Markov kernel $\Kk(X|Y,Z)$
    conditioned on $Y$ given $Z$. In particular $(\Xcal,\Ycal,\Zcal)$ is a disintegration triple, see
    \Cref{def:disintegration-triple}.
\begin{proof}
    Fix $z \in \Zcal$ and consider the probability measure $\Kk(X,Y|Z=z)$ on $\Xcal \times \Ycal$, i.e.\ the
    case $\Zcal = \Asterisk$ of \Cref{reg-cond-trans-kernel-existence-V}. In that case the only place where the
    countable generation of $\Ycal$ was used is the appeal to \Cref{thm:doob-derivation}, which for a single
    parameter value degenerates to the ordinary Radon-Nikodym theorem on $(\Ycal,\Bcal_\Ycal)$ and hence needs no
    assumption on $\Ycal$ at all: the conditional cumulative distribution function $G(q|y)$, $q \in [0,1]\cap\Q$,
    is obtained as a Radon-Nikodym derivative of $B \mapsto \Kk(X \le q, Y \in B|Z=z)$ w.r.t.\
    $B \mapsto \Kk(Y \in B|Z=z)$, and the regularization over the countably many rationals is unchanged.
    This is the classical disintegration theorem for a standard first factor and an arbitrary second factor, see
    also \cite{Kal17} Thm.\ 1.25 and \cite{Kle20} Ch.\ 8.3.
    We thus obtain, for every $z \in \Zcal$, a Markov kernel $\Kk(X|Y,Z=z):\, \Ycal \dshto \Xcal$ with
    $\Kk(X,Y|Z=z)=\Kk(X|Y,Z=z) \otimes \Kk(Y|Z=z)$.
    It remains to check joint measurability. For $A \in \Bcal_\Xcal$ put
    $g(y,z) := \Kk(X \in A|Y=y,Z=z)$. Since $\Zcal$ is countable and discrete we have for every
    $D \in \Bcal_{[0,1]}$:
    \[ g^{-1}(D) = \bigcup_{z \in \Zcal} \lp \lC y \in \Ycal \st \Kk(X\in A|Y=y,Z=z) \in D \rC \times \{z\} \rp
        \;\in\; \Bcal_\Ycal \otimes \Bcal_\Zcal, \]
    a countable union of measurable rectangles. So $\Kk(X|Y,Z)$ is a Markov kernel
    $\Ycal \times \Zcal \dshto \Xcal$, and the desired factorization holds because it holds for every $z$
    separately.
\end{proof}
\end{Cor}

\begin{Cor}[Existence for a discrete second variable]
    \label{reg-cond-trans-kernel-existence-VII}
    Let $\Xcal$ and $\Zcal$ be \emph{arbitrary} measurable spaces and let $\Ycal$ be countable and discrete,
    i.e.\ $\Bcal_\Ycal = \Pws^\Ycal$. Let
    \[\Kk(X,Y|Z): \, \Zcal \dshto  \Xcal \times \Ycal\]
    be a Markov kernel in two variables. Then there exists a conditional Markov kernel $\Kk(X|Y,Z)$
    conditioned on $Y$ given $Z$. In particular $(\Xcal,\Ycal,\Zcal)$ is a disintegration triple, see
    \Cref{def:disintegration-triple}.
\begin{proof}
    If $\Zcal = \emptyset$ there is nothing to show, and if $\Zcal \neq \emptyset$ then $\Xcal \neq \emptyset$,
    since otherwise no Markov kernel $\Zcal \dshto \Xcal \times \Ycal$ exists. So fix a point
    $x_0 \in \Xcal$ and abbreviate, for $y \in \Ycal$, $A \in \Bcal_\Xcal$ and $z \in \Zcal$:
    \[ f_y(z) := \Kk(Y=y|Z=z), \qquad g_{A,y}(z) := \Kk(X \in A, Y=y|Z=z). \]
    Both are measurable in $z$, since $\{y\} \in \Bcal_\Ycal$ and $A \times \{y\} \in \Bcal_\Xcal \otimes
    \Bcal_\Ycal$, and $0 \le g_{A,y} \le f_y$. Now put:
    \[ \Kk(X \in A|Y=y,Z=z) := \begin{cases}
        \frac{g_{A,y}(z)}{f_y(z)} & \text{ if } f_y(z) > 0,\\
        \deltabf_{x_0}(A) & \text{ if } f_y(z) = 0.
    \end{cases}\]
    We check the three required properties.

    \emph{Probability measure in $A$.} Fix $(y,z)$. If $f_y(z)=0$ this is the Dirac measure $\deltabf_{x_0}$.
    If $f_y(z)>0$ then $A \mapsto g_{A,y}(z)$ is a finite measure on $\Bcal_\Xcal$, being the restriction of the
    probability measure $\Kk(X,Y|Z=z)$ to the measurable sets $A \times \{y\}$, and its total mass is
    $g_{\Xcal,y}(z) = f_y(z)$. Dividing by $f_y(z)$ gives a probability measure.

    \emph{Measurability in $(y,z)$.} Since $\Ycal$ is countable and discrete, a map
    $h:\, \Ycal \times \Zcal \to [0,1]$ is measurable as soon as all its sections $h(y,\cdot)$ are, because
    \[ h^{-1}(D) = \bigcup_{y \in \Ycal} \lp \{y\} \times h(y,\cdot)^{-1}(D) \rp \;\in\; \Bcal_\Ycal \otimes
        \Bcal_\Zcal \]
    is then a countable union of measurable rectangles for every $D \in \Bcal_{[0,1]}$. For fixed $y$ and
    $A$ the section $z \mapsto \Kk(X \in A|Y=y,Z=z)$ is measurable, since it agrees with the measurable function
    $g_{A,y}/f_y$ on the measurable set $\lC f_y > 0 \rC$ and is constant on its complement.

    \emph{Factorization.} Let $A \in \Bcal_\Xcal$ and $B \subseteq \Ycal$. Since $\Kk(Y|Z=z)$ is a measure on a
    countable discrete space, integration against it is summation, and:
    \begin{align*}
        \lp \Kk(X|Y,Z) \otimes \Kk(Y|Z) \rp (A \times B|z)
            &= \sum_{y \in B} \Kk(X \in A|Y=y,Z=z) \cdot f_y(z)\\
            &= \sum_{y \in B} g_{A,y}(z) \;=\; \Kk(X \in A, Y \in B|Z=z),
    \end{align*}
    where in the middle step the summands with $f_y(z) = 0$ may be dropped on both sides, because then also
    $g_{A,y}(z) \le f_y(z) = 0$, and the last equality is the countable additivity of
    $\Kk(X,Y|Z=z)$ applied to the disjoint decomposition $A \times B = \biguplus_{y \in B} (A \times \{y\})$.
    The measurable rectangles $A \times B$ form a $\pi$-system generating $\Bcal_\Xcal \otimes \Bcal_\Ycal$, and
    both sides are probability measures on that $\sigma$-algebra, so they agree everywhere by Dynkin's lemma,
    see \cite{Bog07} 1.9.3. This shows
    $\Kk(X,Y|Z) = \Kk(X|Y,Z) \otimes \Kk(Y|Z)$.

    Note that no property of $\Xcal$ or $\Zcal$ entered the argument. In particular, taking $\Ycal = \Asterisk$
    shows that $(\Xcal,\Asterisk,\Zcal)$ is a disintegration triple for arbitrary $\Xcal$ and $\Zcal$.
\end{proof}
\end{Cor}

\begin{Rem}[Weakening the density hypothesis]
    \label{rem:density-hypothesis}
    The corollary below is stated with $\sigma$-finite $\mubf$ and $\nubf$, matching point 4.\ of
    \Cref{thm-regular-conditional-Markov-kernel}, since that is the familiar formulation. Its proof, however,
    uses the absolute continuity only through the iterated identity:
    \[ \Kk(X \in A, Y \in B|Z=z) = \int_B \int_A k(x,y|z)\, \mubf(dx)\,\nubf(dy), \qquad
        A \in \Bcal_\Xcal,\; B \in \Bcal_\Ycal,\; z \in \Zcal, \]
    which under the stated hypotheses is Tonelli's theorem. The product measure $\mubf \otimes \nubf$ is
    therefore never formed and Tonelli is never applied to it; all that is used about $\nubf$ is that
    integration against a measure with $\nubf$-density $k(\cdot|z)$ is integration against $k(y|z)\,\nubf(dy)$,
    which holds for \emph{every} measure. So corollary and proof remain valid verbatim for an
    \emph{arbitrary} measure $\nubf$ on $\Ycal$, once the displayed identity is assumed directly in place of
    the absolute continuity. Note that assuming the density outright also makes the Radon-Nikodym theorem,
and with it the $\sigma$-finiteness needed for the \emph{existence} of a density, superfluous.

    The $\sigma$-finiteness of $\mubf$, in contrast, cannot be dropped. It enters in the appeal to the
    measurability half of Tonelli's theorem, and is not decorative there: for a measure that is not
    $\sigma$-finite the partial integral of a jointly measurable non-negative function need not be measurable
    at all: take the counting measure on $\R$ and the indicator function of a Borel subset of $\R^2$ whose
    projection is not Borel, i.e.\ is analytic but not Borel, see \cite{Kec95} \S14.
\end{Rem}

\begin{Cor}[Existence under a jointly measurable density]
    \label{reg-cond-trans-kernel-existence-VIII}
    Let $\Xcal$, $\Ycal$, $\Zcal$ be \emph{arbitrary} measurable spaces, let $\mubf$ and $\nubf$ be
    $\sigma$-finite measures on $\Xcal$ and $\Ycal$, resp., and let
    \[\Kk(X,Y|Z): \, \Zcal \dshto  \Xcal \times \Ycal\]
    be a Markov kernel with $\Kk(X,Y|Z=z) \ll \mubf \otimes \nubf$ for every $z \in \Zcal$, admitting a density
    $k(x,y|z) \in [0,\infty)$ that is jointly measurable, i.e.\
    $\Bcal_\Xcal \otimes \Bcal_\Ycal \otimes \Bcal_\Zcal$-measurable. Put:
    \[ k(y|z) := \int_\Xcal k(x,y|z)\, \mubf(dx). \]
    Then there exists a conditional Markov kernel $\Kk(X|Y,Z)$ conditioned on $Y$ given $Z$, given by:
    \[ \Kk(X \in A|Y=y,Z=z) = \int_A k(x|y,z)\,\mubf(dx), \qquad
        k(x|y,z) := \frac{k(x,y|z)}{k(y|z)}, \]
    wherever $0 < k(y|z) < \infty$, and by an arbitrary fixed probability measure on $\Xcal$ elsewhere.
\begin{proof}
    We may assume $\Zcal \neq \emptyset$, as otherwise there is nothing to show, and then $\Xcal \neq \emptyset$,
    as otherwise no Markov kernel $\Zcal \dshto \Xcal \times \Ycal$ exists; so fix a point $x_0 \in \Xcal$.

    Since $\mubf$ is $\sigma$-finite and $k$ is jointly measurable, the measurability half of Tonelli's theorem
(in the form that for a $\sigma$-finite measure $\mubf$ and a
    $\Bcal_\Xcal \otimes \Bcal$-measurable $f \ge 0$ the partial integral $\int f(x,\cdot)\,\mubf(dx)$ is
    $\Bcal$-measurable, which requires no measure on the second factor at all) shows that for every
    $A \in \Bcal_\Xcal$ the map:
    \[ (y,z) \; \mapsto \; \int_A k(x,y|z)\,\mubf(dx) \]
    is measurable; the case $A = \Xcal$ gives in particular the measurability of $k(y|z)$. Tonelli's theorem
    also turns the absolute continuity hypothesis into the iterated identity:
    \[ \Kk(X \in A, Y \in B|Z=z) = \int_B \int_A k(x,y|z)\, \mubf(dx)\,\nubf(dy), \]
    for all $A \in \Bcal_\Xcal$, $B \in \Bcal_\Ycal$ and $z \in \Zcal$, and this is the only form in which the
    hypothesis will be used below, see \Cref{rem:density-hypothesis}. Taking $A = \Xcal$ there we get, for
    every $B \in \Bcal_\Ycal$ and $z \in \Zcal$:
    \[ \int_B k(y|z)\,\nubf(dy) = \Kk(Y \in B|Z=z), \]
    i.e.\ $k(\cdot|z)$ is a $\nubf$-density of $\Kk(Y|Z=z)$; in particular $k(\cdot|z)$ is $\nubf$-integrable.
    Hence the measurable set:
    \[ D := \lC (y,z) \in \Ycal \times \Zcal \st 0 < k(y|z) < \infty \rC \]
    satisfies $\Kk(Y \in D_z|Z=z) = 1$ for every $z$, where $D_z$ denotes the $z$-section. Indeed,
    $D_z^\cmpl$ is the union of $\lC k(\cdot|z) = \infty \rC$, which is $\nubf$-null by integrability, and of
    $\lC k(\cdot|z) = 0 \rC$, on which $k(\cdot|z)$ vanishes identically; so
    $\Kk(Y \in D_z^\cmpl|Z=z) = \int_{D_z^\cmpl} k(y|z)\,\nubf(dy) = 0$.
    Now define:
    \[ \Kk(X \in A|Y=y,Z=z) := \begin{cases}
        \frac{1}{k(y|z)}\int_A k(x,y|z)\, \mubf(dx) & \text{ if } (y,z) \in D,\\
        \deltabf_{x_0}(A) & \text{ else.}
    \end{cases}\]

    \emph{Probability measure in $A$.} For $(y,z) \in D$ the map $A \mapsto \int_A k(x,y|z)\,\mubf(dx)$ is a
    measure by monotone convergence, with total mass $k(y|z) \in (0,\infty)$; dividing by $k(y|z)$ makes it a
    probability measure. Otherwise it is the Dirac measure $\deltabf_{x_0}$.

    \emph{Measurability in $(y,z)$.} On the measurable set $D$ the map is the quotient of the two measurable
    functions exhibited above, with a nowhere vanishing finite denominator; on $D^\cmpl$ it is constant.

    \emph{Factorization.} Let $A \in \Bcal_\Xcal$, $B \in \Bcal_\Ycal$ and $z \in \Zcal$. Since $k(\cdot|z)$ is a
    $\nubf$-density of $\Kk(Y|Z=z)$, integration against $\Kk(Y|Z=z)$ is integration against
    $k(y|z)\,\nubf(dy)$. The part of $B$ outside $D_z$ contributes nothing: on $\lC k(\cdot|z)=0\rC$ the factor
    $k(y|z)$ vanishes and so does $\int_A k(x,y|z)\,\mubf(dx) \le k(y|z)$, while the remaining part
    $\lC k(\cdot|z)=\infty\rC$ of $D_z^\cmpl$ is $\nubf$-null. Hence:
    \begin{align*}
        \lp \Kk(X|Y,Z) \otimes \Kk(Y|Z) \rp (A \times B|z)
            &= \int_B \Kk(X \in A|Y=y,Z=z) \cdot k(y|z) \, \nubf(dy)\\
            &= \int_{B} \int_A k(x,y|z)\, \mubf(dx)\, \nubf(dy)\\
            &= \Kk(X \in A, Y \in B|Z=z).
    \end{align*}
    Since the measurable rectangles form a $\pi$-system generating $\Bcal_\Xcal \otimes \Bcal_\Ycal$ and both
    sides are probability measures, Dynkin's lemma, see \cite{Bog07} 1.9.3, gives:
    \[ \Kk(X,Y|Z) = \Kk(X|Y,Z) \otimes \Kk(Y|Z). \]

    Note that no property of $\Xcal$, $\Ycal$ or $\Zcal$ was used. The joint measurability of $k$ in
    $(x,y,z)$, on the other hand, is essential and not automatic: a density for each $z$ separately
    disintegrates every $\Kk(X,Y|Z=z)$ individually, but is not enough to make
    $(y,z) \mapsto \Kk(X \in A|Y=y,Z=z)$ measurable.
\end{proof}
\end{Cor}

\begin{Rem}
    \label{rem:disintegration-triple-scope}
    Being a disintegration triple only depends on the three measurable spaces up to isomorphism: if
    $\iota:\, \Ycal \to \Ycal'$ is an isomorphism of measurable spaces and $(\Xcal,\Ycal,\Zcal)$ is a
    disintegration triple, then so is $(\Xcal,\Ycal',\Zcal)$, since a Markov kernel
    $\Kk(X,Y'|Z)$ pulls back along $\id_\Xcal \times \iota$ to a Markov kernel $\Kk(X,Y|Z)$, and a conditional
    Markov kernel $\Kk(X|Y,Z)$ for the latter pushes forward to $\Kk(X|Y'=y',Z) := \Kk(X|Y=\iota^{-1}(y'),Z)$
    for the former. The same argument applies in the first and the third component.
    Weaker measurability requirements (e.g.\ conditional Markov kernels that are only universally
    measurable, which exist under correspondingly weaker hypotheses on the spaces, see
    \cite{Fad85, Bla63, Bla75, Rao05}) would enlarge the class of disintegration triples further; we do not
    pursue this here.
\end{Rem}

\section{Proofs - Join-Semi-Lattice Rules for Transitional Random Variables}

\label{sec:join-semi-lattice}

In this section we will collect properties of the relation $\ismapof_\Kk$ introduced in the main paper in \Cref{not:ismapof}.
For this let $\lp \Wcal \times \Tcal, \Kk(W|T) \rp$ be a transition probability space and $\Xk:\, \Wcal \times \Tcal \dshto \Xcal$ and
    $\Yk:\, \Wcal \times \Tcal \dshto \Ycal$ and $\Zk:\, \Wcal \times \Tcal \dshto \Zcal$
    and $\Uk:\, \Wcal \times \Tcal \dshto \Ucal$
    be transitional random variables, i.e.\ Markov kernels. We put:
    \[ \Kk(X,Y,Z,U|T) := \lp\Xk(X|W,T) \otimes \Yk(Y|W,T)  \otimes \Zk(Z|W,T)  \otimes \Uk(U|W,T)  \rp \circ \Kk(W|T). \]

 The relation $\ismapof_\Kk$ will be a main ingredient to show that \emph{transitional conditional independence}, see
 \Cref{def:transitional_conditional_independence}, forms a \emph{$T$-$\ast$-separoid}, see \Cref{def:t-k-separoid},
 i.e.\ it satisfies the asymmetric separoid rules of
 \Cref{thm:separoid_axioms-tci}, proven in
 \Cref{sec:proofs-separoid-rules-for-transitional-conditional-independence}.
 We also need to check the compatibility of
 $\ismapof_\Kk$ with the equivalence relation, $\cong$, of isomorphisms of measurable spaces. This will be done in
 \Cref{sec:proofs-separoid-rules-for-transitional-conditional-independence} in \Cref{lem:separoid-compatibility-1}.

\begin{Rem}
    \begin{enumerate}
        \item In general we do \emph{not} have: $\Xk \ismapof_\Kk \Xk$ for arbitrary Markov kernels. It will hold
            for \emph{deterministic} transitional random variables, see \Cref{lem:restricted-reflexivity}.
        \item In general we do \emph{not} have anti-symmetry, i.e.\ that:
            \[\Xk \ismapof_\Kk \Yk \ismapof_\Kk \Xk \implies \Xk = \Yk.\]
    \end{enumerate}
\end{Rem}

\begin{Not} Recall that we write:
     \begin{enumerate}
        \item $\Xk \ismapof_\Kk \Yk$ if there exists a measurable map $\varphi:\, \Ycal \to \Xcal$ such that:
            \[ \Kk(X,Y|T) = \deltabf_{\varphi}(X|Y) \otimes \Kk(Y|T).  \]
    \end{enumerate}
    We further define:
    \begin{enumerate}[resume]
        \item $\Xk \approx_\Kk \Yk \; :\iff \;  \Xk \ismapof_\Kk \Yk \ismapof_\Kk \Xk.$
        \item $\Xk \cong \Yk$ if there exists a measurable isomorphism $\varphi:\, \Ycal \to \Xcal$, i.e.\ a
            bijective measurable map with a measurable inverse, such that $\varphi_*\Yk=\Xk$, see also
            \Cref{sec:proofs-separoid-rules-for-transitional-conditional-independence}.
    \end{enumerate}
\end{Not}

\begin{Rem}[Commutativity and associativity of the join]
    \label{rem:join-commutative}
    By \Cref{lem:separoid-compatibility-1} points 2.\ and 3.\ we have
    $\lp \Xk \otimes \Yk \rp \cong \lp \Yk \otimes \Xk \rp$ and
    $\lp \lp \Xk \otimes \Yk \rp \otimes \Zk \rp \cong \lp \Xk \otimes \lp \Yk \otimes \Zk \rp\rp$,
    via the canonical isomorphisms $\Xcal \times \Ycal \cong \Ycal \times \Xcal$ and
    $(\Xcal \times \Ycal) \times \Zcal \cong \Xcal \times (\Ycal \times \Zcal)$ of measurable spaces.
    By \Cref{lem:separoid-compatibility-1} point 4.\ the relation $\ismapof_\Kk$, and hence also $\approx_\Kk$, is
    invariant under replacing its arguments by $\cong$-equivalent ones. In the following we will therefore reorder and
    regroup the factors of a product of transitional random variables without further mention.
\end{Rem}

The next \Cref{lem:ismapof-null-set} is crucial for most of the following results where $\ismapof_\Kk$ is involved.
Note that a similar result for $\ismapof_\Kk^\ast$ would not hold, i.e.\ where $\deltabf_\varphi(X|Y,T)$ would be
replaced by an arbitrary Markov kernel $\Xk(X|Y,T)$.

\begin{Prp}[Extension lemma for deterministic factorizations]
    \label{lem:ismapof-null-set}
    Consider a Markov kernel:
    \[
        \Kk(X,Y,Z|T):\, (\Tcal,\Bcal_\Tcal) \dshto (\Xcal \times \Ycal \times \Zcal, \Bcal_\Xcal\otimes\Bcal_\Ycal\otimes\Bcal_\Zcal),
    \]
    with the property that the marginal Markov kernel can be written as:
    \[
        \Kk(X,Y|T) = \deltabf_\varphi(X|Y,T) \otimes \Kk(Y|T),
    \]
    for some measurable map $\varphi:\,\Ycal \times \Tcal \to \Xcal$. Then the joint Markov kernel can be written as:
    \[
        \Kk(X,Y,Z|T) = \deltabf_\varphi(X|Y,T) \otimes \Kk(Y,Z|T).
    \]
\begin{proof}
    For $t \in \Tcal$ we abbreviate the following measurable map:
    \begin{align*}
        \varphi_t:\, \Ycal &\to \Xcal, & y & \mapsto \varphi(y,t)=:\varphi_t(y).
    \end{align*}
    Note that for $A \in \Bcal_\Xcal$ we then have:
    \begin{align*}
        \varphi_t^{-1}(A) &= \lC y \in \Ycal \st \varphi(y,t) \in A \rC = \varphi^{-1}(A)_t \in \Bcal_\Ycal.
    \end{align*}
    By assumption we have for every $A \in \Bcal_\Xcal$ and $B \in \Bcal_\Ycal$ and $t \in \Tcal$ the relation:
    \begin{align*}
       \Kk(X \in A, Y \in B, Z \in \Zcal|T=t)
       &= \Kk(X \in A,Y \in B|T=t) \\
       & = \int_B \deltabf_\varphi(X\in A|Y=y,T=t) \, \Kk(Y \in dy|T=t) \\
       & = \int_B \I_{\varphi_t^{-1}(A)}(y) \, \Kk(Y \in dy|T=t) \\
       &= \Kk(Y \in \underbrace{B \cap \varphi_t^{-1}(A)}_{=:D_t \in \Bcal_\Ycal} |T=t) \\
       &= \Kk(X \in \Xcal, Y \in D_t, Z \in \Zcal|T=t).
    \end{align*}
    Similarly, by replacing $B$ with $D_t$ in the above formula and noting that $D_t \cap \varphi_t^{-1}(A) = D_t$, we also get:
    \begin{align*}
        \Kk(X \in A, Y \in D_t, Z \in \Zcal|T=t) &= \Kk(X \in \Xcal, Y \in D_t, Z \in \Zcal|T=t) \\
                                                &= \Kk(X \in A, Y \in B, Z \in \Zcal|T=t).
    \end{align*}
    Note that:
    \begin{align*}
        \lp A \times B \times \Zcal \rp \cap \lp \Xcal \times D_t \times \Zcal \rp &= A \times D_t \times \Zcal,
    \end{align*}
    since $D_t \ins B$.
    So, if we now put:
    \begin{align*}
        N_t &:= \lp A \times B \times \Zcal \rp \syd \lp \Xcal \times D_t \times \Zcal \rp,
    \end{align*}
    then $N_t \in \Bcal_\Xcal \otimes \Bcal_\Ycal \otimes \Bcal_\Zcal$ and we get by the above calculations that:
    \begin{align*}
        &\quad \Kk\lp (X,Y,Z) \in N_t |T=t\rp \\
        &= \Kk(X \in A, Y \in B, Z \in \Zcal|T=t) + \Kk(X \in \Xcal, Y \in D_t, Z \in \Zcal|T=t) \\
        &\quad - 2 \cdot \Kk(X \in A, Y \in D_t, Z \in \Zcal|T=t) \\
        &= 0.
    \end{align*}
    Now let $C \in \Bcal_\Zcal$, recall that $D_t \ins B$ and consider the measurable set:
    \begin{align*}
        M_t &:= \lp A \times B \times C \rp \syd \lp \Xcal \times D_t \times C \rp \\
            &= \lp A \times (B\sm D_t) \times C \rp \dcup \lp A^\cmpl \times D_t \times C \rp \\
            &\ins \lp A \times (B\sm D_t) \times \Zcal \rp \dcup \lp A^\cmpl \times D_t \times \Zcal \rp \\
            &= \lp A \times B \times \Zcal \rp \syd \lp \Xcal \times D_t \times \Zcal \rp  \\
            &= N_t.
    \end{align*}
    This then implies that:
    \begin{align*}
        \Kk\lp (X,Y,Z) \in M_t |T=t\rp & \le \Kk\lp (X,Y,Z) \in N_t |T=t\rp =0.
    \end{align*}
    So $M_t$ is a (measurable) $\Kk(X,Y,Z|T=t)$-null set, which implies the following equality:
    \begin{align*}
       \Kk(X \in A, Y \in B, Z \in C|T=t)
       &= \Kk(X \in \Xcal, Y \in D_t, Z \in C|T=t) \\
       &= \Kk(Y \in B \cap \varphi_t^{-1}(A), Z \in C |T=t) \\
       &= \int \I_B(y) \cdot \I_{\varphi_t^{-1}(A)}(y) \cdot \I_C(z) \, \Kk(Y \in dy, Z \in dz|T=t) \\
       &= \int_{B \times C} \deltabf_\varphi(X\in A|Y=y,T=t) \, \Kk(Y \in dy, Z \in dz|T=t) \\
       &= \lp \deltabf_\varphi(X|Y,T) \otimes \Kk(Y,Z|T) \rp (A \times B \times C|t).
    \end{align*}
    Since this holds for all $A \in \Bcal_\Xcal$, $B \in \Bcal_\Ycal$, $C \in \Bcal_\Zcal$ and $t \in \Tcal$, and since
    the measurable rectangles $A \times B \times C$ form a $\cap$-stable generator of
    $\Bcal_\Xcal\otimes\Bcal_\Ycal\otimes\Bcal_\Zcal$ containing $\Xcal \times \Ycal \times \Zcal$, Dynkin's lemma gives
    the equality:
    \[
        \Kk(X,Y,Z|T) = \deltabf_\varphi(X|Y,T) \otimes \Kk(Y,Z|T).
    \]
    This shows the claim.
\end{proof}
\end{Prp}

\begin{Rem}
    \Cref{lem:ismapof-null-set} is stated for a map $\varphi$ that may also depend on $t \in \Tcal$; the case
    $\varphi:\, \Ycal \to \Xcal$ needed for $\ismapof_\Kk$ is the special case of a $t$-independent $\varphi$.
    Note also that the proof does \emph{not} require the diagonal of $\Xcal$ to be measurable, i.e.\ $\Xcal$ does not
    need to be countably separated: instead of arguing that $\lC (x,y,t) \st x \neq \varphi(y,t)\rC$ is a null set, we
    only use the null sets $M_t$ built from the symmetric differences of measurable rectangles.
\end{Rem}

\begin{Lem}[Exact factorizations are almost-sure factorizations]
    \label{lem:exact-to-as}
    Let $X:\, \Wcal \times \Tcal \to \Xcal$ and $Y:\, \Wcal \times \Tcal \to \Ycal$ be measurable maps,
    considered as the deterministic transitional random variables $\Xk = \deltabf(X|W,T)$ and
    $\Yk = \deltabf(Y|W,T)$ on a transition probability space $\lp \Wcal \times \Tcal, \Kk(W|T) \rp$.
    If $X \ismapof Y$, i.e.\ if $X = \varphi \circ Y$ \emph{pointwise} for a measurable map
    $\varphi:\, \Ycal \to \Xcal$, then $\Xk \ismapof_\Kk \Yk$, i.e.:
    \[ \Kk(X,Y|T) = \deltabf_\varphi(X|Y) \otimes \Kk(Y|T). \]
\begin{proof}
    Fix $t \in \Tcal$, $A \in \Bcal_\Xcal$ and $B \in \Bcal_\Ycal$. Since $X = \varphi \circ Y$ we have,
    pointwise in $w \in \Wcal$:
    \[ \I_A(X(w,t)) \cdot \I_B(Y(w,t)) = \I_{\varphi^{-1}(A)}(Y(w,t)) \cdot \I_B(Y(w,t))
       = \I_{\varphi^{-1}(A) \cap B}(Y(w,t)), \]
    which is exactly the point where determinism enters. Integrating against $\Kk(W|T=t)$ gives:
    \begin{align*}
        \Kk(X \in A, Y \in B|T=t) &= \Kk(Y \in \varphi^{-1}(A) \cap B|T=t) \\
        &= \int \I_B(y) \cdot \deltabf_\varphi(X \in A|Y=y) \, \Kk(Y \in dy|T=t),
    \end{align*}
    because $\deltabf_\varphi(X \in A|Y=y) = \I_{\varphi^{-1}(A)}(y)$.
    The measurable rectangles $A \times B$ form a $\cap$-stable generator of
    $\Bcal_\Xcal \otimes \Bcal_\Ycal$ and both sides are probability measures in $(X,Y)$ for each fixed $t$, so
    Dynkin's lemma (see \cite{Kle20} Thm.\ 1.19) gives the claimed identity of Markov kernels.
\end{proof}
\end{Lem}

\begin{Rem}
    \label{rem:exact-to-as}
    \Cref{lem:exact-to-as} is what allows us to feed a \emph{pointwise} functional relation $X = \varphi \circ Y$
    into every statement below whose hypothesis is the almost-sure relation $\ismapof_\Kk$; we will use it without
    further mention. Note that it is \emph{not} an instance of \Cref{lem:restricted-reflexivity}, which is the
    special case $\varphi = \id$.
\end{Rem}

\begin{Lem}[Product extension]
    \label{lem:join-upper-bound-2}
    We always have the implication:
    \[ \Xk \ismapof_\Kk \Yk \implies \Xk \ismapof_\Kk \Yk \otimes \Zk.  \]
    \begin{proof} By assumption we have for some measurable map $\varphi:\, \Ycal \to \Xcal$ the factorization:
        \begin{align*}
        && \Kk(X,Y|T) &=\deltabf_\varphi(X|Y) \otimes \Kk(Y|T)  \\
           & \xRightarrow{\text{\Cref{lem:ismapof-null-set}}}
            & \Kk(X,Y,Z|T) &= \deltabf_\varphi(X|Y) \otimes \Kk(Y,Z|T).
        \end{align*}
        Since $\varphi$ does not depend on the $\Zcal$-component this exhibits $\Xk$ as a measurable map of
        $\Yk \otimes \Zk$, which shows the claim.
    \end{proof}
\end{Lem}

\begin{Lem}[Bottom element]
    \label{lem:bottom-element}
    We always have:
    \[ \deltabf_\ast \ismapof_\Kk \Xk.\]
    \begin{proof}
        Consider the constant map $\ast:\, \Xcal \to \Asterisk$, which is measurable. With this we get:
        \begin{align*}
            \Kk(\ast,X|T) &=\deltabf_\ast \otimes \Kk(X|T).
        \end{align*}
        This shows the claim.
    \end{proof}
\end{Lem}

\begin{Lem}[Transitivity]
    \label{lem:transitivity}
    We always have the implication:
    \[ \Xk \ismapof_\Kk \Yk \ismapof_\Kk \Zk \implies \Xk \ismapof_\Kk \Zk.\]
    \begin{proof}
        By assumption we have measurable maps $\psi:\,\Zcal \to \Ycal$ and $\varphi:\,\Ycal \to \Xcal$ such that:
        \begin{align*}
           && \Kk(Y,Z|T) &= \deltabf_\psi(Y|Z) \otimes \Kk(Z|T), \\
           && \Kk(X,Y|T) &= \deltabf_\varphi(X|Y) \otimes \Kk(Y|T).
        \end{align*}
        With \Cref{lem:ismapof-null-set} applied to $\varphi$ we then get:
        \begin{align*}
            \Kk(X,Y,Z|T) &= \deltabf_\varphi(X|Y) \otimes \Kk(Y,Z|T) \\
            &=\deltabf_\varphi(X|Y) \otimes \deltabf_\psi(Y|Z) \otimes \Kk(Z|T),
        \end{align*}
        which by marginalizing out $Y$ implies:
        \begin{align*}
            \Kk(X,Z|T) &= \lp \deltabf_{\varphi}(X|Y) \circ \deltabf_{\psi}(Y|Z) \rp \otimes \Kk(Z|T)\\
                       &=  \deltabf_{\varphi \circ \psi}(X|Z) \otimes \Kk(Z|T).
        \end{align*}
        Since the composition $\varphi \circ \psi:\, \Zcal \to \Xcal$ is also measurable, the claim is shown.
    \end{proof}
\end{Lem}

\begin{Lem}[Product stays bounded]
    \label{lem:prod-bounded}
    We always have the implication:
    \[ \Xk \ismapof_\Kk \Zk \quad \land \quad \Yk \ismapof_\Kk \Zk \quad \implies \quad \Xk \otimes \Yk \ismapof_\Kk \Zk. \]
    \begin{proof} By the assumptions we have:
        \begin{align*}
           && \Kk(Y,Z|T) &= \deltabf_\psi(Y|Z) \otimes \Kk(Z|T), \\
           &&  \Kk(X,Z|T) &= \deltabf_\varphi(X|Z) \otimes \Kk(Z|T),
        \end{align*}
        for some measurable maps $\varphi:\, \Zcal \to \Xcal$ and $\psi:\, \Zcal \to \Ycal$.
        \Cref{lem:ismapof-null-set}, applied to $\varphi$ with the roles $(X,Z,Y)$, then implies:
        \begin{align*}
            \Kk(X,Y,Z|T) &= \deltabf_\varphi(X|Z) \otimes \Kk(Y,Z|T) \\
            &=\deltabf_\varphi(X|Z) \otimes \deltabf_\psi(Y|Z) \otimes \Kk(Z|T) \\
            &=\deltabf_{\varphi \times \psi}(X,Y|Z) \otimes \Kk(Z|T).
        \end{align*}
        This shows the claim, as $\varphi \times \psi:\, \Zcal \to \Xcal \times \Ycal$,
        $z \mapsto (\varphi(z),\psi(z))$, is again a measurable map.
    \end{proof}
\end{Lem}

\begin{Lem}[Product compatibility]
    \label{lem:products-bound}
    We always have the implication:
    \[ \Xk \ismapof_\Kk \Zk \quad \land \quad \Yk \ismapof_\Kk \Uk \quad \implies \quad \Xk \otimes \Yk \ismapof_\Kk \Zk \otimes \Uk. \]
    \begin{proof}
        $\Xk \ismapof_\Kk \Zk$ implies $\Xk \ismapof_\Kk \Zk\otimes\Uk$ by \Cref{lem:join-upper-bound-2}.
        Similarly, $\Yk \ismapof_\Kk \Uk$ implies $\Yk \ismapof_\Kk \Uk \otimes \Zk \cong \Zk\otimes\Uk$, using
        \Cref{rem:join-commutative}. By \Cref{lem:prod-bounded} we then get the claim:
        $\quad\Xk \otimes \Yk \ismapof_\Kk \Zk \otimes \Uk.$
    \end{proof}
\end{Lem}

The remaining results of this section concern \emph{deterministic} transitional random variables, i.e.\ those of the form
$\Xk = \deltabf(X|W,T)$ for a measurable map $X:\, \Wcal \times \Tcal \to \Xcal$. Note that for such $\Xk$ and any
$A \in \Bcal_\Xcal$ we have $\Xk(X \in A|W=w,T=t) = \I_A(X(w,t)) \in \{0,1\}$, which is what makes the following
computations work.

\begin{Lem}[Restricted reflexivity]
    \label{lem:restricted-reflexivity}
    If $X:\,\Wcal \times \Tcal \to \Xcal$ is a measurable map and $\Xk=\deltabf(X|W,T)$, then we have:
    \[ \Xk \ismapof_\Kk \Xk.\]
\begin{proof}
    Let $X_1$, $X_2$ denote two copies of $X$ and consider the measurable map $\id:\, \Xcal \to \Xcal$.
    For $A, B \in \Bcal_\Xcal$ and $t \in \Tcal$ we compute:
    \begin{align*}
        & \lp \deltabf_{\id}(X_1|X_2) \otimes \Kk(X_2|T) \rp (A \times B, t) \\
        &= \int_B \deltabf_{\id}(X_1 \in A|X_2=x_2)\, \Kk(X_2 \in dx_2|T=t) \\
        &= \int \I_B(x_2) \cdot \I_A(x_2) \, \Kk(X_2 \in dx_2|T=t) \\
        &= \Kk(X_2 \in A \cap B|T=t) \\
        &= \int \deltabf(X \in A \cap B|W=w,T=t) \, \Kk(W \in dw|T=t) \\
        &= \int \I_{A \cap B}(X(w,t)) \, \Kk(W \in dw|T=t) \\
        &= \int \I_{A}(X(w,t)) \cdot \I_{B}(X(w,t)) \, \Kk(W \in dw|T=t) \\
        &= \int \deltabf(X_1 \in A|W=w,T=t) \cdot \deltabf(X_2 \in B|W=w,T=t) \, \Kk(W \in dw|T=t) \\
        &= \Kk(X_1 \in A, X_2 \in B|T=t),
    \end{align*}
    where the crucial step is the multiplicativity $\Xk(X \in A \cap B|W=w,T=t) = \Xk(X \in A|W=w,T=t) \cdot
    \Xk(X \in B|W=w,T=t)$, which holds precisely because $\Xk$ is a Dirac kernel, i.e.\ deterministic; for a general
    Markov kernel it fails. Since the measurable rectangles $A \times B$ form a $\cap$-stable generator of
    $\Bcal_\Xcal \otimes \Bcal_\Xcal$, Dynkin's lemma gives:
    \[ \Kk(X_1,X_2|T) = \deltabf_{\id}(X_1|X_2) \otimes \Kk(X_2|T), \]
    which is the claim.
\end{proof}
\end{Lem}

\begin{Lem}[Join is upper bound]
    \label{lem:join-upper-bound}
    If $X:\,\Wcal \times \Tcal \to \Xcal$ is a measurable map, $\Xk=\deltabf(X|W,T)$ and $\Zk$ is an arbitrary
    transitional random variable, then we have:
    \[ \Xk \ismapof_\Kk \Xk \otimes \Zk. \]
\begin{proof}
    Consider the projection:
    \begin{align*}
        \pr:\,   \Xcal \times \Zcal &\to \Xcal, & (x,z) & \mapsto x,
    \end{align*}
    which is measurable. With $X_1$, $X_2$ two copies of $X$ we compute for
    $A,B \in \Bcal_\Xcal$, $C \in \Bcal_\Zcal$ and $t \in \Tcal$:
     \begin{align*}
        & \lp \deltabf_{\pr}(X_1|X_2,Z) \otimes \Kk(X_2,Z|T) \rp (A \times B \times C, t) \\
        &= \int_{B \times C} \deltabf_{\pr}(X_1 \in A|X_2=x_2,Z=z)\, \Kk(X_2 \in dx_2,Z \in dz|T=t) \\
        &= \int \I_B(x_2) \cdot \I_C(z) \cdot \I_A(x_2) \, \Kk(X_2 \in dx_2,Z \in dz|T=t) \\
        &= \Kk(X_2 \in A \cap B, Z \in C|T=t) \\
        &= \int \Xk(X \in A \cap B|W=w,T=t) \cdot \Zk(Z \in C|W=w,T=t) \, \Kk(W \in dw|T=t) \\
        &= \int \I_{A \cap B}(X(w,t)) \cdot \Zk(Z \in C|W=w,T=t) \, \Kk(W \in dw|T=t) \\
        &= \int \I_{A}(X(w,t)) \cdot \I_{B}(X(w,t))  \cdot \Zk(Z \in C|W=w,T=t) \, \Kk(W \in dw|T=t) \\
        &= \int \deltabf(X_1 \in A|W=w,T=t) \cdot \deltabf(X_2 \in B|W=w,T=t) \cdot  \\
        &\qquad \qquad \Zk(Z \in C|W=w,T=t) \, \Kk(W \in dw|T=t) \\
        &= \Kk(X_1 \in A, X_2 \in B, Z \in C|T=t).
     \end{align*}
     Again by Dynkin's lemma this implies:
        \[ \Kk(X_1,X_2,Z|T) = \deltabf_{\pr}(X_1|X_2,Z) \otimes \Kk(X_2,Z|T), \]
     and thus $\Xk \ismapof_\Kk \Xk \otimes \Zk$.
\end{proof}
\end{Lem}

\begin{Rem}
    \label{rem:join-upper-bound-both}
    If both $\Xk=\deltabf(X|W,T)$ and $\Zk=\deltabf(Z|W,T)$ are deterministic then \Cref{lem:join-upper-bound},
    applied twice and combined with \Cref{rem:join-commutative}, gives both:
    \[ \Xk \ismapof_\Kk \Xk \otimes \Zk \qquad \text{ and } \qquad \Zk \ismapof_\Kk \Zk \otimes \Xk \cong \Xk \otimes \Zk. \]
\end{Rem}

\begin{Lem}[Bottom element is neutral]
    \label{lem:bot-neutr}
    If $X:\,\Wcal \times \Tcal \to \Xcal$ is a measurable map and $\Xk=\deltabf(X|W,T)$
    then we have:
    \[ \deltabf_\ast \otimes \Xk \approx_\Kk \Xk. \]
\begin{proof}
    Consider the following mutually inverse measurable maps:
    \begin{align*}
      \pr_\Xcal:\,  \Asterisk \times \Xcal &\to \Xcal, & (\ast,x) &\mapsto x, \\
       \id_\Xcal^\ast:\, \Xcal & \to \Asterisk \times \Xcal, & x & \mapsto (\ast,x).
    \end{align*}
    Using \Cref{lem:restricted-reflexivity} in the second step we get:
     \begin{align*}
         \Kk(\ast,X_1,X_2|T) &= \deltabf_\ast \otimes \Kk(X_1,X_2|T) \\
                             &= \deltabf_\ast \otimes \deltabf_{\id}(X_1|X_2) \otimes \Kk(X_2|T) \\
         &= \deltabf_{\id^\ast_\Xcal}(\ast,X_1|X_2) \otimes \Kk(X_2|T),
     \end{align*}
     which shows $\deltabf_\ast \otimes \Xk \ismapof_\Kk \Xk$.
     For the converse, \Cref{lem:join-upper-bound} with $\Zk:=\deltabf_\ast$ gives
     $\Xk \ismapof_\Kk \Xk \otimes \deltabf_\ast$, and
     $\Xk \otimes \deltabf_\ast \cong \deltabf_\ast \otimes \Xk$ by \Cref{rem:join-commutative}.
     Together this shows $\deltabf_\ast \otimes \Xk \approx_\Kk \Xk$.
    Alternatively, this follows from $\deltabf_\ast \otimes \Xk \cong \Xk$, see
    \Cref{lem:separoid-compatibility-1} 8., which holds for \emph{arbitrary} $\Xk$, together with
    \Cref{lem:restricted-reflexivity} and \Cref{lem:separoid-compatibility-1} 4.; the passage from $\cong$ to
    $\approx_\Kk$ is what needs the reflexivity available only for deterministic $\Xk$.
\end{proof}
\end{Lem}

\begin{Lem}[Idempotency]
    \label{lem:idempotency}
    If $X:\,\Wcal \times \Tcal \to \Xcal$ is a measurable map and $\Xk=\deltabf(X|W,T)$ then we have:
    \[ \Xk \otimes \Xk \approx_\Kk \Xk. \]
\begin{proof}
    Consider the measurable diagonal map:
    \begin{align*}
        \Delta:\, \Xcal &\to \Xcal \times \Xcal, & x & \mapsto (x,x).
    \end{align*}
    With $X_1$, $X_2$, $X_3$ three copies of $X$ we compute for $A,B,C \in \Bcal_\Xcal$ and $t \in \Tcal$:
     \begin{align*}
        & \lp \deltabf_{\Delta}(X_1,X_2|X_3) \otimes \Kk(X_3|T) \rp (A \times B \times C, t) \\
        &= \int_C \deltabf_{\Delta}(X_1 \in A,X_2 \in B|X_3=x_3)\,  \Kk(X_3 \in dx_3|T=t) \\
        &= \int \I_C(x_3) \cdot \I_{A \times B}(x_3,x_3)\,  \Kk(X_3 \in dx_3|T=t) \\
        &= \int \I_C(x_3) \cdot \I_{A \cap B}(x_3)\,  \Kk(X_3 \in dx_3|T=t) \\
        &= \Kk(X_3 \in A \cap B \cap C|T=t) \\
        &= \Kk(W \in X^{-1}(A \cap B \cap C)|T=t) \\
        &= \Kk(W \in X^{-1}(A) \cap X^{-1}(B) \cap X^{-1}(C)|T=t) \\
        &= \Kk(X_1 \in A, X_2 \in B, X_3 \in C|T=t),
     \end{align*}
     where the second to last step is the general preimage identity
     $X^{-1}(A \cap B \cap C)=X^{-1}(A)\cap X^{-1}(B)\cap X^{-1}(C)$, and the last step uses that $X$ is deterministic,
     so that $X_1$, $X_2$, $X_3$ are $\Kk$-almost surely equal (here we suppress the dependence of $X$ on
     $t \in \Tcal$ in the notation $X^{-1}$).
     By Dynkin's lemma this implies:
        \[ \Kk(X_1,X_2,X_3|T) = \deltabf_{\Delta}(X_1,X_2|X_3) \otimes \Kk(X_3|T), \]
     and thus $\Xk \otimes \Xk \ismapof_\Kk \Xk$.
     The converse, $\Xk \ismapof_\Kk \Xk \otimes \Xk$, is \Cref{lem:join-upper-bound} with $\Zk:=\Xk$.
\end{proof}
\end{Lem}

\begin{Thm}
    \label{thm:join-semi-lattice}
    Let $\lp \Wcal \times \Tcal, \Kk(W|T) \rp$ be a transition probability space and $X:\, \Wcal \times \Tcal \to \Xcal$ and
    $Y:\, \Wcal \times \Tcal \to \Ycal$ and $Z:\, \Wcal \times \Tcal \to \Zcal$
    be measurable maps.
    We put: $\Xk := \deltabf(X|W,T)$ and $\Yk:=\deltabf(Y|W,T)$ and $\Zk:=\deltabf(Z|W,T)$ and:
    \[ \Kk(X,Y,Z|T) := \lp\Xk(X|W,T) \otimes \Yk(Y|W,T)  \otimes \Zk(Z|W,T)   \rp \circ \Kk(W|T). \]
    We then have:
    \begin{enumerate}
        \item Reflexivity \ref{lem:restricted-reflexivity}: $\Xk \ismapof_\Kk \Xk$.
        \item Transitivity \ref{lem:transitivity}: $\Xk \ismapof_\Kk \Yk \ismapof_\Kk \Zk \implies \Xk \ismapof_\Kk \Zk$.
        \item Almost-sure anti-symmetry (per definition):
            \[ \Xk \ismapof_\Kk \Yk \ismapof_\Kk \Xk \implies \Xk \approx_\Kk \Yk.\]
        \item Closed under join: $ \Xk \otimes \Yk = \deltabf(X,Y|W,T)$.
        \item Join is upper bound \ref{lem:join-upper-bound}: $\Xk \ismapof_\Kk \Xk \otimes \Yk$
             and $\Yk \ismapof_\Kk \Xk \otimes \Yk$.
        \item Join is smallest upper bound \ref{lem:prod-bounded}:
            \[ \Xk \ismapof_\Kk \Zk\quad \land \quad \Yk \ismapof_\Kk \Zk \quad\implies \quad\Xk \otimes \Yk \ismapof_\Kk \Zk. \]
        \item Bottom element \ref{lem:bottom-element}: $\deltabf_\ast \ismapof_\Kk \Xk$.
        \item Bottom element is neutral \ref{lem:bot-neutr}: $\Xk \approx_\Kk \deltabf_\ast \otimes \Xk$.
        \item Idempotent \ref{lem:idempotency}: $\Xk \approx_\Kk\Xk \otimes \Xk$.
    \end{enumerate}
\begin{proof}
    Point 1.\ is \Cref{lem:restricted-reflexivity}, which is applicable since $\Xk$ is deterministic, and point 2.\
    is \Cref{lem:transitivity}, which holds for arbitrary Markov kernels.
    Point 3.\ holds by the definition of $\approx_\Kk$.\\
    Point 4.: for $A \in \Bcal_\Xcal$, $B \in \Bcal_\Ycal$ and $(w,t) \in \Wcal \times \Tcal$ we have:
    \begin{align*}
        \lp \deltabf(X|W,T) \otimes \deltabf(Y|W,T) \rp (A \times B|w,t)
        &= \I_A(X(w,t)) \cdot \I_B(Y(w,t)) \\
        &= \I_{A \times B}\lp (X,Y)(w,t) \rp \\
        &= \deltabf\lp (X,Y)|W=w,T=t\rp (A \times B),
    \end{align*}
    and since the measurable rectangles form a $\cap$-stable generator of $\Bcal_\Xcal \otimes \Bcal_\Ycal$, Dynkin's
    lemma gives $\Xk \otimes \Yk = \deltabf(X,Y|W,T)$. In particular the join of two \emph{deterministic} transitional
    random variables is again deterministic, i.e.\ the class in question is closed under $\otimes$.\\
    Point 5.: $\Xk \ismapof_\Kk \Xk \otimes \Yk$ is \Cref{lem:join-upper-bound} with $\Zk := \Yk$; and
    $\Yk \ismapof_\Kk \Yk \otimes \Xk$ is \Cref{lem:join-upper-bound} with the roles of $X$ and $Y$ exchanged, which by
    \Cref{rem:join-commutative} gives $\Yk \ismapof_\Kk \Xk \otimes \Yk$, see \Cref{rem:join-upper-bound-both}.\\
    Point 6.\ is \Cref{lem:prod-bounded}, point 7.\ is \Cref{lem:bottom-element}, point 8.\ is
    \Cref{lem:bot-neutr} and point 9.\ is \Cref{lem:idempotency}.
\end{proof}
\end{Thm}

\begin{Cor}[The bounded join-semi-lattice of deterministic transitional random variables]
   \label{cor:join-semi-lattice}
    The class of transitional random variables
    of the form $\Xk = \deltabf(X|W,T)$, $\Yk=\deltabf(Y|W,T)$, $\Zk=\deltabf(Z|W,T)$, etc.,
    for some measurable maps $X$, $Y$, $Z$, etc., on the transition probability space
    $(\Wcal \times \Tcal, \Kk(W|T))$ together with the relation $\ismapof_\Kk$, join $\otimes$ and bottom element $\deltabf_\ast$ forms
    a \emph{bounded join-semi-lattice} modulo almost-sure anti-symmetry $\approx_\Kk$ (and up to the fact that such a class might not be a set).
\begin{proof}
    By \Cref{thm:join-semi-lattice} points 1., 2.\ and 3., the relation $\ismapof_\Kk$ is reflexive and transitive on
    this class and thus a pre-order, which induces a partial order on the quotient by $\approx_\Kk$; note that
    $\approx_\Kk$ is an equivalence relation by points 1.\ and 2., and that $\ismapof_\Kk$ is well defined on the
    quotient, again by point 2.
    By point 4.\ the class is closed under $\otimes$. The operation $\otimes$ descends to the quotient: if
    $\Xk \approx_\Kk \Xk'$ and $\Yk \approx_\Kk \Yk'$ then \Cref{lem:products-bound}, applied once in each
    direction, gives $\Xk \otimes \Yk \approx_\Kk \Xk' \otimes \Yk'$. On the quotient $\otimes$ is commutative and
    associative, since by \Cref{rem:join-commutative} the corresponding products are $\cong$-equivalent, and
    $\cong$-equivalent transitional random variables are $\approx_\Kk$-equivalent by point 1.\ together with
    \Cref{lem:separoid-compatibility-1} point 4.
    By points 5.\ and 6.\ the class of $\Xk \otimes \Yk$ is the least upper bound of $\Xk$ and $\Yk$ w.r.t.\
    $\ismapof_\Kk$, so the quotient is a join-semi-lattice with join $\otimes$.
    Finally, by points 7.\ and 8.\ the element $\deltabf_\ast$ is a smallest element and neutral for $\otimes$, so the
    join-semi-lattice is bounded from below. Point 9.\ is the idempotency required of a join.
\end{proof}
\end{Cor}
\section{Proofs - Separoid Rules for Transitional Conditional Independence}
\label{sec:proofs-separoid-rules-for-transitional-conditional-independence}

In this section we want to prove that the class of \emph{transitional random variables} together with
the equivalence relation, $\cong$, isomorphism of measurable spaces, the relation $\ismapof_\Kk$, 
the product $\otimes$ and the ternary relation $\Indep_\Kk$ of \emph{transitional conditional independence}, see 
\Cref{def:transitional_conditional_independence}, 
satisfies all the asymmetric separoid rules of \Cref{thm:separoid_axioms-tci}, at least when restricted to codomains that form disintegration triples, e.g.\ standard measurable spaces; restricted to transitional random variables with standard codomains it forms a \emph{$\Tk$-$\deltabf_\ast$-separoid} (or in different symbols: \emph{$T$-$\ast$-separoid}), see \Cref{def:t-k-separoid} and \Cref{cor:t-star-separoid}.

For this let $\lp \Wcal \times \Tcal, \Kk(W|T) \rp$ be a transition probability space and $\Xk:\, \Wcal \times \Tcal \dshto \Xcal$ and 
    $\Yk:\, \Wcal \times \Tcal \dshto \Ycal$ and $\Zk:\, \Wcal \times \Tcal \dshto \Zcal$ and 
    $\Uk:\, \Wcal \times \Tcal \dshto \Ucal$  be Markov kernels. 
We denote by $T:\, \Wcal \times \Tcal \to \Tcal$ the canonical projection map and 
$\Tk(T|W,T) := \deltabf(T|W,T)$. We also consider the constant map $\ast:\, \Wcal \times \Tcal \to \Asterisk$ and 
$\deltabf_\ast =\deltabf(\ast|W,T)$ the corresponding Markov kernel.
    We put:
    \[ \Kk(X,Y,Z,U|T) := \lp\Xk(X|W,T) \otimes \Yk(Y|W,T)  \otimes \Zk(Z|W,T) \otimes \Uk(U|W,T)  \rp \circ \Kk(W|T), \]
 or similarly if more or other Markov kernels are involved.

 Recall that we say that $\Xk$ is \emph{transitionally independent of $\Yk$ conditioned on $\Zk$ w.r.t.\ $\Kk=\Kk(W|T)$}, in symbols:
\[\Xk \Indep_{\Kk} \Yk \given \Zk, \]
if there exists a Markov kernel $\Qk(X|Z)$ such that: 
\begin{align*} 
    \Kk(X,Y,Z|T) = \Qk(X|Z) \otimes \Kk(Y,Z|T),
\end{align*}
where $\Kk(Y,Z|T)$ is the marginal of $\Kk(X,Y,Z|T)$.

\begin{Not} Recall that we write:
     \begin{enumerate}
        \item $\Xk \ismapof_\Kk \Yk$ if there exists a measurable map $\varphi:\, \Ycal \to \Xcal$ such that:
            \[ \Kk(X,Y|T) = \deltabf_{\varphi}(X|Y) \otimes \Kk(Y|T).  \]
        \item $\Xk \approx_\Kk \Yk \; :\iff \;  \Xk \ismapof_\Kk \Yk \ismapof_\Kk \Xk.$ 
    \end{enumerate}
    We further define:
    \begin{enumerate}[resume]
        \item $ \Xk \cong \Yk$ if  there exists a measurable isomorphism $\varphi:\, \Ycal \to \Xcal$, i.e.\ a 
            bijective measurable map with a measurable inverse,
            such that: $\varphi_*\Yk = \Xk$.
    \end{enumerate}
\end{Not}

We first need to check that $\ismapof_\Kk$, $\otimes$, $\cong$, $\deltabf_\ast$, $\Indep_\Kk$
are all sufficiently compatible with each other. This will be done in the next Lemma.

\begin{Lem}[Compatibility of $\ismapof_\Kk$, $\otimes$, $\cong$, $\deltabf_\ast$, $\Indep_\Kk$]
    \label{lem:separoid-compatibility-1}
    We have the following:
\begin{enumerate}
    \item $(\Xk \cong \Xk') \,\land \, (\Yk \cong \Yk')\, \implies \,(\Xk \otimes \Yk) \cong (\Xk' \otimes \Yk')$.
    \item[] \textit{Proof.} With isomorphisms $\varphi:\,\Xcal \cong \Xcal'$ and $\psi:\, \Ycal \cong \Ycal'$ with $\varphi_*\Xk=\Xk'$ and 
        $\psi_*\Yk=\Yk'$ we get: 
        $\quad (\varphi \times \psi)_*(\Xk\otimes \Yk) = (\varphi_* \Xk) \otimes(\psi_*\Yk) = \Xk'\otimes \Yk'.$
    \item $(\Xk \otimes \Yk) \cong (\Yk \otimes \Xk)$.
    \item[] \textit{Proof.} Use the isomorphism: $\Xcal \times \Ycal \cong \Ycal \times \Xcal$ with $(x,y) \mapsto (y,x)$.
    \item $(\Xk \otimes \Yk ) \otimes \Zk \cong \Xk \otimes (\Yk  \otimes \Zk).$
    \item[] \textit{Proof.} Use the isomorphism: $\id:\, (\Xcal \times \Ycal) \times \Zcal \cong \Xcal \times (\Ycal \times \Zcal)$. 
    \item $(\Xk \ismapof_\Kk \Yk) \,\land\, (\Xk \cong \Xk')\,\land\, (\Yk \cong \Yk')\, \implies\, (\Xk' \ismapof_\Kk \Yk').$    
    \item[] \textit{Proof.} Consider $\Xk'=\xi_*\Xk$ and $\Yk'=\zeta_*\Yk$ with isomorphisms $\xi$, $\zeta$.
    \item[] Let $\varphi:\,\Ycal \to \Xcal$ such that: $\Kk(X,Y|T) = \deltabf_\varphi(X|Y) \otimes \Kk(Y|T)$. Then:
        \begin{align*} 
             \Kk(X',Y'|T) &= \lp \deltabf_\xi(X'|X) \otimes \deltabf_\zeta(Y'|Y) \rp \circ \Kk(X,Y|T)\\ 
              &= \lp \deltabf_\xi(X'|X) \otimes \deltabf_\zeta(Y'|Y) \rp \circ \lp \deltabf_\varphi(X|Y) \otimes \Kk(Y|T) \rp \\
              &= \lp\deltabf_{\xi \circ \varphi}(X'|Y) \otimes  \deltabf_\zeta(Y'|Y) \rp \circ \Kk(Y|T) \\ 
              &= \deltabf_{\xi \circ \varphi \circ \zeta^{-1}}(X'|Y') \otimes  \Kk(Y'|T).
        \end{align*}
    \item $(\Xk \ismapof_\Kk \Yk) \, \implies\, (\Xk \ismapof_\Kk (\Yk \otimes \Zk))$.
    \item[] \textit{Proof.} This is proven in \Cref{lem:ismapof-null-set} and  \ref{lem:join-upper-bound-2}. 
    \item $\lp \Xk \Indep_\Kk \Yk \given \Zk\rp \, \land\,  (\Xk \cong \Xk')\,\land\,(\Yk \cong \Yk') \,\land (\Zk \cong \Zk')\, \implies \, \lp \Xk' \Indep_\Kk \Yk' \given \Zk'\rp.$
    \item[] \textit{Proof.} If $\Xk'=\varphi_*\Xk$ and $\Yk'=\psi_*\Yk$ and $\Zk'=\xi_*\Zk$ and:
        \[ \Kk(X,Y,Z|T) = \Qk(X|Z) \otimes \Kk(Y,Z|T).\]
        Then we get:
        \[ \Kk(X',Y',Z'|T) = \Qk(\varphi(X)|Z=\xi^{-1}(Z')) \otimes \Kk(Y',Z'|T).\]
    \item $\deltabf_\ast \ismapof_\Kk \Xk$.
    \item[] \textit{Proof.} $ \Kk(\ast,X|T) =\deltabf_\ast \otimes \Kk(X|T).$
    \item $\deltabf_\ast \otimes \Xk \cong \Xk$.
    \item[] \textit{Proof.} Use isomorphism: $\Asterisk \times \Xcal \cong \Xcal$.
\end{enumerate}%
\end{Lem}

\subsection{Core Separoid Rules for Transitional Conditional Independence}

 \begin{Lem}[Extended Left Redundancy]
    \label{sep:tci:ext-l-red}
     We have for any $\Yk$ the implication:
    \[ \Uk \ismapof_\Kk \Zk \implies \Uk \Indep_\Kk \Yk \given \Zk.\]
    \begin{proof}
        The assumption implies the existence of a factorization:
        \[ \Kk(U,Z|T) = \deltabf_\varphi(U|Z) \otimes \Kk(Z|T). \]
        \Cref{lem:ismapof-null-set} then shows that this extends to:
        \[ \Kk(U,Y,Z|T) = \deltabf_\varphi(U|Z) \otimes \Kk(Y,Z|T), \]
        which shows the claim.
    \end{proof}
\end{Lem}

\begin{Lem}[Left Redundancy]
    \label{sep:tci:l-red}
    \[ \deltabf_\ast \Indep_\Kk \Yk \given \Zk \qquad \text{always holds.}  \]
\begin{proof}
    \[\Kk(\ast,Y,Z|T) = \deltabf_\ast \otimes \Kk(Y,Z|T).\]
\end{proof}
\end{Lem}

\begin{Lem}[Criteria for a trivial right argument]
    \label{lem:trivial-right-argument}
    Let $\Xk$, $\Zk$ be transitional random variables on $\lp \Wcal \times \Tcal, \Kk(W|T) \rp$ and consider the
    two statements:
    \begin{enumerate}
        \item[(i)] $\Xk \Indep_\Kk \deltabf_\ast \given \Zk$, i.e.\ there is a Markov kernel
            $\Qk(X|Z):\, \Zcal \dshto \Xcal$ with $\Kk(X,Z|T) = \Qk(X|Z) \otimes \Kk(Z|T)$;
        \item[(ii)] the Markov kernel $\Kk(X,Z|T):\, \Tcal \dshto \Xcal \times \Zcal$ admits a conditional Markov
            kernel $\Kk(X|Z,T):\, \Zcal \times \Tcal \dshto \Xcal$, i.e.\
            $\Kk(X,Z|T) = \Kk(X|Z,T) \otimes \Kk(Z|T)$.
    \end{enumerate}
    The only difference between the two is that the kernel in (i) may \emph{not} depend on the input. We have:
    \begin{enumerate}
        \item (i) $\implies$ (ii).
        \item (ii) holds under each of the six conditions of \Cref{thm-regular-conditional-Markov-kernel},
            applied to the Markov kernel $\Kk(X,Z|T)$ with $(\Xcal,\Zcal,\Tcal)$ in the role of
            $(\Xcal,\Ycal,\Zcal)$ there. In particular it holds if $(\Xcal,\Zcal,\Tcal)$ is a disintegration
            triple, but also, on \emph{arbitrary} spaces, if $\Kk(X,Z|T=t) \ll \mubf \otimes \nubf$ for
            $\sigma$-finite $\mubf$, $\nubf$ with a jointly measurable density, or if
            $\Xk \ismapof_\Kk \Zk \otimes \Tk$, or if $\Zk \ismapof_\Kk \Tk$.
        \item (ii) $\implies$ (i) if $\Tk \ismapof_\Kk \Zk$, i.e.\ if the input can be recovered from $\Zk$.
        \item (i) holds outright if $\Xk \ismapof_\Kk \Zk$.
    \end{enumerate}
    In particular (i) and (ii) are equivalent whenever $\Tk \ismapof_\Kk \Zk$, which is the case for
    $\Zk = \Zk' \otimes \Tk$.
\begin{proof}
    1.\ If $\Kk(X,Z|T) = \Qk(X|Z) \otimes \Kk(Z|T)$ then $\Kk(X|Z,T) := \Qk(X|Z)$, read as a Markov kernel
    $\Zcal \times \Tcal \dshto \Xcal$ that happens to be constant in $t$, is a conditional Markov kernel.

    2.\ This is \Cref{thm-regular-conditional-Markov-kernel} verbatim; note that its points 5.\ and 6.\ read as
    $\Xk \ismapof_\Kk \Zk \otimes \Tk$ and $\Zk \ismapof_\Kk \Tk$ here.

    3.\ Let $\varphi:\, \Zcal \to \Tcal$ be measurable with
    $\Kk(T,Z|T) = \deltabf_\varphi(T|Z) \otimes \Kk(Z|T)$, which is what $\Tk \ismapof_\Kk \Zk$ says. Since
    $\Tk = \deltabf(T|T)$ the left hand side is $\deltabf_t \otimes \Kk(Z|T=t)$ at $t$, so the identity states
    that the image of $\Kk(Z|T=t)$ under the measurable map $z \mapsto (\varphi(z),z)$ is
    $\deltabf_t \otimes \Kk(Z|T=t)$. Consequently, for every bounded measurable
    $g:\, \Tcal \times \Zcal \to \R$ and every $t \in \Tcal$:
    \[ \int g(\varphi(z),z)\, \Kk(dz|T=t) \;=\; \int g(t,z)\, \Kk(dz|T=t). \]
    Now put:
    \[ \Qk(X|Z=z) := \Kk(X|Z=z,T=\varphi(z)), \]
    a Markov kernel $\Zcal \dshto \Xcal$, being the composition of the measurable map
    $z \mapsto (z,\varphi(z))$ with $\Kk(X|Z,T)$. Applying the displayed identity to
    $g(s,z) := \I_C(z) \cdot \Kk(X \in A|Z=z,T=s)$ for $A \in \Bcal_\Xcal$ and $C \in \Bcal_\Zcal$ gives:
    \begin{align*}
        \lp \Qk(X|Z) \otimes \Kk(Z|T=t) \rp (A \times C)
            &= \int_C \Kk(X \in A|Z=z,T=t)\,\Kk(dz|T=t)\\
            &= \Kk(X \in A, Z \in C|T=t).
    \end{align*}
    The measurable rectangles form a $\cap$-stable generator containing $\Xcal \times \Zcal$ and both sides are
    probability measures, so Dynkin's lemma, see \cite{Bog07} 1.9.3, gives
    $\Qk(X|Z) \otimes \Kk(Z|T) = \Kk(X,Z|T)$.

    4.\ If $\Xk \ismapof_\Kk \Zk$, say $X = \psi(Z)$ with $\psi:\, \Zcal \to \Xcal$ measurable, then
    $\Kk(X,Z|T) = \deltabf_\psi(X|Z) \otimes \Kk(Z|T)$, so $\Qk(X|Z) := \deltabf_\psi(X|Z)$ witnesses (i).

    For the last sentence note that $\Tk \ismapof_\Kk \Tk$ by \Cref{lem:restricted-reflexivity}, since $\Tk$ is
    deterministic, and hence $\Tk \ismapof_\Kk \Zk' \otimes \Tk$ by \Cref{lem:separoid-compatibility-1}
    points 5., 2.\ and 4.
\end{proof}
\end{Lem}

\begin{Lem}[$\Tk$-Restricted Right Redundancy]
    \label{sep:tci:r-red}
    The statement:
    \[ \Xk \Indep_\Kk \deltabf_\ast \given \Zk\otimes\Tk \]
    is \emph{equivalent} to the existence of a conditional Markov kernel $\Kk(X|Z,T)$ of $\Kk(X,Z|T)$, and
    therefore holds in each of the following cases:
{\setlength{\itemsep}{1pt}\setlength{\parskip}{0pt}
    \begin{enumerate}
        \item $(\Xcal,\Zcal,\Tcal)$ is a disintegration triple, see \Cref{def:disintegration-triple}, which by
            \Cref{thm-regular-conditional-Markov-kernel} is the case if $\Xcal$ is standard and $\Zcal$
            countably generated, if $\Xcal$ is standard and $\Tcal$ is discrete, or if $\Zcal$ is discrete;
        \item $\Kk(X,Z|T=t) \ll \mubf \otimes \nubf$ for every $t \in \Tcal$, for $\sigma$-finite measures
            $\mubf$ on $\Xcal$ and $\nubf$ on $\Zcal$, with a density that is jointly measurable in $(x,z,t)$;
        \item $\Xk \ismapof_\Kk \Zk \otimes \Tk$;
        \item $\Zk \ismapof_\Kk \Tk$.
    \end{enumerate}}
    Only case 1.\ is a condition on the underlying measurable spaces; cases 2.--4.\ hold on arbitrary ones.
\begin{proof}
    \emph{The equivalence.} Let $\Kk(X|Z,T)$ be a conditional Markov kernel of $\Kk(X,Z|T)$. Then
    $\lp (z,s),t \rp \mapsto \Kk(X|Z=z,T=s)$ is a conditional Markov kernel of $\Kk(X,(Z,T)|T)$ conditioned on
    $(Z,T)$: indeed $\Kk(Z,T|T=t) = \Kk(Z|T=t) \otimes \deltabf_t$, and therefore
    \begin{align*}
        \Kk(X,Z,T|T=t) &= \Kk(X,Z|T=t) \otimes \deltabf_t
                        = \Kk(X|Z,T=t) \otimes \Kk(Z|T=t) \otimes \deltabf_t\\
                       &= \Kk(X|Z,T) \otimes \Kk(Z,T|T=t),
    \end{align*}
    the last kernel being evaluated at $(z,s)$ with $s=t$ under $\deltabf_t$. This is statement (ii) of
    \Cref{lem:trivial-right-argument} for $\Zk \otimes \Tk$ in the role of $\Zk$; since
    $\Tk \ismapof_\Kk \Zk\otimes\Tk$, point 3.\ there turns it into (i), which is the claim. Conversely, the
    claim provides a Markov kernel $\Qk(X|Z,T)$ with $\Kk(X,Z,T|T) = \Qk(X|Z,T) \otimes \Kk(Z,T|T)$; evaluating
    at $T=t$, using $\Kk(Z,T|T=t) = \Kk(Z|T=t) \otimes \deltabf_t$ once more and marginalizing the
    $\Tcal$-coordinate out gives $\Kk(X,Z|T=t) = \Qk(X|Z,T=t) \otimes \Kk(Z|T=t)$, so that $\Qk(X|Z,T)$ is a
    conditional Markov kernel of $\Kk(X,Z|T)$.

    \emph{The cases.} Case 1.\ is \Cref{def:disintegration-triple}, and the three space conditions named there
    are points 1.--3.\ of \Cref{thm-regular-conditional-Markov-kernel}, applied to $\Kk(X,Z|T)$ with
    $(\Xcal,\Zcal,\Tcal)$ in the role of $(\Xcal,\Ycal,\Zcal)$ there; cases 2., 3.\ and 4.\ are its points 4.,
    5.\ and 6.\ for the same kernel.
\end{proof}
\end{Lem}

\begin{Lem}[Left Decomposition]
    \label{sep:tci:l-dec}
    \[  \Xk\otimes\Uk \Indep_\Kk  \Yk \given \Zk \implies \Uk \Indep_\Kk \Yk \given \Zk.\]
\begin{proof}
    By assumption we have the factorization:
    \[ \Kk(X,U,Y,Z|T) = \Qk(X,U|Z) \otimes \Kk(Y,Z|T).  \]
    Marginalizing out $X$ gives:
    \[ \Kk(U,Y,Z|T) = \Qk(U|Z) \otimes \Kk(Y,Z|T).  \]
    This shows the claim.
\end{proof}
\end{Lem}

\begin{Lem}[Right Decomposition]
    \label{sep:tci:r-dec}
    \[  \Xk \Indep_\Kk  \Yk\otimes\Uk \given \Zk \implies \Xk \Indep_\Kk \Uk \given \Zk.\]
\begin{proof}
    By assumption we have the factorization:
    \[ \Kk(X,U,Y,Z|T) = \Qk(X|Z) \otimes \Kk(Y,U,Z|T).  \]
    Marginalizing out $Y$ gives:
    \[ \Kk(X,U,Z|T) = \Qk(X|Z) \otimes \Kk(U,Z|T).  \]
    This shows the claim.
\end{proof}
\end{Lem}

\begin{Lem}[$\Tk$-Inverted Right Decomposition]
    \label{sep:tci:inv-r-dec}
    \[ \Xk \Indep_\Kk  \Yk \given \Zk  \implies 
    \Xk \Indep_\Kk \Tk\otimes\Yk \given \Zk.\]
\begin{proof}
    By the assumption $\Xk \Indep_\Kk  \Yk \given \Zk$ we have a factorization:
    \[ \Kk(X,Y,Z|T) = \Qk(X|Z) \otimes \Kk(Y,Z|T).  \]
    Multiplying both sides with $\Tk=\deltabf(T|T)$ gives:
    \[ \Kk(X,T,Y,Z|T) = \Qk(X|Z) \otimes \Kk(T,Y,Z|T).  \]
    This shows the claim.
\end{proof}
\end{Lem}

\begin{Lem}[Left Weak Union]
    \label{sep:tci:l-uni}
    \label{sep:tci:res-l-uni}
    Assume:
    \[ \Xk \otimes \Uk \Indep_\Kk \Yk  \given \Zk. \]
    Then:
    \[ \Xk \Indep_\Kk  \Yk \given \Uk \otimes\Zk \]
    holds in each of the following cases:
    \begin{enumerate}
        \item $(\Xcal,\Ucal,\Zcal)$ is a disintegration triple, see \Cref{def:disintegration-triple};
        \item $\Xk \ismapof_\Kk \Uk \otimes \Zk$;
        \item $\Uk \ismapof_\Kk \Zk$;
        \item there is a probability measure $\Pk(T)$ on $\Tcal$ such that the mixture:
            \[ \Pk(X,U,Z) := \Kk(X,U,Z|T) \circ \Pk(T) \]
            admits a conditional Markov kernel $\Pk(X|U,Z)$ and has a $\Zcal$-marginal $\Pk(Z)$ with
            $\Kk(Z|T=t) \ll \Pk(Z)$ for all $t \in \Tcal$. The latter is automatic for $\Tcal = \Asterisk$, and
            holds for $\Pk(T) := \deltabf_{t_0}$ whenever the family $\lp \Kk(Z|T=t) \rp_{t \in \Tcal}$ is
            \emph{dominated} by one of its own members, i.e.\ whenever $\Kk(Z|T=t) \ll \Kk(Z|T=t_0)$ for all
            $t \in \Tcal$;
        \item $\Kk(X,U,Z|T)$ admits a conditional Markov kernel $\Kk(X|U,Z,T)$ and
            $\Tk \ismapof_\Kk \Uk \otimes \Zk$;
        \item $\Xk \Indep_\Kk \deltabf_\ast \given \Uk \otimes \Zk$.
    \end{enumerate}
    Only case 1.\ is a condition on the underlying measurable spaces; cases 2.--6.\ hold on arbitrary ones.
    Case 6.\ is moreover \emph{necessary}: it follows from the conclusion by Right Decomposition
    \ref{sep:tci:r-dec}, since $\Yk \cong \Yk \otimes \deltabf_\ast$. So, given the premise, 6.\ is
    \emph{equivalent} to the conclusion, and 1.--5.\ are criteria for it.
\begin{proof}
    By assumption there is a Markov kernel $\Qk(X,U|Z)$, a \emph{witness}, with:
    \[\Kk(X,U,Y,Z|T) = \Qk(X,U|Z) \otimes \Kk(Y,Z|T),\]
    and marginalizing out $Y$, resp.\ $X$:
    \begin{align}
        \Kk(X,U,Z|T) &= \Qk(X,U|Z) \otimes \Kk(Z|T), &
        \Kk(U,Y,Z|T) &= \Qk(U|Z) \otimes \Kk(Y,Z|T). \label{eqn:lwu-marg} \tag{$\ast$}
    \end{align}

    \emph{Step 1: a factorizing witness suffices.} Suppose some witness $\tilde \Qk(X,U|Z)$ admits a
    factorization $\tilde\Qk(X,U|Z) = \tilde\Qk(X|U,Z) \otimes \tilde\Qk(U|Z)$. Then, using
    \eqref{eqn:lwu-marg} for $\tilde \Qk$:
    \begin{align*}
        \Kk(X,U,Y,Z|T) &= \tilde\Qk(X,U|Z) \otimes \Kk(Y,Z|T)\\
                     &= \tilde\Qk(X|U,Z) \otimes \tilde\Qk(U|Z) \otimes \Kk(Y,Z|T) \\
                     &= \tilde\Qk(X|U,Z) \otimes  \Kk(U,Y,Z|T),
    \end{align*}
    which is the claim.

    \emph{Step 2: witnesses may be modified on null sets.} Let $\tilde\Qk(X,U|Z)$ be a Markov kernel such that
    for all $A \in \Bcal_\Xcal$, $D \in \Bcal_\Ucal$ and $t \in \Tcal$:
    \[ \tilde\Qk(X \in A, U \in D|Z=z) = \Qk(X \in A, U \in D|Z=z)
       \quad \text{ for } \Kk(Z|T=t)\text{-a.e.\ } z, \]
    where the exceptional null set may depend on $A$, $D$ and $t$. Then $\tilde\Qk$ is again a witness. Indeed,
    for such $A$, $D$, for $B \in \Bcal_\Ycal$, $C \in \Bcal_\Zcal$ and $t \in \Tcal$ the two integrals:
    \begin{align*}
        \int_C \tilde\Qk(A \times D|Z=z)\, \Kk(Y \in B, Z \in dz|T=t), \\
        \int_C \Qk(A \times D|Z=z)\, \Kk(Y \in B, Z \in dz|T=t)
    \end{align*}
    agree, since $\Kk(Y \in B, Z \in \cdot\,|T=t) \le \Kk(Z \in \cdot\,|T=t)$ and the integrands agree
    $\Kk(Z|T=t)$-almost surely. The measurable rectangles $A \times D \times B \times C$ form a $\cap$-stable
    generator of $\Bcal_\Xcal \otimes \Bcal_\Ucal \otimes \Bcal_\Ycal \otimes \Bcal_\Zcal$ containing the whole
    space, and both sides are probability measures, so Dynkin's lemma, see \cite{Bog07} 1.9.3, gives
    $\tilde\Qk(X,U|Z) \otimes \Kk(Y,Z|T) = \Qk(X,U|Z) \otimes \Kk(Y,Z|T) = \Kk(X,U,Y,Z|T)$.

    \emph{Case 1.} Since $(\Xcal,\Ucal,\Zcal)$ is a disintegration triple, the witness $\Qk(X,U|Z)$ itself
    factorizes as $\Qk(X|U,Z) \otimes \Qk(U|Z)$, see \Cref{def:disintegration-triple}. Apply Step 1.

    \emph{Case 2.} By \Cref{lem:trivial-right-argument} 4., applied with $\Uk \otimes \Zk$ in the role of $\Zk$,
    the hypothesis $\Xk \ismapof_\Kk \Uk \otimes \Zk$ gives case 6.

    \emph{Case 3.} Let $\psi:\, \Zcal \to \Ucal$ be measurable with
    $\Kk(U,Z|T) = \deltabf_\psi(U|Z) \otimes \Kk(Z|T)$, which is what $\Uk \ismapof_\Kk \Zk$ says. By the
    extension \Cref{lem:ismapof-null-set}, applied with $\Uk$, $\Zk$ and $\Xk \otimes \Yk$ in the roles of
    $\Xk$, $\Yk$, $\Zk$ there, this extends to:
    \[ \Kk(X,U,Y,Z|T) = \deltabf_\psi(U|Z) \otimes \Kk(X,Y,Z|T), \qquad
       \Kk(U,Y,Z|T) = \deltabf_\psi(U|Z) \otimes \Kk(Y,Z|T), \]
    the second identity by marginalizing out $X$. Marginalizing out $U$ in the witness identity gives
    $\Kk(X,Y,Z|T) = \Qk(X|Z) \otimes \Kk(Y,Z|T)$, and putting the three together:
    \begin{align*}
        \Kk(X,U,Y,Z|T) &= \deltabf_\psi(U|Z) \otimes \Qk(X|Z) \otimes \Kk(Y,Z|T)\\
                       &= \Qk(X|Z) \otimes \deltabf_\psi(U|Z) \otimes \Kk(Y,Z|T)
                        = \Qk(X|Z) \otimes \Kk(U,Y,Z|T),
    \end{align*}
    which is the claim, with $\Qk(X|Z)$ read as a kernel $\Ucal \times \Zcal \dshto \Xcal$ that is constant in
    $u$. Here Steps 1 and 2 are not needed.

    \emph{Case 4.} Put $\tilde \Qk(X,U|Z) := \Pk(X|U,Z) \otimes \Qk(U|Z)$, which factorizes by construction.
    Multiplying the first identity of \eqref{eqn:lwu-marg} by $\Pk(T)$ from the right and marginalizing $\Tcal$
    out gives $\Pk(X,U,Z) = \Qk(X,U|Z) \otimes \Pk(Z)$, by the computation in Step 1 of the proof of
    \Cref{thm:supp:tci-prior} applied with $\Xk \otimes \Uk$ in the role of $\Xk$ and $\deltabf_\ast$ in the role
    of $\Yk$; marginalizing $\Xcal$ out of it gives $\Pk(U,Z) = \Qk(U|Z) \otimes \Pk(Z)$. Hence:
    \begin{align*}
        \tilde\Qk(X,U|Z) \otimes \Pk(Z) &= \Pk(X|U,Z) \otimes \Qk(U|Z) \otimes \Pk(Z)
                                         = \Pk(X|U,Z) \otimes \Pk(U,Z) \\
                                        &= \Pk(X,U,Z) = \Qk(X,U|Z) \otimes \Pk(Z),
    \end{align*}
    so $\tilde \Qk$ and $\Qk$ are both conditional Markov kernels of $\Pk(X,U,Z)$ conditioned on $Z$. By the
    essential uniqueness \Cref{ess-unique} the set
    $N_A := \lC z \st \tilde\Qk(A|Z=z) \neq \Qk(A|Z=z) \rC$ is measurable and $\Pk(Z)$-null for every
    $A \in \Bcal_{\Xcal\times\Ucal}$, hence $\Kk(Z|T=t)$-null for every $t$ by the domination hypothesis. The
    argument of Step 2, run for one $A$ at a time, therefore shows that $\tilde\Qk$ is a witness; and it
    factorizes. Apply Step 1.

    \emph{Case 5.} By \Cref{lem:trivial-right-argument} 3., applied with $\Uk\otimes\Zk$ in the role of $\Zk$,
    the two hypotheses give case 6.

    \emph{Case 6.} By assumption there is a Markov kernel $\Pk(X|U,Z)$ with:
    \[ \Kk(X,U,Z|T) = \Pk(X|U,Z) \otimes \Kk(U,Z|T). \]
    Together with \eqref{eqn:lwu-marg} this gives:
    \[\Kk(X,U,Z|T) = \Pk(X|U,Z) \otimes \Qk(U|Z) \otimes \Kk(Z|T)
        = \Qk(X,U|Z) \otimes \Kk(Z|T),\]
    so by the essential uniqueness \Cref{ess-unique} we get, for every $A \in \Bcal_\Xcal$ and
    $D \in \Bcal_\Ucal$:
    \[ \Qk(X \in A,U \in D|Z)  =  \int_D \Pk(X \in A|U=u,Z) \, \Qk(U \in du|Z) \quad  \Kk(Z|T)\text{-a.s.},\]
    i.e.\ $\tilde \Qk(X,U|Z) := \Pk(X|U,Z) \otimes \Qk(U|Z)$ satisfies the hypothesis of Step 2 and is thus a
    witness. It factorizes, so Step 1 applies.

    \emph{Necessity of 6.} If $\Xk \Indep_\Kk \Yk \given \Uk\otimes\Zk$ then, since
    $\Yk \cong \Yk \otimes \deltabf_\ast$ by \Cref{lem:separoid-compatibility-1}, Right Decomposition
    \ref{sep:tci:r-dec} gives $\Xk \Indep_\Kk \deltabf_\ast \given \Uk\otimes\Zk$.
\end{proof}
\end{Lem}

\begin{Lem}[Right Weak Union]
    \label{sep:tci:r-uni}
    \[\Xk \Indep_\Kk  \Yk\otimes\Uk \given \Zk \implies \Xk \Indep_\Kk  \Yk \given \Uk \otimes\Zk.\]
\begin{proof}
    We have the factorization:
    \[\Kk(X,Y,U,Z|T) = \Qk(X|Z) \otimes \Kk(Y,U,Z|T),\]
    with some Markov kernel $\Qk(X|Z)$. If we view $\Qk(X|Z)$ as a function in $(u,z)$ via:
    \[ (u,z) \mapsto \Qk(X|Z=z), \]
    by just ignoring the argument $u$ then the claim follows from the same factorization above.
\end{proof}
\end{Lem}

\begin{Lem}[Left Contraction]
    \label{sep:tci:l-con}
  \[ (\Xk \Indep_\Kk  \Yk \given \Uk \otimes\Zk) \land (\Uk \Indep_\Kk  \Yk \given \Zk) \implies \Xk\otimes\Uk \Indep_\Kk\Yk \given \Zk.\]
\begin{proof}
    By assumption we have the two factorizations:
    \begin{align*}
        \Kk(X,Y,U,Z|T) & = \Qk(X|U,Z) \otimes \Kk(Y,U,Z|T),\\
        \Kk(Y,U,Z|T) & = \Pk(U|Z) \otimes \Kk(Y,Z|T),
    \end{align*}
    with some Markov kernels $\Qk(X|U,Z)$, $\Pk(U|Z)$. Putting these equations together using $\Qk(X|U,Z) \otimes \Pk(U|Z)$ we get:
    \[ \Kk(X,Y,U,Z|T) = \lp \Qk(X|U,Z) \otimes \Pk(U|Z)\rp \otimes \Kk(Y,Z|T).\]
    This shows the claim.
\end{proof}
\end{Lem}

\begin{Lem}[Right Contraction]
    \label{sep:tci:r-con}
    \[ (\Xk \Indep_\Kk \Yk \given \Uk \otimes\Zk) \land (\Xk \Indep_\Kk \Uk \given \Zk)  \implies \Xk \Indep_\Kk \Yk\otimes\Uk 
    \given \Zk.\]
\begin{proof}
 By assumption we have the two factorizations:
    \begin{align*}
        \Kk(X,Y,U,Z|T) & = \Qk(X|U,Z) \otimes \Kk(Y,U,Z|T),\\
        \Kk(X,U,Z|T) & = \Pk(X|Z) \otimes \Kk(U,Z|T),
    \end{align*}
    with some Markov kernels $\Qk(X|U,Z)$, $\Pk(X|Z)$.\\
    Marginalizing out $Y$ we get the equalities:
    \begin{align*}
        \Kk(X,U,Z|T) & = \Qk(X|U,Z) \otimes \Kk(U,Z|T),\\
        \Kk(X,U,Z|T) & = \Pk(X|Z) \otimes \Kk(U,Z|T).
    \end{align*}
    By the essential uniqueness (see \Cref{ess-unique}) of such factorization we get that for every $A \in \Bcal_\Xcal$:
\[ \Qk(X \in A|U,Z) = \Pk(X \in A|Z) \qquad \Kk(U,Z|T)\text{-a.s.}\]
The same equation then holds also $\Kk(Y,U,Z|T)$-a.s., since the null set $N_A$ of \Cref{lem:ess-unique} does not depend on $y$ and the $(U,Z)$-marginal of $\Kk(Y,U,Z|T=t)$ is $\Kk(U,Z|T=t)$.
Plugging that back into the first equation gives:
\[ \Kk(X,Y,U,Z|T) = \Pk(X|Z) \otimes \Kk(Y,U,Z|T).\]
    This shows the claim.
\end{proof}
\end{Lem}

\begin{Lem}[Right Cross Contraction]
    \label{sep:tci:rc-con}
    \[ (\Xk\Indep_\Kk\Yk\given\Uk\otimes\Zk) \land (\Uk\Indep_\Kk\Xk\given\Zk) \implies \Xk\Indep_\Kk\Yk\otimes\Uk\given\Zk.\]
\begin{proof}
     By assumption we have the two factorizations:
    \begin{align}
        \Kk(X,Y,U,Z|T) & = \Qk(X|U,Z) \otimes \Kk(Y,U,Z|T), \label{eq:rcc-a}\\
        \Kk(X,U,Z|T) & = \Pk(U|Z) \otimes \Kk(X,Z|T), \label{eq:rcc-b}
    \end{align}
     with some Markov kernels $\Qk(X|U,Z)$, $\Pk(U|Z)$.\\
     We then define the Markov kernel:
    \begin{align}
        \Rk(X,U|Z) & := \Qk(X|U,Z) \otimes \Pk(U|Z). \label{eq:rcc-c}
    \end{align}
    We will now show that its marginal: 
    \begin{align}
        \Rk(X|Z) & = \Qk(X|U,Z) \circ \Pk(U|Z). \label{eq:rcc-c2}
    \end{align} 
    will satisfy the claim.\\
    If we marginalize out $Y$ from equation \ref{eq:rcc-a} we get:
    \begin{align}
        \Kk(X,U,Z|T) & = \Qk(X|U,Z) \otimes \Kk(U,Z|T). \label{eq:rcc-d}
    \end{align}
    Equating equations \ref{eq:rcc-b} and \ref{eq:rcc-d} gives:
    \begin{align}
        \Pk(U|Z) \otimes \Kk(X,Z|T) &=      \Kk(X,U,Z|T)  
                                 = \Qk(X|U,Z) \otimes \Kk(U,Z|T). \label{eq:rcc-q}
    \end{align}
    Marginalizing out $X$ in equation \ref{eq:rcc-q} on both sides gives:
    \begin{align}
        \Kk(U,Z|T) &=  \Pk(U|Z) \otimes \Kk(Z|T).  \label{eq:rcc-r}
    \end{align}
    If we now plug equation \ref{eq:rcc-r} into \ref{eq:rcc-d} then we get:
    \begin{align}
        \Kk(X,U,Z|T) & = \Qk(X|U,Z) \otimes \Pk(U|Z) \otimes \Kk(Z|T) \\
            &\stackrel{\ref{eq:rcc-c}}{=} \Rk(X,U|Z) \otimes \Kk(Z|T). \label{eq:rcc-e}
    \end{align}
    If we marginalize out $U$ in equation \ref{eq:rcc-e} and use definition \ref{eq:rcc-c2}  we arrive at:
    \begin{align}
        \Kk(X,Z|T) & = \Rk(X|Z) \otimes \Kk(Z|T). \label{eq:rcc-f}
    \end{align}
    We now get:
    \begin{align}
        \Qk(X|U,Z) \otimes \Kk(U,Z|T) 
        & \stackrel{\ref{eq:rcc-d}}{=}  \Kk(X,U,Z|T) \\
        & \stackrel{\ref{eq:rcc-b}}{=} \Pk(U|Z) \otimes \Kk(X,Z|T) \\
        &\stackrel{\ref{eq:rcc-f}}{=} \Pk(U|Z) \otimes \Rk(X|Z) \otimes \Kk(Z|T) \\
        &=  \Rk(X|Z)  \otimes \Pk(U|Z) \otimes \Kk(Z|T) \\
        &\stackrel{\ref{eq:rcc-r}}{=} \Rk(X|Z)  \otimes  \Kk(U,Z|T).
    \end{align}
    By the essential uniqueness (see \Cref{ess-unique}) of such a factorization we get that for every $A \in \Bcal_\Xcal$:
    \begin{align}
        \Qk(X \in A|U,Z) &= \Rk(X \in A|Z) \qquad \Kk(U,Z|T)\text{-a.s.}  \label{eq:ess-un}
    \end{align}
The same equation then holds also $\Kk(Y,U,Z|T)$-a.s. (by ignoring the non-occurring argument $y$).
Plugging \ref{eq:ess-un} back into the equation \ref{eq:rcc-a} we get:
    \begin{align}
        \Kk(X,Y,U,Z|T) & = \Qk(X|U,Z) \otimes \Kk(Y,U,Z|T),\\
                     & = \Rk(X|Z)  \otimes \Kk(Y,U,Z|T).
    \end{align}
This shows the claim.
\end{proof}
\end{Lem}

\begin{Lem}[Flipped Left Cross Contraction]
    \label{sep:tci:flc-con}
    \[ (\Xk \Indep_\Kk \Yk \given \Uk\otimes\Zk) \land (\Yk \Indep_\Kk \Uk \given \Zk) \implies \Yk \Indep_\Kk \Xk \otimes \Uk 
    \given \Zk.\]
\begin{proof}
 By assumption we have the two factorizations:
    \begin{align*}
        \Kk(X,Y,U,Z|T) & = \Qk(X|U,Z) \otimes \Kk(Y,U,Z|T),\\
        \Kk(Y,U,Z|T) & = \Pk(Y|Z) \otimes \Kk(U,Z|T),
    \end{align*}
     with some Markov kernels $\Qk(X|U,Z)$, $\Pk(Y|Z)$.\\
     Marginalizing out $Y$ in the first equation we get the equality:
    \begin{align*}
        \Kk(X,U,Z|T) & = \Qk(X|U,Z) \otimes \Kk(U,Z|T).
    \end{align*}
Plugging all three equations into each other we get:
    \begin{align*}
        \Kk(X,Y,U,Z|T) & = \Qk(X|U,Z) \otimes \Kk(Y,U,Z|T)\\
                     & = \Qk(X|U,Z) \otimes  \Pk(Y|Z) \otimes \Kk(U,Z|T)\\
                     & = \Pk(Y|Z) \otimes  \Qk(X|U,Z)   \otimes \Kk(U,Z|T) \\
                     & =  \Pk(Y|Z)  \otimes   \Kk(X,U,Z|T).
    \end{align*}
    This shows the claim.
\end{proof}
\end{Lem}

\begin{Cor}[The $T$-$\ast$-separoid of transitional random variables]
    \label{cor:t-star-separoid}
    Consider, on a transition probability space $\lp \Wcal \times \Tcal, \Kk(W|T) \rp$ with $\Tcal$ standard,
    the class of \emph{all} transitional random variables $\Xk:\, \Wcal \times \Tcal \dshto \Xcal$ whose
    codomains $\Xcal$ are standard measurable spaces. This class is closed under $\otimes$, it contains
    $\Tk = \deltabf(T|W,T)$ and $\deltabf_\ast$, and all triples of codomains occurring in it are disintegration
    triples by \Cref{thm-regular-conditional-Markov-kernel} 1., since standard measurable spaces are countably
    generated.
    Together with the product of Markov kernels $\otimes$, the equivalence of isomorphisms of measurable spaces
    $\cong$, the relation $\ismapof_\Kk$, the one-point Markov kernel $\deltabf_\ast$ and transitional
    conditional independence $\Indep_\Kk$ it forms a $\Tk$-$\deltabf_\ast$-separoid (or in simpler symbols,
    a $T$-$\ast$-separoid), see \Cref{def:t-k-separoid}.
    Restricted to the \emph{deterministic} transitional random variables $\deltabf(X|W,T)$ it is in addition a
    bounded join-semi-lattice with join $\otimes$ and bottom element $\deltabf_\ast$ up to almost-sure
    anti-symmetry $\approx_\Kk$, see \Cref{cor:join-semi-lattice}.
\end{Cor}

\begin{Rem}
    \label{rem:tau-idempotency-free}
    \Cref{def:t-k-separoid} imposes two conditions on the distinguished element $\tau$, namely reflexivity of
    $\ll$ at $\tau$ and $\tau \lor \tau \approx \tau$. Here they come for free and impose \emph{no}
    restriction on the class: the input variable $\Tk = \deltabf(T|W,T)$ is by construction the Dirac kernel of the
    canonical projection $T:\, \Wcal \times \Tcal \to \Tcal$ and hence \emph{deterministic}, so
    \Cref{lem:restricted-reflexivity} gives $\Tk \ismapof_\Kk \Tk$ and \Cref{lem:idempotency} gives
    \[ \Tk \otimes \Tk \approx_\Kk \Tk \]
    on \emph{every} transition probability space, whatever the ambient class of transitional random variables and
    whatever the space $\Tcal$. The same applies to $\kappa = \deltabf_\ast$, the Dirac kernel of the constant
    map. Note that $\Tk \otimes \Tk \cong \Tk$ would be \emph{false}, since $\Tcal \times \Tcal$ and $\Tcal$
    need not be measurably isomorphic; this is exactly why \Cref{def:t-k-separoid} is formulated with the coarser
    equivalence $\approx$, whose invariance for $\Indep_\Kk$ is Full Equivalent Exchange
    \ref{sep:tci:full-eq-ex}, again for arbitrary transitional random variables.
    What determinism \emph{is} needed for is the join-semi-lattice statement, i.e.\ reflexivity and idempotency of
    the \emph{elements} of the class: for a genuinely stochastic $\Xk$ one does not even have
    $\Xk \ismapof_\Kk \Xk$, see \Cref{rem:ismapof-properties} item \ref{rem:no-reflexivity}.
    Finally note that the global Markov property of \Cref{sec:applications-causal_models} uses only the rules
    a)--k), whose individual hypotheses are listed in \Cref{tab:rules-hypotheses}, and therefore needs neither the
    full separoid structure nor standardness of the input spaces.
\end{Rem}

\subsection{Derived Separoid Rules for Transitional Conditional Independence}

Most of the following rules follow directly from the $T$-$\ast$-separoid rules proven in the last subsection.
Since we have to track which of the spaces form disintegration triples, we go through the proofs carefully.

\begin{Lem}[Extended $\Tk$-Restricted Right Redundancy]
    \label{sep:tci:ext-r-red}
    Let $(\Xcal,\Zcal,\Tcal)$ be a disintegration triple, see \Cref{def:disintegration-triple}. Then:
  \[ \Tk \ismapof_\Kk \Zk \implies \Xk \Indep_\Kk \deltabf_\ast \given \Zk.  \]
\begin{proof}
Extended $\Tk$-Restricted Right Redundancy \ref{sep:tci:ext-r-red} can  be proven using 
$\Tk$-Restricted Right Redundancy \ref{sep:tci:r-red}
($\implies\Xk \Indep_\Kk \deltabf_\ast \given\Zk\otimes\Tk$) together with Extended Left Redundancy \ref{sep:tci:ext-l-red} ($\implies \Tk \Indep_\Kk\Xk\given\Zk$)
and Right Cross Contraction \ref{sep:tci:rc-con} ($\implies \Xk \Indep_\Kk \Tk\otimes\deltabf_\ast \given \Zk$) and 
then Right Decomposition \ref{sep:tci:r-dec} 
    ($\implies \Xk \Indep_\Kk \deltabf_\ast \given \Zk$).
    It is also the instance $\Zk$ of \Cref{lem:trivial-right-argument} 2.\ and 3.\ read together, and the
    disintegration triple hypothesis may accordingly be replaced by any of the cases of \Cref{sep:tci:r-red}.
\end{proof}
\end{Lem}

\begin{Lem}[Restricted Symmetry]
    \label{sep:tci:res-sym}
    \[ (\Xk\Indep_\Kk\Yk\given\Zk) \land (\Yk\Indep_\Kk\deltabf_\ast\given\Zk) \implies \Yk\Indep_\Kk\Xk\given\Zk.\]
    \begin{proof}
        This follows from Flipped Left Cross Contraction \ref{sep:tci:flc-con} with $\Uk=\deltabf_\ast$.
    \end{proof}
\end{Lem}

\begin{Lem}[$\Tk$-Restricted Symmetry]
    \label{sep:tci:t-res-sym}
    Assume $\Xk\Indep_\Kk\Yk\given\Zk\otimes\Tk$. Then:
    \[ \Yk \Indep_\Kk \Xk\given \Zk\otimes\Tk \]
    holds in each of the four cases of \Cref{sep:tci:r-red}, read with $\Yk$ in the role of $\Xk$: if
    $(\Ycal,\Zcal,\Tcal)$ is a disintegration triple, if $\Kk(Y,Z|T)$ has a jointly measurable density w.r.t.\ a
    product of $\sigma$-finite measures, if $\Yk \ismapof_\Kk \Zk \otimes \Tk$, or if $\Zk \ismapof_\Kk \Tk$.
    Only the first is a condition on the underlying measurable spaces. What these criteria establish,
    $\Yk \Indep_\Kk \deltabf_\ast \given \Zk\otimes\Tk$, is moreover \emph{necessary}, so that, given the
    premise, it is equivalent to the conclusion.
    \begin{proof}
        By \Cref{sep:tci:r-red}, applied with $\Yk$ in the role of $\Xk$, each of the four cases gives:
    \[ \Yk \Indep_\Kk  \deltabf_\ast \given \Zk\otimes\Tk. \]
    Together with $\Xk \Indep_\Kk \Yk \given \Zk\otimes\Tk$ and Restricted Symmetry \ref{sep:tci:res-sym} we get:
    \[ \Yk \Indep_\Kk  \Xk \given \Zk\otimes\Tk.\]
    For the necessity, write $\Xk \cong \Xk \otimes \deltabf_\ast$, see \Cref{lem:separoid-compatibility-1}, and
    apply Right Decomposition \ref{sep:tci:r-dec} to the conclusion, which drops the first right factor.
    \end{proof}
\end{Lem}

\begin{Lem}[Symmetry]
    \label{sep:tci:sym}
    Let $\Tcal=\Asterisk$ and assume $\Xk \Indep_\Kk \Yk \given \Zk$. Then:
    \[ \Yk \Indep_\Kk\Xk \given \Zk \]
    holds in each of the four cases of \Cref{sep:tci:r-red}, read with $\Yk$ in the role of $\Xk$ and
    $\Tcal = \Asterisk$, and $\Yk \Indep_\Kk \deltabf_\ast \given \Zk$ is again necessary.
\begin{proof}
    This is $\Tk$-Restricted Symmetry \ref{sep:tci:t-res-sym} for $\Tcal = \Asterisk$, where
    $\Zk \otimes \Tk \cong \Zk$.
\end{proof}
\end{Lem}

\begin{Lem}[Inverted Left Decomposition]
    \label{sep:tci:inv-l-dec}
    \[ \lp \Xk \Indep_\Kk  \Yk \given \Zk \rp \land \lp \Uk \ismapof_\Kk \Xk \otimes \Zk  \rp  \implies 
    \Xk\otimes\Uk \Indep_\Kk \Yk \given \Zk.\]
\begin{proof}
Inverted Left Decomposition \ref{sep:tci:inv-l-dec} can be proven using Extended 
Left Redundancy \ref{sep:tci:ext-l-red} ($\implies \Uk \Indep_\Kk\Yk\given\Xk\otimes\Zk$)
together with the assumption ($\Xk\Indep_\Kk\Yk\given\Zk$) and Left Contraction \ref{sep:tci:l-con} 
    ($\implies \Xk\otimes\Uk \Indep_\Kk \Yk \given \Zk$).
\end{proof}
\end{Lem}

\begin{Lem}[$\Tk$-Extended Inverted Right Decomposition]
    \label{sep:tci:ext-inv-r-dec}
    \[ \lp \Xk \Indep_\Kk  \Yk \given \Zk \rp \land \lp \Uk \ismapof_\Kk \Tk \otimes \Yk \otimes \Zk  \rp  \implies 
    \Xk \Indep_\Kk \Tk\otimes\Yk \otimes \Uk \given \Zk.\]
\begin{proof}
    $\Tk$-Extended Inverted Right Decomposition \ref{sep:tci:ext-inv-r-dec} can be proven using
    $\Tk$-Inverted Right Decomposition \ref{sep:tci:inv-r-dec} ($\implies \Xk\Indep_\Kk\Tk \otimes\Yk \given \Zk$) 
    in combination with Extended Left Redundancy \ref{sep:tci:ext-l-red} ($\implies \Uk\Indep_\Kk\Xk\given\Tk\otimes\Yk\otimes\Zk$)
    and Flipped Left Cross Contraction \ref{sep:tci:flc-con} ($\implies \Xk \Indep_\Kk \Tk\otimes\Yk \otimes \Uk \given \Zk$).
\end{proof}
\end{Lem}

\begin{Lem}[Equivalent Exchange]
    \label{sep:tci:eq-ex}
    \[ \lp \Xk \Indep_\Kk  \Yk \given \Zk \rp \land \lp \Zk \approx_\Kk \Zk'  
        \rp  \implies 
    \Xk \Indep_\Kk\Yk \given \Zk'.\]
\begin{proof}
    We get:
    \begin{align*}
\Zk' \ismapof_\Kk \Zk  & \qquad\xRightarrow{\text{\Cref{lem:join-upper-bound-2}}} 
            && \Zk' \ismapof_\Kk \Tk \otimes \Yk\otimes \Zk,\\
\Xk \Indep_\Kk  \Yk \given \Zk
             & \qquad\xRightarrow{\Tk\text{-Ext. Inv. Right Decomposition \ref{sep:tci:ext-inv-r-dec}}}
            &&\Xk \Indep_\Kk \Tk\otimes\Yk\otimes\Zk' \given \Zk\\
            &\qquad\xRightarrow{\text{Right Decomposition \ref{sep:tci:r-dec}}}
            &&\Xk \Indep_\Kk \Yk\otimes\Zk' \given \Zk\\       
            &\qquad\xRightarrow{\text{Right Weak Union \ref{sep:tci:r-uni}}}
            &&\Xk \Indep_\Kk \Yk \given \Zk'\otimes\Zk, \tag{a} \label{eq:int-a}\\
            \Zk \ismapof_\Kk \Zk' & \qquad\xRightarrow{\text{Extended Left Redundancy \ref{sep:tci:ext-l-red}}} &&\Zk \Indep_\Kk \Yk \given \Zk',\tag{b} \label{eq:int-b}\\\
        \eqref{eq:int-a}  \land  \eqref{eq:int-b}& 
        \qquad \xRightarrow{\text{Left Contraction \ref{sep:tci:l-con}}}&& \Xk \otimes \Zk \Indep_\Kk \Yk \given \Zk'\\
         &\qquad \xRightarrow{\text{Left Decomposition \ref{sep:tci:l-dec}}}&& \Xk \Indep_\Kk \Yk \given \Zk'.
    \end{align*}
\end{proof}
\end{Lem}

\begin{Lem}[Full Equivalent Exchange]
    \label{sep:tci:full-eq-ex}
   If  $\Xk' \approx_\Kk \Xk$ and $\Yk' \approx_\Kk \Yk$ and $\Zk' \approx_\Kk \Zk$ 
then we have the equivalence:
\[ \Xk \Indep_\Kk  \Yk \given \Zk \qquad \iff \qquad 
        \Xk' \Indep_\Kk \Yk' \given \Zk'.  \]
\begin{proof}
    \begin{align*}
    \Xk' \ismapof_\Kk \Xk  & \qquad\xRightarrow{\text{\Cref{lem:join-upper-bound-2}}} 
            && \Xk' \ismapof_\Kk \Xk \otimes \Zk\\
             & \qquad\xRightarrow{\text{Inverted Left Decomposition \ref{sep:tci:inv-l-dec}}} 
            && \Xk \otimes \Xk' \Indep_\Kk  \Yk \given \Zk \\
           & \qquad\xRightarrow{\text{Left Decomposition \ref{sep:tci:l-dec}}} 
            && \Xk' \Indep_\Kk  \Yk \given \Zk, \\
  \Yk' \ismapof_\Kk \Yk  & \qquad\xRightarrow{\text{\Cref{lem:join-upper-bound-2}}} 
            && \Yk' \ismapof_\Kk \Tk \otimes \Yk \otimes \Zk\\
             & \qquad\xRightarrow{\Tk\text{-Ext. Inv. Right Decomposition \ref{sep:tci:ext-inv-r-dec}}} 
            && \Xk'  \Indep_\Kk \Tk \otimes \Yk \otimes \Yk' \given \Zk \\
           & \qquad\xRightarrow{\text{Right Decomposition \ref{sep:tci:r-dec}}} 
            && \Xk' \Indep_\Kk  \Yk' \given \Zk, \\
   \Zk' \approx_\Kk \Zk  & \qquad\xRightarrow{\text{Equivalent Exchange \ref{sep:tci:eq-ex}}} 
            && \Xk' \Indep_\Kk  \Yk' \given \Zk'.
    \end{align*}
    The other direction works similarly.
\end{proof}
\end{Lem}

\section{Proofs - Applications to Statistical Theory}
\label{sec:supp:appl-statistics}

For the reader's convenience we restate the results of \Cref{sec:applications-statistics} before proving them; the hypotheses are identical to the ones given there.

Next we will give a proof that the classical Fisher-Neyman factorization criterion (see \cite{Fis22, Ney35, HS49, Bur61})
 is equivalent to sufficiency reformulated as transitional conditional independence.

\begin{Thm}[Fisher-Neyman]
    \label{thm:supp:fisher-neyman}
    Let $\Xcal$, $\Scal$, $\Thetacal$ be measurable spaces with $\Xcal$ standard. 
    Let $\mubf$ be a $\sigma$-finite measure on $\Xcal$ and $S:\, \Xcal \to \Scal$ a measurable map. 
    Let $\Pk(X|\Theta):\, \Thetacal \dshto \Xcal$ be a statistical model that is absolutely continuous w.r.t.\ $\mubf$: $\Pk(X|\Theta) \ll \mubf$. Then the following two statements are equivalent:
    \begin{enumerate}
        \item $\Pk(X|\Theta)$ has a Radon-Nikodym derivative\footnote{It is not necessary to assume joint measurability for the equivalence to hold; the proof \emph{produces} versions of $p_\theta$, $g_\theta$ and $f$ for which the maps $(x,\theta) \mapsto p_\theta(x)$, $(s,\theta) \mapsto g_\theta(s)$ and $x \mapsto f(x)$ are jointly measurable, so that one may always assume this w.l.o.g.; see the last paragraph of the proof of \Cref{thm:supp:fisher-neyman}. This joint measurability is what the likelihood principle, \Cref{thm:proper-likelihood-principle}, uses.} $p_\theta$ w.r.t.\ $\mubf$ of the form:
            \begin{align*}  p_\theta(x) = h(x) \cdot g_\theta(S(x)), \end{align*} 
            with measurable maps $h:\, \Xcal \to \R_{\ge 0}$ and $g_\theta:\, \Scal \to \R_{\ge0}$ for $\theta \in \Thetacal$.
        \item $S$ is a sufficient statistic for $\Pk(X|\Theta)$, i.e.\ we have the transitional conditional independence: \begin{align*}  X \Indep_{\Pk(X|\Theta)} \Theta \given S. \end{align*} 
    \end{enumerate}
\begin{proof}
    First note that $\mubf \neq 0$, since $\Pk(X|\Theta=\theta) \ll \mubf$ is a probability measure. We claim that
    there is a probability measure $\Qk(X)$ on $\Xcal$ such that $\mubf$ has a density $m$ w.r.t.\ $\Qk(X)$ with
    values in $(0,\infty)$ \emph{everywhere}, i.e.\ in particular $\mubf \ll \Qk(X) \ll \mubf$. Indeed, write
    $\Xcal = \biguplus_{n \in \N} \Xcal_n$ with $\mubf(\Xcal_n) < \infty$, put
    $I := \lC n \,|\, \mubf(\Xcal_n) > 0 \rC$, which is non-empty, and set
    \[ \Qk(X) := c \cdot \sum_{n \in I} 2^{-n} \cdot \frac{\mubf(\,\cdot \cap \Xcal_n)}{\mubf(\Xcal_n)},
       \qquad c := \lp \sum_{n \in I} 2^{-n} \rp^{-1}. \]
    Then $m := \frac{2^{n}\,\mubf(\Xcal_n)}{c}$ on $\Xcal_n$ for $n \in I$, and $m := 1$ on the remaining part
    $\biguplus_{n \notin I} \Xcal_n$, which is both $\mubf$- and $\Qk(X)$-null, is such a density. So with the relation:
    \begin{align*}  f(x) := h(x) \cdot m(x), \end{align*} 
    we can equivalently replace the first statement with the existence of a Radon-Nikodym derivative $p_\theta$ for $\Pk(X|\Theta)$ w.r.t.\ $\Qk(X)$ of the form:
            \begin{align}  
                p_\theta(x) = f(x) \cdot g_\theta(S(x)). \label{eq:supp:fisher-neyman}
            \end{align} 
   Further note, that the joint distribution $\Qk(X,S)$ has a regular conditional probability distribution $\Qk(X|S)$ by \Cref{reg-cond-trans-kernel-existence-V} and by the assumptions that $\Xcal$ is standard. 

``1.$\implies$2.'': We assume that we have a density like in \Cref{eq:supp:fisher-neyman}.
    For $s \in \Scal$ we then put:
    \begin{align*}  F(s) := \int f(x) \, \Qk(X \in dx|S=s) \in [0,\infty]. \end{align*} 
  Then the so defined $F:\Scal \to [0,\infty]$ is measurable.
  With this we get for $B \in \Bcal_\Scal$:
    \begin{align*} 
     &\Pk(S \in B|\Theta=\theta)\\
      & = \Pk(X \in S^{-1}(B)|\Theta=\theta)  \\ 
      &= \int \I_B(S(x)) \cdot g_\theta(S(x)) \cdot f(x) \,\Qk(X \in dx)\\
      &= \int \int\I_B(s)  \cdot g_\theta(s) \cdot f(x) \, \Qk(X \in dx, S \in ds)\\
      &= \int \I_B(s)  \cdot \lp \int f(x) \, \Qk(X \in dx|S=s) \rp \cdot  g_\theta(s) \, \Qk(S \in ds)\\
      &= \int \I_B(s)  \cdot F(s)  \cdot g_\theta(s) \, \Qk(S \in ds).
    \end{align*}
This means for each $\theta$ separately we have a Radon-Nikodym derivative:
\begin{align*}  \frac{\Pk(S \in ds|\Theta=\theta)}{\Qk(S \in ds)\phantom{|\Theta=\theta} }(s) =  F(s)  \cdot g_\theta(s).\end{align*} 
In particular, taking $B=\Scal$ shows $\int F(s) \cdot g_\theta(s) \, \Qk(S \in ds) = 1$ for every $\theta \in \Thetacal$, and thus:
\begin{align*}
    F(s) \cdot g_\theta(s) < \infty \qquad \text{ for } \Qk(S)\text{-almost-all } s \in \Scal.
\end{align*}
Note, however, that $F$ itself may well attain the values $0$ and $\infty$ on sets of strictly positive $\Qk(S)$-measure.
We thus put:
\begin{align*}
    N_0 &:= \lC s \in \Scal \st F(s) \in \{0,\infty\} \rC \, \in \Bcal_\Scal,
\end{align*}
which is measurable since $F$ is, and claim that $N_0$ is a $\Pk(S|\Theta)$-null set.
Indeed, fix $\theta \in \Thetacal$. For $s$ with $F(s)=0$ we clearly have $F(s) \cdot g_\theta(s) = 0$, and 
the above finiteness implies $g_\theta(s)=0$ for $\Qk(S)$-almost-all $s$ with $F(s)=\infty$, so that also there 
$F(s) \cdot g_\theta(s) = 0$ (with the usual convention $0 \cdot \infty =0$). So $F \cdot g_\theta$ vanishes 
$\Qk(S)$-almost-everywhere on $N_0$ and we get for every $\theta \in \Thetacal$:
\begin{align*}
    \Pk(S \in N_0|\Theta=\theta) = \int \I_{N_0}(s) \cdot F(s) \cdot g_\theta(s) \, \Qk(S \in ds) = 0.
\end{align*}
We now define the Markov kernel $\Kk(X|S):\, \Scal \dshto \Xcal$ via:
\begin{align*}
    \Kk(X \in A|S=s) &:= 
    \begin{cases}
        \int_A \frac{f(x)}{F(s)} \, \Qk(X \in dx|S=s), & \text{ if } s \notin N_0, \\
        \Qk(X \in A|S=s), & \text{ if } s \in N_0,
    \end{cases}
\end{align*}
i.e.\ on the exceptional (and, as just seen, $\Pk(S|\Theta)$-null) set $N_0$ we can choose $\Kk(X|S=s)$ arbitrarily.
Note that $\Kk(X|S)$ is a well-defined Markov kernel: it is measurable in $s$, and for $s \notin N_0$ we have 
$F(s) \in (0,\infty)$ and thus $\Kk(X \in \Xcal|S=s) = \frac{F(s)}{F(s)} = 1$.
Now consider the joint distribution:
    \begin{align*} 
     &\Pk(X \in A, S \in B|\Theta=\theta)\\
      & = \Pk(X \in A \cap S^{-1}(B)|\Theta=\theta)  \\ 
      &= \int \I_A(x) \cdot\I_B(S(x)) \cdot g_\theta(S(x)) \cdot f(x) \,\Qk(X \in dx)\\
      &= \int \int \I_A(x) \cdot\I_B(s)  \cdot g_\theta(s) \cdot f(x) \, \Qk(X \in dx, S \in ds)\\
      &= \int \lp \int \I_A(x) \cdot f(x) \, \Qk(X \in dx|S=s) \rp \cdot \I_B(s)  \cdot g_\theta(s) \, \Qk(S \in ds)\\
      &= \int \Kk(X \in A|S=s) \cdot \I_B(s) \cdot F(s)  \cdot g_\theta(s) \, \Qk(S \in ds)\\
      &= \int \Kk(X \in A|S=s) \cdot \I_B(s) \cdot \Pk(S \in ds|\Theta=\theta)\\
      &= \lp \Kk(X|S) \otimes \Pk(S|\Theta)\rp(A \times B,\theta).
    \end{align*}
Here the fifth equality holds because for $s \notin N_0$ we have, by the very definition of $\Kk(X|S=s)$:
\begin{align*}
    \int \I_A(x) \cdot f(x) \, \Qk(X \in dx|S=s) = \Kk(X \in A|S=s) \cdot F(s),
\end{align*}
while for $\Qk(S)$-almost-all $s \in N_0$ both integrands vanish: if $F(s)=0$ then 
$\int \I_A(x) \cdot f(x) \, \Qk(X \in dx|S=s) \le F(s)=0$, and for $\Qk(S)$-almost-all $s$ with $F(s)=\infty$ 
we have $g_\theta(s)=0$, as shown above.
So we get the factorization:
\begin{align*}  \Pk(X,S|\Theta) = \Kk(X|S) \otimes \Pk(S|\Theta).  \end{align*} 
This shows the transitional conditional independence: 
\begin{align*}  X \Indep_{\Pk(X|\Theta)} \Theta \given S, \end{align*} 
and thus the claim.

``2.$\implies$1.'': Assume $X \Indep_{\Pk(X|\Theta)} \Theta \given S$ and let $p_\theta$ be Radon-Nikodym derivative of $\Pk(X|\Theta=\theta)$ w.r.t.\ $\Qk(X)$ that is jointly measurable as a map $(x,\theta) \mapsto p_\theta(x)$, which exists by \Cref{thm:doob-derivation} with the assumption that $\Xcal$ is standard and thus countably generated.

    Then note, since $p_\theta$ is a density for $\Pk(X|\Theta=\theta)$ w.r.t.\ $\Qk(X)$, we also have that $p_\theta$ is a density for $\Pk(X,S|\Theta=\theta)$ w.r.t.\ $\Qk(X,S)$.

We then define the measurable maps $g_\theta:\, \Scal \to \R_{\ge0}$ for $\theta \in \Thetacal$ via:
\begin{align*}  g_\theta(s) := \int p_\theta(x) \, \Qk(X \in dx|S=s). \end{align*} 
It is then clear that $g_\theta$ is a density of the marginal $\Pk(S|\Theta=\theta)$ w.r.t.\ $\Qk(S)$. Indeed:
\begin{align*}
    \int_B g_\theta(s) \, \Qk(S \in ds) &= \int_B \int  p_\theta(x) \, \Qk(X \in dx|S=s) \, \Qk(S \in ds) \\
                                        &= \int \I_B(s) \cdot  p_\theta(x) \, \Qk(X \in dx, S \in ds) \\
                                        &= \int \I_B(S(x)) \cdot  p_\theta(x) \, \Qk(X \in dx) \\
                                        &= \int \I_{S^{-1}(B)}(x) \cdot  p_\theta(x) \, \Qk(X \in dx) \\
                                        &= \Pk(X \in S^{-1}(B)|\Theta = \theta) \\
                                        &= \Pk(S \in B|\Theta=\theta).
\end{align*}
Note that, since $(x,\theta) \mapsto p_\theta(x)$ is jointly measurable, so is $(s,\theta)\mapsto g_\theta(s)$.
Furthermore, we can then define the following map:
\begin{align*}
    \tilde p_\theta(x|s) &:=
    \begin{cases}
         \frac{p_\theta(x)}{g_\theta(s)},  & \text{ if } 0 < g_\theta(s) < \infty, \\
        1, & \text{ if } g_\theta(s) \in \{0,\infty\},
\end{cases}
\end{align*}
which is jointly measurable in the arguments $(x,s,\theta)$. The case $g_\theta(s)=\infty$ has to be excluded
explicitly, since the construction so far only produces $g_\theta$ with values in $[0,\infty]$. It occurs on a
$\Qk(S)$-null set only: by the computation above $g_\theta$ is a density of $\Pk(S|\Theta=\theta)$ with respect to
$\Qk(S)$, so $\int g_\theta \, d\Qk(S) = 1$ and hence $g_\theta < \infty$ $\Qk(S)$-almost surely. Redefining
$g_\theta := 0$ on $\lC g_\theta = \infty \rC$ therefore changes $q_\theta$ only on a $\Qk(X)$-null set and yields
the required $g_\theta:\, \Scal \to \R_{\geq 0}$; with that convention all three branches give a probability
measure and none of the identities below is affected.
With this we then define the following Markov kernel $\Pk(X|S,\Theta)$ via:
\begin{align*}
    \Pk(X \in A|S=s,\Theta=\theta) &:= 
    \int_A \tilde p_\theta(x|s)  \, \Qk(X \in dx|S=s)
\end{align*}
Note that:
\begin{align*}
    &\lp\Pk(X|S,\Theta) \otimes \Pk(S|\Theta)\rp (A \times B,\theta) \\
   &= \int_B \lp \int_A \tilde p_\theta(x|s) \, \Qk(X \in dx|S=s) \rp \cdot g_\theta(s) \, \Qk(S \in ds) \\
   &= \int_B \int_A p_\theta(x) \, \Qk(X \in dx|S=s) \, \Qk(S \in ds) \\
   &= \int \I_B(s) \cdot \I_A(x) \cdot p_\theta(x) \, \Qk(X \in dx,S\in ds) \\
   &= \int \I_B(S(x)) \cdot \I_A(x) \cdot p_\theta(x) \, \Qk(X \in dx) \\
   &= \Pk(X \in A, S \in B|\Theta=\theta).
\end{align*}
So $\Pk(X|S,\Theta)$ is a conditional Markov kernel of $\Pk(X,S|\Theta)$:
\begin{align*}  \Pk(X,S|\Theta) = \Pk(X|S,\Theta) \otimes \Pk(S|\Theta). \end{align*} 
On the other hand, $X \Indep_{\Pk(X|\Theta)} \Theta \given S$ implies that there exists a Markov kernel $\tilde\Kk(X|S)$ such that:
            \begin{align*}  \Pk(X,S|\Theta) = \tilde\Kk(X|S) \otimes \Pk(S|\Theta).  \end{align*} 
By the essential uniqueness, see \Cref{ess-unique}, we know that the set:
\begin{align*}
    N&:=\lC (s,\theta) \in \Scal \times \Thetacal \st \Pk(X|S=s,\Theta=\theta) \neq \tilde\Kk(X|S=s)   \rC
\end{align*}
is a measurable $\Pk(S|\Theta)$-null set of $\Scal \times \Thetacal$. So for $(s,\theta) \in N^\cmpl$ we get:
\begin{align*}
    \tilde\Kk(X \in A|S=s) &= \Pk(X \in A|S=s,\Theta=\theta) \\
                           &= \int_A \tilde p_\theta(x|s)  \, \Qk(X \in dx|S=s).
\end{align*}
Now consider the set $\tilde \Scal := \pr_\Scal(N^\cmpl) \ins \Scal$, which we endow with the subspace-$\sigma$-algebra.
Note that $\tilde\Scal$ need \emph{not} be a measurable subset of $\Scal$ (a projection of a measurable set is in
general not measurable, see \Cref{rem:no-measurable-selection}), but this is immaterial here. Indeed,
\Cref{thm:doob-derivation} allows an arbitrary measurable space in its parameter slot, and
$\lp \tilde\Scal, \Bcal_{\Scal|\tilde\Scal} \rp$ is one; and \Cref{thm-kuratowski-extension} is stated for an
arbitrary subset with its subspace-$\sigma$-algebra, which applies to
$\Xcal \times \tilde\Scal \ins \Xcal \times \Scal$ because the generating rectangles
$A \times (B \cap \tilde\Scal)$ lie in the trace, so that
$\Bcal_\Xcal \otimes \Bcal_{\Scal|\tilde\Scal} \ins \lp \Bcal_\Xcal \otimes \Bcal_\Scal \rp_{|\Xcal \times \tilde\Scal}$.
Then for $s \in \tilde \Scal$, by definition, there exists a $\theta \in \Thetacal$ with 
$(s,\theta) \in N^\cmpl$ and thus:
\begin{align*}
    \tilde\Kk(X \in A|S=s) &= \int_A \tilde p_\theta(x|s)  \, \Qk(X \in dx|S=s).
\end{align*}
In particular, we have for $s \in \tilde \Scal$: $\tilde\Kk(X|S=s) \ll \Qk(X|S=s)$.
By \Cref{thm:doob-derivation} with $\Xcal$ standard there exists a measurable map:
\begin{align*}  
    \tilde k:\, \Xcal \times \tilde \Scal &\to \R_{\ge0},  & (x,s) & \mapsto \tilde k(x|s),
\end{align*} 
such that for all $A \in \Bcal_\Xcal$ and $s \in \tilde \Scal$ we have:
\begin{align*}  \tilde\Kk(X \in A|S=s) = \int_A \tilde k(x|s) \, \Qk(X \in dx|S=s).  \end{align*} 
By Kuratowski's extension theorem, \Cref{thm-kuratowski-extension}, $\tilde k$ can be extended to a measurable map:
\begin{align*}  \tilde k:\, \Xcal \times \Scal \to \R_{\ge0}.  \end{align*} 
We can then normalize $\tilde k$ as:
\begin{align*}  
    k:\, \Xcal \times \Scal &\to \R_{\ge0},  \\ 
    k(x|s) & := 
    \begin{cases}
        \frac{\tilde k(x|s)}{\int \tilde k(\tilde x|s)\, \Qk(X \in d\tilde x|S=s)}, & \text{ if } \int \tilde k(\tilde x|s)\, \Qk(X \in d\tilde x|S=s) \neq 0, \\
        1, & \text{ if } \int \tilde k(\tilde x|s)\, \Qk(X \in d\tilde x|S=s) =0.
    \end{cases}
\end{align*}
This then defines a Markov kernel $\Kk(X|S)$ via:
\begin{align*}
    \Kk(X \in A |S=s) &:= \int_A k( x|s)\, \Qk(X \in dx|S=s).
\end{align*}
Note that for $(s,\theta) \in N^\cmpl$ we still have for every $A \in \Bcal_\Xcal$:
\begin{align*}
    \Kk(X \in A|S=s) 
    &= \tilde \Kk(X \in A|S=s) 
    = \Pk(X \in A|S=s,\Theta=\theta),
\end{align*}
and thus:
        \begin{align*}  \Pk(X,S|\Theta) &= \Kk(X|S) \otimes \Pk(S|\Theta).  \end{align*} 
This implies:
\begin{align*}
    \Pk(X\in A|\Theta=\theta) &= \Pk(X \in A,S \in \Scal|\Theta=\theta) \\
                              &= \int \int \I_A(x) \, \Kk(X \in dx|S=s) \, \Pk(S \in ds|\Theta=\theta) \\
                              &= \int \int \I_A(x) \cdot k(x|s) \, \Qk(X \in dx|S=s) \cdot g_\theta(s) \, \Qk(S \in ds) \\
                              &= \int \I_A(x) \cdot k(x|s) \cdot g_\theta(s) \, \Qk(X \in dx,S\in ds) \\
                              &= \int \I_A(x) \cdot \underbrace{k(x|S(x))}_{=: f(x)} \cdot g_\theta(S(x)) \, \Qk(X \in dx) \\
                              &= \int \I_A(x) \cdot f(x) \cdot g_\theta(S(x)) \, \Qk(X \in dx),
\end{align*}
where we defined the measurable map $f$ as:
\begin{align*}
    f:\, \Xcal &\to \R_{\ge 0}, &  f(x) &:= k(x|S(x)).  
\end{align*}
This shows that $\Pk(X|\Theta)$ has a Radon-Nikodym derivative w.r.t.\ $\Qk(X)$ of the form:
\begin{align*}
    q_\theta(x) :=    f(x) \cdot g_\theta(S(x)),
\end{align*}
which shows the claim.

Finally, note that the maps $(x,\theta) \mapsto p_\theta(x)$, $(s,\theta) \mapsto g_\theta(s)$ and $f$ constructed 
in this direction are jointly measurable, and, since $h=f/m$ with the strictly positive measurable density $m$, 
the same then holds w.r.t.\ the original reference measure $\mubf$. Together with 
``1.$\implies$2.'' this shows that in the first statement we can always w.l.o.g.\ assume this joint measurability, 
as claimed in the footnote there.
\end{proof}
\end{Thm}

\begin{Thm}[Conditional Basu]
    \label{thm:supp:basu}
    Let $\Xk$, $\Yk$, $\Uk$, $\Zk$ be transitional random variables on a transition probability space
    $\lp \Wcal \times \Tcal, \Kk(W|T) \rp$ such that:
{\setlength{\itemsep}{1pt}\setlength{\parskip}{0pt}
    \begin{enumerate}
        \item $\Xk \Indep_\Kk \Yk \given \Uk \otimes \Zk$;
        \item $\Xk \Indep_\Kk \Yk \given \Zk$;
        \item $\Uk$ is boundedly complete given $\Zk$ w.r.t.\ $\Kk$, see \Cref{def:boundedly-complete}.
    \end{enumerate}}
    Then we have:
    \[ \Xk \Indep_\Kk \Yk \otimes \Uk \given \Zk. \]
\begin{proof}
    \emph{Step 1: the two witnesses.}
    By 1.\ and \Cref{def:transitional_conditional_independence} there is a Markov kernel
    $\Qk_1(X|U,Z):\, \Ucal \times \Zcal \dshto \Xcal$ with:
    \begin{align}
        \Kk(X,Y,U,Z|T) &= \Qk_1(X|U,Z) \otimes \Kk(Y,U,Z|T), \label{eq:cbasu-1}
    \end{align}
    and by 2.\ a Markov kernel $\Qk_2(X|Z):\, \Zcal \dshto \Xcal$ with:
    \begin{align}
        \Kk(X,Y,Z|T) &= \Qk_2(X|Z) \otimes \Kk(Y,Z|T). \label{eq:cbasu-2}
    \end{align}
    Marginalizing out $Y$ gives:
    \begin{align}
        \Kk(X,U,Z|T) &= \Qk_1(X|U,Z) \otimes \Kk(U,Z|T), &
        \Kk(X,Z|T) &= \Qk_2(X|Z) \otimes \Kk(Z|T). \label{eq:cbasu-3}
    \end{align}

    \emph{Step 2: comparing the two on $\Zk$.}
    Fix $A \in \Bcal_\Xcal$ and define:
    \[ g_A:\, \Ucal \times \Zcal \to \R, \qquad
       g_A(u,z) := \Qk_1(X \in A|U=u,Z=z) - \Qk_2(X \in A|Z=z), \]
    which is measurable with values in $[-1,1]$, hence bounded. Marginalizing $U$ out of the first identity of
    \eqref{eq:cbasu-3} and comparing with the second gives, for every $C \in \Bcal_\Zcal$ and $t \in \Tcal$:
    \begin{align*}
        \int_{\Ucal \times C} \Qk_1(X \in A|U,Z)\, \Kk(du,dz|T=t)
            &= \Kk(X \in A, Z \in C|T=t)\\
            &= \int_C \Qk_2(X \in A|Z)\, \Kk(dz|T=t)\\
            &= \int_{\Ucal \times C} \Qk_2(X \in A|Z)\, \Kk(du,dz|T=t),
    \end{align*}
    the last step because the integrand does not depend on $u$ and $\Kk(Z|T=t)$ is the $\Zk$-marginal of
    $\Kk(U,Z|T=t)$. Subtracting, $\int_{\Ucal \times C} g_A \, \Kk(du,dz|T=t) = 0$ for all $C$ and $t$, so
    hypothesis 3.\ yields:
    \begin{align}
        \Qk_1(X \in A|U=u,Z=z) = \Qk_2(X \in A|Z=z) \qquad \Kk(U,Z|T=t)\text{-a.s., for every } t \in \Tcal.
        \label{eq:cbasu-ae}
    \end{align}
    The exceptional null set may depend on $A$ and on $t$; this does no harm, since we integrate
    \eqref{eq:cbasu-ae} for one fixed $A$ at a time.

    \emph{Step 3: the factorization.}
    Let $A \in \Bcal_\Xcal$, $B \in \Bcal_\Ycal$, $D \in \Bcal_\Ucal$, $C \in \Bcal_\Zcal$ and $t \in \Tcal$.
    Since $\Kk(U,Z|T=t)$ is the $(\Uk,\Zk)$-marginal of $\Kk(Y,U,Z|T=t)$, the identity \eqref{eq:cbasu-ae} also
    holds $\Kk(Y,U,Z|T=t)$-almost surely, so with \eqref{eq:cbasu-1}:
    \begin{align*}
        &\Kk(X \in A, Y \in B, U \in D, Z \in C|T=t)\\
        &\qquad = \int_{B \times D \times C} \Qk_1(X \in A|U,Z)\, \Kk(dy,du,dz|T=t)\\
        &\qquad = \int_{B \times D \times C} \Qk_2(X \in A|Z)\, \Kk(dy,du,dz|T=t)\\
        &\qquad = \lp \Qk_2(X|Z) \otimes \Kk(Y,U,Z|T=t) \rp (A \times B \times D \times C).
    \end{align*}
    The measurable rectangles $A \times B \times D \times C$ form a $\cap$-stable generator of
    $\Bcal_\Xcal \otimes \Bcal_\Ycal \otimes \Bcal_\Ucal \otimes \Bcal_\Zcal$ containing the whole space and both
    sides are probability measures, so Dynkin's lemma, see \cite{Bog07} 1.9.3, gives:
    \[ \Kk(X,Y,U,Z|T) = \Qk_2(X|Z) \otimes \Kk(Y,U,Z|T), \]
    which by \Cref{def:transitional_conditional_independence} is exactly
    $\Xk \Indep_\Kk \Yk \otimes \Uk \given \Zk$.
    Note that no property of the measurable spaces was used: the only conditional Markov kernels in the argument,
    $\Qk_1$ and $\Qk_2$, are supplied by the hypotheses and not constructed.
\end{proof}
\end{Thm}

\begin{Cor}[Basu]
    \label{cor:supp:basu-classical}
    In the situation of \Cref{cor:basu-classical} we have $R \Indep_{\Pk(X|\Theta)} \Theta, S$.
\begin{proof}
    Apply \Cref{thm:supp:basu} with $R$, $\Theta$, $S$, $\deltabf_\ast$ in the roles of $\Xk$, $\Yk$, $\Uk$,
    $\Zk$. Its first hypothesis is obtained from sufficiency as in \Cref{cor:basu-classical}, its second is
    ancillarity, and its third is bounded completeness, which is \Cref{def:boundedly-complete} with
    $\Zk = \deltabf_\ast$.
\end{proof}
\end{Cor}

\begin{Thm}[Blackwell's order is a transitional conditional independence]
    \label{thm:supp:blackwell}
    Let $\Xcal_1$, $\Xcal_2$, $\Thetacal$ be measurable spaces and let
    $\Pk_i(X_i|\Theta):\, \Thetacal \dshto \Xcal_i$, $i=1,2$, be experiments, see \Cref{def:experiment}. Then the
    following are equivalent:
    \begin{enumerate}
        \item $\Ek_1 \succeq \Ek_2$;
        \item there is a Markov kernel $\Kk(X_1,X_2|\Theta):\, \Thetacal \dshto \Xcal_1 \times \Xcal_2$ with
            $\Kk(X_i|\Theta) = \Pk_i(X_i|\Theta)$ for $i=1,2$ such that
            $X_2 \Indep_{\Kk(X_1,X_2|\Theta)} \Theta \given X_1$.
    \end{enumerate}
\begin{proof}
    Throughout, $X_1$, $X_2$ and $\Theta$ denote the coordinate projections of
    $\lp \Xcal_1 \times \Xcal_2 \rp \times \Thetacal$, so that they are \emph{deterministic} transitional random
    variables and $\Theta$ is the input.

    \emph{1.\ $\implies$ 2.} Let $\Qk(X_2|X_1)$ be a garbling as in \Cref{def:experiment} and put:
    \[ \Kk(X_1,X_2|\Theta) := \Qk(X_2|X_1) \otimes \Pk_1(X_1|\Theta). \]
    This is a Markov kernel $\Thetacal \dshto \Xcal_1 \times \Xcal_2$. Its $\Xcal_1$-marginal is
    $\Pk_1(X_1|\Theta)$, since $\Qk(X_2|X_1)$ is a probability kernel, and its $\Xcal_2$-marginal is
    $\Qk(X_2|X_1) \circ \Pk_1(X_1|\Theta) = \Pk_2(X_2|\Theta)$ by hypothesis, so $\Kk$ is a coupling of the two
    experiments. Since $\Theta$ is the input, its value is $\deltabf_\theta$ under $\Kk(X_2,\Theta,X_1|\Theta=\theta)$, so
    that $\Kk(X_2,\Theta,X_1|\Theta=\theta) = \deltabf_\theta(\Theta) \otimes \Kk(X_1,X_2|\Theta=\theta)$ and
    $\Kk(\Theta,X_1|\Theta=\theta) = \deltabf_\theta(\Theta) \otimes \Kk(X_1|\Theta=\theta)$ for every
    $\theta \in \Thetacal$. Therefore:
    \begin{align*}
        \Kk(X_2,\Theta,X_1|\Theta=\theta)
        &= \deltabf_\theta(\Theta) \otimes \Kk(X_1,X_2|\Theta=\theta)\\
        &= \deltabf_\theta(\Theta) \otimes \Qk(X_2|X_1) \otimes \Pk_1(X_1|\Theta=\theta)\\
        &= \Qk(X_2|X_1) \otimes \deltabf_\theta(\Theta) \otimes \Kk(X_1|\Theta=\theta)\\
        &= \Qk(X_2|X_1) \otimes \Kk(\Theta,X_1|\Theta=\theta),
    \end{align*}
    where the constant factor $\deltabf_\theta(\Theta)$ may be moved past $\Qk(X_2|X_1)$, which carries no
    $\Thetacal$-argument, by \Cref{rem:markov-kernels-products} 2.\ and 3.\ in the generality of
    \Cref{rem:product-wirings}. By
    \Cref{def:transitional_conditional_independence} this is $X_2 \Indep_\Kk \Theta \given X_1$.

    \emph{2.\ $\implies$ 1.} Let $\Kk$ be a coupling with $X_2 \Indep_\Kk \Theta \given X_1$. By
    \Cref{def:transitional_conditional_independence} there is a Markov kernel $\Qk(X_2|X_1)$ with
    $\Kk(X_2,\Theta,X_1|\Theta) = \Qk(X_2|X_1) \otimes \Kk(\Theta,X_1|\Theta)$. Evaluating at $\Theta=\theta$ and
    marginalizing the $\Thetacal$-component out of both sides gives:
    \[ \Kk(X_1,X_2|\Theta=\theta) = \Qk(X_2|X_1) \otimes \Kk(X_1|\Theta=\theta), \]
    and marginalizing $X_1$ out as well gives
    $\Kk(X_2|\Theta=\theta) = \Qk(X_2|X_1) \circ \Kk(X_1|\Theta=\theta)$. Since $\Kk$ is a coupling this reads
    $\Pk_2(X_2|\Theta=\theta) = \Qk(X_2|X_1) \circ \Pk_1(X_1|\Theta=\theta)$, for every $\theta \in \Thetacal$, so
    $\Qk(X_2|X_1)$ is a garbling and $\Ek_1 \succeq \Ek_2$.
    Note that no property of the spaces was used: as in \Cref{thm:supp:basu} the only kernel in the argument is
    supplied by the hypothesis and not constructed.
\end{proof}
\end{Thm}

Note that the Fisher-Neyman factorization theorem for sufficiency (see \cite{Fis22, Ney35, HS49}) requires the existence of a Radon-Nikodym derivative w.r.t.\ a reference measure.
Our definition of conditional independence generalizes the factorization theorem to Markov kernels (per definition) without the necessity of densities and/or reference measures.\\

A direct application of the Fisher-Neyman factorization theorem, \Cref{thm:supp:fisher-neyman}, is the \emph{likelihood principle}, \Cref{thm:proper-likelihood-principle} from \Cref{sec:likelihood}, which we can formalise and prove here in the following:

\begin{Thm}[The likelihood principle]
    \label{thm:supp:proper-likelihood-principle}
    Let $\Xcal$, $\Thetacal$ be measurable spaces with $\Xcal$ standard and let $\mubf$ be a $\sigma$-finite measure on $\Xcal$.
    Consider a statistical model, written as the Markov kernel: $\Pk(X|\Theta): \, \Thetacal \dshto \Xcal$.
    For each $\theta \in \Thetacal$ assume that the Radon-Nikodym derivative $p_\theta$ exists ($\Pk(X|\Theta)\ll \mubf$):
    \begin{align*}  p_\theta(x) := \frac{\Pk(X \in dx|\Theta=\theta)}{\mubf(dx)}(x).\end{align*} 
    Then consider the \emph{likelihood function}\footnote{We endow $\Rcal$ with the smallest $\sigma$-algebra $\Bcal_\Rcal$ such that the evaluation map $\ev_\theta:\, \Rcal \to \R_{\ge0}$, $r \mapsto r(\theta)$, is measurable for every $\theta \in \Thetacal$.}:
    \begin{align*}  L_\mubf:\, \Xcal \to \Rcal:=\R_{\ge 0}^\Thetacal, \qquad x \mapsto \lp \theta \mapsto p_\theta(x) \rp. \end{align*} 
    Then $L_\mubf$ is measurable, satisfies the Fisher-Neyman criterion, \Cref{eq:fisher-neyman}, w.r.t.\ $\mubf$ and is thus a sufficient statistic for $\Pk(X|\Theta)$, i.e.\ 
    we have the transitional conditional independence:
     \begin{align*}  X \Indep_{\Pk(X|\Theta) } \Theta \given L_\mubf. \end{align*} 
   Furthermore, let $S$ be any other measurable map of $X$, i.e.\ $S \ismapof X$. Then we have:
   \begin{enumerate}
       \item \emph{Sufficiency}: If $L_\mubf \ismapof S$ then also: $\displaystyle X \Indep_{\Pk(X|\Theta) } \Theta \given S$.
       \item \emph{Quasi-minimality}: If $S$ satisfies: $\displaystyle X \Indep_{\Pk(X|\Theta) } \Theta \given S$, then there exists a measure $\nubf$ such that $\nubf$ has a density w.r.t.\ $\mubf$, $\Pk(X|\Theta)$ has a density w.r.t.\ $\nubf$ and the corresponding likelihood function $L_\nubf$ satisfies: $L_\nubf \ismapof S$.
     \end{enumerate}
\begin{proof}
    $\displaystyle X \Indep_{\Pk(X|\Theta) } \Theta \given L_\mubf$ follows from Fisher-Neyman, \Cref{thm:supp:fisher-neyman}, via: 
    \begin{align*}
        p_\theta(x) = \ev_\theta(L_\mubf(x)).
    \end{align*}

    1.) Now assume $L_\mubf \ismapof S \ismapof X$. So there exists a measurable map $G$ such that: $L_\mubf(x)=G(S(x))$. With this we get: 
    \begin{align*}  p_\theta(x) = \ev_\theta(L_\mubf(x)) = \ev_\theta(G(S(x))), \end{align*} 
    which again satisfies Fisher-Neyman, \Cref{thm:supp:fisher-neyman}, w.r.t.\ $\mubf$, and thus: $\displaystyle X \Indep_{\Pk(X|\Theta) } \Theta \given S$.
    
    2.) For the reverse, apply Fisher-Neyman, \Cref{thm:supp:fisher-neyman} w.r.t.\ $\mubf$:
    \begin{align*}  
        p_\theta(x) = h(x) \cdot g_\theta(S(x)).  
    \end{align*} 
    So w.r.t.\ the measure $\nubf$ given by:
     \begin{align*}
         \nubf(A) &:= \int_A h(x) \, \mubf(dx), 
     \end{align*}
     we get the Radon-Nikodym derivative: 
    \begin{align*}  \tilde p_\theta(x) = g_\theta(S(x)).  \end{align*} 
    So we  can define the measurable map:
    \begin{align*} 
        G:\, \Scal &\to \Rcal=\R_{\ge0}^\Thetacal, & s &\mapsto (\theta \mapsto g_\theta(s)).
    \end{align*}
    This shows for every $\theta \in \Thetacal$ and $x \in \Xcal$ the equation:
    \begin{align*}  \ev_\theta(G(S(x))) = g_\theta(S(x)) = \tilde p_\theta(x) = \ev_\theta(L_\nubf(x)),\end{align*} 
    which implies:
    \begin{align*}  L_\nubf = G(S),  \end{align*} 
    and thus the claim: $L_\nubf \ismapof S$.
\end{proof}
\end{Thm}

\begin{Cor}[The likelihood ratio principle]
        \label{thm:supp:likelihood-ratio-principle}
    Let $\Xcal$, $\Thetacal$, $\Scal$ be measurable spaces with $\Xcal$ standard.
    Consider a statistical model, written as the Markov kernel: $\Pk(X|\Theta): \, \Thetacal \dshto \Xcal$.
    Let $S:\, \Xcal \to \Scal$ be a measurable map.
    Let $\Pk(\Theta)$ be any ``prior'' probability distribution on $\Thetacal$ such that the following absolute continuity holds:
    \begin{align*}
        \Pk(X|\Theta) \ll \Pk(X|\Theta) \circ \Pk(\Theta) =: \Qk(X).
    \end{align*}
Then the following statements are equivalent:
    \begin{enumerate}
        \item $\Pk(X|\Theta)$ has a Radon-Nikodym derivative $(x,\theta) \mapsto p_\theta(x)$ w.r.t.\ $\Qk(X)$ such that the corresponding likelihood function $L_\Qk$ is a measurable map in $S$: $L_\Qk \ismapof S$.
        \item There exist measurable maps $g_\theta:\, \Scal \to \R_{\ge 0}$ for $\theta \in \Thetacal$ such that: 
            \begin{align*}
                p_\theta(x) &:= g_\theta(S(x)),
            \end{align*}
            defines a Radon-Nikodym derivative of  $\Pk(X|\Theta)$ w.r.t.\ $\Qk(X)$.
        \item The sufficiency condition holds:
            $\displaystyle X \Indep_{\Pk(X|\Theta) } \Theta \given S$.
    \end{enumerate}
    In particular, $L_\Qk$ is an (a.s.) minimal sufficient statistic for $\Pk(X|\Theta)$.
\begin{proof}
    By \Cref{thm:supp:proper-likelihood-principle} and \Cref{thm:supp:fisher-neyman} the equivalence is clear up to a multiplicative factor $h(x)$ in:
            \begin{align*}
                p_\theta(x) &= h(x) \cdot g_\theta(S(x)),
            \end{align*}
            where, by the footnote in \Cref{thm:supp:fisher-neyman}, we can w.l.o.g.\ assume that the maps 
            $(x,\theta) \mapsto p_\theta(x)$ and $(s,\theta) \mapsto g_\theta(s)$ are jointly measurable, 
            which justifies the following applications of Fubini's theorem.
            First note that the constant function $1$ is a density of $\Qk(X)$ w.r.t.\ $\Qk(X)$.
            Then integrating the above equation over $\Pk(\Theta)$ gives:
  \begin{align*}
      1 &\stackrel{\Qk(X)\text{-a.s.}}{=} \int p_{\tilde \theta}(x)\,\Pk(\Theta \in d\tilde \theta) \\
        &= \int h(x) \cdot g_{\tilde\theta}(S(x)) \,\Pk(\Theta \in d\tilde\theta) \\
        &= h(x) \cdot \int g_{\tilde\theta}(S(x)) \,\Pk(\Theta \in d\tilde\theta),
  \end{align*}
  which implies:
  \begin{align*}
      h(x) &\stackrel{\Qk(X)\text{-a.s.}}{=} \frac{1}{\int g_{\tilde \theta}(S(x)) \,\Pk(\Theta \in d\tilde \theta)},
  \end{align*}
  and, in particular, that $\int g_{\tilde \theta}(s) \,\Pk(\Theta \in d\tilde \theta) \in (0,\infty)$ for 
  $\Qk(S)$-almost-all $s \in \Scal$, so that the following quotient is well-defined up to a $\Qk(S)$-null set, 
  on which we can put it to $1$:
      \begin{align*}
          p_\theta(x) &\stackrel{\Qk(X)\text{-a.s.}}{=} \frac{g_\theta(S(x))}{\int g_{\tilde \theta}(S(x)) \,\Pk(\Theta \in d\tilde \theta)} =: \tilde g_\theta(S(x)).
      \end{align*}
  This implies all claims.
\end{proof}
\end{Cor}

\begin{Rem}[Savage-Dickey e-posterior]
    The density from \Cref{thm:supp:likelihood-ratio-principle}, in arguments given by $p_\theta(x)$, can be written as the ratio of densities $\frac{q_\theta(x)}{q(x)}$ when using a different reference measure $\mubf \gg \Qk(X)$.
    This ratio was called the \emph{Savage-Dickey e-posterior/density ratio} in \cite{Gru23},
    also see \cite{GHK19, WR21, PLW23}. 
    So, \Cref{thm:supp:likelihood-ratio-principle} can be reformulated as saying that the Savage-Dickey e-posterior is an (a.s.) minimal sufficient statistic for $\Pk(X|\Theta)$, when considered as a function in $\theta$ and evaluated at the variable $X$. So all inference about $\Theta$ based on $X$ can be done by processing the Savage-Dickey e-posterior only.
\end{Rem}

The next theorem contains the proof for the propensity score of \Cref{sec:propensity} and the 
dual likelihood principle of \Cref{sec:likelihood}.

\begin{Thm}
\label{thm:supp:propensity}
Let $\Pk(Y|X)$ be a Markov kernel.
For $x \in \Xcal$ we put:
\begin{align*} E(x):= \Pk(Y|X=x) \; \in \; \Pcal(\Ycal).\end{align*} 
Note that the map $E: \,\Xcal \to \Pcal(\Ycal)=:\Ecal$ is measurable and $E \ismapof X$.
Now let $S: \, \Xcal \to \Scal$ be another measurable map ($S \ismapof X$). 
Then we have the equivalence:
\begin{align*}  Y \Indep_{\Pk(Y|X) } X \given S \qquad \iff \qquad E \ismapof S.\end{align*} 
In particular, as $E \ismapof E$, we have:
\begin{align*}  Y \Indep_{\Pk(Y|X) } X \given E.\end{align*} 
\begin{proof}
First assume the conditional independence:
    \begin{align*}  Y \Indep_{\Pk(Y|X) } X \given S.\end{align*} 
This then implies that there exists a measurable function: 
\begin{align*}  \Pk(Y|S,\cancel{X}):\,  \Scal \to \Pcal(\Ycal),\end{align*} 
such that:
\begin{align*}  \Pk(Y,S|X) = \Pk(Y|S,\cancel{X}) \otimes \Pk(S|X). \end{align*} 
Noting that $\Pk(S|X) = \deltabf(S|X)$ and marginalizing $S$ out we get:
\begin{align*}  E(x) = \Pk(Y|X=x) = \Pk(Y|S = S(x),\cancel{X}),  \end{align*} 
for every $x \in \Xcal$. This shows: $E \ismapof S$.

For the reverse, now assume: $E \ismapof S$. Then there exists a measurable map $g:\,\Scal \to \Ecal=\Pcal(\Ycal)$ such that $E=g(S)$. We then define a Markov kernel via:
\begin{align*}  \Kk(Y|S):\, \Scal \dshto \Ycal, \quad \Kk(Y \in B|S=s) := g(s)(B). \end{align*} 
With this we then get for $B \in \Bcal_\Ycal$, $D \in \Bcal_\Scal$ and $x \in \Xcal$:
\begin{align*}
   & \lp \Kk(Y|S) \otimes \Pk(S|X) \rp (B \times D,x) \\
   &= \int_D \Kk(Y \in B|S=s) \, \Pk(S \in ds|X=x) \\
   &= \int_D \Kk(Y \in B|S=s) \, \deltabf(S \in ds|X=x) \\
   &= \Kk(Y \in B|S=S(x)) \cdot \I_D(S(x))  \\
   &= g(S(x))(B) \cdot \I_D(S(x))  \\
   &= E(x)(B) \cdot \I_D(S(x))  \\
   &= \Pk(Y \in B|X=x) \cdot \deltabf(S \in D|X=x)  \\
   &= \Pk(Y \in B, S \in D|X=x).
\end{align*}
Since this holds for all measurable sets $B \in \Bcal_\Ycal$, $D \in \Bcal_\Scal$ and $x \in \Xcal$
this shows the equality:
\begin{align*}  \Pk(Y,S|X) = \Kk(Y|S) \otimes \Pk(S|X),  \end{align*} 
which implies:
    \begin{align*}  Y \Indep_{\Pk(Y|X) } X \given S.\end{align*} 
This shows the claim.
\end{proof}
\end{Thm}

We now turn to Bayesian statistics.

\begin{Thm} [Bayesian statistics]
    \label{thm:supp:bayesian-statistics}
Let $\Pk(X|\Theta)$ be a Markov kernel between standard measurable spaces
    and $\Pk(\Theta|\Pi)$ be another Markov kernel.
Then put:
\begin{align*} \Pk(X,\Theta|\Pi) := \Pk(X|\Theta) \otimes \Pk(\Theta|\Pi). \end{align*} 
Then by \Cref{reg-cond-trans-kernel-existence-V} we have a conditional Markov kernel:
\begin{align*} \Pk(\Theta|X,\Pi),\end{align*} 
which is unique up to a $\Pk(X|\Pi)$-null set. 
We now define the transitional random variable: 
\begin{align*}  Z(x,\pi):= \Pk(\Theta|X=x,\Pi=\pi), \end{align*}  
which gives us a joint (transition) probability distribution: $\Pk(X,\Theta,Z|\Pi)$.

With the above notations we then have the conditional independence:
\begin{align*}  \Theta \Indep_{\Pk(X,\Theta|\Pi) } X \given Z.\end{align*} 
Now let $S$ be another deterministic measurable function in $(X,\Pi)$. 
Then we have the equivalence:
\begin{align*}  \Theta \Indep_{\Pk(X,\Theta|\Pi) } X \given S \qquad \iff \qquad Z \ismapof_{\Pk(X,\Theta|\Pi)} S  \end{align*} 
\begin{proof}
    For the first statement consider the Markov kernel given by:
\begin{align*}  \Kk(\Theta \in D|Z=z) := z(D). \end{align*} 
Then we get:
\begin{align*}
& \left(\Kk(\Theta| Z) \otimes \Pk(Z,X|\Pi)\right)(B\times C \times A, \pi) \\
&= \int_{C \times A} \Kk(\Theta \in B|Z=z) \, \Pk(Z\in dz,X \in dx|\Pi=\pi) \\ 
&= \int_A \int_C \Kk(\Theta \in B|Z=z) \,\deltabf(Z \in dz|X=x,\Pi=\pi) \,\Pk(X \in dx|\Pi=\pi) \\
&= \int_A \I_C(Z(x,\pi)) \cdot \Pk(\Theta \in B|X=x,\Pi=\pi) \,\Pk(X \in dx|\Pi=\pi) \\
&= \int_A \int_C \Pk(\Theta \in B|X=x,\Pi=\pi)\,\deltabf(Z \in dz|X=x,\Pi=\pi) \,\Pk(X \in dx|\Pi=\pi) \\
&= \Pk(\Theta \in B, Z \in C, X \in A |\Pi=\pi).
\end{align*}
This shows the first claim. 

For the second claim now first assume:
\begin{align*} \Theta \Indep_{\Pk(X,\Theta|\Pi) } X \given S.\end{align*}  
This conditional independence gives us the factorization:
\begin{align*}
    \Pk(\Theta, S, X|\Pi) &= \Pk(\Theta|S,\cancel{X,\Pi}) \otimes \Pk(S,X|\Pi)\\
    &= \Pk(\Theta|S,\cancel{X,\Pi}) \otimes \deltabf(S|X,\Pi) \otimes \Pk(X|\Pi).
\end{align*}
On the other hand we have:
\begin{align*}
    \Pk(\Theta, S, X|\Pi) &= \Pk(\Theta|X,\Pi) \otimes \deltabf(S|X,\Pi) \otimes \Pk(X|\Pi).
\end{align*}
Marginalizing $S$ out in those equations gives:
\begin{align*}
    \Pk(\Theta|S = S(X,\Pi),\cancel{X,\Pi}) \otimes \Pk(X|\Pi) &= \Pk(\Theta|X,\Pi) \otimes \Pk(X|\Pi).
\end{align*}
Because conditional Markov kernels are essentially unique by \Cref{ess-unique} we get:
\begin{align*}  \Pk(\Theta|S=S(X,\Pi),\cancel{X,\Pi})) = \Pk(\Theta|X,\Pi) = Z(X,\Pi) \qquad \Pk(X|\Pi)\text{-a.s.}.\end{align*} 
This shows:
\begin{align*}   Z \ismapof_{\Pk(X,\Theta|\Pi)} S  \end{align*} 

For the reverse let $S$ be a measurable deterministic map in  $(X,\Pi)$ such that:
\begin{align*}   Z \ismapof_{\Pk(X,\Theta|\Pi)} S  \end{align*} 
So there exists a measurable map $g:\, \Scal \to \Zcal=\Pcal(\Thetacal)$ such that:
\begin{align*}   Z =g(S) \quad \Pk(X|\Pi)\text{-a.s.}  \end{align*} 
We now define a Markov kernel via:
\begin{align*}  \Qk(\Theta|S):\, \Scal \dshto \Thetacal,\quad \Qk(\Theta \in B|S=s):=g(s)(B).\end{align*} 
Then we get:
\begin{align*}
& \left(\Qk(\Theta| S) \otimes \Pk(S,X|\Pi)\right)(B\times D \times A, \pi) \\
&= \int_{D \times A} \Qk(\Theta \in B|S=s) \, \Pk(S\in ds,X \in dx|\Pi=\pi) \\ 
&= \int_A \int_D \Qk(\Theta \in B|S=s) \,\deltabf(S \in ds|X=x,\Pi=\pi) \,\Pk(X \in dx|\Pi=\pi) \\
&= \int_A  \Qk(\Theta \in B|S=S(x,\pi)) \cdot \I_D(S(x,\pi))  \,\Pk(X \in dx|\Pi=\pi) \\
&= \int_A  g(S(x,\pi))(B) \cdot \I_D(S(x,\pi))  \,\Pk(X \in dx|\Pi=\pi) \\
&= \int_A  Z(x,\pi)(B) \cdot \I_D(S(x,\pi))  \,\Pk(X \in dx|\Pi=\pi) \\
&= \int_A  \Pk(\Theta \in B|X=x,\Pi=\pi) \cdot \deltabf(S\in D|X=x,\Pi=\pi)  \,\Pk(X \in dx|\Pi=\pi) \\
&= \int_A  \Pk(\Theta \in B,S \in D|X=x,\Pi=\pi) \,\Pk(X \in dx|\Pi=\pi) \\
&= \Pk(\Theta \in B, S \in D, X \in A |\Pi=\pi),
\end{align*}
where we used that $Z=g(S)$ holds $\Pk(X|\Pi)$-almost-surely, which is enough for the above equalities 
of the occurring integrals.
Since this holds for all inputs and measurable sets we get:
\begin{align*} \Pk(\Theta, S, X|\Pi) = \Qk(\Theta| S) \otimes \Pk(S,X|\Pi).  \end{align*} 
This implies the conditional independence:
\begin{align*}  \Theta \Indep_{\Pk(X,\Theta|\Pi) } X \given S.\end{align*} 
So all claims are shown.
\end{proof}
\end{Thm}

\begin{Thm}[A dual likelihood principle for Bayesian statistics]
    \label{thm:supp:bayesian-likelihood}
    We also have the transitional conditional independence with $R(\theta):=\Pk(X|\Theta=\theta)$:
    \begin{align*}  X \Indep_{\Pk(X,\Theta|\Pi) } \Theta,\Pi \given R. \end{align*} 
    Let $\Xcal$ be countably generated and $S:\,\Thetacal \to \Scal$ be another measurable map in $\Theta$.
    Then we get the equivalence:
    \begin{align*} X \Indep_{\Pk(X,\Theta|\Pi) } \Theta,\Pi \given S \qquad \iff \qquad R \ismapof_{\Pk(X,\Theta|\Pi)} S\end{align*}  
   \begin{proof} We have:
        \begin{align*}  \Pk(X,R,\Theta,\Pi|\Pi) = \Qk(X|R) \otimes \Pk(R,\Theta,\Pi|\Pi), \end{align*} 
        with $\Qk(X \in A|R=r):=r(A)$ for $r \in \Pcal(\Xcal)$, which implies the first claim:
     \begin{align*}  X \Indep_{\Pk(X,\Theta|\Pi) } \Theta,\Pi \given R. \end{align*} 
        For the second claim first assume:
        \begin{align*} X \Indep_{\Pk(X,\Theta|\Pi) } \Theta,\Pi \given S.\end{align*}  
        Then we get a factorization:
        \begin{align*} 
            \Pk(S,X,\Theta|\Pi) &= \Kk(X|S) \otimes \Pk(S,\Theta|\Pi)\\
                                       &= \Kk(X|S) \otimes \deltabf(S|\Theta) \otimes \Pk(\Theta|\Pi).
        \end{align*}
        Marginalizing $S$ out on both sides gives:
        \begin{align*}  \Pk(X|\Theta) \otimes \Pk(\Theta|\Pi) = \Kk(X|S=S(\Theta)) \otimes \Pk(\Theta|\Pi). \end{align*} 
        Since such factorizations are essentially unique by \Cref{ess-unique} and $\Xcal$ is countably generated, we have that
        for $\Pk(\Theta|\Pi)$-almost-all $(\theta,\pi)$ we get:
        \begin{align*}  R(\theta)=\Pk(X|\Theta=\theta) = \Kk(X|S=S(\theta)).  \end{align*} 
        This shows:
    \begin{align*}   R \ismapof_{\Pk(X,\Theta|\Pi)} S  \end{align*} 
    
    For the reverse direction now assume:
    \begin{align*}   R \ismapof_{\Pk(X,\Theta|\Pi)} S  \end{align*} 
    Then there exists a measurable map $g:\,\Scal \to \Pcal(\Xcal)$ such that:
    \begin{align*}  R=g(S) \quad \Pk(\Theta|\Pi)\text{-a.s.} \end{align*} 
    We then get the factorization:
    \begin{align*}  \Pk(X,R,S,\Theta,\Pi|\Pi) = \Qk(X|R) \otimes \deltabf_g(R|S) \otimes \Pk(S,\Theta,\Pi|\Pi) 
        \qquad \Pk(\Theta|\Pi)\text{-a.s.}, \end{align*} 
    where $\deltabf_g(R|S)$ is the deterministic Markov kernel given by $g$ and where the identity 
    $R=g(S)$ only holds $\Pk(\Theta|\Pi)$-almost-surely, which is enough for the equality of the two 
    (transition) probability distributions, as the exceptional set is a $\Pk(\Theta|\Pi)$-null set.
    This implies the conditional independence:
    \begin{align*} X,R \Indep_{\Pk(X,\Theta|\Pi) } \Theta,\Pi \given S,\end{align*} 
    and thus:
    \begin{align*} X\Indep_{\Pk(X,\Theta|\Pi) } \Theta,\Pi \given S.\end{align*} 
    So all claims are proven.
    \end{proof}
\end{Thm}
\section{Proofs - Reparameterization of Transitional Random Variables}
\label{sec:proofs-reparameterizing}
In this section we lift two classical one-dimensional facts (that $F(X)$ is uniform for a continuous cumulative
distribution function $F$, see \cite{Dar53}, and that $X$ can be written as a measurable function of $Z$ and
independent uniform noise, see \cite{Cen82}) from random variables to \emph{transitional} random variables.
What is new here is the parametrized version: all constructions are carried out simultaneously for all values $z$ of
the parameter, with joint measurability in $z$, and for distributions with atoms.
The one-dimensional building block is \Cref{unif-cdf} below: it is the classical quantile/inverse-transform construction
(see e.g.\ \cite{Kle20} Ch.\ 1), extended to distributions with atoms by interpolating the cumulative distribution function with
independent uniform noise; for kernel/randomization versions also see \cite{Kal17}.

\begin{Lem}[Interpolated cdf and quantile transform]
    \label{unif-cdf}
    Let $\bar\R := [-\infty,\infty]$ be endowed with the usual ordering and Borel $\sigma$-algebra,
    let $\Pk$ be a probability measure on $\bar\R$ and $\lambdabf$ the uniform distribution on $[0,1]$.
    For $x \in \bar\R$, $u \in [0,1]$ and $t \in [0,1]$ define:
    \begin{align*}
        F_0(x) &:= \Pk([-\infty,x)),  &   p(x) &:= \Pk(\{x\}),   &   F(x)=F_1(x) &:= \Pk([-\infty,x]),\\
        E(x;u)=F_u(x) &:= F_0(x) + u \cdot p(x),  &&&  R(t) &:= \inf F^{-1}([t,1]).
    \end{align*}
    Then:
    \begin{enumerate}
        \item $F$ is non-decreasing, right-continuous with $F(\infty)=1$, and $F_0(x)=\sup_{\tilde x < x} F(\tilde x)$
            (with $\sup \emptyset :=0$), so $F_0(-\infty)=0$ and:
            \[ F_u(\tilde x) \le F(\tilde x) \le F_0(x) \le F_u(x), \qquad \tilde x < x,\; u \in [0,1].\]
            Furthermore, $R:\,[0,1] \to \bar\R$ is non-decreasing with $R(0)=-\infty$, the maps $F_0,F,R$ are measurable and
            $E:\, \bar\R \times [0,1] \to [0,1]$ is (jointly) measurable.
        \item For all $x \in \bar\R$ and $t \in [0,1]$ we have the equivalence:
            \[ t \le F(x) \iff R(t) \le x. \]
            In particular, $F(R(t)) \ge t$, i.e.\ $R(t) = \min F^{-1}([t,1])$, and $R(F(x)) \le x$.
            Moreover, $R$ is a reflexive generalized inverse of $F$, i.e.:
            \[ F \circ R \circ F = F, \qquad R \circ F \circ R = R,\]
            and $R_*\lambdabf = \Pk$.
        \item Let $\bar\Pk := \Pk \otimes \lambdabf$ be the product distribution on $\bar\R \times [0,1]$ and let
            $X(x,u):=x$ and $U(x,u):=u$ be the two projections, so that $X$ has distribution $\Pk$ and $U$ is uniformly
            distributed under $\bar\Pk$. Then the random variable:
            \[ \begin{array}{ccccl}
                E &:& \bar\R\times[0,1] &\to& [0,1], \\
                  && (x,u) &\mapsto& \Pk([-\infty,x))+ u \cdot \Pk(\{x\}),
            \end{array} \]
            is uniformly distributed under $\bar\Pk$, i.e.\ $\bar\Pk(E \le e) = e$ for every $e \in [0,1]$, and:
            \[ R(E) = X \qquad \bar\Pk\text{-a.s.} \]
    \end{enumerate}
    \begin{proof}
        1.) Monotonicity, $F(\infty)=1$ and the right-continuity of $F$ ($[-\infty,x_n] \downarrow [-\infty,x]$ for
        $x_n \downarrow x$) follow from the continuity from above of $\Pk$, and
        $F_0(x)=\sup_{\tilde x<x}F(\tilde x)$ from the continuity from below applied to
        $[-\infty,x) = \bigcup_{\tilde x < x}[-\infty,\tilde x]$. The displayed chain of inequalities is immediate from this
        and from $F_0 \le F_u \le F_1=F$. Monotone maps between (subsets of) $\bar\R$ are Borel measurable, so $F_0,F,R$ are
        measurable, and $E(x;u)=F_0(x)+u \cdot (F(x)-F_0(x))$ is jointly measurable as the composition of the measurable map
        $(x,u) \mapsto (F_0(x),F(x),u)$ with the continuous map $(a,b,u)\mapsto a + u\cdot(b-a)$. Finally
        $R(0)=\inf \bar\R = -\infty$.

        2.) The set $A_t := F^{-1}([t,1])$ is non-empty (as $F(\infty)=1 \ge t$), it is an up-set (as $F$ is non-decreasing)
        and it is closed under non-increasing limits (as $F$ is right-continuous). So $A_t = [R(t),\infty]$ with
        $R(t)= \min A_t$, which is precisely the claimed equivalence:
        \[ t \le F(x) \iff x \in A_t \iff R(t) \le x.\]
        Putting $x:=R(t)$ resp.\ $t:=F(x)$ gives $F(R(t)) \ge t$ resp.\ $R(F(x)) \le x$.
        Applying the non-decreasing map $R$ to $t \le F(R(t))$ gives $R(t) \le R(F(R(t)))$, while $R(F(x))\le x$ with
        $x:=R(t)$ gives the reverse inequality, so $R\circ F \circ R = R$. For the identity $F \circ R \circ F = F$ apply
        the non-decreasing map $F$ to $R(F(x)) \le x$ to get $F(R(F(x))) \le F(x)$, while $F(R(t)) \ge t$ with
        $t:=F(x)$ gives the reverse inequality. For $R_*\lambdabf$ we get for every $x \in \bar\R$:
        \begin{align*}
            (R_*\lambdabf)([-\infty,x]) &= \lambdabf(\lC t \in [0,1]\,|\, R(t) \le x \rC) \\
            &= \lambdabf(\lC t \in [0,1]\,|\, t \le F(x) \rC) = F(x) = \Pk([-\infty,x]).
        \end{align*}
        Since the sets $[-\infty,x]$, $x \in \bar\R$, form a $\cap$-stable generator of the Borel $\sigma$-algebra of
        $\bar\R$, this implies $R_*\lambdabf = \Pk$.

        3.) Fix $e \in [0,1]$, put $x_e := R(e)$ and $\delta := e - F_0(x_e)$. By 2.) we have $F(x_e)=F(R(e)) \ge e$, and by
        2.) again every $x < x_e=R(e)$ satisfies $F(x) < e$, thus $F_0(x_e)=\sup_{x<x_e}F(x) \le e$. Together:
        \[ 0 \le \delta \le F(x_e)-F_0(x_e) = p(x_e). \]
        Note also that $x > x_e$ implies $F_0(x) \ge F(x_e) \ge e$ by 1.).

        We first record the following elementary computation. Since $E(x_e;u)=F_0(x_e) + u\cdot p(x_e)$ we get by Fubini:
        \[ \bar\Pk\lp X=x_e,\; U \cdot p(x_e) \le \delta \rp = p(x_e) \cdot \lambdabf(\lC u \in [0,1]\,|\, u \cdot p(x_e) \le \delta\rC) = \delta,\]
        because for $p(x_e)>0$ the right hand side equals $p(x_e) \cdot \min\lp 1, \delta/p(x_e)\rp = \delta$ by
        $\delta \le p(x_e)$, and for $p(x_e)=0$ both sides vanish by $0 \le \delta \le p(x_e)=0$. The same computation with
        ``$\le \delta$'' replaced by ``$< \delta$'' gives a value $\le \delta$ as well.

        \emph{Lower bound.} If $x < x_e$ then $E(x;u) \le F(x) < e$ for every $u$, and if $x=x_e$ and
        $u \cdot p(x_e) \le \delta$ then $E(x_e;u) = F_0(x_e)+ u\cdot p(x_e) \le F_0(x_e)+\delta = e$. So:
        \[ \lC X < x_e \rC \;\cup\; \lC X=x_e,\, U \cdot p(x_e) \le \delta \rC \;\ins\; \lC E \le e\rC, \]
        and since the two events on the left are disjoint with $\bar\Pk(X<x_e) = F_0(x_e)$ we get:
        \[ \bar\Pk(E \le e) \ge F_0(x_e) + \delta = e. \]

        \emph{Upper bound.} If $x > x_e$ then $E(x;u) \ge F_0(x) \ge e$, so $\lC E < e \rC \ins \lC X \le x_e \rC$ and thus:
        \[ \lC E<e \rC \;\ins\; \lC X < x_e \rC \;\cup\; \lC X=x_e,\, U \cdot p(x_e) < \delta \rC, \]
        which gives $\bar\Pk(E<e) \le F_0(x_e)+\delta = e$ for every $e \in [0,1]$. For $e<1$ and any $e' \in (e,1]$ we then
        get $\bar\Pk(E \le e) \le \bar\Pk(E<e') \le e'$, and thus $\bar\Pk(E \le e)\le e$; for $e=1$ this holds trivially.
        Together with the lower bound this shows $\bar\Pk(E \le e)=e$ for all $e \in [0,1]$, i.e.\ $E$ is uniformly
        distributed under $\bar\Pk$.

        \emph{The identity $R(E)=X$.} Since $E(x;u) \le F(x)$ for all $(x,u)$, the equivalence in 2.) gives the pointwise
        inequality $R(E) \le X$. Furthermore, by 2.) and the uniformity of $E$ just shown we have for every $x \in \bar\R$:
        \[ \bar\Pk( R(E) \le x ) = \bar\Pk( E \le F(x)) = F(x) = \Pk([-\infty,x]),\]
        so $R(E)$ has the same distribution $\Pk$ as $X$ under $\bar\Pk$ (again by uniqueness of measures on the
        $\cap$-stable generator $\{[-\infty,x]\,|\,x \in \bar\R\}$). Now fix any strictly increasing map
        $\varphi:\, \bar\R \to [0,1]$, e.g.\ $\varphi(x):=\frac{1}{2}+\frac{1}{\pi}\arctan(x)$ with
        $\varphi(\pm\infty):=\frac{1}{2}\pm\frac{1}{2}$. Then $\varphi(X)-\varphi(R(E)) \ge 0$ pointwise and:
        \[ \int \varphi(X) \,d\bar\Pk = \int \varphi \,d\Pk = \int \varphi(R(E))\,d\bar\Pk \in [0,1], \]
        so the non-negative integrand $\varphi(X)-\varphi(R(E))$ has vanishing integral and thus
        $\varphi(X)=\varphi(R(E))$ $\bar\Pk$-a.s. Since $\varphi$ is injective this shows:
        \[ R(E) = X \qquad \bar\Pk\text{-a.s.} \]
    \end{proof}
\end{Lem}

\begin{Thm}
    \label{ccdf2}
    Let $\Zcal$ be any measurable space and $\Xcal$ be a standard measurable space with a fixed embedding $\iota: \Xcal \inj \bar\R=[-\infty,\infty]$ onto a Borel subset (which always exists; we therefore identify $\Xcal$ with that Borel subset and argue on
    $\bar\R$ throughout, which is harmless since $\Kk(X|Z=z)$ is carried by $\Xcal$).
    Let $\Kk(X|Z) : \Zcal \dshto \Xcal$ be a Markov kernel.
    Furthermore, let $\Ucal:=[0,1]$ and $\Kk(U)$ be the uniform distribution/Markov kernel on $\Ucal$. We write:
    \[ \Kk(U,X|Z) := \Kk(U)\otimes \Kk(X|Z).\]
Also put:
\begin{align*}
    F(x;u|z) &:= \Kk(X<x|Z=z) + u \cdot \Kk(X=x|Z=z)\\
	R(e|z) &:= \inf \left\{ \tilde{x} \in \bar\R\,| F(\tilde{x};1|z) \ge e \right\}.
\end{align*}
Let $E:= F(X;U|Z)$. We consider $X,U,Z,E$  as the measurable maps:
\[\begin{array}{cccl}
    X: & \Xcal \times \Ucal \times \Zcal &\to& \Xcal,\\
     & (x,u,z) &\mapsto& x,\\
    U: & \Xcal \times \Ucal \times \Zcal &\to& \Ucal,\\
     & (x,u,z) &\mapsto& u,\\
    Z: & \Xcal \times \Ucal \times \Zcal &\to& \Zcal,\\
     & (x,u,z) &\mapsto& z,\\
    E: & \Xcal \times \Ucal \times \Zcal &\to& \Ecal:=[0,1],\\
     & (x,u,z) &\mapsto& F(x;u|z).
\end{array}\]
Then for all $e \in \Ecal$ and $z \in \Zcal$ we have:
\[ \Kk( E \le e | Z=z ) = e, \]
implying $E \Indep_{\Kk(U,X|Z)} Z$.
Furthermore, we have:
\[ X = R(E|Z) \quad \Kk(U,X|Z)\text{-a.s.}.\]
\begin{proof}
    Once the measurabilities are established, the two displayed statements follow from \Cref{unif-cdf} 3.\ applied
    to $\Pk := \Kk(X|Z=z)$ for every $z \in \Zcal$ separately, since for fixed $z$ the objects
    $F(\cdot\,;\cdot|z)$ and $R(\cdot|z)$ are precisely the maps $E$ and $R$ of \Cref{unif-cdf} for that $\Pk$.
    So it remains to prove that:
    \[ (x,u,z) \mapsto F(x;u|z) \qquad\text{and}\qquad (e,z) \mapsto R(e|z) \]
    are jointly measurable, which is the only point where the parameter $z$ causes work.\\
    \emph{Step 1: $(x,z) \mapsto \Kk(X \le x|Z=z)$ is jointly measurable.}
    For each fixed $q \in \bar\R$ the map $z \mapsto \Kk(X \le q|Z=z)$ is measurable, being the composition of
    the measurable map $\Kk(X|Z):\, \Zcal \to \Pcal(\bar\R)$ with the evaluation map $j_{[-\infty,q]}$.
    Hence for every $q \in \Q$ the map
    \[ h_q:\, \bar\R \times \Zcal \to [0,1], \qquad
       h_q(x,z) := \I_{[-\infty,q]}(x) \cdot \Kk(X \le q|Z=z) + \I_{(q,\infty]}(x) \]
    is jointly measurable. Since $\Kk(X|Z=z)$ is a probability measure, continuity from above gives
    $\Kk(X \le x|Z=z) = \inf_{q \in \Q} h_q(x,z)$ for every $x \in \bar\R$: the terms with $q \ge x$ contribute
    $\Kk(X \le q|Z=z) \downarrow \Kk(X \le x|Z=z)$, and the terms with $q < x$ contribute the value $1$, which is
    never smaller; for $x = +\infty$ there is no rational $q \ge x$ and all terms equal $1 = \Kk(X \le \infty|Z=z)$,
    and for $x=-\infty$ the infimum is $\Kk(X = -\infty|Z=z)$ by continuity from above.
    A countable infimum of jointly measurable maps is jointly measurable.\\
    \emph{Step 2: $(x,z) \mapsto \Kk(X<x|Z=z)$ and $(x,z)\mapsto \Kk(X=x|Z=z)$ are jointly measurable.}
    Analogously, continuity from below gives
    \[ \Kk(X<x|Z=z) = \sup_{q \in \Q} \lp \I_{(q,\infty]}(x) \cdot \Kk(X \le q|Z=z) \rp, \]
    a countable supremum of jointly measurable maps, and
    $\Kk(X=x|Z=z) = \Kk(X \le x|Z=z) - \Kk(X<x|Z=z)$.
    Consequently $F(x;u|z) = \Kk(X<x|Z=z) + u \cdot \Kk(X=x|Z=z)$ is jointly measurable in $(x,u,z)$.\\
    \emph{Step 3: $(e,z) \mapsto R(e|z)$ is jointly measurable.}
    By \Cref{unif-cdf} 2., applied for each fixed $z$, we have for all $e \in [0,1]$ and $x \in \bar\R$ the
    equivalence $R(e|z) \le x \iff e \le F(x;1|z)$. Hence for every $x \in \bar\R$:
    \[ \lC (e,z) \in [0,1] \times \Zcal \,|\, R(e|z) \le x \rC
       = \lC (e,z) \,|\, e \le \Kk(X \le x|Z=z) \rC \in \Bcal_{[0,1]} \otimes \Bcal_\Zcal, \]
    by Step 1. Since the sets $[-\infty,x]$, $x \in \bar\R$, generate $\Bcal_{\bar\R}$, the map
    $(e,z) \mapsto R(e|z)$ is jointly measurable.\\
    Finally, $E = F(X;U|Z)$ is measurable as the composition of the coordinate projections with $F$, and
    $R(E|Z)$ is measurable as the composition of $(E,Z)$ with $R$.
\end{proof}
\end{Thm}

\begin{Rem}[Transfer to an arbitrary transition probability space]
    \label{rem:reparam-transfer}
    \Cref{ccdf2} is formulated on the canonical space $\Xcal \times \Ucal \times \Zcal$, whereas
    \Cref{thm:reparamererizing-trv} of the main paper is formulated on an arbitrary transition probability space
    $\lp \bar\Wcal \times \Zcal, \Kk(\bar W|Z)\rp$ with $\bar\Wcal = \Ucal \times \Wcal$. The passage from
    the former to the latter is immediate: apply \Cref{ccdf2} to the push-forward Markov kernel $\Kk(X|Z)$ of
    $\Kk(W|Z)$ along $X$. Both conclusions only refer to the joint Markov kernel
    $\Kk(X,U|Z) = \Kk(U) \otimes \Kk(X|Z)$: the first states $\Kk(E \le e|Z=z)=e$ for the measurable function
    $E = F(X;U|Z)$ of $(X,U,Z)$, and the second states that the two measurable functions
    $(x,u,z) \mapsto x$ and $(x,u,z)\mapsto R(F(x;u|z)|z)$ of $(X,U,Z)$ agree $\Kk(X,U|Z=z)$-almost surely for every
    $z \in \Zcal$. Since $\Kk(X,U|Z=z)$ is exactly the push-forward of $\Kk(\bar W|Z=z)$ along $(X,U)$, both
    statements pull back to $\bar\Wcal \times \Zcal$, which proves \Cref{thm:reparamererizing-trv}.
\end{Rem}

\begin{Cor}
    Let $X$ and $Z$ be random variables with values in any standard measurable spaces $\Xcal$ and $\Zcal$, resp., and with a joint distribution $\Pk(X,Z)$.
    Then there exists a uniformly distributed random variable $E$ on $[0,1]$ that is $\Pk$-independent of $Z$ and a measurable function $g$ such that $X = g(E,Z)$ $\Pk$-almost-surely.
    Furthermore, $E$ can be constructed via a deterministic measurable function in $X$ and $Z$ and (uniformly distributed) independent noise $U$ (on $[0,1]$).
    \begin{proof}
        The regular conditional probability distribution $\Pk(X|Z)$ exists for standard measurable spaces (and is unique up to a $\Pk(Z)$-zero-set), and is a Markov kernel. Then apply the result from above for $\Kk(X|Z):=\Pk(X|Z)$ to get $g(e,z):=R(e|z)$ and $E$.
    \end{proof}
\end{Cor}

\begin{Rem}
    Any Polish space, i.e.\ any completely metrizable topological space with a countable dense subset (separable), is a standard measurable space in its Borel $\sigma$-algebra. These are fundamental theorems in classical descriptive set theory, see
    \cite{Bog07, Fremlin, Kec95}.
    Examples of Polish and thus standard measurable spaces are $[0,1]$, $\R$, $\R^d$, $\N$, any (discrete) finite or countable set, any topological or smooth manifold $\Xcal$, any finite (or even countable) CW-complex $\Ycal$, etc., (in its usual Borel $\sigma$-algebra). So these are all measurably isomorphic to a Borel subset of $[0,1]$ (or $\bar\R$), and measurably isomorphic to $[0,1]$ (or $\bar\R$) itself if non-countable (excluding finite and countable sets).
\end{Rem}
\section{Proofs - Separoid Rules for d-Separation and id-Separation}
\label{sec:supp:separoid-rules-sigma}

In the following let $\Gk=(J,V,E)$ be a CDAG, i.e.\ an acyclic conditional directed graph with input nodes $J$,
output nodes $V$ and directed edges $E$, and let $A,B,C,D \ins J \cup V$ (not necessarily disjoint) be subsets of nodes.
We abbreviate the two ternary relations in the following:
$\displaystyle \Perp^d \; := \; \Perp^d_{\Gk}$ and $\displaystyle \Perp^{\id} \; := \; \Perp^{\id}_{\Gk}$.\\

Recall from \Cref{def:d-blocked} and \Cref{def:d-separation} that $A \Perp^d B \given C$ holds if and only if
every walk from a node in $A$ to a node in $B$ is blocked by $C$, and that a walk
$\pi = \lp v_0 \sus \cdots \sus v_n \rp$, $n \ge 0$, is \emph{blocked by $C$} (or \emph{$C$-blocked}) if it has a
non-collider in $C$ or a collider outside of $\Anc^\Gk(C)$.
Here the node $v_k$ of $\pi$ is a \emph{collider of $\pi$} if two arrow heads of $\pi$ point at it, i.e.\ if
$0<k<n$ and $v_{k-1} \tuh v_k \hut v_{k+1}$, and a \emph{non-collider of $\pi$} if at most one arrow head of $\pi$
points at it, which is the case for the two \emph{end nodes} $v_0$, $v_n$ and for those \emph{inner nodes} $v_k$,
$0<k<n$, that form a chain ($v_{k-1} \tuh v_k \tuh v_{k+1}$ or $v_{k-1} \hut v_k \hut v_{k+1}$) or a fork
($v_{k-1} \hut v_k \tuh v_{k+1}$).
A walk that is not $C$-blocked is called \emph{$C$-open}.
Note that the same node of $\Gk$ may occur at several positions of $\pi$ and may be a collider at some of these
positions and a non-collider at others.
Recall further that id-separation is defined via d-separation by:
\[ A \Perp^{\id} B \given C \qquad :\iff \qquad A \Perp^{d} B \cup J \given C. \]

Every walk that we will actually have to inspect can be taken to have all its colliders in $C$ itself. We therefore
call $\pi$ \emph{strictly $C$-open} if every non-collider of $\pi$ lies outside of $C$ and every collider of $\pi$
lies \emph{in $C$}, and \emph{strictly $C$-open at its inner nodes} if every collider of $\pi$ lies in $C$ and every
inner non-collider of $\pi$ lies outside of $C$. With this wording:
\[ \pi \text{ is strictly } C\text{-open} \quad \iff \quad
\lp \pi \text{ is strictly } C\text{-open at its inner nodes} \rp \;\land\; v_0,v_n \notin C. \]
Since $C \ins \Anc^\Gk(C)$, every strictly $C$-open walk is $C$-open. The converse fails for an individual walk, but
not for the existence of one, and that is all we shall need:

\begin{Lem}[Colliders may be assumed to lie in $C$]
    \label{lem:strictly-open}
    Let $C \ins J \cup V$ and let $\pi$ be a $C$-open walk from a node $a$ to a node $b$ in $\Gk$. Then there is a
    \emph{strictly} $C$-open walk from $a$ to $b$. Consequently, for all $A,B,C \ins J \cup V$:
    \[ A \Perp^d B \given C \quad \iff \quad
       \text{no strictly } C\text{-open walk runs from a node in } A \text{ to one in } B. \]
\begin{proof}
    Let $\pi = \lp v_0 \sus \cdots \sus v_n \rp$ be $C$-open and suppose some collider $v_k$ of $\pi$ satisfies
    $v_k \notin C$. Being a collider of a $C$-open walk, $v_k \in \Anc^\Gk(C)$, so there is a directed walk
    $v_k \tuh z_1 \tuh \cdots \tuh z_r$ in $\Gk$ with $z_r \in C$; choosing $r \ge 1$ minimal we may assume
    $z_1,\dots,z_{r-1} \notin C$. Replace the occurrence of $v_k$ in $\pi$ by the detour:
    \[ v_{k-1} \tuh v_k \tuh z_1 \tuh \cdots \tuh z_r \hut z_{r-1} \hut \cdots \hut z_1 \hut v_k \hut v_{k+1}. \]
    In the resulting walk the two occurrences of $v_k$ are chains, hence non-colliders, and $v_k \notin C$; the two
    occurrences of each $z_i$, $i<r$, are chains, hence non-colliders, and $z_i \notin C$; the single occurrence of
    $z_r$ is a collider, and $z_r \in C$. All other positions keep their collider/non-collider status and their
    node, and the two end nodes are unchanged. So the new walk is again a $C$-open walk from $a$ to $b$, and it has
    strictly fewer collider positions outside of $C$ than $\pi$, since the detour creates none. Iterating removes
    them all and yields a strictly $C$-open walk from $a$ to $b$.

    For the displayed equivalence, ``$\implies$'' holds because every strictly $C$-open walk is $C$-open, and
    ``$\Longleftarrow$'' is the statement just proven.
\end{proof}
\end{Lem}

Whenever a collider has to be inspected below we will work with strictly $C$-open walks.

\begin{Rem}[Properties of strictly $C$-open walks]
    \label{rem:d-open-walks}
    Let $\pi=\lp v_0 \sus \cdots \sus v_n \rp$ be a walk in $\Gk$ and $C,D \ins J \cup V$.
    \begin{enumerate}
        \item[a)] \emph{(Locality.)} Whether $\pi$ is strictly $C$-open depends only on the node sequence of
            $\pi$, on the orientations of its edges and on the set $C$; no further reference to the ambient
            graph $\Gk$ occurs. (This is the reason for working with strict openness: $C$-openness itself refers to
            $\Anc^\Gk(C)$, and hence to $\Gk$.)
        \item[b)] \emph{(Reversal.)} The reversed walk $\pi^{-}:=\lp v_n \sus \cdots \sus v_0 \rp$ has the same
            colliders and non-colliders as $\pi$ and the same end nodes. So $\pi^{-}$ is strictly $C$-open iff $\pi$
            is, and $\pi^-$ is $C$-open iff $\pi$ is.
        \item[c)] \emph{(Sub-walks.)} For $0 \le m \le l \le n$ put $\pi_{[m,l]}:=\lp v_m \sus \cdots \sus v_l \rp$.
            Every inner node $v_k$, $m<k<l$, of $\pi_{[m,l]}$ has the same two adjacent edges in $\pi_{[m,l]}$ as in
            $\pi$ and thus is a collider of $\pi_{[m,l]}$ iff it is a collider of $\pi$ (at that position). Consequently,
            if $\pi$ is strictly $C$-open at its inner nodes then so is $\pi_{[m,l]}$.
        \item[d)] \emph{(Changing the conditioning set.)} If no inner node of $\pi$ lies in $D \sm C$ then $\pi$ is
            strictly $C$-open at its inner nodes if and only if it is strictly $(D\cup C)$-open at its inner nodes.
            If, in addition, no end node of $\pi$ lies in $D\sm C$, then:
            \[ \pi \text{ is strictly } C\text{-open} \quad \iff \quad \pi \text{ is strictly } (D\cup C)\text{-open}. \]
        \item[e)] \emph{(End nodes.)} Every walk with an end node in $C$ is $C$-blocked, since end nodes are
            non-colliders; in particular it is not strictly $C$-open. So $A \Perp^d B \given C$ holds iff no strictly
            $C$-open walk runs from a node in $A \sm C$ to a node in $B \sm C$.
    \end{enumerate}
\begin{proof}
    Only d) needs an argument. For every node $x \in J \cup V$ with $x \notin D \sm C$ we have the equivalence:
    $x \in C \iff x \in D \cup C$. Applying this to all inner nodes of $\pi$ (resp.\ to all inner nodes and both end
    nodes of $\pi$) gives the two claims.
\end{proof}
\end{Rem}

\begin{Lem}[Minimal open walks]
    \label{lem:minimal-open-walk}
    Let $X,Y,C \ins J \cup V$ and assume that there exists a strictly $C$-open walk from a node in $X$ to a node in
    $Y$. Among all such walks let $\pi = \lp v_0 \sus \cdots \sus v_n \rp$, $v_0 \in X$, $v_n \in Y$, be one of
    minimal length $n$. Then $v_0 \in X \sm C$, $v_n \in Y \sm C$ and for all $0 \le k \le n$:
    \[ v_k \in X\sm C \implies k=0, \qquad\qquad v_k \in Y \sm C \implies k=n. \]
\begin{proof}
    Since $\pi$ is strictly $C$-open its end nodes lie outside of $C$, which gives $v_0 \in X\sm C$ and
    $v_n \in Y \sm C$.
    Now assume $v_k \in X\sm C$ for some $k \ge 1$. By \Cref{rem:d-open-walks} c) the sub-walk $\pi_{[k,n]}$ is
    strictly $C$-open at its inner nodes, and its end nodes $v_k$ and $v_n$ lie outside of $C$. So $\pi_{[k,n]}$ is a
    strictly $C$-open walk from $v_k \in X$ to $v_n \in Y$ of length $n-k<n$, contradicting the minimality of $n$. The second implication
    follows in the same way, using the sub-walk $\pi_{[0,k]}$ for $v_k \in Y \sm C$ with $k \le n-1$.
\end{proof}
\end{Lem}

\subsection{Operations on Graphs}
\label{sec:operations-on-graphs}

The only operation on graphs that we will need in the following is the removal of a childless node.

\begin{Lem}[Removing a childless node preserves d-separation]
    \label{lem:remove-childless}
    Let $\Gk=(J,V,E)$ be a CDAG and $v \in V$ with $\Ch^\Gk(v)=\emptyset$. Let
    \[ \Gk_{-v}:=\lp J,\, V \sm \{v\},\, E_{-v} \rp, \qquad
       E_{-v}:= \lC w \tuh u \in E \,|\, u,w \in (J \cup V)\sm\{v\} \rC, \]
    be the induced subgraph of $\Gk$ on the nodes $(J, V\sm\{v\})$. Then $\Gk_{-v}$ is again a CDAG and for all subsets
    $A,B,C \ins (J\cup V)\sm\{v\}$ we have the equivalence:
    \[ A \Perp^d_{\Gk} B \given C \qquad \iff \qquad A \Perp^d_{\Gk_{-v}} B \given C. \]
\begin{proof}
    First note that $E_{-v} \ins E$. So $\Gk_{-v}$ contains no directed cycle, no edge of $\Gk_{-v}$ has an arrow head
    pointing to a node of $J$ and every edge of $\Gk_{-v}$ points to a node of $V \sm\{v\}$. So $\Gk_{-v}$ is a CDAG.
    Furthermore, $\Ch^\Gk(v)=\emptyset$ means that no edge of $E$ starts at $v$, so every edge of $E$ that is incident to
    $v$ is of the form $w \tuh v$ with $w \in \Pa^\Gk(v)$ and thus:
    \[ E_{-v} = E \sm \lC w \tuh v \,|\, w \in \Pa^\Gk(v) \rC. \]

    ``$\implies$'': Assume $A \Perp^d_\Gk B \given C$ and let $\pi$ be a walk in $\Gk_{-v}$ from a node in $A$ to a node
    in $B$. Since $E_{-v} \ins E$, the same sequence of nodes and edges is a walk in $\Gk$ from a node in $A$ to a
    node in $B$, hence, by \Cref{lem:strictly-open}, it is not strictly $C$-open. By \Cref{rem:d-open-walks} a)
    strict $C$-openness does not depend on which of the two graphs the walk is considered in, so $\pi$ is not
    strictly $C$-open in $\Gk_{-v}$ either. Since $\pi$ was arbitrary, \Cref{lem:strictly-open} applied in
    $\Gk_{-v}$ gives $A \Perp^d_{\Gk_{-v}} B \given C$.\\

    ``$\Longleftarrow$'': Assume $A \Perp^d_{\Gk_{-v}} B \given C$ and let $\pi=\lp v_0 \sus \cdots \sus v_n\rp$ be a
    walk in $\Gk$ with $v_0 \in A$ and $v_n \in B$. We distinguish two cases.\\
    Case 1: $v_k \neq v$ for all $k=0,\dots,n$. Then every edge of $\pi$ joins two nodes different from $v$ and thus lies
    in $E_{-v}$. So $\pi$ is a walk in $\Gk_{-v}$ from a node in $A$ to a node in $B$ and hence, by
    \Cref{lem:strictly-open} applied in $\Gk_{-v}$, not strictly $C$-open; by \Cref{rem:d-open-walks} a) it is not
    strictly $C$-open as a walk in $\Gk$ either.\\
    Case 2: $v_k = v$ for some $k \in \{0,\dots,n\}$. Since $A \ins (J\cup V)\sm\{v\}$ and $B \ins (J\cup
    V)\sm\{v\}$ we have $v_0 \neq v$ and $v_n \neq v$, so $0<k<n$ and $v_k$ is an inner node of $\pi$. Since
    $\Ch^\Gk(v)=\emptyset$ there is no edge in $\Gk$ that starts at $v$, so both edges of $\pi$ adjacent to the position
    $k$ point towards $v_k=v$:
    \[ v_{k-1} \tuh v_k \hut v_{k+1}. \]
    So $v_k$ is a collider of $\pi$. Since $C \ins (J\cup V)\sm\{v\}$ we
    have $v_k = v \notin C$, and thus $\pi$ is not strictly $C$-open.\\
    In both cases $\pi$ is not strictly $C$-open, which by \Cref{lem:strictly-open} shows
    $A \Perp^d_\Gk B \given C$.
\end{proof}
\end{Lem}

\begin{Rem}
    \label{rem:remove-childless-id}
    The corresponding statement for id-separation follows immediately: the two graphs $\Gk$ and $\Gk_{-v}$ have the
    same set of input nodes $J$, and $v \in V$ implies $v \notin J$, so for $B \ins (J\cup V)\sm\{v\}$ we also have
    $B \cup J \ins (J \cup V)\sm\{v\}$. Applying \Cref{lem:remove-childless} to the set $B \cup J$ thus gives, for all
    $A,B,C \ins (J\cup V)\sm\{v\}$:
    \[ A \Perp^{\id}_{\Gk} B \given C \qquad \iff \qquad A \Perp^{\id}_{\Gk_{-v}} B \given C. \]
\end{Rem}

\subsection{Symmetric Separoid Rules for d-Separation}

In this subsection we prove the classical, symmetric separoid rules for d-separation, see
\Cref{thm:d-sep-sym-rules} and \cite{PP85, Spo94, Daw01, Gei90, Ver93, Lau90, Lau96, SGS2000, Pearl09}.
These are the only statements about walks that will be needed: all rules for id-separation in the following
subsections are then derived from them by purely formal arguments.

\begin{Lem}[Symmetry]
    \label{sep:d:sym}
    \[ A \Perp^d B \given C \implies B \Perp^d A \given C. \]
\begin{proof}
    Let $\pi$ be a walk from a node in $B$ to a node in $A$. By \Cref{rem:d-open-walks} b) the reversed walk $\pi^-$
    is a walk from a node in $A$ to a node in $B$, and $\pi^-$ is $C$-blocked iff $\pi$ is. By assumption $\pi^-$ is
    $C$-blocked, hence so is $\pi$.
\end{proof}
\end{Lem}

\begin{Lem}[Redundancy]
    \label{sep:d:red}
    \[ A \ins C \implies A \Perp^d B \given C. \]
\begin{proof}
    If $\pi$ is a walk from a node $v$ in $A$ to a node $w$ in $B$ then its first end node $v$ lies in $A \ins C$.
    End nodes are non-colliders, so $\pi$ is blocked by $C$, see \Cref{rem:d-open-walks} e).
\end{proof}
\end{Lem}

\begin{Lem}[Decomposition]
    \label{sep:d:dec}
    \[ A \Perp^d B \cup D \given C \implies A \Perp^d B \given C. \]
\begin{proof}
    If $\pi$ is a walk from a node in $A$ to a node in $B$ then $\pi$ is a walk from a node in $A$ to a node in
    $B \cup D$, which by assumption is blocked by $C$.
\end{proof}
\end{Lem}

\begin{Lem}[Weak Union]
    \label{sep:d:uni}
    \[ A \Perp^d B \cup D \given C \implies A \Perp^d B \given D \cup C. \]
\begin{proof}
    Assume the contrary. Then, by \Cref{lem:strictly-open}, there is a strictly $(D \cup C)$-open walk
    $\pi = \lp v_0 \sus \cdots \sus v_n\rp$ with $v_0 \in A$ and $v_n \in B$. In particular $v_0, v_n \notin D \cup C$.
    Put $m:=n$ if there is no index $k$ with $v_k \in D \sm C$, and otherwise let $m$ be the smallest such index, in
    which case $1 \le m \le n-1$. In both cases:
    \[ v_m \in B \cup D, \qquad v_0 \notin C, \qquad v_m \notin C. \]
    None of the inner nodes $v_k$, $0<k<m$, of $\tilde \pi := \pi_{[0,m]}$ lies in $D \sm C$, so
    \Cref{rem:d-open-walks} c) and d) show that $\tilde \pi$ is strictly $C$-open at its inner nodes, and thus
    strictly $C$-open.
    So $\tilde \pi$ is a strictly $C$-open walk from a node in $A$ to a node in $B \cup D$, which by
    \Cref{lem:strictly-open} contradicts the assumption:
    $A \Perp^d B \cup D \given C$.
\end{proof}
\end{Lem}

\begin{Lem}[Contraction]
    \label{sep:d:con}
    \[ (A \Perp^d B \given D \cup C) \land (A \Perp^d D \given C) \implies A \Perp^d B \cup D \given C. \]
\begin{proof}
    Assume the contrary. Then, by \Cref{lem:strictly-open}, we can pick a strictly $C$-open walk
    $\pi = \lp v_0 \sus \cdots \sus v_n \rp$ from a node
    $v_0 \in A$ to a node $v_n \in B \cup D$ of minimal length. By \Cref{lem:minimal-open-walk} (with $Y = B \cup D$)
    we have $v_n \in (B \cup D)\sm C$ and:
    \[ v_k \in (B \cup D) \sm C \implies k = n. \tag{$\ast$} \]
    If $v_n \in D$ then $\pi$ is a strictly $C$-open walk from a node in $A$ to a node in $D$, contradicting the
    assumption:
    $A \Perp^d D \given C$. So $v_n \in B \sm C$ and $v_n \notin D$.
    Then no node of $\pi$ lies in $D \sm C$: if $v_k \in D \sm C \ins (B \cup D)\sm C$ then $k=n$ by $(\ast)$ and thus
    $v_n \in D$, which we just excluded. By \Cref{rem:d-open-walks} d) the walk $\pi$ is therefore a strictly
    $(D\cup C)$-open walk from a node in $A$ to a node in $B$. This contradicts the other assumption:
    $A \Perp^d B \given D \cup C$.
\end{proof}
\end{Lem}

\begin{Lem}[Composition]
    \label{sep:d:com}
    \[ (A \Perp^d B \given C) \land (A \Perp^d D \given C) \implies A \Perp^d B \cup D \given C. \]
\begin{proof}
    Let $\pi$ be a walk from a node in $A$ to a node $w \in B \cup D$. If $w \in B$ then $\pi$ is blocked by $C$ by
    assumption: $A \Perp^d B \given C$. If $w \in D$ then $\pi$ is blocked by $C$ by assumption:
    $A \Perp^d D \given C$.
\end{proof}
\end{Lem}

\begin{Lem}[Intersection]
    \label{sep:d:int}
    Assume that $B \cap D = \emptyset$, then:
    \[ (A \Perp^d B \given D \cup C) \land (A \Perp^d D \given B \cup C) \implies A \Perp^d B \cup D \given C. \]
\begin{proof}
    Assume the contrary. Then, by \Cref{lem:strictly-open}, we can pick a strictly $C$-open walk
    $\pi = \lp v_0 \sus \cdots \sus v_n \rp$ from a node
    $v_0 \in A$ to a node $v_n \in B \cup D$ of minimal length. By \Cref{lem:minimal-open-walk} (with $Y = B \cup D$)
    we have $v_n \in (B \cup D)\sm C$ and:
    \[ v_k \in (B \cup D)\sm C \implies k=n. \tag{$\ast$} \]
    Since $B \cap D = \emptyset$ we have $v_n \notin D$ or $v_n \notin B$.\\
    If $v_n \notin D$ then $v_n \in B \sm C$, and no node of $\pi$ lies in $D \sm C$, because
    $v_k \in D \sm C \ins (B \cup D)\sm C$ would force $k=n$ by $(\ast)$ and hence $v_n \in D$. So by
    \Cref{rem:d-open-walks} d) the walk $\pi$ is a strictly $(D \cup C)$-open walk from a node in $A$ to a node in
    $B$. This contradicts the assumption: $A \Perp^d B \given D \cup C$.\\
    If $v_n \notin B$ then $v_n \in D \sm C$, and in the same way no node of $\pi$ lies in $B \sm C$, so $\pi$ is
    a strictly $(B \cup C)$-open walk from a node in $A$ to a node in $D$. This contradicts the assumption:
    $A \Perp^d D \given B \cup C$.
\end{proof}
\end{Lem}

\begin{Lem}[More Redundancies]
    \label{sep:d:m-red}
    Let $E_1,E_2 \ins C$. Then:
    \[ A \Perp^d B \given C \quad\iff\quad (A \cup E_1) \Perp^d (B \cup E_2) \given C \quad\iff\quad
       (A \sm C) \Perp^d (B \sm C) \given C. \]
    In particular, the validity of $A \Perp^d B \given C$ only depends on the three sets $A \sm C$, $B \sm C$ and $C$.
\begin{proof}
    We first show for $E_2 \ins C$: $\; A \Perp^d B \given C \iff A \Perp^d B \cup E_2 \given C$.
    ``$\Longleftarrow$'' is Decomposition \ref{sep:d:dec}. For ``$\implies$'' note that $E_2 \ins C$ gives
    $E_2 \Perp^d A \given C$ by Redundancy \ref{sep:d:red}, hence $A \Perp^d E_2 \given C$ by Symmetry \ref{sep:d:sym},
    and then Composition \ref{sep:d:com} yields $A \Perp^d B \cup E_2 \given C$.
    Applying Symmetry \ref{sep:d:sym} twice this also gives, for $E_1 \ins C$:
    $\; A \Perp^d B \given C \iff (A \cup E_1) \Perp^d B \given C$. Combining the two gives the first equivalence.
    For the second one apply the first equivalence with $E_1 := A \cap C \ins C$ and $E_2 := B \cap C \ins C$ to the sets
    $A \sm C$ and $B \sm C$, and use $(A\sm C)\cup(A \cap C) = A$ and $(B \sm C)\cup (B \cap C) = B$.
\end{proof}
\end{Lem}

\subsection{Core Separoid Rules for id-Separation}

From now on no walk will be inspected anymore. All the rules of this subsection are derived from the five
\emph{core} symmetric rules Symmetry \ref{sep:d:sym}, Redundancy \ref{sep:d:red}, Decomposition
\ref{sep:d:dec}, Weak Union \ref{sep:d:uni} and Contraction \ref{sep:d:con} alone, i.e.\ from the (symmetric)
semi-graphoid axioms; the two extra rules Composition \ref{sep:d:com} and Intersection \ref{sep:d:int}, which are
special to d-separation, are only needed for the rules o)--s), see \Cref{sec:supp:further-id-rules} and the
derived rules thereafter.
The derivations are purely formal, using over and over again the definition:
\[ A \Perp^{\id} B \given C \qquad \iff \qquad A \Perp^{d} B \cup J \given C. \]
Note that ``left'' versions of the symmetric rules (dropping or moving sets on the left hand side) are available as
well, by an application of Symmetry \ref{sep:d:sym} before and after the respective rule.

\begin{Lem}[Extended Left Redundancy]
    \label{sep:sig:l-red}
    \[ A \ins C \implies  A \Perp^{\id}  B \given C.\]
\begin{proof}
    Unfolding the definition the claim reads: $A \Perp^d B \cup J \given C$. Since $A \ins C$ this is an instance of
    Redundancy \ref{sep:d:red}.
\end{proof}
\end{Lem}

\begin{Lem}[$J$-Restricted Right Redundancy]
    \label{sep:sig:r-red}
    \[A \Perp^{\id}  \emptyset \given C \cup J \qquad \text{always holds.}\]
\begin{proof}
    Unfolding the definition the claim reads: $A \Perp^d \emptyset \cup J \given C \cup J$, i.e.\
    $A \Perp^d J \given C \cup J$. Since $J \ins C \cup J$, Redundancy \ref{sep:d:red} gives
    $J \Perp^d A \given C \cup J$, and Symmetry \ref{sep:d:sym} turns this into $A \Perp^d J \given C \cup J$.
\end{proof}
\end{Lem}

\begin{Lem}[Left Decomposition]
    \label{sep:sig:l-dec}
    \[ A \cup D \Perp^{\id}  B \given C \implies D \Perp^{\id}  B \given C.\]
\begin{proof}
    Unfolding the definition this reads: $A \cup D \Perp^d B \cup J \given C \implies D \Perp^d B \cup J \given C$.
    Apply Symmetry \ref{sep:d:sym}, then Decomposition \ref{sep:d:dec} (dropping $A$), then Symmetry \ref{sep:d:sym}
    again.
\end{proof}
\end{Lem}

\begin{Lem}[Right Decomposition]
    \label{sep:sig:r-dec}
    \[ A \Perp^{\id}  B \cup D\given C \implies A \Perp^{\id}  D \given C.\]
\begin{proof}
    Unfolding the definition this reads:
    $A \Perp^d B \cup D \cup J \given C \implies A \Perp^d D \cup J \given C$.
    Since $B \cup D \cup J = (D \cup J) \cup B$ this is Decomposition \ref{sep:d:dec} (dropping $B$).
\end{proof}
\end{Lem}

\begin{Lem}[$J$-Inverted Right Decomposition]
    \label{sep:sig:inv-r-dec}
    \[A \Perp^{\id}  B \given C \implies A \Perp^{\id}  J \cup B \given C.\]
\begin{proof}
    Unfolding the definition both sides read $A \Perp^d B \cup J \given C$, because $(J \cup B) \cup J = B \cup J$.
    So the two statements are even equivalent.
\end{proof}
\end{Lem}

\begin{Lem}[Left Weak Union]
    \label{sep:sig:l-uni}
    \[A \cup D \Perp^{\id}  B  \given C \implies A \Perp^{\id}  B \given D\cup C.\]
\begin{proof}
    Unfolding the definition this reads:
    $A \cup D \Perp^d B \cup J \given C \implies A \Perp^d B \cup J \given D \cup C$.
    By Symmetry \ref{sep:d:sym} the assumption reads $B \cup J \Perp^d A \cup D \given C$, so Weak Union
    \ref{sep:d:uni} gives $B \cup J \Perp^d A \given D \cup C$, and Symmetry \ref{sep:d:sym} gives the claim.
\end{proof}
\end{Lem}

\begin{Lem}[Right Weak Union]
    \label{sep:sig:r-uni}
    \[A \Perp^{\id}  B \cup  D \given C \implies A \Perp^{\id}  B \given D \cup  C.\]
\begin{proof}
    Unfolding the definition this reads:
    $A \Perp^d (B \cup J) \cup D \given C \implies A \Perp^d B \cup J \given D \cup C$,
    which is Weak Union \ref{sep:d:uni}, applied with the sets $B \cup J$ and $D$.
\end{proof}
\end{Lem}

\begin{Lem}[Left Contraction]
    \label{sep:sig:l-con}
    \[ (A \Perp^{\id}  B \given D \cup  C) \land (D \Perp^{\id}  B \given C) \implies A \cup  D \Perp^{\id}  B \given C.\]
\begin{proof}
    Unfolding the definition and applying Symmetry \ref{sep:d:sym} the two assumptions read:
    \[ B \cup J \Perp^d A \given D \cup C \qquad \text{ and } \qquad B \cup J \Perp^d D \given C. \]
    Contraction \ref{sep:d:con} (with $B \cup J$ in the left slot) gives $B \cup J \Perp^d A \cup D \given C$, and
    Symmetry \ref{sep:d:sym} turns this into the claim $A \cup D \Perp^d B \cup J \given C$.
\end{proof}
\end{Lem}

\begin{Lem}[Right Contraction]
    \label{sep:sig:r-con}
    \[ (A \Perp^{\id}  B \given D \cup  C) \land (A \Perp^{\id}  D \given C)  \implies A \Perp^{\id}  B \cup  D \given C.\]
\begin{proof}
    Unfolding the definition the two assumptions read:
    \[ A \Perp^d B \cup J \given D \cup C \qquad \text{ and } \qquad A \Perp^d D \cup J \given C, \]
    and the claim reads $A \Perp^d B \cup D \cup J \given C$.
    Decomposition \ref{sep:d:dec} turns the second assumption into $A \Perp^d D \given C$. Now Contraction
    \ref{sep:d:con}, applied to $A \Perp^d B \cup J \given D \cup C$ and $A \Perp^d D \given C$, gives
    $A \Perp^d (B \cup J) \cup D \given C$, which is the claim, since $(B \cup J) \cup D = B \cup D \cup J$.
\end{proof}
\end{Lem}

\begin{Lem}[Right Cross Contraction]
    \label{sep:sig:rc-con}
    \[ (A \Perp^{\id}  B \given D \cup  C) \land (D \Perp^{\id}  A \given C) \implies A  \Perp^{\id}  B \cup  D \given C.\]
\begin{proof}
    Unfolded, the second assumption reads $D \Perp^d A \cup J \given C$, which by Decomposition \ref{sep:d:dec} and
    Symmetry \ref{sep:d:sym} gives $A \Perp^d D \given C$. Together with the first assumption
    $A \Perp^d B \cup J \given D \cup C$ we can now conclude exactly as in Right Contraction \ref{sep:sig:r-con}.
\end{proof}
\end{Lem}

\begin{Lem}[Flipped Left Cross Contraction]
    \label{sep:sig:flc-con}
    \[ (A \Perp^{\id}  B \given D \cup  C) \land (B \Perp^{\id}  D \given C) \implies B \Perp^{\id}  A \cup  D \given C.\]
\begin{proof}
    Unfolding the definition the two assumptions read:
    \[ \text{(i) } A \Perp^d B \cup J \given D \cup C, \qquad\qquad \text{(ii) } B \Perp^d D \cup J \given C, \]
    and the claim reads: $B \Perp^d A \cup D \cup J \given C$.
    From (i) we get $A \Perp^d B \given J \cup D \cup C$ by Weak Union \ref{sep:d:uni} (moving $J$, not $B$, into the
    conditioning set) and thus $B \Perp^d A \given (D \cup J) \cup C$ by Symmetry \ref{sep:d:sym}.
    Contraction \ref{sep:d:con}, applied to this and to (ii), i.e.\ to $B \Perp^d D \cup J \given C$, gives
    $B \Perp^d A \cup D \cup J \given C$, which is the claim.\\
    Note that the input nodes are carried along inside the \emph{set} $D \cup J$ of assumption (ii) and are never
    split off from $D$. In particular only Symmetry \ref{sep:d:sym}, Weak Union \ref{sep:d:uni} and
    Contraction \ref{sep:d:con} are used; neither Composition \ref{sep:d:com} nor Intersection \ref{sep:d:int}
    is needed here.
\end{proof}
\end{Lem}

\subsection{Further Separoid Rules for id-Separation}
\label{sec:supp:further-id-rules}

The rules of this subsection and the id-version of More Redundancies below are the ones that additionally use
Composition \ref{sep:d:com} and Intersection \ref{sep:d:int}. They are thus special to d-separation and have no
counterpart for general transitional conditional independence.

\begin{Lem}[Left Composition]
    \label{sep:sig:l-com}
    \[ (A \Perp^{\id}  B \given  C) \land (D \Perp^{\id}  B \given C) \implies A \cup D \Perp^{\id}  B \given C.\]
\begin{proof}
    Unfolding the definition and applying Symmetry \ref{sep:d:sym} the two assumptions read
    $B \cup J \Perp^d A \given C$ and $B \cup J \Perp^d D \given C$. Composition \ref{sep:d:com} gives
    $B \cup J \Perp^d A \cup D \given C$ and Symmetry \ref{sep:d:sym} the claim.
\end{proof}
\end{Lem}

\begin{Lem}[Right Composition]
    \label{sep:sig:r-com}
    \[ (A \Perp^{\id}  B \given  C) \land (A \Perp^{\id}  D \given C) \implies A \Perp^{\id}  B \cup  D \given C.\]
\begin{proof}
    Unfolding the definition the two assumptions read $A \Perp^d B \cup J \given C$ and $A \Perp^d D \cup J \given C$.
    Composition \ref{sep:d:com} gives $A \Perp^d (B \cup J) \cup (D \cup J) \given C$, which is the claim, since
    $(B \cup J) \cup (D \cup J) = B \cup D \cup J$.
\end{proof}
\end{Lem}

\begin{Lem}[Left Intersection]
    \label{sep:sig:l-int}
    Assume that $A \cap D = \emptyset$, then:
    \[ (A \Perp^{\id}  B \given D \cup  C) \land (D \Perp^{\id}  B \given A \cup C) \implies A \cup D \Perp^{\id}  B \given   C.\]
\begin{proof}
    Unfolding the definition and applying Symmetry \ref{sep:d:sym} the two assumptions read:
    \[ B \cup J \Perp^d A \given D \cup C \qquad \text{ and } \qquad B \cup J \Perp^d D \given A \cup C. \]
    Since $A \cap D = \emptyset$, Intersection \ref{sep:d:int} applies (with $B \cup J$ in the left slot and the two
    disjoint sets $A$ and $D$) and gives $B \cup J \Perp^d A \cup D \given C$. Symmetry \ref{sep:d:sym} gives the
    claim.
\end{proof}
\end{Lem}

\begin{Lem}[Right Intersection]
    \label{sep:sig:r-int}
    Assume that $B \cap D = \emptyset$, then:
    \[ (A \Perp^{\id}  B \given D \cup  C) \land (A \Perp^{\id}  D \given B \cup C) \implies A  \Perp^{\id}  B \cup D \given   C.\]
\begin{proof}
    Unfolding the definition the two assumptions read:
    \[ \text{(i) } A \Perp^d B \cup J \given D \cup C, \qquad\qquad \text{(ii) } A \Perp^d D \cup J \given B \cup C, \]
    and the claim reads: $A \Perp^d B \cup D \cup J \given C$.
    Note that we cannot apply Intersection \ref{sep:d:int} to the pair of sets $B \cup J$ and $D$ directly, since
    these two sets need not be disjoint: $D$ may contain input nodes. So we put:
    \[ B^\ast \; := \; (B \cup J) \sm D, \]
    and record, using $B \cap D = \emptyset$, the following four identities:
    \[ B^\ast = B \cup (J \sm D), \qquad B^\ast \cap D = \emptyset, \qquad B^\ast \cup D = B \cup D \cup J,
       \qquad B \cup J = B^\ast \cup (J \cap D). \]
    By the last identity, Decomposition \ref{sep:d:dec} turns (i) into: $A \Perp^d B^\ast \given D \cup C$.
    (No information is lost here: $J \cap D$ is contained in the conditioning set $D \cup C$, so by More Redundancies
    \ref{sep:d:m-red} the two statements are even equivalent.)
    Writing $D \cup J = D \cup (J \sm D)$, Weak Union \ref{sep:d:uni} turns (ii) into:
    $A \Perp^d D \given (J \sm D) \cup B \cup C$, and by the first identity we have
    $(J \sm D) \cup B \cup C = B^\ast \cup C$.
    Since $B^\ast \cap D = \emptyset$, Intersection \ref{sep:d:int} now applies to
    \[ A \Perp^d B^\ast \given D \cup C \qquad \text{ and } \qquad A \Perp^d D \given B^\ast \cup C \]
    and gives $A \Perp^d B^\ast \cup D \given C$, which by the third identity is the claim.
\end{proof}
\end{Lem}

\subsection{Derived Separoid Rules for id-Separation}

\begin{Lem}[Restricted Symmetry]
    \label{sep:sig:res-sym}
    \[(A \Perp^{\id}  B \given C) \land (B \Perp^{\id}  \emptyset \given C)  \implies B \Perp^{\id}  A \given C.\]
\begin{proof}
    Follows from Flipped Left Cross Contraction \ref{sep:sig:flc-con} with $D=\emptyset$.
\end{proof}
\end{Lem}

\begin{Lem}[$J$-Restricted Symmetry]
    \label{sep:sig:j-res-sym}
    \[A \Perp^{\id}  B \given C \cup J \implies B \Perp^{\id}  A \given C \cup J.\]
\begin{proof}
    Follows from Restricted Symmetry \ref{sep:sig:res-sym}, applied with $C \cup J$ in place of $C$, together with
    $J$-Restricted Right Redundancy \ref{sep:sig:r-red}, which provides $B \Perp^{\id} \emptyset \given C \cup J$.
\end{proof}
\end{Lem}

\begin{Lem}[Symmetry]
    \label{sep:sig:sym}
    If $J=\emptyset$ then we have:
    \[A \Perp^{\id}  B \given C \implies B \Perp^{\id} A \given C.\]
\begin{proof}
    Follows directly from $J$-Restricted Symmetry \ref{sep:sig:j-res-sym}.
\end{proof}
\end{Lem}

\begin{Lem}[More Redundancies]
    \label{sep:sig:m-red}
    \[A \Perp^{\id} B \given C \iff (A\sm C) \Perp^{\id} (B \sm C) \given C  \iff A \cup C \Perp^{\id} J \cup B \cup C \given C.\]
\begin{proof}
    Unfolding the definition the three statements read:
    \[ A \Perp^d B \cup J \given C, \qquad (A \sm C) \Perp^d (B\sm C) \cup J \given C, \qquad
       (A \cup C) \Perp^d J \cup B \cup C \given C, \]
    where for the third one we used $(J \cup B \cup C) \cup J = J \cup B \cup C$.
    By More Redundancies \ref{sep:d:m-red}, applied to the sets $A$ and $B \cup J$ (in the second case with
    $E_1:=E_2:=C$), the first and the third statement are both equivalent to:
    \[ (A \sm C) \Perp^d \lp (B \cup J) \sm C \rp \given C. \]
    By More Redundancies \ref{sep:d:m-red}, applied to the sets $A \sm C$ and $(B \sm C) \cup J$, so is the second
    statement, since $(A\sm C) \sm C = A \sm C$ and:
    \[ \lp (B\sm C) \cup J \rp \sm C \; = \;  (B \sm C) \cup (J \sm C) \; = \; (B \cup J) \sm C. \]
    So all three statements are equivalent.
\end{proof}
\end{Lem}
\section{Proofs - Global Markov Property}
\label{sec:supp:markov_property}

For the reader's convenience we restate \Cref{thm-gmp-mI-CBN} before proving it.

The proof of the global Markov property follows similar arguments as used in \cite{Lau90, Ver93, Ric03, FM17, FM18, RERS17},
namely chaining the separoid rules together in an inductive way.
The main difference here is that we never rely on the Symmetry property but instead use
the left and right versions of the separoid rules separately.

Throughout this section all transitional conditional independences are understood w.r.t.\ the transition probability space
$\lp \Xcal_V \times \Xcal_J, \Pk(X_V\Vert X_J) \rp$, i.e.\ $\Wcal = \Xcal_V$, $\Tcal = \Xcal_J$ and $\Tk = X_J$,
and for $S \ins J \dcup V$ we write $X_S$ for the \emph{deterministic} transitional random variable
$\deltabf(X_S|X_V,X_J)$ given by the coordinate projection $\Xcal_V \times \Xcal_J \to \Xcal_S$.
By \Cref{lem:separoid-compatibility-1} the relations $\ismapof_\Kk$ and $\Indep_\Kk$ are invariant under
isomorphisms of the codomains, so we may and will freely identify $X_S \otimes X_{S'}$ with $X_{S \cup S'}$ for
\emph{disjoint} $S$, $S'$ and reorder the factors of a product, all without further mention.
We will also drop \emph{repeated} factors. This is \emph{not} an isomorphism of codomains (in general
$\Xcal_S \times \Xcal_{S'} \not\cong \Xcal_{S \cup S'}$ when $S \cap S' \neq \emptyset$), so we record it
separately.

\begin{Lem}[Merging index sets]
    \label{lem:merge-index-sets}
    For all $S, S' \ins J \dcup V$ we have $X_S \otimes X_{S'} \approx_\Kk X_{S \cup S'}$, and consequently any
    occurrence of $X_S \otimes X_{S'}$ in any of the three arguments of $\Indep_\Kk$ may be replaced by
    $X_{S \cup S'}$, and conversely.
\begin{proof}
    All the $X_S$ are \emph{deterministic}, so $X_S \otimes X_S \approx_\Kk X_S$ by idempotency, see
    \Cref{lem:idempotency}, and $X_S \otimes X_{S'} \approx_\Kk X_{S \cup S'}$ follows with
    \Cref{lem:products-bound}. The second claim holds since $\Indep_\Kk$ is invariant under $\approx_\Kk$ in all
    three arguments, by Full Equivalent Exchange \ref{sep:tci:full-eq-ex}.
\end{proof}
\end{Lem}

\noindent
We will use \Cref{lem:merge-index-sets} silently from here on.

Furthermore, $\Gk=(J,V,E)$ will be a (finite) conditional directed acyclic graph (CDAG), i.e.\ every edge of $E$
is a directed edge $w \tuh u$ with $u \in V$, and $\Gk$ contains no directed cycles. In particular there are no
bi-directed edges and no latent confounders.
Recall from \Cref{def:d-separation} that for $A,B,C \ins J \dcup V$ (not necessarily disjoint) we write:
\[ A \Perp^{\id}_\Gk B \given C \qquad :\iff \qquad A \Perp^{d}_\Gk B \cup J \given C, \]
where the symmetric relation on the right holds if every walk in $\Gk$ from a node in $A$ to a node in $J \cup B$
is d-blocked by $C$.
Here a walk $\pi = \lp v_0 \sus \cdots \sus v_n \rp$ is \emph{d-blocked by $C$}, or \emph{$C$-blocked}, if either one of its
end nodes $v_0$, $v_n$ lies in $C$, or if $\pi$ contains a non-collider (a fork or a left/right chain) $v_k \in C$,
or if $\pi$ contains a collider $v_k \notin \Anc^\Gk(C)$.
By \Cref{lem:strictly-open} we may throughout replace $\Anc^\Gk(C)$ by $C$ here, i.e.\ work with \emph{strictly
$C$-open} walks: $\pi$ is strictly $C$-open if and only if $v_0,v_n \notin C$, all inner non-colliders of $\pi$ lie
outside of $C$ and all colliders of $\pi$ lie \emph{in $C$}; and $A \Perp^d_\Gk B \given C$ holds if and only if no
strictly $C$-open walk runs from a node in $A$ to a node in $B$. Every strictly $C$-open walk is in particular
$C$-open, i.e.\ not $C$-blocked.
Note that, per definition, we always have the equivalence:
\[ A \Perp^{\id}_\Gk B \given C \qquad \iff \qquad A \Perp^{\id}_\Gk J \cup B \given C. \]
We will also use \Cref{lem:remove-childless} through its id-version, see \Cref{rem:remove-childless-id}.

Left Decomposition \ref{sep:tci:l-dec} is stated so that it drops the \emph{first} factor of a product; whenever we
use it below to drop the \emph{second} one, this is legitimate by the commutativity of $\otimes$ up to isomorphism
and the $\cong$-invariance of $\Indep_\Kk$, see \Cref{lem:separoid-compatibility-1} points 2.\ and 6.; we will not
mention this again. The same applies to Left Decomposition \ref{sep:sig:l-dec} for id-separation, where $\cup$ is
commutative on the nose.

For a childless node $v$ we will write $\Gk_{-v}$ for the induced subgraph of $\Gk$ on $(J,V \sm \{v\})$ and repeatedly use
\Cref{lem:remove-childless}, which states that removing a childless node changes neither d- nor id-separation statements
among the remaining nodes.

\begin{Thm}[Global Markov property for Bayesian networks with input nodes]
\label{GMP-CCBN}
Consider a Bayesian network $\Mk$ with input nodes, see \Cref{def:cbn}, with CDAG $\Gk=(J,V,E)$
and joint Markov kernel $\Pk(X_V\Vert X_J)$.
Then for all $A, B, C \ins J \dcup V$ (not-necessarily disjoint) we have the implication:
\[ A \Perp^{\id}_\Gk B \given C \qquad \implies \qquad X_A \Indep_{\Pk(X_V\Vert X_J)} X_B \given X_C.   \]
If one wants to make the implicit dependence on $J$ more explicit one can equivalently also write:
\[ A \Perp^{\id}_\Gk J \cup B \given C \qquad \implies \qquad X_A \Indep_{\Pk(X_V\Vert X_J)} X_J,X_{B} \given X_C.   \]
Indeed, the two left hand sides are equivalent by $J$-Inverted Right Decomposition \ref{sep:sig:inv-r-dec} together with
Right Decomposition \ref{sep:sig:r-dec}, and the two right hand sides are equivalent by \Cref{add-T}, using
$\Tk = X_J$. We will prove the first form.
\begin{proof}
We do induction by $\#V$. \\

0.) Induction start:  $V= \emptyset$. This means that $A,B,C \ins J$.
The assumption:
\[A \Perp^{\id}_\Gk J \cup B \given C,\]
implies that we must have that $A \ins C$. Otherwise a trivial walk from a node in $A \ins J$ to the same node in $J \cup B$
would be $C$-open.
Since $A \ins C$ the coordinate projection $\pr:\, \Xcal_C \to \Xcal_A$ satisfies $X_A = \pr \circ X_C$
pointwise, i.e.\ $X_A \ismapof X_C$, and therefore $X_A \ismapof_\Kk X_C$ by \Cref{lem:exact-to-as}.
Extended Left Redundancy \ref{sep:tci:ext-l-red} then already gives:
\[ X_A \Indep_{\Pk(X_V\Vert X_J)} X_B \given X_C.\]

(IND): Induction assumption:
The global Markov property holds for all Bayesian networks with input nodes
with $\#V < n$ (and arbitrary $J$).\\

1.) Now assume: $\#V =n>0$ and $A \Perp^{\id}_\Gk J \cup B \given C$.

Since $\Gk$ is acyclic it has a topological order, and since $\Pa^\Gk(j)=\emptyset$ for every input node $j \in J$,
see \Cref{not-cdag}, we can choose a topological order $<$ of $J \dcup V$ in which the elements of $J$ are ordered first.
Let $v \in V$ be its last element. Since $<$ is a topological order and $v$ is the maximal element, no node can be a child
of $v$, i.e.\ $\Ch^\Gk(v) = \emptyset$, so $v$ is childless. \\
Recall from \Cref{not-cdag} that $\Pred^\Gk_<(u) := \lC w \in J \dcup V \,|\, w < u \rC$ denotes the set of
\emph{predecessors} of $u$ w.r.t.\ $<$. For our maximal element $v$ we thus have:
\[ \Pred^\Gk_<(v) = (J \dcup V) \sm \{v\}. \]

Since $\Ch^\Gk(v)=\emptyset$ no node $w \in V \sm \{v\}$ has $v$ as a parent, i.e.:
\[ \Pa^{\Gk_{-v}}(w) = \Pa^{\Gk}(w) \qquad \text{for all } w \in V \sm \{v\}. \]
As an induced subgraph of an acyclic graph, $\Gk_{-v}=(J,V\sm\{v\},E_{-v})$ is again a CDAG, and
by the above the restricted family of Markov kernels:
\[ \Mk_{-v} := \lp \Gk_{-v}, \lp \Pk_w\lp X_w\,\Vert\, X_{\Pa^{\Gk}(w)}\rp \rp_{w \in V \sm \{v\}} \rp \]
is again a Bayesian network with input nodes $J$ and now with $\#(V \sm \{v\}) = n-1 <n$ output nodes.
Its joint Markov kernel is:
\[ \Pk(X_{V \sm \{v\}}\Vert X_J) = \bigotimes^>_{w \in V \sm \{v\}} \Pk_w\lp X_w\,\Vert\, X_{\Pa^{\Gk}(w)}\rp. \]
Furthermore, since $v$ is the last element of the topological order $<$, we have the factorization:
\begin{align*} \Pk(X_V\Vert X_J) &= \Pk_v\lp X_v\,\Vert\, X_{\Pa^\Gk(v)}\rp \otimes
    \underbrace{\bigotimes^>_{w \in V \sm \{v\}} \Pk_w\lp X_w\,\Vert\, X_{\Pa^\Gk(w)}\rp}_{=\Pk(X_{V \sm \{v\}}\Vert X_J)}.
\end{align*}
Here we used that by \Cref{def:joint-markov-kernel-bn} and \Cref{rem:markov-kernels-products} the joint Markov kernel
does not depend on the chosen topological order, so we may compute both joint kernels w.r.t.\ our $<$.
Since $\Pk_v$ is a probability (and not merely a sub-probability) kernel, marginalizing out the leftmost factor gives
that $\Pk(X_{V\sm\{v\}}\Vert X_J)$ is the marginal of $\Pk(X_V\Vert X_J)$ onto $\Xcal_{V \sm \{v\}}$.
Since transitional conditional independence only depends on the joint Markov kernel of the involved (transitional) random
variables, see \Cref{def:transitional_conditional_independence}, we get for all
$\tilde A,\tilde B,\tilde C \ins (J \dcup V) \sm \{v\}$ the equivalence:
\[ X_{\tilde A} \Indep_{\Pk(X_{V \sm \{v\}}\Vert X_J)} X_{\tilde B} \given X_{\tilde C}
\qquad \iff \qquad
X_{\tilde A} \Indep_{\Pk(X_{V}\Vert X_J)} X_{\tilde B} \given X_{\tilde C}. \]
So we may and will apply the induction assumption (IND) to $\Mk_{-v}$ and state its conclusions directly in terms of
$\Pk(X_V\Vert X_J)$.\\
The above factorization also gives us the conditional independence:
\begin{align}
    X_v \Indep_{\Pk(X_V\Vert X_J) } X_{\Pred^\Gk_<(v)} \given X_D, \label{eqn:pred} \tag{$\dagger$}
\end{align}
where we put $D:=\Pa^\Gk(v) \ins \Pred^\Gk_<(v)$. Note that $v \notin D$, as $\Gk$ has no self-loops.
Indeed, put $\Qk(X_v|X_D) := \Pk_v\lp X_v\,\Vert\, X_{\Pa^\Gk(v)}\rp$, which by \Cref{def:cbn} is a Markov
kernel $\Xcal_D \dshto \Xcal_v$; in particular it depends on $T=X_J$ only through the coordinates $X_{D \cap J}$, as
required by \Cref{def:transitional_conditional_independence}. Since $X_{\Pred^\Gk_<(v)}$ and $X_D$ are coordinate
projections of $(X_{V \sm \{v\}},X_J)$, pushing the displayed factorization forward along them yields:
\[ \Kk\lp X_v,X_{\Pred^\Gk_<(v)}\big|T \rp = \Qk(X_v|X_D) \otimes \Kk\lp X_{\Pred^\Gk_<(v)}\big|T\rp, \]
which is exactly \eqref{eqn:pred}.\\

In the following we will distinguish between 4 cases:
\begin{enumerate}[label = \Alph*.)]
    \item $v \in A \sm C$,
    \item $v \in B \sm C$,
    \item $v \in C$,
    \item $ v \notin A \cup J \cup B \cup C$,
\end{enumerate}
Note that $v \in V$, thus $v \notin J$, which shows that the above cover all possible cases.\\
Further note that:
\[A \Perp^{\id}_\Gk J \cup B \given C,\]
implies that:
\[ A \cap ( J \cup B) \ins C.  \]
Otherwise a trivial walk from $A$ to $J \cup B$ would be $C$-open.
This shows that $A\sm C$, $(J \cup B) \sm C$ and $C$ are pairwise disjoint.
It also shows $A \cap J \ins C$ and thus:
\begin{align}
    A \sm C \ins V. \label{eqn:ainv} \tag{$\ddagger$}
\end{align}
So all nodes of $A \sm C$ are output nodes and their measurable spaces are standard. This will be used in case C.\ to
justify the use of Left Weak Union \ref{sep:tci:l-uni}, which requires a disintegration triple.\\

Case D.):  $v \notin A \cup J \cup B \cup C$. Then we can remove the childless node $v$ and use the equivalence of
\Cref{lem:remove-childless}:
\[A \Perp^{\id}_\Gk J \cup B \given C\quad \iff \quad A \Perp^{\id}_{\Gk_{-v}} J \cup B \given C.\]
With $\#(V \sm \{v\})<n$ and induction (IND) applied to $\Mk_{-v}$ we then get:
\[X_A \Indep_{\Pk(X_V\Vert X_J)} X_B \given X_C.\]
This shows the claim in case D.\\

Case A.):  $v \in A\sm C$. Then we can write:
\begin{align*}
    A &= A' \dcup (A \cap C) \dcup \{v\}, \\
    B &= B' \dcup (B \cap C),
\end{align*}
with some disjoint $A' \ins A \sm C$ and $B' \ins B \sm C$.
Note that $v \notin A' \cup J \cup B' \cup C$, since $A'$ and $\{v\}$ are disjoint, $v \in V$, $v \notin C$ and
$v \in A \sm C$ is disjoint from $(J \cup B) \sm C \ni$ every element of $B'$.
We then have the implications:
\begin{align*}
    A \Perp^{\id}_\Gk  J \cup B \given C & \xRightarrow{\text{Right Decomposition \ref{sep:sig:r-dec}}}  & A \Perp^{\id}_\Gk  J \cup B' \given C \\
    & \xRightarrow{\text{Left Decomposition \ref{sep:sig:l-dec}}}  & A' \Perp^{\id}_\Gk  J \cup B' \given C \\
    & \xRightarrow{\text{\Cref{lem:remove-childless}, } v \notin A'\cup J \cup B' \cup C}  & A' \Perp^{\id}_{\Gk_{-v}} J \cup B' \given C \\
 & \xRightarrow{\text{induction (IND)}}  & X_{A'} \Indep_{\Pk(X_V\Vert X_J)}  X_{B'} \given X_C. \label{eqn:bli} \tag{\#1}
\end{align*}
On the other hand we have with $D = \Pa^\Gk(v)$:
\begin{align*}
    A \Perp^{\id}_\Gk  J \cup B \given C & \xRightarrow{\text{Right Decomposition \ref{sep:sig:r-dec}, }B' \ins B}  & A \Perp^{\id}_\Gk  J \cup  B' \given C   \\
  & \xRightarrow{\substack{\text{Left Weak Union \ref{sep:sig:l-uni} and Left Decomposition \ref{sep:sig:l-dec},}\\ A= A' \dcup (A \cap C) \dcup \{v\}}}  & \{v\} \Perp^{\id}_\Gk  J \cup  B' \given A' \dcup C.    \\
 & \xRightarrow{(*), \text{ see below}}  & D  \Perp^{\id}_\Gk  J \cup  B' \given A' \dcup C \\
 & \xRightarrow{\text{\Cref{lem:remove-childless}, } v \notin D \cup J \cup B' \cup A' \cup C}
 & D \Perp^{\id}_{\Gk_{-v}} J \cup B' \given A' \dcup C \\
& \xRightarrow{\text{induction (IND)}}  &   X_D  \Indep_{\Pk(X_V\Vert X_J)}  X_{B'} \given X_{A' \dcup C}\\
& \xRightarrow{X_{A'\dcup C} = (X_{A'},X_C)}  &   X_D  \Indep_{\Pk(X_V\Vert X_J)}  X_{B'} \given X_{A'}, X_C. \label{eqn:bla} \tag{\#2}
\end{align*}
$(*)$ holds by contraposition. So assume $D \nPerp^{\id}_\Gk J \cup B' \given A' \dcup C$. Then, by
\Cref{lem:strictly-open}, there is a strictly $(A' \dcup C)$-open walk:
\[ \pi:\quad D \ni w = w_0 \sus \cdots \sus w_l \in J \cup B'.  \]
Since $w$ is an end node of the strictly $(A' \dcup C)$-open walk $\pi$ we have $w \notin A' \dcup C$.
Because $w \in D = \Pa^\Gk(v)$ the edge $v \hut w$ exists in $\Gk$ and we can prepend it to $\pi$:
\[ \tilde \pi:\quad v \hut w_0 \sus \cdots \sus w_l \in J \cup B'.  \]
The edge $w \tuh v$ has its tail at $w$. So $w$ is a non-collider of $\tilde \pi$ (an end node of $\tilde \pi$ in case $l=0$),
and we already know $w \notin A' \dcup C$.
All other nodes of $\tilde \pi$ are nodes of $\pi$ and keep their collider/non-collider status, and the other end node
$w_l \in J \cup B'$ is unchanged.
Finally, the new end node $v$ satisfies $v \notin A' \dcup C$, since $v \notin C$ (case A.) and $v \notin A'$.
Hence $\tilde \pi$ is a strictly $(A' \dcup C)$-open walk from $v$ to $J \cup B'$, which by
\Cref{lem:strictly-open} contradicts:
\[\{v\} \Perp^{\id}_\Gk J \cup B' \given A' \dcup C.\]
This shows $(*)$.\\

As discussed above we also already have the conditional independence \eqref{eqn:pred}:
\[ X_v \Indep_{\Pk(X_V\Vert X_J)}  X_{\Pred^\Gk_<(v)} \given X_D. \]
With this and $A' \dcup B' \dcup C \ins (J \dcup V) \sm \{v\} = \Pred^\Gk_<(v)$, where $A'$, $B'$ and $C$ are
pairwise disjoint, so that $\lp X_{A'},X_{B'},X_C\rp$ really is a marginal of $X_{\Pred^\Gk_<(v)}$ and Right
Decomposition is applicable, we get the implications:
\begin{align*}
&&X_v \Indep_{\Pk(X_V\Vert X_J)}  X_{\Pred^\Gk_<(v)} \given X_D\\
&\xRightarrow{\text{Right Decomposition \ref{sep:tci:r-dec}}} &
X_v \Indep_{\Pk(X_V\Vert X_J)}  X_{A'}, X_{B'}, X_C \given X_D \\
& \xRightarrow{\text{Right Weak Union \ref{sep:tci:r-uni}}}  &
X_v \Indep_{\Pk(X_V\Vert X_J)} X_{B'} \given X_{A'}, X_C, X_D \\
& \xRightarrow{\text{Left Contraction \ref{sep:tci:l-con}, }\eqref{eqn:bla}}  & X_v, X_D  \Indep_{\Pk(X_V\Vert X_J)}  X_{B'} \given X_{A'}, X_C \\
& \xRightarrow{\text{Left Decomposition \ref{sep:tci:l-dec}}}  & X_v \Indep_{\Pk(X_V\Vert X_J)}  X_{B'} \given X_{A'}, X_C\\
& \xRightarrow{\text{Left Contraction \ref{sep:tci:l-con}, } \eqref{eqn:bli}}  &  X_{A'},X_v \Indep_{\Pk(X_V\Vert X_J)}  X_{B'} \given X_C\\
& \xRightarrow{X_J\text{-Ext. Inv. Right Decomposition \ref{sep:tci:ext-inv-r-dec}, }X_C \ismapof_\Kk X_J\otimes X_{B'} \otimes X_C}  & X_{A'}, X_v \Indep_{\Pk(X_V\Vert X_J)}  X_J,X_{B'},X_C \given X_C\\
& \xRightarrow{\text{Right Decomposition \ref{sep:tci:r-dec}, } B \ins B' \dcup C}  & X_{A'}, X_v \Indep_{\Pk(X_V\Vert X_J)}  X_B \given X_C.
 \label{eqn:bloop} \tag{\#3}
\end{align*}
By Extended Left Redundancy \ref{sep:tci:ext-l-red} we have:
\[ X_{A'},X_v, X_C \Indep_{\Pk(X_V\Vert X_J)}  X_B \given X_{A'},X_v,X_C. \]
With this we get the implications:
\begin{align*}
&&X_{A'},X_v, X_C \Indep_{\Pk(X_V\Vert X_J)}  X_B \given X_{A'},X_v,X_C\\
& \xRightarrow{\text{Left Contraction \ref{sep:tci:l-con}, } \eqref{eqn:bloop}}  &
X_{A'},X_v,X_{A'},X_v, X_C \Indep_{\Pk(X_V\Vert X_J)}  X_B \given X_C\\
& \xRightarrow{\text{Left Decomposition \ref{sep:tci:l-dec}, } A \ins A' \dcup \{v\} \dcup C}  &
X_A \Indep_{\Pk(X_V\Vert X_J)}  X_B \given X_C.
\end{align*}
This shows the claim in case A.\\

Case B.):  $v \in B\sm C$. Then we can write:
\begin{align*}
    A &= A' \dcup (A \cap C), \\
    B &= B' \dcup (B \cap C) \dcup \{v\},
\end{align*}
with some disjoint $A' \ins A \sm C$ and $B' \ins B \sm C$.
Note that $v \notin A' \cup J \cup B' \cup C \cup D$, since $v \in B \sm C$ is disjoint from $A \sm C \supseteq A'$,
$v \in V$, $B'$ and $\{v\}$ are disjoint, $v \notin C$ and $v \notin D=\Pa^\Gk(v)$.
We then have the implications:
\begin{align*}
A  \Perp^{\id}_\Gk  J \cup B \given C
& \xRightarrow{\text{Left Decomposition \ref{sep:sig:l-dec}}}  & A' \Perp^{\id}_\Gk  J \cup B \given C \\
& \xRightarrow{\text{Right Decomposition \ref{sep:sig:r-dec}}}  & A' \Perp^{\id}_\Gk  J \cup B' \given C \\
 & \xRightarrow{\text{\Cref{lem:remove-childless}, } v \notin A'\cup J \cup B' \cup C}  & A' \Perp^{\id}_{\Gk_{-v}} J \cup B' \given C \\
 & \xRightarrow{\text{induction (IND)}}  & X_{A'} \Indep_{\Pk(X_V\Vert X_J)}  X_{B'} \given X_C. \label{eqn:blione} \tag{\#1'}
\end{align*}
Again with $D = \Pa^\Gk(v)$ we get:
\begin{align*}
    A  \Perp^{\id}_\Gk  J \cup B \given C  & \xRightarrow{\text{Left Decomposition \ref{sep:sig:l-dec}}}  & A' \Perp^{\id}_\Gk  J \cup B \given C \\
 & \xRightarrow{\text{Right Decomposition \ref{sep:sig:r-dec}}}  & A' \Perp^{\id}_\Gk  J \cup B' \cup \{v\} \given C \\
& \xRightarrow{\text{Right Weak Union \ref{sep:sig:r-uni}}}  & A' \Perp^{\id}_\Gk  J \cup \{v\} \given B' \dcup C  \\
 & \xRightarrow{(\bullet), \text{ see below}}  & A'  \Perp^{\id}_\Gk  J \cup D \given B' \dcup C \\
 & \xRightarrow{\text{\Cref{lem:remove-childless}, } v \notin A' \cup J \cup D \cup B' \cup C}  & A'  \Perp^{\id}_{\Gk_{-v}}  J \cup D \given B' \dcup C \\
& \xRightarrow{\text{induction (IND)}}  &   X_{A'}  \Indep_{\Pk(X_V\Vert X_J)}  X_D \given X_{B' \dcup C}\\
& \xRightarrow{X_{B' \dcup C} = (X_{B'},X_C)}  &   X_{A'}  \Indep_{\Pk(X_V\Vert X_J)}  X_D \given X_{B'}, X_C. \label{eqn:blaone} \tag{\#2'}
\end{align*}
$(\bullet)$ holds by contraposition. So assume $A' \nPerp^{\id}_\Gk J \cup D \given B' \dcup C$. Then, by
\Cref{lem:strictly-open}, there is a strictly $(B' \dcup C)$-open walk:
\[ \pi:\quad A' \ni w_0 \sus \cdots \sus w_l = w \in J \cup D.  \]
If $w \in J$ then $\pi$ is already a strictly $(B' \dcup C)$-open walk from $A'$ to $J \cup \{v\}$.
Otherwise $w \in D = \Pa^\Gk(v)$, the edge $w \tuh v$ exists in $\Gk$ and we can append it to $\pi$:
\[ \tilde \pi:\quad A' \ni w_0 \sus \cdots \sus w_l \tuh v.  \]
The edge $w_l \tuh v$ has its tail at $w_l$. So $w_l$ is a non-collider of $\tilde \pi$ (an end node of $\tilde \pi$ in case $l=0$),
and $w_l \notin B' \dcup C$, since $w_l$ was an end node of the strictly $(B' \dcup C)$-open walk $\pi$.
All other nodes of $\tilde \pi$ are nodes of $\pi$ and keep their collider/non-collider status, and the other end node
$w_0 \in A'$ is unchanged.
Finally, the new end node $v$ satisfies $v \notin B' \dcup C$, since $v \notin C$ (case B.) and $v \notin B'$.
So in both cases we found a strictly $(B' \dcup C)$-open walk from $A'$ to $J \cup \{v\}$, which by
\Cref{lem:strictly-open} contradicts:
\[A' \Perp^{\id}_\Gk J \cup \{v\} \given B' \dcup C.\]
This shows $(\bullet)$.\\

As before we will use the conditional independence \eqref{eqn:pred}:
\[ X_v \Indep_{\Pk(X_V\Vert X_J)}  X_{\Pred^\Gk_<(v)} \given X_D. \]
With this and $A' \dcup B' \dcup C \ins (J \dcup V) \sm \{v\} = \Pred^\Gk_<(v)$, where again $A'$, $B'$ and $C$ are
pairwise disjoint, so that $\lp X_{A'},X_{B'},X_C\rp$ is a marginal of $X_{\Pred^\Gk_<(v)}$, we get the implications:
\begin{align*}
   && X_v \Indep_{\Pk(X_V\Vert X_J)} X_{\Pred^\Gk_<(v)} \given X_D \\
    &\xRightarrow{\text{Right Decomposition \ref{sep:tci:r-dec}}} & X_v \Indep_{\Pk(X_V\Vert X_J)} X_{A'},X_{B'},X_C \given X_D\\
& \xRightarrow{\text{Right Weak Union \ref{sep:tci:r-uni}}}  &  X_v \Indep_{\Pk(X_V\Vert X_J)} X_{A'} \given X_{B'},X_C,X_D\\
& \xRightarrow{\text{Flipped Left Cross Contraction \ref{sep:tci:flc-con}, }\eqref{eqn:blaone}}
& X_{A'}  \Indep_{\Pk(X_V\Vert X_J)} X_D,X_v \given X_{B'},X_C \\
& \xRightarrow{\text{Right Decomposition \ref{sep:tci:r-dec}}}
& X_{A'}  \Indep_{\Pk(X_V\Vert X_J)} X_v \given X_{B'},X_C \\
& \xRightarrow{\text{Right Contraction \ref{sep:tci:r-con}, }\eqref{eqn:blione}}
& X_{A'}  \Indep_{\Pk(X_V\Vert X_J)} X_{B'},X_v \given X_C \\
& \xRightarrow{\text{$X_J$-Ext. Inv. Right Dec. \ref{sep:tci:ext-inv-r-dec}, }X_C \ismapof_\Kk X_J\otimes X_{B'}\otimes X_v \otimes X_C}  & X_{A'}  \Indep_{\Pk(X_V\Vert X_J)} X_J,X_{B'},X_v,X_C \given X_C \\
& \xRightarrow{\text{Right Decomposition \ref{sep:tci:r-dec}, } B \ins B'\dcup \{v\} \dcup C}  & X_{A'}  \Indep_{\Pk(X_V\Vert X_J)} X_B \given X_C.  \label{eqn:bloopone} \tag{\#3'}
\end{align*}
By Extended Left Redundancy \ref{sep:tci:ext-l-red} we have:
\[ X_{A'}, X_C \Indep_{\Pk(X_V\Vert X_J)}  X_B \given X_{A'},X_C. \]
With this we get the implications:
\begin{align*}
&&X_{A'}, X_C \Indep_{\Pk(X_V\Vert X_J)}  X_B \given X_{A'},X_C\\
& \xRightarrow{\text{Left Contraction \ref{sep:tci:l-con}, } \eqref{eqn:bloopone}}  &
X_{A'},X_{A'}, X_C \Indep_{\Pk(X_V\Vert X_J)}  X_B \given X_C\\
& \xRightarrow{\text{Left Decomposition \ref{sep:tci:l-dec}, } A \ins A' \dcup C}  &
X_A \Indep_{\Pk(X_V\Vert X_J)}  X_B \given X_C.
\end{align*}
This shows the claim in case B.\\

Case C.): $v \in C$.
This case will be \emph{reduced} to Cases A.\ and B., applied to a \emph{different} triple of node sets, but
with the same graph $\Gk$, the same node $v$, the same factorization \eqref{eqn:pred} and the same induction
assumption (IND). It is therefore not a circular argument; see the note after \eqref{eqn:caseC-disjunction} below.
 Then we can write:
\begin{align*}
    A &= A' \dcup (A \cap C), \\
    B &= B' \dcup (B \cap C),\\
    C &= C' \dcup \{v\},
\end{align*}
with some pairwise disjoint $A' \ins A \sm C$,
$B' \ins B \sm C$ and $C' \ins C$.
In particular $v \notin A'$, $v \notin B'$ and $v \notin C'$.

We then get the implications:
\begin{align*}
A  \Perp^{\id}_\Gk  J \cup B \given C
& \xRightarrow{\text{Left Decomposition \ref{sep:sig:l-dec}}}  & A' \Perp^{\id}_\Gk  J \cup B \given C \\
& \xRightarrow{\text{Right Decomposition \ref{sep:sig:r-dec}}}  & A' \Perp^{\id}_\Gk  J \cup B' \given C \\
& \xRightarrow{C = C' \dcup \{v\}}  & A' \Perp^{\id}_\Gk  J \cup B' \given C'\dcup\{v\}.
\end{align*}

We now claim that:
\[A' \Perp^{\id}_\Gk  J \cup B' \given C'\dcup\{v\}\]
implies that one of the following statements holds:
\[  A' \dcup \{v\} \Perp^{\id}_\Gk  J \cup B' \given C' \qquad \lor \qquad A' \Perp^{\id}_\Gk  J \cup (B' \dcup \{v\}) \given C'. \]
Assume the contrary:
\[  A' \dcup \{v\} \nPerp^{\id}_\Gk  J \cup B' \given C' \qquad \land \qquad A' \nPerp^{\id}_\Gk  J \cup (B' \dcup \{v\}) \given C'. \]
So, by \Cref{lem:strictly-open}, there exist strictly $C'$-open walks $\pi_1$ and $\pi_2$ in $\Gk$:
\[ \pi_1:\quad  A' \cup \{v\}\ni u_0 \sus \cdots \sus u_k \in J \cup B',  \]
and:
\[ \pi_2:\quad  A' \ni w_0 \sus \cdots \sus w_m \in J \cup (B' \dcup \{v\}).  \]
So all colliders of $\pi_1$ and $\pi_2$ lie in $C'$, all their inner non-colliders lie outside of $C'$ and all their end
nodes lie outside of $C'$.\\
Since $v$ is childless, every edge of $\Gk$ incident to $v$ has its arrowhead at $v$.
So if $v$ occurred as an \emph{inner} node of a walk $\pi$, then both edges of $\pi$ adjacent to that occurrence would point into
$v$, i.e.\ $v$ would be a collider of $\pi$ there. As $v \notin C'$ such a walk $\pi$ would not be strictly
$C'$-open.
Since $\pi_1$ and $\pi_2$ are strictly $C'$-open, the node $v$ can therefore occur at most as an end node of $\pi_1$
and $\pi_2$.\\
Then note that $v \notin A'$ and $v \notin J \cup B'$, thus: $u_k \neq v$ and $w_0 \neq v$.\\
If now $\pi_i$ does not contain $v$ at all, then $\pi_i$ would also be strictly $(C'\dcup \{v\})$-open, since strict
openness w.r.t.\ $C'$ and w.r.t.\ $C' \dcup \{v\}$ can only differ at occurrences of $v$. For $i=1$ we would then
have $u_0 \in A'$ and thus a strictly $(C'\dcup\{v\})$-open walk from $A'$ to $J \cup B'$; for $i=2$ we would have
$w_m \in J \cup B'$ and again such a walk. Both contradict, via \Cref{lem:strictly-open}, the assumption:
\[A' \Perp^{\id}_\Gk  J \cup B' \given C'\dcup\{v\}.\]
So we can assume that the other end nodes equal $v$, i.e.: $u_0=v$ and $w_m =v$.\\
Furthermore, both $\pi_1$ and $\pi_2$ are non-trivial walks, i.e.\ $k,m \ge 1$, since $u_0 =v \neq u_k$ and $w_0 \neq v= w_m$.
Since $v$ is childless the first edge of $\pi_1$ and the last edge of $\pi_2$ point into $v$, so the $\pi_i$ are of the forms:
\[ \pi_1:\quad   v \hut u_1 \sus \cdots \sus u_k,  \]
and:
\[ \pi_2:\quad  w_0 \sus \cdots \sus w_{m-1} \tuh v,  \]
with $u_1, w_{m-1} \in D=\Pa^\Gk(v)$. Then consider the concatenated walk:
\[ \pi:\quad A' \ni  w_0 \sus \cdots \sus w_{m-1} \tuh v \hut u_1 \sus \cdots \sus u_k \in J \cup B'.  \]
Its end nodes are $w_0 \in A'$ and $u_k \in J \cup B'$, which lie outside of $C'$ (as end nodes of strictly
$C'$-open walks) and
differ from $v$, hence lie outside of $C' \dcup \{v\}$.
The single occurrence of $v$ in $\pi$ is a collider and $v \in C' \dcup \{v\}$, so it does not block $\pi$.
Every other node of $\pi$ is a node of $\pi_1$ or $\pi_2$ with the same two adjacent edges, hence with unchanged
collider/non-collider status, and it differs from $v$ (as $v$ occurs in $\pi_1$, $\pi_2$ only as the end nodes $u_0$, $w_m$).
So all remaining colliders of $\pi$ lie in $C' \ins C' \dcup \{v\}$ and all remaining non-colliders lie outside of
$C' \dcup \{v\}$.
Therefore $\pi$ is a strictly $(C' \dcup \{v\})$-open walk from $A'$ to $J \cup B'$, in contradiction, via
\Cref{lem:strictly-open}, to:
\[A' \Perp^{\id}_\Gk  J \cup B' \given C'\dcup\{v\}.\]
So the claim:
\begin{align}
    \lp A' \dcup \{v\} \Perp^{\id}_\Gk  J \cup B' \given C' \rp \qquad \lor \qquad
    \lp A' \Perp^{\id}_\Gk  J \cup (B' \dcup \{v\}) \given C' \rp \label{eqn:caseC-disjunction} \tag{$\ast$}
\end{align}
must be true.
Since $v \in (A' \dcup \{v\}) \sm C'$ in the first case and $v \in (B' \dcup \{v\}) \sm C'$ in the second case,
we have reduced case C.\ to case A.\ (applied to the triple $A' \dcup \{v\}$, $B'$, $C'$) or to case B.\
(applied to the triple $A'$, $B' \dcup \{v\}$, $C'$). Note that cases A.\ and B.\ were shown for arbitrary triples satisfying the hypothesis of the theorem, using only the
induction assumption (IND), while $\Gk$, $v$, the factorization \eqref{eqn:pred} and (IND) itself are unchanged; in
particular the consequences derived above, pairwise disjointness and \eqref{eqn:ainv}, are re-derived for the new
triple. So this is not circular. They imply:
\[ X_{A'},X_v \Indep_{\Pk(X_V\Vert X_J)}  X_{B'} \given X_{C'} \qquad \lor \qquad X_{A'} \Indep_{\Pk(X_V\Vert X_J)}  X_{B'},X_v \given X_{C'}.\]
By \eqref{eqn:ainv} we have $A' \ins A \sm C \ins V$, so $\Xcal_{A'}=\prod_{w \in A'}\Xcal_w$ is a finite product of
standard measurable spaces and hence itself standard, see \cite{Fremlin} 424B (for $A'=\emptyset$ we have
$\Xcal_{A'}=\Asterisk$, which is standard as well). Since $v \in V$, the space $\Xcal_v$ is standard and thus countably
generated. So $(\Xcal_{A'},\Xcal_v,\Xcal_{C'})$ is a disintegration triple by \Cref{thm-regular-conditional-Markov-kernel} point 1.,
for an arbitrary $\Xcal_{C'}$, which matters here, since $C'$ may contain input nodes.
So we may apply Left Weak Union \ref{sep:tci:l-uni} to the left statement and Right Weak Union \ref{sep:tci:r-uni} to the
right statement and get in both cases:
\[ X_{A'} \Indep_{\Pk(X_V\Vert X_J)}  X_{B'} \given X_{C'},X_v = X_C.\]
With $X_J$-Extended Inverted Right Decomposition \ref{sep:tci:ext-inv-r-dec}, applied with
$X_C \ismapof_\Kk X_J \otimes X_{B'} \otimes X_C$ (by \Cref{lem:restricted-reflexivity} and
\Cref{lem:join-upper-bound-2}, using \Cref{rem:join-commutative} to reorder the factors), and
Right Decomposition \ref{sep:tci:r-dec}, using $B \ins B' \dcup C$, this gives:
\[ X_{A'} \Indep_{\Pk(X_V\Vert X_J)}  X_B \given X_C.\]
By Extended Left Redundancy \ref{sep:tci:ext-l-red} we have:
\[ X_{A'}, X_C \Indep_{\Pk(X_V\Vert X_J)}  X_B \given X_{A'},X_C, \]
so Left Contraction \ref{sep:tci:l-con} and then Left Decomposition \ref{sep:tci:l-dec}, using $A \ins A' \dcup C$,
finally imply:
\[ X_A \Indep_{\Pk(X_V\Vert X_J)}  X_B \given X_C.\]
This shows the claim in case C.
\end{proof}

\end{Thm}

\begin{Rem}
    \begin{enumerate}
        \item Note that we only needed to use Left Weak Union \ref{sep:tci:l-uni} for the nodes $\{v\} \cup A'$ in case C.,
            which all lie in $V$ by \eqref{eqn:ainv}.
            So no assumptions about standard measurable spaces for $\Xcal_j$, $j \in J$, were needed.
    \end{enumerate}
\end{Rem}
\section{From Symmetric to Asymmetric Separoid Rules}
\label{sec:sym-sep-asym-sep}

Several of the notions compared in \Cref{sec:supp:comparison} are \emph{symmetric} ternary relations, while
transitional conditional independence and id-separation are \emph{asymmetric}. The passage between the two is always
the same and purely formal, so we record it once here, together with the terminology that goes with it.

Let $\Omega$ be a class equipped with an associative and commutative operation $\lor$ (up to a fixed notion of
isomorphism $\cong$), a neutral element $\varnothing$, a transitive relation $\ll$ that is compatible with $\cong$ and
satisfies \emph{product extension}, i.e.\ $\alpha \ll \beta \implies \alpha \ll \beta \lor \gamma$, and a ternary
relation $\Indep$ on $\Omega$. We write
\[ \alpha \approx \beta \qquad :\iff \qquad \lp \alpha \ll \beta \rp \;\land\; \lp \beta \ll \alpha \rp \]
for the equivalence induced by $\ll$ (on the sub-class where $\ll$ is reflexive) and we always assume that $\Indep$ is
invariant under $\cong$ and under $\approx$ in each of its three arguments.
Note that $\alpha \ll \alpha \lor \beta$ is \emph{not} assumed for general $\alpha$: together with
$\alpha \lor \varnothing \cong \alpha$ it would amount to reflexivity of $\ll$ everywhere, which fails in the
transitional instance, see \Cref{rem:ismapof-properties} item \ref{rem:no-reflexivity}. Product extension is the
weaker property that does hold there, for arbitrary transitional random variables.
We say that $\Indep$ satisfies the \emph{symmetric separoid rules} if for all $\alpha,\beta,\gamma,\delta \in \Omega$:
\begin{enumerate}[label=(S\arabic*)]
    \item \emph{Symmetry}: $\alpha \Indep \beta \given \gamma \implies \beta \Indep \alpha \given \gamma$;
    \item \emph{Redundancy}: $\alpha \ll \gamma \implies \alpha \Indep \beta \given \gamma$;
    \item \emph{Decomposition}: $\alpha \Indep \beta \lor \delta \given \gamma \implies \alpha \Indep \beta \given \gamma$;
    \item \emph{Weak Union}: $\alpha \Indep \beta \lor \delta \given \gamma \implies \alpha \Indep \beta \given \delta \lor \gamma$;
    \item \emph{Contraction}: $\lp \alpha \Indep \beta \given \delta \lor \gamma \rp \land \lp \alpha \Indep \delta \given \gamma\rp \implies \alpha \Indep \beta \lor \delta \given \gamma$.
\end{enumerate}
Some relations satisfy the two further rules
\begin{enumerate}[label=(S\arabic*),resume]
    \item \emph{Composition}: $\lp\alpha \Indep \beta \given \gamma\rp \land \lp \alpha \Indep \delta \given \gamma\rp \implies \alpha \Indep \beta \lor \delta \given \gamma$;
    \item \emph{Intersection}: $\lp\alpha \Indep \beta \given \delta \lor \gamma\rp \land \lp \alpha \Indep \delta \given \beta \lor \gamma\rp \implies \alpha \Indep \beta \lor \delta \given \gamma$,
        under a suitable disjointness assumption on $\beta$ and $\delta$;
\end{enumerate}
in which case we speak of a \emph{compositional} symmetric separoid, resp.\ of a \emph{graphoid}. These two rules
are, however, \emph{not} needed for the transfer below.

The asymmetric counterpart of these rules is the notion of a $\tau$-$\kappa$-separoid, which was set up
abstractly in \Cref{def:t-k-separoid} of the main text and which applies verbatim in the present setting.

\begin{Thm}[Asymmetric rules from symmetric ones]
    \label{thm:plus-t-k-separoid}
    Let $\Indep$ satisfy the symmetric separoid rules (S1)--(S5) and fix an element $\tau \in \Omega$ such that
    $\ll$ is reflexive at $\tau$ and $\tau \lor \tau \approx \tau$. Recall that $\ll$ satisfies product extension
(so that $\tau \ll \tau \lor \gamma$ for every $\gamma$, by reflexivity at $\tau$), and that $\Indep$ is
    $\cong$- and $\approx$-invariant.
    Define the \emph{$\tau$-shifted} ternary relation:
    \[ \alpha \Indep^\tau \beta \given \gamma \qquad :\iff \qquad \alpha \Indep \tau \lor \beta \given \gamma. \]
    Then $\Indep^\tau$ is a $\tau$-$\varnothing$-separoid in the sense of \Cref{def:t-k-separoid}, i.e.\ it
    satisfies the eleven rules a)--k) of \Cref{thm:separoid_axioms-tci} under the dictionary given there. It then
    also satisfies the symmetry rules l) and m) of \Cref{cor:sep-ci:symmetry} and, in the case
    $\tau \cong \varnothing$, the rule n). These are not part of the definition, see \Cref{rem:symmetry-derived};
    we include their derivations below to make visible that they cost nothing beyond a)--k).
    If $\Indep$ in addition satisfies (S6) then $\Indep^\tau$ also satisfies Left and Right Composition.
\begin{proof}
    Everything is obtained by unfolding the definition on both sides and applying (S1)--(S5) to the enlarged
    elements; not a single property of $\Indep$ beyond (S1)--(S5), the associativity, commutativity and neutrality of
    $\lor$ up to $\cong$, product extension for $\ll$, reflexivity of $\ll$ at $\tau$ and the $\cong$- and
    $\approx$-invariance of $\Indep$ is used. In particular no monotonicity of $\ll$ under $\lor$ is invoked
    anywhere, and the hypothesis $\tau \lor \tau \approx \tau$ enters only through \Cref{def:t-k-separoid}, not
    through any of the derivations. We spell the derivations out.

    \emph{Invariance of $\Indep^\tau$ under $\approx$.} In the first and third argument this is inherited from
    $\Indep$. In the second it needs an argument, since the wrapper $\tau \lor {-}$ is not known to be monotone:
    let $\beta \approx \beta'$ and assume $\alpha \Indep \tau \lor \beta \given \gamma$. From $\beta' \ll \beta$ and
    product extension we get $\beta' \ll \beta \lor (\tau \lor \gamma) \cong (\tau \lor \beta) \lor \gamma$, so
    Redundancy (S2) and Symmetry (S1) give $\alpha \Indep \beta' \given (\tau \lor \beta) \lor \gamma$; Contraction
    (S5) with the single element $\tau \lor \beta$ gives $\alpha \Indep \beta' \lor \tau \lor \beta \given \gamma$,
    and Decomposition (S3) dropping $\beta$ gives $\alpha \Indep \tau \lor \beta' \given \gamma$.

    \emph{a) Extended Left Redundancy.} The claim
    $\alpha \ll \gamma \implies \alpha \Indep \tau \lor \beta \given \gamma$
    is an instance of Redundancy (S2).

    \emph{b) $\tau$-Restricted Right Redundancy.} Unfolded, the claim reads
    $\alpha \Indep \tau \lor \varnothing \given \gamma \lor \tau$, i.e.\ $\alpha \Indep \tau \given \tau \lor \gamma$.
    Since $\tau \ll \tau \lor \gamma$, Redundancy (S2) gives $\tau \Indep \alpha \given \tau \lor \gamma$, and
    Symmetry (S1) turns this into the claim.

    \emph{c) Left Decomposition.} Unfolded:
    $\alpha \lor \delta \Indep \tau \lor \beta \given \gamma \implies \delta \Indep \tau \lor \beta \given \gamma$.
    Apply (S1), then Decomposition (S3) dropping $\alpha$, then (S1) again.

    \emph{d) Right Decomposition.} Unfolded:
    $\alpha \Indep \tau \lor \beta \lor \delta \given \gamma \implies \alpha \Indep \tau \lor \delta \given \gamma$.
    Since $\tau \lor \beta \lor \delta \cong (\tau \lor \delta) \lor \beta$, this is (S3) dropping $\beta$.

    \emph{e) $\tau$-Inverted Right Decomposition.} Unfolded, the rule reads
    $\alpha \Indep \tau \lor \beta \given \gamma \implies \alpha \Indep \tau \lor (\tau \lor \beta) \given \gamma$,
    and it is in fact an equivalence.
    ``$\Longleftarrow$'': since $\tau \lor \tau \lor \beta \cong (\tau \lor \beta) \lor \tau$, Decomposition (S3)
    dropping $\tau$ gives the claim.
    ``$\Longrightarrow$'': by product extension $\tau \ll \tau \lor \beta \lor \gamma$, so Redundancy (S2) and
    Symmetry (S1) give $\alpha \Indep \tau \given (\tau \lor \beta) \lor \gamma$; Contraction (S5) applied to this
    and to the hypothesis, with the \emph{single} element $\tau \lor \beta$, gives
    $\alpha \Indep \tau \lor (\tau \lor \beta) \given \gamma$.
    Note that the tempting shortcut ``$\tau \lor (\tau \lor \beta) \approx (\tau \lor \tau) \lor \beta \approx
    \tau \lor \beta$'' is \emph{not} available: passing from $\tau \lor \tau \approx \tau$ to
    $(\tau \lor \tau) \lor \beta \approx \tau \lor \beta$ would need $\ll$ to be monotone under $\lor$, which is not
    assumed and fails in the transitional instance.

    \emph{f) Left Weak Union.} Unfolded:
    $\alpha \lor \delta \Indep \tau \lor \beta \given \gamma \implies \alpha \Indep \tau \lor \beta \given \delta \lor \gamma$.
    By (S1) the assumption reads $\tau \lor \beta \Indep \alpha \lor \delta \given \gamma$, so Weak Union (S4) gives
    $\tau \lor \beta \Indep \alpha \given \delta \lor \gamma$, and (S1) gives the claim.

    \emph{g) Right Weak Union.} Unfolded:
    $\alpha \Indep (\tau \lor \beta) \lor \delta \given \gamma \implies \alpha \Indep \tau \lor \beta \given \delta \lor \gamma$,
    which is Weak Union (S4) applied to the two elements $\tau \lor \beta$ and $\delta$.

    \emph{h) Left Contraction.} Unfolding and applying (S1), the two assumptions read
    $\tau \lor \beta \Indep \alpha \given \delta \lor \gamma$ and $\tau \lor \beta \Indep \delta \given \gamma$.
    Contraction (S5), with $\tau \lor \beta$ in the left slot, gives
    $\tau \lor \beta \Indep \alpha \lor \delta \given \gamma$, and (S1) gives the claim.

    \emph{i) Right Contraction.} Unfolded, the assumptions read
    $\alpha \Indep \tau \lor \beta \given \delta \lor \gamma$ and $\alpha \Indep \tau \lor \delta \given \gamma$, and
    the claim reads $\alpha \Indep \tau \lor \beta \lor \delta \given \gamma$. Decomposition (S3) turns the second
    assumption into $\alpha \Indep \delta \given \gamma$; now (S5), applied to
    $\alpha \Indep \tau \lor \beta \given \delta \lor \gamma$ and $\alpha \Indep \delta \given \gamma$, gives
    $\alpha \Indep (\tau \lor \beta) \lor \delta \given \gamma$, which is the claim.

    \emph{j) Right Cross Contraction.} Unfolded, the second assumption reads
    $\delta \Indep \tau \lor \alpha \given \gamma$, which by (S3) and (S1) gives $\alpha \Indep \delta \given \gamma$.
    Together with the first assumption we conclude exactly as in i).

    \emph{k) Flipped Left Cross Contraction.} Unfolded, the assumptions read
    \[ \text{(i) } \alpha \Indep \tau \lor \beta \given \delta \lor \gamma, \qquad\qquad
       \text{(ii) } \beta \Indep \tau \lor \delta \given \gamma, \]
    and the claim reads $\beta \Indep \tau \lor \alpha \lor \delta \given \gamma$.
    From (i), Weak Union (S4) (moving $\tau$, \emph{not} $\beta$, into the conditioning position) gives
    $\alpha \Indep \beta \given \tau \lor \delta \lor \gamma$, hence
    $\beta \Indep \alpha \given (\delta \lor \tau) \lor \gamma$ by (S1).
    Contraction (S5), applied to this and to (ii) with the \emph{single} element $\delta \lor \tau$, gives
    $\beta \Indep \alpha \lor \delta \lor \tau \given \gamma$, which is the claim. In particular no Composition is
    needed here.

    \emph{l) Restricted Symmetry.} The premise $\beta \Indep^\tau \kappa \given \gamma$ gives
    $\beta \Indep^\tau \varnothing \given \gamma$ by d) Right Decomposition, since
    $\kappa \cong \kappa \lor \varnothing$. Now apply k) with $\delta := \varnothing$, using
    $\varnothing \lor \gamma \cong \gamma$ and $\alpha \lor \varnothing \cong \alpha$. No relation between
    $\kappa$ and $\varnothing$ is needed.

    \emph{m) $\tau$-Restricted Symmetry.} Apply l) with $\gamma \lor \tau$ in place of $\gamma$; its second premise
    $\beta \Indep^\tau \varnothing \given \gamma \lor \tau$ is b).

    \emph{n) Symmetry.} If $\tau \cong \varnothing$ then $\gamma \lor \tau \cong \gamma$ by neutrality, so m)
    together with the $\cong$-invariance of $\Indep$ gives the claim. Note that $\tau \approx \varnothing$ would
    not suffice here: nothing in \Cref{def:t-k-separoid} makes $\ll$ monotone in $\lor$, and in the transitional
    instance $\Zk \otimes \Tk \approx_\Kk \Zk$ fails for stochastic $\Zk$, see
    \Cref{rem:ismapof-properties} item \ref{rem:no-reflexivity}.

    \emph{Left and Right Composition under (S6).} For Left Composition, unfolding and (S1) turn the assumptions into
    $\tau \lor \beta \Indep \alpha \given \gamma$ and $\tau \lor \beta \Indep \delta \given \gamma$; (S6) and (S1)
    give the claim. For Right Composition, the assumptions read $\alpha \Indep \tau \lor \beta \given \gamma$ and
    $\alpha \Indep \tau \lor \delta \given \gamma$; Decomposition (S3) turns the second one into
    $\alpha \Indep \delta \given \gamma$, and (S6) applied to the first assumption and to this gives
    $\alpha \Indep (\tau \lor \beta) \lor \delta \given \gamma$, which is the claim.
\end{proof}
\end{Thm}

\begin{Rem}
    The two \emph{intersection} rules are deliberately absent from \Cref{thm:plus-t-k-separoid}: both need a notion
    of disjointness, which a general $\lp\Omega,\lor\rp$ does not provide, and the derivation of Right Intersection
    \ref{sep:sig:r-int} for id-separation in addition passes to the set $(B \cup J)\sm D$, i.e.\ it uses set
    differences.
\end{Rem}

\begin{Rem}
    \Cref{thm:plus-t-k-separoid} is the abstract reason why id-separation, see \Cref{def:d-separation}, satisfies
    exactly the same rules as transitional conditional independence, see \Cref{thm:separoid_axioms-tci}:
    id-separation is the $J$-shift of ordinary d-separation, and \Cref{sec:supp:separoid-rules-sigma} is nothing but
    the instance $\Indep = \Perp^{d}_\Gk$, $\Indep^\tau = \Perp^{\id}_\Gk$, $\tau = J$, $\lor = \cup$,
    $\ll \,=\, \ins$ of the proof above, spelled out for sets of nodes.
    The theorem also applies to the symmetric notions of extended conditional independence
    discussed in \Cref{sec:supp:comparison}, e.g.\ to variation conditional independence and to
    $\Qcal$-extended conditional independence.
\end{Rem}

Conversely, one can symmetrize an asymmetric relation with the logical ``or''.

\begin{Thm}[Symmetrized transitional conditional independence]
    \label{thm:sym-t-k-separoid}
    Let the setting be as in \Cref{thm:separoid_axioms-tci} and assume that all occurring codomains form disintegration
    triples, see \Cref{def:disintegration-triple}. Define the symmetrized ternary relation:
    \[ \Xk \Indep^\lor_{\Kk(W|T)} \Yk \given \Zk \qquad :\iff \qquad
    \lp \Xk \Indep_{\Kk(W|T)} \Yk \given \Zk \rp \; \lor \; \lp \Yk \Indep_{\Kk(W|T)} \Xk \given \Zk \rp, \]
    where the $\lor$ in the middle of the right hand side is the logical OR.
    Then $\Indep^\lor_{\Kk(W|T)}$ satisfies the symmetric separoid rules (S1)--(S5), i.e.\ Symmetry, Redundancy,
    Decomposition, Weak Union and Contraction, where $\Xk \ll \Zk$ is read as $\Xk \ismapof_\Kk \Zk$.
\begin{proof}
    Symmetry holds by construction, and Redundancy is exactly Extended Left Redundancy \ref{sep:tci:ext-l-red}.
    Decomposition: if $\Xk \Indep \Yk \otimes \Uk \given \Zk$ then Right Decomposition \ref{sep:tci:r-dec} gives
    $\Xk \Indep \Yk \given \Zk$; and if $\Yk \otimes \Uk \Indep \Xk \given \Zk$ then Left Decomposition
    \ref{sep:tci:l-dec} gives $\Yk \Indep \Xk \given \Zk$. In both cases the disjunction holds.
    Weak Union: analogously from Right Weak Union \ref{sep:tci:r-uni} and Left Weak Union \ref{sep:tci:l-uni}.
    Contraction: unfolding the two disjunctions gives four cases, and these are handled by
    Right Contraction \ref{sep:tci:r-con}, Right Cross Contraction \ref{sep:tci:rc-con}, Flipped Left Cross
    Contraction \ref{sep:tci:flc-con} and Left Contraction \ref{sep:tci:l-con}, one case each. This is precisely the reason why all
    four contraction rules were included in \Cref{thm:separoid_axioms-tci}.
\end{proof}
\end{Thm}

\begin{Rem}
    Note that $\Indep^\lor$ is genuinely weaker than $\Indep$: we always have the implication
    $\lp \Xk \Indep_\Kk \Yk \given \Zk \rp \implies \lp \Xk \Indep^\lor_\Kk \Yk \given \Zk \rp$, but not conversely, so
    the symmetrized version may have lost information about the interplay between $\Xk$, $\Yk$, $\Zk$ and $\Tk$.
    In particular it can no longer express the asymmetric statistical concepts of
    \Cref{sec:applications-statistics}.
\end{Rem}
\section{Comparison to Other Notions of Conditional Independence}
\label{sec:supp:comparison}

In this section we want to look at other notions of conditional independence and compare them
to transitional conditional independence.

Recall that for transition probability space
  $\lp \Wcal \times \Tcal, \Kk(W|T) \rp$  and transitional random variables $\Xk:\, \Wcal \times \Tcal \dshto \Xcal$ and 
    $\Yk:\, \Wcal \times \Tcal \dshto \Ycal$ and $\Zk:\, \Wcal \times \Tcal \dshto \Zcal$ 
    with joint Markov kernel:
    \[ \Kk(X,Y,Z|T) := \lp\Xk(X|W,T) \otimes \Yk(Y|W,T)  \otimes \Zk(Z|W,T)  \rp \circ \Kk(W|T), \]
    we define the \emph{transitional conditional independence} of $\Xk$ from $\Yk$ given $\Zk$:
    \[ \Xk \Indep_{\Kk(W|T) } \Yk \given \Zk \quad :\iff \quad 
    \exists \Qk(X|Z):\; \Kk(X,Y,Z|T) = \Qk(X|Z) \otimes \Kk(Y,Z|T).  \]

\subsection{Variation Conditional Independence}
\label{sec:supp:variation-ci}

We follow \cite{CD17, Daw01b} and their supplementary material \cite{CD17-supp} to review variation conditional independence and then comment on some possible generalizations.

For this let $\Wcal$ be a  set and $X:\, \Wcal \to \Xcal$, $Y:\, \Wcal \to \Ycal$, $Z:\, \Wcal \to \Zcal$, $U:\,\Wcal \to \Ucal$ 
be maps.

\begin{Not}[See \cite{CD17} \S2.2]
    We define:
    \[ \Rcal(X) := X(\Wcal) = \{ X(w) \,|\, w \in \Wcal\} \in 2^\Xcal,\]
    and for $z \in \Zcal$:
    \[ \Rcal(X|Z=z) := X\lp Z^{-1}(z)\rp = \{ X(w)\,|\, w \in \Wcal, Z(w)=z \} \in 2^\Xcal. \]
    We then define the map:
    \[ \Rcal(X|Z):\, \Zcal \to 2^\Xcal, \quad z \mapsto \Rcal(X|Z=z).   \]
    In this sense we then can also make sense of:
    \[ \Rcal(X,Y|Z,U) : \Zcal \times \Ucal \to 2^{\Xcal \times \Ycal},\] 
   \[ (z,u) \mapsto \Rcal(X,Y|Z=z,U=u) := \{ (X(w),Y(w))\,|\, w\in \Wcal, Z(w)=z, U(w)=u \}.\]
\end{Not}

\begin{Def}[Variation conditional independence]
    We will say that \emph{$X$ is variation conditionally independent of $Y$ given $Z$}
    if:
   \[ \forall (y,z) \in \Rcal(Y,Z) :\, 
   \Rcal(X|Y=y,Z=z) = \Rcal(X|Z=z).\]
   In symbols we will write then:
   \[X \Indep_v Y \given Z.\]
\end{Def}

\begin{Not}
    We will write:
    \[ X \ismapof_v Y \]
    if there exists a map $\varphi:\, \Ycal \to \Xcal$ such that $\varphi \circ Y = X$.
    Note that we use a slightly simpler, but equivalent relation, than \cite{CD17} \S2.2, Prop. 2.6.
\end{Not}

\begin{Rem}[See \cite{CD17} Thm. 2.7., \cite{CD17-supp}]
    \label{rem:supp:var-tci-sep}
    The ternary relation $\Indep_v$ together with $\ismapof_v$ and $X \lor Y:= (X,Y)$
    is a (symmetric) separoid; \cite{CD17} Thm.\ 2.7 shows that it is even a \emph{strong} separoid, which
    stochastic conditional independence is not.
\end{Rem}

\begin{Rem}
    \label{rem:supp:var-tci}
    The relation between variation conditional independence 
    and stochastic conditional independence for random variables seems rather on
    the structural side, i.e.\ both follow similar functorial
    relations. In short, if one wants to go from variation to stochastic conditional independence 
    one could start by replacing $2^\Xcal$ with $\Pcal(\Xcal)$ and maps with measurable maps, etc.
    Then maps $\Zcal \to 2^\Xcal$ become measurable maps $\Zcal \to \Pcal(\Xcal)$, which are nothing
    else but Markov kernels, reflecting the approach we went down for 
    transitional conditional independence.
    So $\Rcal(X,Y,Z)$ can be represented as the (constant) map:
    \[ \Rcal(X,Y,Z):\,\Asterisk \to 2^{\Xcal \times \Ycal \times \Zcal}, 
    \quad \ast \mapsto \Rcal(X,Y,Z). \]
    We also have the ``marginal'':
    \[ \Rcal(Y,Z):\,\Asterisk \to 2^{\Ycal \times \Zcal}, \quad \ast \mapsto \Rcal(Y,Z)
    = \pr_{\Ycal \times \Zcal} \lp \Rcal(X,Y,Z) \rp,\]
    which is received by ignoring the $X$-entries.\\
    It is then easily seen that  we have: $X \Indep_v Y\given Z $ iff there exists a 
    map: 
    \[\Qcal(X|Z): \, \Zcal \to 2^\Xcal,\] 
    such that:
    \[  \Rcal(X,Y,Z) = \Qcal(X|Z) \otimes_v \Rcal(Y,Z),\]
    where we put:
    \[  \Qcal(X|Z) \otimes_v \Rcal(Y,Z) := 
    \bigcup_{(y,z) \in \Rcal(Y,Z)} \bigcup_{x \in \Qcal(X|Z=z)} \{ (x,y,z) \}. \]
\end{Rem}

\begin{Eg}
    \label{eg:supp:var-tci}
    We can now apply the above to
    $X,Y,Z$ with common domain $\Wcal \times \Tcal$ and $T:\,\Wcal \times \Tcal$ the canonical projection map.
    Then we get:
    \begin{align*}
       &\qquad\qquad\; X \Indep_{v,T} Y \given Z\\ 
       & :\iff \quad 
       X \Indep_{v} T,Y \given Z \\ 
       &\iff \quad  \exists \Qcal(X|Z):\, \Zcal \to 2^{\Xcal}:\quad
     \Rcal(X,Y,Z|T) = \Qcal(X|Z) \otimes_v \Rcal(Y,Z|T),
     \end{align*}
    where we again now have:
    \[\Rcal(X,Y,Z|T=t) = \{(x,y,z) \,|\, \exists w \in \Wcal: X(w,t)=x, Y(w,t)=y, Z(w,t)=z\}.  \]
    This shows the close formal relationship between variation conditional independence and
    transitional conditional independence $\Indep_\Kk$. 
    It then follows from \Cref{rem:supp:var-tci-sep} and 
    the general theory in \Cref{sec:sym-sep-asym-sep} with \Cref{thm:plus-t-k-separoid} 
    that $\Indep_{v,T}$ forms a $T$-$\ast$-separoid.
    One can thus combine $\Indep_\Kk$ and $\Indep_{v,T}$ with a logical ``and'', while still preserving the
    $T$-$\ast$-separoid rules.
\end{Eg}

    It seems, more generally,
    that one can formulate a (transitional) conditional independence relation 
    in any monad with products and some extra structure.
    We leave this for future research.

\subsection{Transitional Conditional Independence for Random Variables}

If we wanted to re-define the notion of \emph{independent} random variables $X$ and $Y$ on a probability space $(\Wcal,\Pk(W))$ 
we would have a hard time coming up with something else than
the classical definition of:
\begin{align}
    \Pk(X,Y) &= \Pk(X) \otimes \Pk(Y), \label{eq:independence}
\end{align}
where $\Pk(X)$ and $\Pk(Y)$ are the marginals. This is in contrast to \emph{conditionally independent} 
random variables $X$ and $Y$ given a third $Z$, where many nuances can play a role.
For instance, the direct analogue of relation \ref{eq:independence} would read like:
\begin{align}
    \Pk(X,Y|Z) &= \Pk(X|Z) \otimes \Pk(Y|Z) &\Pk(Z)\text{-a.s.}  \label{eq:cond-ind-fact}
\end{align}
The problem with definition \ref{eq:cond-ind-fact} is that the conditional probability distributions, like $\Pk(X,Y|Z)$, 
may not exist on general measurable spaces $\Xcal$, $\Ycal$, $\Zcal$, in contrast to the marginals $\Pk(X)$, $\Pk(Y)$ 
in equation \ref{eq:independence}. This then forces one to restrict oneself to only work 
with measurable spaces 
where regular conditional probability distributions exist, like standard measurable spaces. 
But even if the existence were guaranteed, they would only be unique up to some null sets.
One then either ends up with a notion of conditional independence that would depend on the choices made 
or, as the better alternative, one would work with almost-sure equations like we
already indicated in equation \ref{eq:cond-ind-fact}.\\
If one wanted to work with more general measurable spaces one could demand that equation \ref{eq:cond-ind-fact}
only holds for every $A \in \Bcal_\Xcal$ and $B \in \Bcal_\Ycal$ individually. Furthermore, one then could replace 
$\Pk(X \in A, Y \in B|Z)$ with conditional expectations $\E[\I_A(X) \cdot \I_B(Y)|Z]$, etc., which exist on all measurable spaces. 
One then arrives at the most general 
and \emph{weak form of conditional independence} for random variables:
\begin{align} 
&\forall A \in \Bcal_\Xcal\, \forall B \in \Bcal_\Ycal:\, 
    \E[\I_A(X) \cdot \I_B(Y)|Z] = \E[\I_A(X)|Z] \cdot \E[\I_B(Y)|Z]& \Pk(W)\text{-a.s.}, \label{eq:cond-ind-exp-1}
\end{align}
which can, equivalently, but more compactly, also be written as:
\begin{align} 
&\forall A \in \Bcal_\Xcal:\; 
    \E[\I_A(X)|Y,Z] = \E[\I_A(X)|Z] & \Pk(W)\text{-a.s.} \label{eq:cond-ind-exp-2}
\end{align}
We will use the following symbols for weak conditional independence:
\[ X \Indep_{\Pk(W)}^\omega Y \given Z. \]
Furthermore, if we used definition \ref{eq:cond-ind-exp-1} or \ref{eq:cond-ind-exp-2} on standard measurable spaces, 
where regular conditional probability distributions like $\Pk(X,Y|Z)$ exist, the equation \ref{eq:cond-ind-fact} would 
automatically be implied. So the equations \ref{eq:cond-ind-exp-1} or \ref{eq:cond-ind-exp-2}  seem to be the way to go, 
as one does not need to bother with existence questions, and when existence is secured the above versions 
are equivalent anyways. The only downside is that this definition does not provide one with a meaningful factorization.
Furthermore, the conditional expectations, like $\E[\I_A(X)|Z]$, are only defined for each event $A$ separately 
and thus might not be countably additive in $A$. So we are not given an object like a conditional distribution that we could use
to further work with.\\
In contrast, our definition of \emph{transitional conditional independence} $X \Indep_{\Pk(W)} Y \given Z$, 
when restricted to random variables, would read like:
\begin{align}
    \exists \Pk(X|Z):\quad  \Pk(X,Y,Z) &= \Pk(X|Z) \otimes \Pk(Y,Z), \label{eq:trans-cond-ind}
\end{align}
where clearly $\Pk(X|Z)$ would be a regular conditional probability distribution of $X$ given $Z$.
So the existence of one of the regular conditional probability distributions $\Pk(X|Z)$ 
and a proper factorization of the joint distribution
are directly built into the definition of transitional conditional independence. 
This makes this notion also meaningful for general measurable spaces, with the tendency that random variables are 
declared conditional dependent if such a regular conditional probability distribution does not even exist.
Definitions \ref{eq:trans-cond-ind} and \ref{eq:cond-ind-exp-2} are equivalent as soon as $\Xcal$ is standard and
$\Zcal$ is countably generated, with $\Ycal$ \emph{arbitrary}, by \Cref{thm:supp:ci-equivalences} 4.\ with
$\Tcal = \Asterisk$ together with \Cref{thm-regular-conditional-Markov-kernel}. Definition \ref{eq:cond-ind-fact}
presupposes in addition the existence of $\Pk(X,Y|Z)$ and $\Pk(Y|Z)$, so all three are equivalent
on standard measurable spaces. Note that transitional conditional independence \ref{eq:trans-cond-ind} is asymmetric in nature,
which at this level might look like a flaw, but which allows one to generalize the definition of transitional conditional independence
to \emph{transitional random variables}, where dependencies are asymmetric from the start.

\subsection{Transitional Conditional Independence for Deterministic Variables}
\label{sec:tci-det}

\begin{Thm}[Transitional conditional independence for deterministic variables]
    \label{thm:supp:variation-ci}
    Let $F:\, \Tcal \to \Fcal$ and $H:\,\Tcal \to \Hcal$ be measurable maps with $\Fcal$ standard.
    We now consider them as (deterministic) transitional random variables on the transition probability space $(\Wcal \times \Tcal,\Kk(W|T))$.
    Let $\Yk:\, \Wcal \times \Tcal \dshto \Ycal$ be another transitional random variable.\\
    Then the following statements are equivalent:
    \begin{enumerate}
        \item $\displaystyle F \Indep_{\Kk(W|T)} \Yk \given H$.
        \item There exists a measurable function $\varphi:\, \Hcal \to \Fcal$ such that $F=\varphi \circ H$.
    \end{enumerate}
\begin{proof}
    ``$\Longleftarrow$'': This direction follows from Extended Left Redundancy \ref{sep:tci:ext-l-red}.\\
  ``$\implies$'': Since $F$ and $H$ are deterministic and only dependent on $T$ we get that:
    \[ \Kk(F,Y,H|T) = \deltabf(F|T) \otimes \deltabf(H|T) \otimes \Kk(Y|T).  \]
    By the conditional independence we now have a Markov kernel $\Qk(F|H)$ such that we have the factorization:
    \[ \Kk(F,Y,H|T) = \Qk(F|H) \otimes \Kk(Y,H|T) = \Qk(F|H) \otimes \deltabf(H|T) \otimes \Kk(Y|T).  \]
    Marginalizing out $Y$, $H$ and taking $T=t$ we get from these equations:
    \[ \deltabf_{F(t)} = \deltabf(F|T=t) = \Qk(F|H(t)),\]
    which is a Dirac measure centered at $F(t)$.
    We can now define the mapping:
    \[\varphi:\, H(\Tcal) \to \Fcal, \quad H(t) \mapsto F(t), \]
    which is well-defined, because $h:=H(t_1)=H(t_2)$ implies that $\Qk(F|H=h)$ is a Dirac measure centered at $F(t_1)$ and $F(t_2)$.
    Since $\Bcal_\Fcal$ separates points ($\Fcal$ is standard) we get: $F(t_1)=F(t_2)$.
    $\varphi$ is measurable.
    Indeed, its composition with $\delta:\,\Fcal \to \Pcal(\Fcal), \, z \mapsto \deltabf_z$ equals $\Qk(F|H)$, which is measurable.
    Since $j_A \circ \delta = \I_A$ for every $A \in \Bcal_\Fcal$ and since the evaluation maps $j_A$ generate 
    $\Bcal_{\Pcal(\Fcal)}$, we have $\Bcal_\Fcal = \delta^*\Bcal_{\Pcal(\Fcal)}$, so also $\varphi$ is measurable.
    Since $\Fcal$ is a standard measurable space, $\varphi$ extends to a measurable mapping $\varphi:\, \Hcal \to \Fcal$ 
    by Kuratowski's extension theorem for standard measurable spaces (see \cite{Kec95} 12.2 and \Cref{thm-kuratowski-extension}). 
    Finally, note that we have $F(t) = \varphi(H(t))$ for all $t \in \Tcal$, which shows the claim.
    \[\xymatrix{ \Tcal \ar^-F[r] \ar_-H[d] & \Fcal \ar@{^(->}^-\delta[d] \\
            \Hcal \ar_-\Qk[r] \ar^-{\exists \varphi}[ur] & \Pcal(\Fcal)
    }\]
\end{proof}
\end{Thm}

\subsection{Equivalent Formulations of Transitional Conditional Independence}
\label{sec:tci-equiv}

Our groundwork of developing the framework of transition probability spaces, transitional random variables 
and transitional conditional independence
now allows us to rigorously compare different notions of extended conditional independence in the literature.
We will compare to three of them, namely the one from \cite{CD17, RERS17, FM20}.

To relate transitional conditional independence to other notion of conditional independence it is useful to reformulate
transitional conditional independence in other terms. The main result for this will be the next theorem.

\begin{Thm}
    \label{thm:supp:ci-equivalences}
Let $\lp \Wcal \times \Tcal, \Kk(W|T) \rp$ be a transition probability space and $\Xk:\, \Wcal \times \Tcal \dshto \Xcal$ and 
    $\Yk:\, \Wcal \times \Tcal \dshto \Ycal$ and $\Zk:\, \Wcal \times \Tcal \dshto \Zcal$ transitional random variables.
    We put:
    \[ \Kk(X,Y,Z|T) := \lp\Xk(X|W,T) \otimes \Yk(Y|W,T)  \otimes \Zk(Z|W,T)  \rp \circ \Kk(W|T), \]
We will write $\Kk(X|Z,\cancel{T})$ 
for any version of the Markov kernel appearing in the conditional independence $\Xk \Indep_{\Kk(W|T)} \deltabf_\ast \given \Zk$
(only in case it holds). Statement 3.\ below does not depend on that choice: two versions agree
$\Kk(Z|T=t)$-almost surely for every $t$, see \Cref{ess-unique}, the identity is tested on rectangles
$A \times B \times C$ against $\Kk(Y,Z|T=t)$, whose $\Zcal$-marginal is $\Kk(Z|T=t)$, and only one event $A$ is
used at a time, so no countable generation of $\Bcal_\Xcal$ is needed either.\\
With these notations, the following are equivalent:
\begin{enumerate}
    \item $\displaystyle \Xk \Indep_{\Kk(W|T)} \Yk \given \Zk,$
    \item $\displaystyle \Xk \Indep_{\Kk(W|T)} \Tk\otimes\Yk \given \Zk,$
    \item $\displaystyle \Xk \Indep_{\Kk(W|T)} \deltabf_\ast \given \Zk$ and 
        $\displaystyle \Kk(X,Y,Z|T) = \Kk(X|Z,\cancel{T}) \otimes \Kk(Y,Z|T)$.
     \item  $\displaystyle  \Xk \ismapof_\Kk^\ast \Zk$ and for every $t \in \Tcal$ we have:
         $\displaystyle X_t \Indep_{\Kk(X,Y,Z|T=t)}^\omega Y_t \given Z_t$ (in the weak sense).
\end{enumerate}
Furthermore, any of those points implies
the following:
\begin{enumerate}[resume]
    \item For every probability distribution $\Qk(T) \in \Pcal(\Tcal)$ we have the conditional independence:
        \[ \Xk \Indep_{\Kk(W|T)\otimes \Qk(T) } \Tk\otimes\Yk \given \Zk.\]
\end{enumerate}
\begin{proof}
    3. $\implies$ 1. is clear by definition.\\
    1. $\implies$ 2.: by $\Tk$-Inverted Right Decomposition \ref{sep:tci:inv-r-dec}.\\
    2. $\implies$ 4.,5.:
    By assumption we have the factorization:
    \[\Kk(X,Y,Z,T|T) = \Kk(X|Z) \otimes \Kk(Y,Z,T|T),\]
    for some Markov kernel $\Kk(X|Z)$. Via marginalization and multiplication this implies the two equations:
    \begin{align*}
        \Kk(X,Z,T|T)            &= \Kk(X|Z) \otimes \Kk(Z,T|T),\\
        \underbrace{\Kk(X,Y,Z|T)\otimes \Qk(T)}_{=:\Qk(X,Y,Z,T)} &= \Kk(X|Z) \otimes \underbrace{\Kk(Y,Z|T) \otimes \Qk(T)}_{=\Qk(Y,Z,T)},
    \end{align*}
    for every $\Qk(T) \in \Pcal(\Tcal)$. The last equation shows 5.\\ 
    If we take $\Qk(T)=\deltabf_t$ we get:
    \[\Kk(X,Y,Z|T=t) = \Kk(X|Z) \otimes \Kk(Y,Z|T=t). \]
    Together with the first of the above equations this shows 4.\\
    4. $\implies$ 3.:  By $\Xk \ismapof_\Kk^\ast \Zk$ we have a factorization:
    \[\Kk(X,Z|T) = \Pk(X|Z) \otimes \Kk(Z|T).\]
    This means that for every $t \in \Tcal$ and every measurable $A \ins \Xcal$, $C \ins \Zcal$ we have:
    \[ \E_t\lB \I_A(X_t) \cdot \I_C(Z_t)\rB = \E_t\lB \Pk(X \in A|Z_t) \cdot \I_C(Z_t)\rB,  \]
    where the expectation $\E_t$ is w.r.t.\ $\Kk(X,Y,Z|T=t)$. This shows that $\Pk(X \in A|Z_t)$ is a version of $\E_t[\I_A(X_t)|Z_t]$ for every $t \in \Tcal$, 
    by the defining properties of conditional expectation.\\
    By the assumption $ X_t \Indep_{\Kk(X,Y,Z|T=t)}^\omega Y_t \given Z_t$ we then have for every fixed $t\in \Tcal$ and measurable $A \ins \Xcal$:
    \[ \E_t[\I_A(X_t)|Y_t,Z_t] = \E_t[\I_A(X_t)|Z_t] = \Pk(X \in A|Z_t)\qquad \Kk(X,Y,Z|T=t)\text{-a.s.}\]
    By the defining properties of conditional expectation for $\E_t[\I_A(X_t)|Y_t,Z_t]$  we then get that for every measurable  $A \ins \Xcal$, $B \ins \Ycal$, $C \ins \Zcal$:
    \[ \E_t\lB \I_A(X_t) \cdot \I_B(Y_t) \cdot \I_C(Z_t)\rB = \E_t\lB \Pk(X \in A|Z_t) \cdot \I_B(Y_t) \cdot \I_C(Z_t)\rB. \]
    Since this holds for every $t \in \Tcal$ we get:
    \[ \Kk(X,Y,Z|T) = \Pk(X|Z) \otimes \Kk(Y,Z|T),\]
    which shows the claim.
\end{proof}
\end{Thm}

\begin{Cor}
    If $\Xcal$ is standard and $\Zcal$ countably generated (e.g.\ also standard) 
then we have the equivalence:
\[ \Xk \Indep_{\Kk(W|T)} \Yk \given \Zk\otimes\Tk \qquad \iff \qquad \forall t \in \Tcal:\quad X_t \Indep_{\Kk(X,Y,Z|T=t)}^\omega Y_t \given Z_t.\]
\begin{proof}
    This directly follows from \Cref{thm:supp:ci-equivalences} 4.\ with $(Z,T)$ in the role of $Z$
    and \Cref{conditional-markov-kernel-as-ci} to get the first part of 4.
    One step deserves to be spelled out, since it is used in both directions: applied in that form, point 4.\ gives
    the weak conditional independence of $X_t$ and $Y_t$ given the \emph{pair} $(Z_t,T_t)$. But under
    $\Kk(X,Y,Z|T=t)$ the transitional random variable $T_t$ is almost surely constant, equal to $t$, so
    $\sigma(T_t)$ is trivial modulo $\Kk(X,Y,Z|T=t)$-null sets and therefore
    $\sigma(Z_t,T_t) = \sigma(Z_t)$ modulo null sets. The two conditional expectations agree almost surely, and
    the statement given $(Z_t,T_t)$ is the statement given $Z_t$.
\end{proof}
\end{Cor}

\subsection{Criteria for Checking Transitional Conditional Independence}
\label{sec:supp:checking}

In this subsection we supply the proofs for the criteria stated in \Cref{sec:checking-tci}. Throughout, the
underlying transition probability space is $\lp \Wcal \times \Tcal, \Kk(W|T) \rp$ with
$\Wcal := \Xcal \times \Ycal \times \Zcal$ and $\Kk(W|T) := \Kk(X,Y,Z|T)$, taken with its three coordinate
variables, so that
\Cref{def:transitional_conditional_independence} reads: $\Xk \Indep_\Kk \Yk \given \Zk$ holds if and only if there
is a Markov kernel $\Qk(X|Z):\, \Zcal \dshto \Xcal$ with:
\[ \Kk(X,Y,Z|T=t) = \Qk(X|Z) \otimes \Kk(Y,Z|T=t) \qquad \text{ for every } t \in \Tcal. \]

We first record how the marginal densities are obtained.

\begin{Rem}[Marginal densities]
    \label{rem:supp:marginal-densities}
    Let $\mubf_\Xcal$, $\mubf_\Ycal$, $\mubf_\Zcal$ be $\sigma$-finite measures on $\Xcal$, $\Ycal$, $\Zcal$, put
    $\mubf := \mubf_\Xcal \otimes \mubf_\Ycal \otimes \mubf_\Zcal$. Let $\Nk$ be a probability measure on
    $\Xcal \times \Ycal \times \Zcal$ with $\Nk \ll \mubf$ and let $n(x,y,z) \in [0,\infty]$ be a density of
    $\Nk$ w.r.t.\ $\mubf$. Put:
    \begin{align*}
        n(x,z) &:= \int_\Ycal n(x,y,z)\, \mubf_\Ycal(dy), &
        n(y,z) &:= \int_\Xcal n(x,y,z)\, \mubf_\Xcal(dx), \\
        n(z)   &:= \int_{\Xcal \times \Ycal} n(x,y,z)\, \lp\mubf_\Xcal \otimes \mubf_\Ycal\rp(d(x,y)). & &
    \end{align*}
    Since all the measures involved are $\sigma$-finite and $n \ge 0$ is jointly measurable, Tonelli's theorem
    shows that these functions are again jointly measurable in their arguments and that they are densities of the
    marginals $\Nk(X,Z)$, $\Nk(Y,Z)$ and $\Nk(Z)$ w.r.t.\ $\mubf_\Xcal \otimes \mubf_\Zcal$,
    $\mubf_\Ycal \otimes \mubf_\Zcal$ and $\mubf_\Zcal$, respectively. In particular $\int_\Zcal n(z)\, \mubf_\Zcal(dz) = 1$, so $n(z) < \infty$ for
    $\mubf_\Zcal$-almost all $z$, and $\int_\Xcal n(x,z)\,\mubf_\Xcal(dx) = n(z)$ for \emph{every} $z \in \Zcal$.
    Note also that $\mubf_\Xcal \neq 0$, since otherwise $\mubf = 0$ and no probability measure could be dominated
    by it; as $\mubf_\Xcal$ is $\sigma$-finite there thus exists a probability density $w:\, \Xcal \to [0,\infty)$
    w.r.t.\ $\mubf_\Xcal$, a fact we will use to normalize on null sets. We apply all of this to
    $\Nk := \Kk(X,Y,Z|T=t)$ with $n := k(x,y,z|t)$, and to the mixture $\Nk := \Pk(X,Y,Z)$ with
    $n := p(x,y,z)$.
\end{Rem}

We first prove the criterion of \Cref{thm:tci-prior}, which rests on no density assumption at all.

\begin{Thm}[A witness from a prior]
    \label{thm:supp:tci-prior}
    Let $\Pk(T)$ be a probability measure on $\Tcal$, put
    $\Pk(X,Y,Z) := \Kk(X,Y,Z|T) \circ \Pk(T)$ and assume that:
    \begin{enumerate}
        \item $\Kk(Z|T=t) \ll \Pk(Z)$ for every $t \in \Tcal$;
        \item $\Pk(X,Z) = \Pk(X|Z) \otimes \Pk(Z)$ for a Markov kernel $\Pk(X|Z):\, \Zcal \dshto \Xcal$.
    \end{enumerate}
    Then, for every such $\Pk(X|Z)$:
    \[ \Xk \Indep_\Kk \Yk \given \Zk \qquad \iff \qquad
       \Kk(X,Y,Z|T) = \Pk(X|Z) \otimes \Kk(Y,Z|T). \]
\begin{proof}
    ``$\Longleftarrow$'' is \Cref{def:transitional_conditional_independence}, since $\Pk(X|Z)$ is a Markov
    kernel. For ``$\implies$'' let $\Qk(X|Z)$ be a witness, so that:
    \begin{align}
        \Kk(X,Y,Z|T=t) = \Qk(X|Z) \otimes \Kk(Y,Z|T=t) \qquad \text{ for every } t \in \Tcal.
        \label{eq:supp-prior-witness}
    \end{align}

    \emph{Step 1: averaging over the input.} Fix measurable $A \ins \Xcal$, $B \ins \Ycal$, $C \ins \Zcal$ and
    put:
    \[ \varphi(t) := \int_{B \times C} \Qk(X \in A|Z=z)\, \Kk(Y \in dy, Z \in dz|T=t). \]
    The function $(y,z) \mapsto \Qk(X \in A|Z=z) \cdot \I_{B \times C}(y,z)$ is measurable and bounded by $1$, so
    $\varphi$ is measurable and bounded by $1$ as well, by the defining property of the Markov kernel
    $\Kk(Y,Z|T)$ together with a monotone class argument, which extends the measurability of
    $t \mapsto \Kk(E|T=t)$ from indicators to bounded measurable integrands. Integrating \eqref{eq:supp-prior-witness} against $\Pk(T)$ and using that the integral of a
    bounded measurable function against a mixture is the mixture of its integrals, which for indicators is the
    definition of $\circ$ and in general follows by a monotone class argument, we get:
    \begin{align*}
        \Pk(X \in A, Y \in B, Z \in C)
            &= \int_\Tcal \varphi(t)\, \Pk(T \in dt) \\
            &= \int_{B \times C} \Qk(X \in A|Z=z)\, \Pk(Y \in dy, Z \in dz),
    \end{align*}
    where $\Pk(Y,Z) = \Kk(Y,Z|T) \circ \Pk(T)$. The measurable rectangles form a $\pi$-system generating
    $\Bcal_\Xcal \otimes \Bcal_\Ycal \otimes \Bcal_\Zcal$ and both sides are probability measures, so Dynkin's
    lemma gives $\Pk(X,Y,Z) = \Qk(X|Z) \otimes \Pk(Y,Z)$ and, taking $B := \Ycal$,
    $\Pk(X,Z) = \Qk(X|Z) \otimes \Pk(Z)$.

    \emph{Step 2: identification of the two conditionals.} By hypothesis 2.\ and Step 1 the Markov kernels
    $\Pk(X|Z)$ and $\Qk(X|Z)$ are both conditional Markov kernels of $\Pk(X,Z)$ given $Z$. By essential
    uniqueness, \Cref{lem:ess-unique} applied with $\Asterisk$ in the input slot, we get for every fixed
    measurable $A \ins \Xcal$ a $\Pk(Z)$-null set $N_A \ins \Zcal$ with:
    \[ \Pk(X \in A|Z=z) = \Qk(X \in A|Z=z) \qquad \text{ for every } z \notin N_A. \]
    Only the single event $A$ is involved, so no countable generation of $\Bcal_\Xcal$ is needed.

    \emph{Step 3: transfer to every input value.} Fix $t \in \Tcal$ and measurable $A \ins \Xcal$,
    $B \ins \Ycal$, $C \ins \Zcal$. By hypothesis 1.\ the set $N_A$ is $\Kk(Z|T=t)$-null, hence the cylinder
    $\Ycal \times N_A$ is $\Kk(Y,Z|T=t)$-null, since $\Kk(Z|T=t)$ is the $\Zcal$-marginal of $\Kk(Y,Z|T=t)$.
    Therefore:
    \begin{align*}
        \lp \Pk(X|Z) \otimes \Kk(Y,Z|T=t) \rp &(A \times B \times C) \\
            &= \int_{B \times C} \Pk(X \in A|Z=z)\, \Kk(Y \in dy, Z \in dz|T=t) \\
            &= \int_{B \times C} \Qk(X \in A|Z=z)\, \Kk(Y \in dy, Z \in dz|T=t) \\
            &= \Kk(X \in A, Y \in B, Z \in C|T=t),
    \end{align*}
    by \eqref{eq:supp-prior-witness}. Dynkin's lemma again gives
    $\Kk(X,Y,Z|T=t) = \Pk(X|Z) \otimes \Kk(Y,Z|T=t)$, for every $t \in \Tcal$.
\end{proof}
\end{Thm}

We next restate the lemma of \cite{HS49} that produces a qualifying prior, see \Cref{rem:halmos-savage}. It is
stated there for families of finite measures; the proof below is the classical one, and we record that it needs no finiteness assumption on the
members of the family, since the argument only ever compares their null sets.

\begin{Lem}[Halmos-Savage, see \cite{HS49} Lemma 7]
    \label{lem:supp:halmos-savage}
    Let $\lp \Wcal, \Bcal_\Wcal, \mubf \rp$ be a $\sigma$-finite measure space and let $\Mcal \neq \emptyset$ be
    any family of measures $\nubf$ on $\lp \Wcal, \Bcal_\Wcal \rp$ with $\nubf \ll \mubf$. Then there exist a
    countable subfamily $\lC \nubf_1, \nubf_2, \dots \rC \ins \Mcal$ and weights $p_i > 0$ with
    $\sum_i p_i = 1$ such that:
    \[ \rhobf := \sum_i \nubf_i \cdot p_i \ll \mubf \qquad \text{ and } \qquad
       \nubf \ll \rhobf \quad \text{ for every } \nubf \in \Mcal. \]
    Note that $\rhobf$ need not be $\sigma$-finite in general, so $\nubf \ll \rhobf$ need not provide a density
    of $\nubf$ w.r.t.\ $\rhobf$; this can happen already for a family of finite measures, since the weights are
    not adapted to their masses. If $\Mcal$ consists of probability measures, as it does in every application
    below, then $\rhobf$ is a probability measure and the difficulty disappears.
\begin{proof}
    Since $\mubf$ is $\sigma$-finite it is equivalent to a \emph{finite} measure on $\Bcal_\Wcal$ (write
    $\Wcal = \biguplus_n \Wcal_n$ with $\mubf(\Wcal_n)<\infty$ and replace $\mubf$ by
    $\sum_n 2^{-n}\mubf(\,\cdot \cap \Wcal_n)/(1+\mubf(\Wcal_n))$), and absolute continuity w.r.t.\ two
    equivalent measures is the same relation, so we may and do assume that $\mubf$ is finite. (We cannot in
    general normalize to a probability measure, since $\mubf$ may be the zero measure.)

    \emph{Step 1: supports.} Every measure $\nubf \ll \mubf$ admits a \emph{support}, i.e.\ a set
    $S(\nubf) \in \Bcal_\Wcal$ with:
    \[ \nubf(A) = 0 \quad \iff \quad \mubf\lp A \cap S(\nubf) \rp = 0
       \qquad \text{ for all } A \in \Bcal_\Wcal. \]
    Indeed, put $s := \sup \lC \mubf(A) \mid A \in \Bcal_\Wcal, \, \nubf(A) = 0 \rC$ and choose $A_n$ with
    $\nubf(A_n)=0$ and $\mubf(A_n) \to s$. Then $A_\infty := \bigcup_n A_n$ satisfies $\nubf(A_\infty)=0$ by
    countable subadditivity and $\mubf(A_\infty)=s$; put $S(\nubf) := A_\infty^\cmpl$. If $\nubf(A)=0$ then
    $\nubf(A \cup A_\infty)=0$, so $\mubf(A \cup A_\infty) \le s = \mubf(A_\infty)$ and therefore
    $\mubf(A \cap S(\nubf)) = \mubf(A \sm A_\infty) = 0$. Conversely, if $\mubf(A \cap S(\nubf))=0$ then
    $\nubf(A \cap S(\nubf))=0$ because $\nubf \ll \mubf$, while $\nubf(A \cap A_\infty) \le \nubf(A_\infty)=0$;
    adding the two gives $\nubf(A)=0$. Fix one such $S(\nubf)$ for every $\nubf \in \Mcal$.

    \emph{Step 2: a maximal countable subfamily.} For a countable $\Dcal \ins \Mcal$ we put
    $S(\Dcal) := \bigcup_{\nubf \in \Dcal} S(\nubf) \in \Bcal_\Wcal$, and:
    \[ s^\ast := \sup \lC \mubf\lp S(\Dcal) \rp \mid \Dcal \ins \Mcal
       \text{ countable, non-empty} \rC \le \mubf(\Wcal) < \infty. \]
    Choose countable non-empty $\Dcal_n \ins \Mcal$ with $\mubf(S(\Dcal_n)) \to s^\ast$ and put
    $\Ccal := \bigcup_n \Dcal_n$, which is countable and non-empty and satisfies
    $\mubf(S(\Ccal)) \ge \mubf(S(\Dcal_n))$ for every $n$, hence $\mubf(S(\Ccal)) = s^\ast$. Enumerate
    $\Ccal = \lC \nubf_1, \nubf_2, \dots \rC$, choose any weights $p_i>0$ with $\sum_i p_i = 1$ and put
    $\rhobf := \sum_i \nubf_i \cdot p_i$.

    \emph{Step 3: conclusion.} If $\mubf(A)=0$ then $\nubf_i(A)=0$ for every $i$ and hence $\rhobf(A)=0$, so
    $\rhobf \ll \mubf$. Since all weights are strictly positive, $\rhobf(A)=0$ holds if and only if
    $\nubf_i(A)=0$ for every $i$, which by Step 1 holds if and only if $\mubf(A \cap S(\nubf_i))=0$ for every $i$,
    i.e.\ if and only if $\mubf(A \cap S(\Ccal))=0$. Now let $\nubf \in \Mcal$ be arbitrary. Applying the
    definition of $s^\ast$ to the countable family $\Ccal \cup \lC \nubf \rC$ gives
    $\mubf(S(\Ccal) \cup S(\nubf)) \le s^\ast = \mubf(S(\Ccal))$ and therefore
    $\mubf(S(\nubf) \sm S(\Ccal)) = 0$. So if $\rhobf(A)=0$, i.e.\ $\mubf(A \cap S(\Ccal))=0$, then:
    \[ \mubf\lp A \cap S(\nubf) \rp \;\le\; \mubf\lp A \cap S(\Ccal) \rp
       + \mubf\lp S(\nubf) \sm S(\Ccal) \rp \;=\; 0, \]
    and thus $\nubf(A)=0$ by Step 1. This is $\nubf \ll \rhobf$.
\end{proof}
\end{Lem}

\begin{Cor}[A witness from a prior in terms of densities]
    \label{cor:supp:tci-density}
    Let $\mubf_\Xcal$, $\mubf_\Ycal$, $\mubf_\Zcal$ be $\sigma$-finite measures on
    $\Xcal$, $\Ycal$, $\Zcal$, put $\mubf := \mubf_\Xcal \otimes \mubf_\Ycal \otimes \mubf_\Zcal$ and assume that
    $\Kk(X,Y,Z|T=t) \ll \mubf$ for every $t \in \Tcal$. For every $t \in \Tcal$ fix a density $k(x,y,z|t)$ of
    $\Kk(X,Y,Z|T=t)$ w.r.t.\ $\mubf$, with marginal densities $k(y,z|t)$ and $k(z|t)$ formed as in
    \Cref{rem:supp:marginal-densities}; no joint measurability in $t$ is required. Let $\Pk(T)$ be a probability
    measure on $\Tcal$ such that $\Kk(Y,Z|T=t) \ll \Pk(Y,Z)$ for every $t \in \Tcal$, where we have put
    $\Pk(X,Y,Z) := \Kk(X,Y,Z|T) \circ \Pk(T)$, and let $p(x,y,z)$ be a density of $\Pk(X,Y,Z)$ w.r.t.\ $\mubf$,
    with marginal densities $p(x,z)$ and $p(z)$ formed as in \Cref{rem:supp:marginal-densities}. Fix a
    probability density $w$ w.r.t.\ $\mubf_\Xcal$, which exists by \Cref{rem:supp:marginal-densities}, and put:
    \[ p(x|z) := \begin{cases} p(x,z) \big/ p(z), & \text{ if } 0 < p(z) < \infty,\\
                               w(x), & \text{ otherwise,}\end{cases} \]
    \[ \Pk(X \in A|Z=z) := \int_A p(x|z)\, \mubf_\Xcal(dx). \]
    Then $\Pk(X|Z):\, \Zcal \dshto \Xcal$ is a Markov kernel with $\Pk(X,Z) = \Pk(X|Z) \otimes \Pk(Z)$ and the
    following are equivalent:
    \begin{enumerate}
        \item $\Xk \Indep_\Kk \Yk \given \Zk$;
        \item $\Kk(X,Y,Z|T) = \Pk(X|Z) \otimes \Kk(Y,Z|T)$;
        \item for every $t \in \Tcal$: $k(x,y,z|t) \cdot p(z) = p(x,z) \cdot k(y,z|t)$ for $\mubf$-almost all
            $(x,y,z)$, where the exceptional set may depend on $t$;
        \item $X \Indep^\omega_{\Pk(X,Y,Z)} Y \given Z$ and $\Xk \Indep_\Kk \Tk \given \Yk \otimes \Zk$.
    \end{enumerate}
\begin{proof}
    \emph{Step 0: the mixture and its densities.} If $\mubf(E)=0$ then $\Kk(E|T=t)=0$ for every $t \in \Tcal$,
    hence $\Pk(E) = \int_\Tcal \Kk(E|T=t)\, \Pk(T \in dt) = 0$; so $\Pk(X,Y,Z) \ll \mubf$ and a density
    $p(x,y,z)$ exists by the Radon-Nikodym theorem, $\mubf$ being $\sigma$-finite and $\Pk(X,Y,Z)$ a probability
    measure. \Cref{rem:supp:marginal-densities}, applied to $\Pk(X,Y,Z)$, provides
    the jointly measurable marginal densities $p(x,z)$ and $p(z)$ of $\Pk(X,Z)$ and $\Pk(Z)$ and
    gives $p(z) < \infty$ for $\mubf_\Zcal$-almost all $z$ as well as
    $\int_\Xcal p(x,z)\,\mubf_\Xcal(dx) = p(z)$ for every $z$. The map $p(x|z)$ is jointly measurable, since
    $\lC 0 < p(z) < \infty \rC$ is a measurable set, and $\int_\Xcal p(x|z)\,\mubf_\Xcal(dx) = 1$ for every
    $z \in \Zcal$: for $z$ with $0 < p(z)<\infty$ this is $p(z)/p(z)$, and otherwise it holds by the choice of
    $w$. So $\Pk(X \in A|Z=z)$ is a probability measure in $A$ for every $z$, by monotone convergence, and
    measurable in $z$ for every $A$, by Tonelli's theorem; that is, $\Pk(X|Z)$ is a Markov kernel. Note that
    nothing so far presupposes any conditional independence, and that altering the chosen densities on
    $\mubf$-null sets alters $p(x|z)$ and $\Pk(X|Z)$ only on null sets of the corresponding reference measures,
    so that none of 1.--4.\ is affected.

    \emph{Step 1: $\Pk(X,Z) = \Pk(X|Z) \otimes \Pk(Z)$.} Let $A \ins \Xcal$ and $C \ins \Zcal$ be measurable and
    split $\Zcal$ into $N_0 := \lC p(z)=0\rC$, $N_\infty := \lC p(z)=\infty \rC$ and
    $G := \lC 0 < p(z)<\infty\rC$. The set $N_\infty$ is $\mubf_\Zcal$-null, hence also $\Pk(Z)$-null, and
    contributes nothing to either side. For $z \in N_0$ we have
    $\int_\Xcal p(x,z)\,\mubf_\Xcal(dx) = p(z) = 0$, so $p(x,z)=0$ for $\mubf_\Xcal$-almost all $x$, and both
    $\Pk(X \in A, Z \in C \cap N_0)$ and $\int_{C \cap N_0} \Pk(X \in A|Z=z)\,p(z)\,\mubf_\Zcal(dz)$ vanish. For
    $z \in G$ we have $\Pk(X \in A|Z=z)\cdot p(z) = \int_A p(x,z)\,\mubf_\Xcal(dx)$ by the definition of
    $p(x|z)$. Integrating over $C \cap G$ and adding up the three parts gives
    $\lp \Pk(X|Z) \otimes \Pk(Z) \rp (A \times C) = \Pk(X \in A, Z \in C)$. The measurable rectangles form a
    $\pi$-system generating $\Bcal_\Xcal \otimes \Bcal_\Zcal$ and both sides are probability measures, so
    Dynkin's lemma gives the claim.

    \emph{Step 2: absolute continuity of the marginals.} Absolute continuity passes to marginals, since the
    preimage of a null set under a projection is a null cylinder, so the hypothesis
    $\Kk(Y,Z|T=t) \ll \Pk(Y,Z)$ gives $\Kk(Z|T=t) \ll \Pk(Z)$, for every $t \in \Tcal$.

    \emph{1.\ $\iff$ 2.} Steps 1 and 2 are precisely the two hypotheses of \Cref{thm:supp:tci-prior} for the
    prior $\Pk(T)$, so that theorem gives the equivalence.

    \emph{Step 3: mixing.} If $\Rk(X|Z)$ is any Markov kernel with
    $\Kk(X,Y,Z|T) = \Rk(X|Z) \otimes \Kk(Y,Z|T)$, then $\Pk(X,Y,Z) = \Rk(X|Z) \otimes \Pk(Y,Z)$, by Step 1 of the
    proof of \Cref{thm:supp:tci-prior}. That computation goes through unchanged for a Markov kernel
    $\Rk(X|Y,Z)$ in place of $\Rk(X|Z)$, since $(y,z) \mapsto \Rk(X \in A|Y=y,Z=z)$ is again measurable and
    bounded by $1$, and yields $\Pk(X,Y,Z) = \Rk(X|Y,Z) \otimes \Pk(Y,Z)$; only its last step, the
    marginalization of $\Ycal$, is special to $\Rk(X|Z)$ and is not needed here.

    \emph{2.\ $\iff$ 3.} Fix $t \in \Tcal$. By Tonelli's theorem the measure
    $\Pk(X|Z) \otimes \Kk(Y,Z|T=t)$ has density $p(x|z) \cdot k(y,z|t)$ w.r.t.\ $\mubf$, and it is a probability
    measure, as is $\Kk(X,Y,Z|T=t)$ with density $k(x,y,z|t)$. Since $\mubf$ is $\sigma$-finite, 2.\ holds for
    this $t$ if and only if:
    \begin{align}
        k(x,y,z|t) = p(x|z) \cdot k(y,z|t) \qquad \mubf\text{-almost everywhere.} \label{eq:supp-hs-star}
    \end{align}
    On $G$ we have $p(x|z) = p(x,z)/p(z)$ with $p(z) \in (0,\infty)$, so there \eqref{eq:supp-hs-star} and 3.\
    are obtained from one another by multiplying, respectively dividing, by $p(z)$. On $\lC p(z)=\infty\rC$ both
    are silent, the set being $\mubf$-null. On $N_0 = \lC p(z)=0 \rC$ both hold automatically. Indeed, $p(z)$ is
    a density of $\Pk(Z)$ w.r.t.\ $\mubf_\Zcal$, so $N_0$ is $\Pk(Z)$-null and hence, by Step 2,
    $\Kk(Z|T=t)$-null; that is, $k(z|t)=0$ for $\mubf_\Zcal$-almost all $z \in N_0$. Since
    $k(z|t) = \int_{\Xcal \times \Ycal} k(x,y,z|t)\, \lp\mubf_\Xcal \otimes \mubf_\Ycal\rp(d(x,y))$ and
    $k(y,z|t) = \int_\Xcal k(x,y,z|t)\, \mubf_\Xcal(dx)$, Tonelli's theorem turns this into $k(x,y,z|t)=0$ for
    $\mubf$-almost all $(x,y,z)$ with $z \in N_0$ and $k(y,z|t)=0$ for
    $\mubf_\Ycal \otimes \mubf_\Zcal$-almost all $(y,z)$ with $z \in N_0$; and $p(x,z)=0$ for
    $\mubf_\Xcal \otimes \mubf_\Zcal$-almost all such $(x,z)$, as in Step 1. So over $N_0$ the identity
    \eqref{eq:supp-hs-star} reads $0 = w(x) \cdot 0$ and both sides of 3.\ vanish, $\mubf$-almost everywhere.
    Note that only the hypothesis $\Kk(Y,Z|T=t) \ll \Pk(Y,Z)$, through Step 2, has been used.

    \emph{2.\ $\implies$ 4.} By Step 3, $\Pk(X,Y,Z) = \Pk(X|Z) \otimes \Pk(Y,Z)$, so $\Pk(X \in A|Z)$ is a
    version of $\E_\Pk\lB \I_A(X) \,\middle|\, Y,Z \rB$; being $\sigma(Z)$-measurable it is also a version of
    $\E_\Pk\lB \I_A(X) \,\middle|\, Z \rB$, by Step 1 and the defining property of conditional expectation. This
    is $X \Indep^\omega_{\Pk(X,Y,Z)} Y \given Z$. Reading $\Pk(X|Z)$ as a Markov kernel
    $\Ycal \times \Zcal \dshto \Xcal$ that ignores its first argument, 2.\ also exhibits a conditional Markov
    kernel of $\Kk(X,Y,Z|T)$ given $(Y,Z)$ that does not depend on the input, which is
    $\Xk \Indep_\Kk \deltabf_\ast \given \Yk \otimes \Zk$ and hence
    $\Xk \Indep_\Kk \Tk \given \Yk \otimes \Zk$ by \Cref{add-T}.

    \emph{4.\ $\implies$ 2.} By \Cref{add-T} the second half of 4.\ provides a Markov kernel $\Rk(X|Y,Z)$ with
    $\Kk(X,Y,Z|T=t) = \Rk(X|Y,Z) \otimes \Kk(Y,Z|T=t)$ for every $t \in \Tcal$. By Step 3,
    $\Pk(X,Y,Z) = \Rk(X|Y,Z) \otimes \Pk(Y,Z)$, so for a fixed measurable $A \ins \Xcal$ the function
    $\Rk(X \in A|Y=y,Z=z)$ is a version of $\E_\Pk\lB \I_A(X) \,\middle|\, Y,Z \rB$, while $\Pk(X \in A|Z=z)$ is
    a version of $\E_\Pk\lB \I_A(X) \,\middle|\, Z \rB$ by Step 1. The first half of 4.\ says that these two
    conditional expectations agree $\Pk(X,Y,Z)$-almost surely; both are functions of $(y,z)$ alone, so
    $\Rk(X \in A|Y=y,Z=z) = \Pk(X \in A|Z=z)$ for $\Pk(Y,Z)$-almost all $(y,z)$. That exceptional set is a
    $\Pk(Y,Z)$-null subset of $\Ycal \times \Zcal$ and need not be a cylinder over $\Zcal$, so here the
    hypothesis $\Kk(Y,Z|T=t) \ll \Pk(Y,Z)$ is used at full strength, rather than its marginal consequence from
    Step 2: the exceptional set is $\Kk(Y,Z|T=t)$-null for every $t \in \Tcal$, so for all measurable
    $B \ins \Ycal$, $C \ins \Zcal$ and all $t \in \Tcal$:
    \begin{align*}
        \Kk(X \in A, &Y \in B, Z \in C|T=t) \\
            &= \int_{B \times C} \Rk(X \in A|Y=y,Z=z)\, \Kk(Y \in dy, Z \in dz|T=t) \\
            &= \lp \Pk(X|Z) \otimes \Kk(Y,Z|T=t) \rp(A \times B \times C).
    \end{align*}
    Again only one event $A$ is involved at a time, and Dynkin's lemma gives 2.
\end{proof}
\end{Cor}

We now turn to the linear Gaussian case. Recall that for a positive definite block matrix $M$ with blocks
$\Sigma_{YY}$, $\Sigma_{YZ}$, $\Sigma_{ZY}$, $\Sigma_{ZZ}$ and Schur complement
$S := \Sigma_{YY} - \Sigma_{YZ}\Sigma_{ZZ}^{-1}\Sigma_{ZY}$ of $\Sigma_{ZZ}$ we have:
\[ M^{-1} = \begin{bmatrix} S^{-1} & -S^{-1}\Sigma_{YZ}\Sigma_{ZZ}^{-1} \\
    -\Sigma_{ZZ}^{-1}\Sigma_{ZY}S^{-1} &
    \Sigma_{ZZ}^{-1} + \Sigma_{ZZ}^{-1}\Sigma_{ZY}S^{-1}\Sigma_{YZ}\Sigma_{ZZ}^{-1} \end{bmatrix}. \]

\begin{Thm}[Linear Gaussian Markov kernels]
    \label{thm:supp:tci-gaussian}
    Let $\Tcal = \R^{d_T}$, $\Xcal = \R^{d_X}$, $\Ycal = \R^{d_Y}$, $\Zcal = \R^{d_Z}$, let
    $\Gamma_X \in \R^{d_X \times d_T}$, $\Gamma_Y \in \R^{d_Y \times d_T}$, $\Gamma_Z \in \R^{d_Z \times d_T}$ and
    $\gamma_X \in \R^{d_X}$, $\gamma_Y \in \R^{d_Y}$, $\gamma_Z \in \R^{d_Z}$, and let $\Pk(X,Y,Z|T)$ be the linear
    Gaussian Markov kernel:
    \[ \Pk(X,Y,Z|T=t) = \Ncal \lp
        \begin{bmatrix} \Gamma_X \\ \Gamma_Y \\ \Gamma_Z \end{bmatrix} t +
        \begin{bmatrix} \gamma_X \\ \gamma_Y \\ \gamma_Z \end{bmatrix} ,\;
        \begin{bmatrix} \Sigma_{XX} & \Sigma_{XY} & \Sigma_{XZ} \\
                        \Sigma_{YX} & \Sigma_{YY} & \Sigma_{YZ} \\
                        \Sigma_{ZX} & \Sigma_{ZY} & \Sigma_{ZZ} \end{bmatrix} \rp,  \]
    with a positive definite covariance matrix $\Sigma$. Then:
    \[ X \Indep_{\Pk(X,Y,Z|T)} Y \given Z \qquad \iff \qquad
       \Sigma_{XY} = \Sigma_{XZ}\Sigma_{ZZ}^{-1}\Sigma_{ZY} \quad \land \quad
       \Gamma_X = \Sigma_{XZ}\Sigma_{ZZ}^{-1}\Gamma_Z, \]
    and in that case a witness is given by the linear Gaussian Markov kernel:
    \begin{align*}
        \Qk(X|Z=z) &= \Pk(X|Z=z,T=0) \\
        &= \Ncal \lp \gamma_X + \Sigma_{XZ}\Sigma_{ZZ}^{-1}(z-\gamma_Z),\;
          \Sigma_{XX} - \Sigma_{XZ}\Sigma_{ZZ}^{-1}\Sigma_{ZX} \rp.
    \end{align*}
\begin{proof}
    Write $\mu_X(t) := \Gamma_X t + \gamma_X$ and similarly $\mu_Y(t)$, $\mu_Z(t)$, and put:
    \[ M := \begin{bmatrix} \Sigma_{YY} & \Sigma_{YZ} \\ \Sigma_{ZY} & \Sigma_{ZZ}\end{bmatrix}, \qquad
       S := \Sigma_{YY} - \Sigma_{YZ}\Sigma_{ZZ}^{-1}\Sigma_{ZY}, \qquad
       D := \Sigma_{XY} - \Sigma_{XZ}\Sigma_{ZZ}^{-1}\Sigma_{ZY}. \]
    As principal submatrices of the positive definite $\Sigma$ both $M$ and $\Sigma_{ZZ}$ are positive definite,
    and hence so is the Schur complement $S$. In particular $M$, $\Sigma_{ZZ}$ and $S$ are invertible.

    \emph{Step 1: the conditional distribution of $X$ given $(Y,Z)$.} For each fixed $t \in \Tcal$ the vector
    $(X,Y,Z)$ is Gaussian under $\Pk(X,Y,Z|T=t)$ with mean $(\mu_X(t),\mu_Y(t),\mu_Z(t))$ and covariance $\Sigma$.
    Since $M$ is invertible the usual Gaussian conditioning formula provides a regular conditional distribution:
    \[ \Pk(X|Y=y,Z=z,T=t) = \Ncal \lp m_t(y,z),\, V \rp, \]
    with:
    \begin{align*}
        m_t(y,z) &:= \mu_X(t) + \begin{bmatrix} \Sigma_{XY} & \Sigma_{XZ}\end{bmatrix} M^{-1}
            \lp \begin{bmatrix} y \\ z \end{bmatrix} -
                \begin{bmatrix} \mu_Y(t) \\ \mu_Z(t) \end{bmatrix} \rp, \\
        V &:= \Sigma_{XX} - \begin{bmatrix} \Sigma_{XY} & \Sigma_{XZ}\end{bmatrix} M^{-1}
            \begin{bmatrix} \Sigma_{YX} \\ \Sigma_{ZX} \end{bmatrix}.
    \end{align*}
    Note that the right hand side is defined, and is a Markov kernel in $(y,z)$, for \emph{all} $(y,z)$, not
    merely almost all.

    \emph{Step 2: block inversion.} Inserting the displayed formula for $M^{-1}$ and multiplying out gives:
    \[ \begin{bmatrix} \Sigma_{XY} & \Sigma_{XZ}\end{bmatrix} M^{-1}
       = \begin{bmatrix} A & B \end{bmatrix}, \qquad
       A := D S^{-1}, \qquad
       B := \Sigma_{XZ}\Sigma_{ZZ}^{-1} - D S^{-1} \Sigma_{YZ}\Sigma_{ZZ}^{-1}. \]
    Substituting this into $m_t$ and collecting the terms containing $t$ we obtain, for all $y,z,t$:
    \begin{align}
        m_t(y,z) = \gamma_X + A\,(y - \gamma_Y) + B\,(z-\gamma_Z)
                   + \lp \Gamma_X - A\Gamma_Y - B\Gamma_Z \rp t. \label{eq:supp-gauss-mean}
    \end{align}

    \emph{Step 3: ``$\Longleftarrow$''.} Assume both conditions. The first says $D = 0$, whence $A = 0$ and
    $B = \Sigma_{XZ}\Sigma_{ZZ}^{-1}$; the second then says $\Gamma_X - A\Gamma_Y - B\Gamma_Z
    = \Gamma_X - \Sigma_{XZ}\Sigma_{ZZ}^{-1}\Gamma_Z = 0$. So \eqref{eq:supp-gauss-mean} reduces to
    $m_t(y,z) = \gamma_X + \Sigma_{XZ}\Sigma_{ZZ}^{-1}(z - \gamma_Z)$, which is free of $y$ and $t$, and
    $V = \Sigma_{XX} - A\Sigma_{YX} - B\Sigma_{ZX} = \Sigma_{XX} - \Sigma_{XZ}\Sigma_{ZZ}^{-1}\Sigma_{ZX}$, which
    is free of $(y,z,t)$ anyway. Defining the Markov kernel:
    \[ \Qk(X|Z=z) := \Ncal \lp \gamma_X + \Sigma_{XZ}\Sigma_{ZZ}^{-1}(z-\gamma_Z),\;
        \Sigma_{XX} - \Sigma_{XZ}\Sigma_{ZZ}^{-1}\Sigma_{ZX} \rp, \]
    Step 1 thus reads $\Pk(X|Y=y,Z=z,T=t) = \Qk(X|Z=z)$ for all $(y,z,t)$, and therefore for every $t \in \Tcal$:
    \[ \Pk(X,Y,Z|T=t) = \Pk(X|Y,Z,T=t) \otimes \Pk(Y,Z|T=t) = \Qk(X|Z) \otimes \Pk(Y,Z|T=t), \]
    which is $X \Indep_{\Pk(X,Y,Z|T)} Y \given Z$. Finally, applying the Gaussian conditioning formula to the pair
    $(X,Z)$ under $\Pk(X,Y,Z|T=0)$ gives exactly:
    \begin{align*}
        \Pk(X|Z=z,T=0) &= \Ncal\lp \gamma_X + \Sigma_{XZ}\Sigma_{ZZ}^{-1}(z-\gamma_Z),\,
        \Sigma_{XX} - \Sigma_{XZ}\Sigma_{ZZ}^{-1}\Sigma_{ZX}\rp \\
        &= \Qk(X|Z=z).
    \end{align*}

    \emph{Step 4: ``$\implies$''.} Assume $X \Indep_{\Pk(X,Y,Z|T)} Y \given Z$ with witness $\Qk(X|Z)$. For each
    $t \in \Tcal$ the identity $\Pk(X,Y,Z|T=t) = \Qk(X|Z) \otimes \Pk(Y,Z|T=t)$ says precisely that
    $(y,z) \mapsto \Qk(X|Z=z)$ is a regular conditional distribution of $X$ given $(Y,Z)$ under
    $\Pk(X,Y,Z|T=t)$; by Step 1 so is $(y,z) \mapsto \Ncal(m_t(y,z),V)$. Since $\Xcal = \R^{d_X}$ is standard,
    two such versions agree for $\Pk(Y,Z|T=t)$-almost all $(y,z)$, see \Cref{ess-unique}. As $M$ is positive
    definite, $\Pk(Y,Z|T=t)$ is a non-degenerate Gaussian and hence equivalent to the Lebesgue measure $\lambda$
    on $\R^{d_Y} \times \R^{d_Z}$. So for every $t \in \Tcal$ there is a set $E_t \ins \R^{d_Y} \times \R^{d_Z}$
    with $\lambda\lp E_t^c \rp = 0$ such that:
    \[ \Qk(X|Z=z) = \Ncal\lp m_t(y,z), V \rp \qquad \text{ for all } (y,z) \in E_t. \]
    Comparing the means of two normal distributions with the same covariance $V$, this gives in particular that
    $m_t(y,z)$ equals the mean vector of $\Qk(X|Z=z)$ for all $(y,z) \in E_t$.

    First we show $A = 0$. By Fubini's theorem applied to $\I_{E_0^c}$ there is a $z \in \R^{d_Z}$ whose section
    $(E_0)_z := \lC y \mid (y,z) \in E_0 \rC$ has full Lebesgue measure in $\R^{d_Y}$. For $y, y' \in (E_0)_z$ both
    $m_0(y,z)$ and $m_0(y',z)$ equal the mean of $\Qk(X|Z=z)$, so $A\,(y-y') = 0$ by
    \eqref{eq:supp-gauss-mean}. Since $(E_0)_z$ has full measure, for every $v \in \R^{d_Y}$ the set
    $(E_0)_z \cap \lp (E_0)_z + v \rp$ has full measure and is in particular non-empty, i.e.\
    $(E_0)_z - (E_0)_z = \R^{d_Y}$. Hence $A = 0$, i.e.\ $D S^{-1} = 0$, and as $S$ is invertible $D = 0$, which is
    the first condition $\Sigma_{XY} = \Sigma_{XZ}\Sigma_{ZZ}^{-1}\Sigma_{ZY}$. Consequently
    $B = \Sigma_{XZ}\Sigma_{ZZ}^{-1}$.

    Now let $t \in \Tcal$ be arbitrary. The set $E_0 \cap E_t$ has full Lebesgue measure and is in particular
    non-empty; picking $(y,z)$ in it, both $m_0(y,z)$ and $m_t(y,z)$ equal the mean of $\Qk(X|Z=z)$, so
    \eqref{eq:supp-gauss-mean} gives $\lp \Gamma_X - A\Gamma_Y - B\Gamma_Z \rp t = 0$. With $A = 0$ and
    $B = \Sigma_{XZ}\Sigma_{ZZ}^{-1}$ this reads
    $\lp \Gamma_X - \Sigma_{XZ}\Sigma_{ZZ}^{-1}\Gamma_Z \rp t = 0$, and since $t \in \R^{d_T}$ was arbitrary we
    get the second condition $\Gamma_X = \Sigma_{XZ}\Sigma_{ZZ}^{-1}\Gamma_Z$.
\end{proof}
\end{Thm}

\subsection{The Extended Conditional Independence}

We shortly review the definition of \emph{extended conditional independence} introduced in \cite{CD17}.

\begin{Def}[Extended conditional independence, see \cite{CD17} Def. 3.2]
    \label{def:ext-ci}
    Let $\Wcal$ and $\Tcal$ be measurable spaces
    and $\Ecal=(\Pk_t(W))_{t \in \Tcal}$ be a family of probability measures on $\Wcal$.
    Let $X,Y,Z$ be measurable maps on $\Wcal$ and $\Phi,\Theta$ measurable maps on $\Tcal$ such that
    the joint map $(\Phi,\Theta)$ is injective. For these cases \emph{extended conditional independence} was defined as:
    \[ X \Indep_\Ecal (Y,\Theta) \given (Z,\Phi)  \]
    if for all $\phi \in \Phi(\Tcal)$ and all real bounded measurable $h$
    there exists a function $g_{h,\phi}$ such that for all $t \in \Phi^{-1}(\phi)$ we have that:
    \[ \E_t\lB h(X) |Y,Z] \rB = g_{h,\phi}(Z)\quad \Pk_t(W)\text{-a.s.},  \]
    where the conditional expectation $\E_t$ is w.r.t.\ $\Pk_t(W)$.
\end{Def}

We now show that when $(Y,\Theta)$ and $(Z,\Phi)$ are considered as transitional random variables
on transition probability space $\lp \Wcal \times \Tcal, \Pk(W|T) \rp$, where we put $\Pk(W|T=t):=\Pk_t(W)$,
then transitional conditional independence implies extended conditional independence.

\begin{Lem}
    \label{lem:tci-implies-eci}
We have the implication:
\[ X \Indep_{\Pk(W|T)} (Y,\Theta) \given (Z,\Phi) \qquad \implies \qquad X \Indep_\Ecal (Y,\Theta) \given (Z,\Phi).   \]
\begin{proof}
Indeed, by the assumption we get a Markov kernel $\Qk(X|Z,\Phi)$ such that:
\begin{align*}
    \Pk(X,Y,\Theta,Z,\Phi|T) &= \Qk(X|Z,\Phi) \otimes \Pk(Y,\Theta,Z,\Phi|T) \\
                             &= \Qk(X|Z,\Phi) \otimes \deltabf(\Theta|T) \otimes \deltabf(\Phi|T) \otimes \Pk(Y,Z|T).
\end{align*}
Marginalizing out $\Theta$ and $\Phi$ gives:
   \[ \Pk(X,Y,Z|T=t) = \Qk(X|Z,\Phi=\Phi(t)) \otimes \Pk(Y,Z|T=t).\]
For any $\phi$ and function $h$ we then define:
\[ g_{h,\phi}(z) := \int_\Xcal h(x) \, \Qk(X \in dx|Z=z,\Phi=\phi).  \]
Then for each $t \in \Phi^{-1}(\phi)$ and $B \in \Bcal_\Ycal$ and $C \in \Bcal_\Zcal$ we get:
\begin{align*}
&  \int_C \int_B  \int_\Xcal h(x) \, \Pk(X \in dx, Y \in dy, Z \in dz|T=t) \\
&= \int_C \int_B  \int_\Xcal h(x) \, \Qk(X \in dx|Z=z,\Phi=\Phi(t)) \, \Pk(Y \in dy, Z \in dz|T=t) \\
&= \int_C \int_B  g_{h,\phi}(z)\, \Pk(Y \in dy, Z \in dz|T=t).
\end{align*}
Since this is the defining equation for the conditional expectation we get the claim:
    \[ \E_t\lB h(X) |Y,Z \rB = g_{h,\phi}(Z)\quad \Pk_t(W)\text{-a.s.}.  \]
\end{proof}
\end{Lem}

\begin{Rem}
    So transitional conditional independence is the stronger notion and implies extended conditional independence, 
    but it works for all transitional random variables, 
    not just of the restricted type in \Cref{def:ext-ci}, where the stochastic variables $X,Y,Z$ and the decision
    variables $\Phi,\Theta$ live on different spaces and the first argument has to be fully stochastic.
    Transitional conditional independence satisfies all the (asymmetric) separoid rules,
    \Cref{thm:separoid_axioms-tci}, for all measurable spaces,
    except Left Weak Union, $T$-Restricted Right Redundancy and $T$-Restricted Symmetry, which hold when one can ensure the existence of
    conditional Markov kernels, e.g.\ on standard measurable spaces, see \Cref{thm-regular-conditional-Markov-kernel}.
    For \cite{CD17} the corresponding situation is that all their axioms hold in general except the mirror weak
    union P4$''$, which they leave open and recover under a discrete regime space, discrete stochastic variables or
    a dominating regime, and except symmetry, which they have in a restricted form only; see
    \Cref{sec:main:comparison:eci} for the comparison and for how those hypotheses reappear here as criteria.
    This makes transitional conditional independence a better fit for the use in graphical models, 
    see, for instance, the global Markov property, 
    \Cref{thm-gmp-mI-CBN}. 
    Note that both notions only have a restricted direct relation to variation conditional independence, 
    see \Cref{thm:variation-ci} and \cite{CD17, CD17-supp}.
    A formal analogy between variation conditional independence and transitional conditional independence 
    was discussed in \Cref{rem:supp:var-tci} and \Cref{eg:supp:var-tci}.
\end{Rem}

\subsection{Symmetric Extended Conditional Independence}

If we wanted to arrive at a symmetric version of extended conditional independence that satisfies all (symmetric)
separoid rules (at least when restricted to codomains forming disintegration triples)
we could just use \emph{symmetrized transitional conditional independence}:
\[ \Xk \Indep_{\Kk(W|T)}^\lor \Yk \given \Zk \qquad :\iff \qquad 
\Xk \Indep_{\Kk(W|T)} \Yk \given \Zk \quad \lor \quad \Yk \Indep_{\Kk(W|T)} \Xk \given \Zk. \]
Since $\Indep_{\Kk(W|T)}$ forms a $T$-$\ast$-separoid it is immediate that $\Indep_{\Kk(W|T)}^\lor$ is a symmetric separoid
by the general theory of $\tau$-$\kappa$-separoids, see \Cref{thm:sym-t-k-separoid}.
We clearly have the implication:
\[ \Xk \Indep_{\Kk(W|T)} \Yk \given \Zk \qquad \implies \qquad 
\Xk \Indep_{\Kk(W|T)}^\lor \Yk \given \Zk, \]
showing that the asymmetric version is stronger than the symmetrized version, where the latter might have lost some information about 
the interplay between $\Xk$, $\Yk$, $\Zk$ and $\Tk$. Furthermore, it is not invariant under the equivalences that $\Indep_{\Kk(W|T)}$ itself enjoys, see \Cref{add-T}:
in the equivalent spelling without $\Tk$ in the second argument it is strictly weaker, and with $\deltabf_\ast$
there it is vacuous. So classical statistical concepts like ancillarity, sufficiency and adequacy are expressible
by it in one particular spelling only, see \Cref{eg:symmetrization-loses}.
A symmetric notion of extended conditional independence was introduced in \cite{RERS17}.
Also in \cite{Cho17, Fri20} a symmetric version of conditional independence for categorical probability theory was
proposed. The later \cite{FK23}, in contrast, uses an \emph{asymmetric} conditional independence for morphisms
with inputs, which the authors introduce as the categorical generalization of transitional conditional
independence, see \Cref{sec:main:comparison:cat}.

It is worth spelling out what a disjunctive symmetrization costs, since $\Indep^\lor$ is the obvious candidate for
a symmetric notion. Left Redundancy \ref{sep:tci:l-red} gives
$\deltabf_\ast \Indep_{\Kk(W|T)} \Yk \given \Zk$ for \emph{all} $\Yk$ and $\Zk$. Hence
\[ \Xk \Indep^\lor_{\Kk(W|T)} \deltabf_\ast \given \Zk \qquad \text{holds for all } \Xk, \Zk, \]
i.e.\ the symmetrized relation is \emph{vacuous} on every statement that has $\deltabf_\ast$ in the
second slot. But among these are the statements that carry much of the content of the theory: the existence of a
conditional Markov kernel, see \Cref{conditional-markov-kernel-as-ci}, and the invariance of a predictor across
environments, see \Cref{prp:invariance-tci}. More generally, $\Indep^\lor$ records that \emph{one} of the two
kernels exists without recording \emph{which}, so it cannot be used to produce the kernel one is after. This is the
price of symmetry, and it is not a defect of our particular symmetrization: any disjunctive definition pays it.
An explicit model in which the symmetrized relation holds while the statistical property it is meant to express
fails is given in \Cref{eg:symmetrization-loses}.

It is also worth mentioning that id-separation becomes symmetric as soon as one conditions on all input nodes, see $J$-Restricted Symmetry \ref{sep:sig:j-res-sym}. Together with the global Markov property, \Cref{thm-gmp-mI-CBN}, this immediately implies the following symmetrized version:
\[ A \Perp^{\id}_\Gk B \given C \cup J \qquad \implies \qquad X_A \Indep_{\Pk(X_V\Vert X_J)}^\lor X_B \given X_C, X_J.   \]
This recovers, generalizes and strengthens the corresponding results from \cite{RERS17} and \cite{FM20} Appendix C.

\subsection{Extended Conditional Independence for Families of Probability Distributions}

In this subsection we will introduce a strikingly simple and powerful form of extended conditional independence
that works for all measurable spaces and satisfies all the separoid rules.
For this consider a transition probability space  $\lp \Wcal \times \Tcal, \Kk(W|T) \rp$ and transitional random variables
$X:\, \Wcal \times \Tcal \to \Xcal$, $Y:\, \Wcal \times \Tcal \to \Ycal$, $Z:\,\Wcal \times \Tcal \to \Zcal$.
Furthermore, fix a set $\Qcal \ins \Pcal(\Tcal)$ of probability measures on $\Tcal$, e.g.\ $\Qcal= \Pcal(\Tcal)$ or 
$\Qcal=\lC\deltabf_t\,|\,t \in \Tcal \rC$.
Then we can define \emph{$\Qcal$-extended conditional independence} as:
\[ X \Indep_{\Kk(W|T) \otimes \Qcal}^\omega  Y \given Z \qquad :\iff \qquad \forall \Qk(T) \in \Qcal:\quad 
X \Indep_{\Kk(W|T) \otimes \Qk(T)}^\omega Y \given Z.  \]
It is easily seen that the usual weak conditional independence $\Indep^\omega$ satisfies all separoid axioms \cite{Daw01, CD17}
for arbitrary measurable spaces. Furthermore, if one combines several separoids by conjunction
$\land$ then one gets another separoid, see \cite{Daw01}. 
So $\Indep_{\Kk(W|T) \otimes \Qcal}^\omega $ clearly satisfies all (symmetric) separoid axioms for arbitrary measurable spaces.
By \Cref{thm:supp:ci-equivalences} we have the implications:
\[ X \Indep_{\Kk(W|T)} Y \given Z  \quad \implies \quad X \Indep_{\Kk(W|T) \otimes \Qcal}^\omega T,Y \given Z  
\quad \implies \quad X \Indep_{\Kk(W|T)  \otimes \Qcal}^\omega Y \given Z.\]
The middle ternary relation in $X,Y,Z$ satisfy the asymmetric separoid rules
from \Cref{thm:separoid_axioms-tci}, but without any requirement on the underlying measurable spaces,
in contrast to transitional conditional independence on the left. The asymmetric separoid rules for the middle relation follow from the right relation and \Cref{thm:plus-t-k-separoid}.

It seems that $\Qcal$-extended conditional independence checks all boxes that one would like to have from 
a notion of extended conditional independence. It is certainly simpler than most other notions.
It just comes with one drawback: it does not provide one with the existence or factorization of certain Markov kernels.
When the reverse implication holds is stated in \Cref{thm:supp:ci-equivalences}, e.g.\ if $\deltabf_t \in \Qcal$ for all $t\in \Tcal$
and $X \ismapof_\Kk^\ast Z$, where the latter encodes the existence of a certain Markov kernel, which is thus the main obstruction
to arrive at transitional conditional independence.

To elaborate further, a specialized version of this $\Qcal$-extended conditional independence was first introduced in \cite{FM20}, where it was used to
derive the causal do-calculus rules, see \cite{Pearl09, FM20}, for certain structural causal models.
In their proofs they had to construct certain Markov kernels and then check for $\Qcal$-extended conditional independence. 
Since the construction of such Markov kernels became complicated many corner cases have not been proved.
The main ingredient of their proof was a \emph{global Markov property} for $\Qcal$-extended conditional independence:
\[ A \Perp^{\id}_\Gk J \cup B \given C \qquad \implies \qquad X_A \Indep_{\Pk(X_V\Vert X_J)\otimes \Qcal }^\omega X_J,X_B \given X_C.   \]
Here, $\Qcal$-extended conditional independence was not strong enough to produce the needed Markov kernels. 
This is in contrast to transitional conditional independence,
whose global Markov property, \Cref{thm-gmp-mI-CBN} now gives:
\[ A \Perp^{\id}_\Gk B \given C \qquad \implies \qquad X_A \Indep_{\Pk(X_V\Vert X_J)} X_B \given X_C,   \]
which is a stronger conclusion and provides us with the needed Markov kernels for free. 
This was one of the core motivation for developing transitional conditional independence.
The derivation of the causal do-calculus rules from a global Markov property of this kind is carried out in
\cite{FM20} for $\Qcal$-extended conditional independence. The same strategy applies to transitional
conditional independence, once the global Markov property, \Cref{thm-gmp-mI-CBN}, is established for the
relevant class of graphs; we do not carry this out here.

\end{document}